\documentclass[10pt]{amsart}

\usepackage{xypic,amscd,amssymb,latexsym}
\usepackage{amsmath}	
\usepackage[breaklinks=true]{hyperref}

\usepackage{mathrsfs}

\usepackage{url}

\usepackage{pdfsync}

\usepackage{cancel}

\begin{document}
\pagenumbering{arabic}
\newtheorem{theorem}{Theorem}
\newtheorem{proposition}[theorem]{Proposition}
\newtheorem{lemma}[theorem]{Lemma}
\newtheorem{corollary}[theorem]{Corollary}
\newtheorem{remark}[theorem]{Remark}
\newtheorem{definition}[theorem]{Definition}
\newtheorem{question}[theorem]{Question}
\newtheorem{claim}[theorem]{Claim}
\newtheorem{conjecture}[theorem]{Conjecture}
\newtheorem{defprop}[theorem]{Definition and Proposition}
\newtheorem{example}[theorem]{Example}
\newtheorem{deflem}[theorem]{Definition and Lemma}

\def\qed{{\quad \vrule height 8pt width 8pt depth 0pt}}

\newcommand{\cplx}[0]{\mathbb{C}}

\newcommand{\vs}[0]{\vspace{2mm}}

\newcommand{\til}[1]{\widetilde{#1}}

\newcommand{\mcal}[1]{\mathcal{#1}}

\newcommand{\ul}[1]{\underline{#1}}

\newcommand{\ol}[1]{\overline{#1}}

\newcommand{\wh}[1]{\widehat{#1}}

\newcommand{\mut}[1]{\stackrel{#1}{\to}}

\address{School of Mathematics, Korea Institute for Advanced Study (KIAS), 85 Hoegiro Dongdaemun-gu, Seoul 130-722, Republic of Korea}

\email[H.~Kim]{hkim@kias.re.kr, hyunkyu87@gmail.com}

\author{Hyun Kyu Kim}

\numberwithin{equation}{section}

\title[Phases in quantum cluster varieties: long version]{Phase constants in the Fock-Goncharov quantization of cluster varieties: long version}

\begin{abstract}
A cluster variety of Fock and Goncharov is a scheme constructed from the data related to the cluster algebras of Fomin and Zelevinsky. A seed is a combinatorial data which can be encoded as an $n\times n$ matrix with integer entries, or as a quiver in special cases, together with $n$ formal variables. A mutation is a certain rule for transforming a seed into another seed; the new variables are related to the previous variables by some rational expressions. To each seed one attaches an $n$-dimensional torus, and by gluing the tori along the birational maps defined by the mutation formulas, one constructs a cluster variety. Quantization of a cluster variety assigns to each seed a non-commutative ring which deforms the classical ring of functions on the torus attached to the seed, as well as to each mutation an isomorphism of skew fields of fractions of these non-commutative rings. A representation realizes the non-commutative rings as algebras of operators on Hilbert spaces, and the quantum mutation isomorphisms as unitary maps between the Hilbert spaces that intertwine the operators for the rings. These unitary intertwiners are one of the major results of the Fock-Goncharov quantization of cluster varieties, and are given by the special function called the quantum dilogarithm. The classical mutations satisfy certain algebraic relations, which were known to be satisfied also by the corresponding intertwiners up to complex constants of modulus $1$. The present paper shows by computation that these constants are all $1$. One implication is that the mapping class group representations resulting from the application of the Fock-Goncharov quantization to the quantum Teichm\"uller theory are genuine, not projective.
\end{abstract}

\maketitle


\tableofcontents

\section{Introduction}

A `cluster variety' of Fock and Goncharov \cite{FG09} \cite{FG09b} can be viewed as a algebro-geometric space whose ring of regular functions is an `upper cluster algebra' of Fomin and Zelevinsky \cite{FoZ02}. More precisely, it is a scheme over $\mathbb{Z}$ constructed by gluing split algebraic tori $(\mathbb{G}_m)^n$ along certain birational maps, where the tori are enumerated by the data called `seeds' in the theory of cluster algebras, while the birational maps correspond to the `mutations' of seeds. There are three different types of cluster varieties, namely $\mcal{A}$, $\mcal{X}$, and $\mcal{D}$, according to how these birational maps are given. In fact, only the cluster $\mcal{A}$-variety provides a genuine example of a cluster algebra in the sense of \cite{FoZ02}, while the $\mcal{X}$- and the $\mcal{D}$-varieties are certain generalizations. 

\vs

Fix a `rank' $n \in \mathbb{Z}_{>0}$. Underlying a \emph{seed} is a combinatorial data, which essentially can be encoded as an $n\times n$ integer matrix $\varepsilon$, called the \emph{exchange matrix}, which is required to be `skew-symmetrizable', i.e. must be skew-symmetric when multiplied by some diagonal matrix from the left. In particular, when $\varepsilon$ is skew-symmetric, it can be realized as the adjacency matrix of a `quiver', i.e. a graph with oriented edges, without cycles of length $2$. A seed is also equipped with $n$ algebraically independent (formal) variables, which can be thought of as attached to the vertices of a quiver. If one picks any $k \in \{1,\ldots,n\}$, or a vertex of a quiver, one can apply the \emph{mutation} $\mu_k$ to the seed to obtain another seed; the quiver or the exchange matrix changes according to a certain combinatorial rule, while the new $n$ variables are related to the former $n$ variables by certain rational formulas which depend on the exchange matrix and are `positive', i.e. can be written as quotients of polynomials with positive coefficients. The mutation formulas for the exchange matrices and the attached variables might seem quite ad hoc at a first glance. So it is remarkable that such a structure appears in very many areas of mathematics.

\vs

What is even more remarkable is a Poisson-like structure on a cluster variety, which is written in a simple form on each torus in terms of the exchange matrix and is compatible under the mutation formulas. More precisely, the $\mcal{A}$-, the $\mcal{X}$-, and the $\mcal{D}$-varieties are equipped with a $2$-form (of `$K_2$ class'), a Poisson structure, and a symplectic form, respectively. So it is sensible to investigate whether there exist `deformation quantizations' of the $\mcal{X}$- and the $\mcal{D}$-varieties, which I explain now. 

\vs

For a (classical) cluster variety we associate to each seed a torus $(\mathbb{G}_m)^n$, which can be viewed as being defined as the space whose `ring of regular functions' is the ring of Laurent polynomials in $n$ variables over $\mathbb{Z}$. So, one can think of this situation as attaching such a commutative ring to each seed. To each mutation we associate a birational map between the two tori, that is, an isomorphism of the fields of fractions of the attached Laurent polynomials rings. 

\vs

To quantize a cluster variety means that we would like to associate to each seed a one-parameter family of non-commutative rings which deforms the classical commutative ring of Laurent polynomials, while to each mutation we associate an isomorphism of the skew fields of fractions of the non-commutative rings associated to the relevant seeds. First, this `quantum' isomorphism of skew fields must be a deformation of the corresponding `classical' isomorphism of fields of fractions of the commutative Laurent polynomial rings. Second, more importantly, there must be a `consistency'. Namely, there exists sequences of mutations that return the same exchange matrices after their application and whose induced classical isomorphisms of (commutative) fields of fractions are the identity maps. It is natural to require that the quantum isomorphisms of the skew fields of fractions induced by these sequences are also identity maps. However, such sequences are not classified in the classical setting. At this point one can only say that there are some known examples: the `rank $1$ identity' $\mu_k\circ \mu_k = \mathrm{id}$ and `rank $2$ identities', including the `pentagon identity' $\mu_i \mu_j \mu_i \mu_j \mu_i = (i \, j)$ in case $\varepsilon_{ij}=-\varepsilon_{ji} \in \{1,-1\}$, where $(i \, j)$ stands for the relabling $i\leftrightarrow j$. Fock and Goncharov \cite{FG09} constructed quantum isomorphisms of skew fields associated to mutations and checked that they satisfy this consistency for these known sequences of mutations.

\vs

As is often the case for a `physical' quantization, one would like to `represent' the non-commutative ring associated to a seed as an algebra of operators, i.e. to find a representation of this ring on a Hilbert space, satisfying certain desirable analytic properties. Then, to each mutation one would like to associate a unitary map between the Hilbert spaces assigned to the relevant seeds that  `intertwines' the representations of the non-commutative rings which are related by the quantum isomorphism. In the Fock-Goncharov quantization \cite{FG09}, one advantage of having this intertwining map is that its formula is more transparent than the quantum isomorphism formulas for the non-commutative algebras. In particular, the intertwiner is written quite neatly in terms of a famous special function called the `quantum dilogarithm' of Faddeev and Kashaev \cite{F95} \cite{FK94}, while the formula for the quantum isomorphism of algebras is not as enlightening.

\vs

A priori, the `consistency' of these unitary intertwiners for mutations for the above mentioned sequences of mutations is not guaranteed, and must be verified separately if one wants to have it. However, Fock and Goncharov \cite{FG09} deduce this consistency of intertwiners essentially by taking advantage of the `irreducibility' of the representations of rings which are being intertwined. A rough idea is like the famous Schur's lemma, which says that a vector space endomorphism that commutes with all the operators for an irreducible representation of a (finite) group on that vector space is a scalar operator. In particular, for each sequence of mutations inducing the identity classical isomorphism of fields of fractions, the composition of corresponding unitary intertwiners is a scalar operator on a relevant Hilbert space. One may interpret this result as having a unitary projective representation  on a Hilbert space of the groupoid of mutations, where the word `projective' carries the connotation that the relations are satisfied only up to complex scalars. As all operators are unitary, these scalars are of modulus $1$, hence can be viewed as \emph{phases}. The existence of these phase scalars, which can be thought of as a certain `anomaly', is not at all a defect of the program of quantization of cluster varieties. As a matter of fact, this (projective) anomaly alone may already contain a very interesting bit of information. One of the prominent examples of cluster varieties are the \emph{Teichm\"uller spaces} of Riemann surfaces, where the \emph{mapping class groups} can be embedded into the groupoid of mutations. Quantization of cluster varieties is thus one approach to the `quantum Teichm\"uller theory', one of the main results of which is the unitary projective representations of mapping class groups on Hilbert spaces. The anomaly can then be interpreted as group $2$-cocycles of mapping class groups, or equivalently, central extensions of these groups, and this led to the works of Funar-Sergiescu \cite{FS10} and Xu \cite{X14}.

\vs

However, these phase scalars for the projective anomaly have not been determined as precise complex numbers. The main original result of the present paper is that these scalars are all $1$ (Propositions \ref{prop:c_A1_is_1}, \ref{prop:c_A1_times_A1_is_1}, \ref{prop:c_A2_is_1}, \ref{prop:c_B2_is_1}, \ref{prop:c_G2_is_1}). So, we get \emph{genuine} representations of groupoid of mutations and mapping class groups on Hilbert spaces, instead of projective representations. To prove this I perform explicit computation involving unitary operators, instead of resorting just to the irreducibility of representations and the ``Schur's lemma"-type argument. While doing so I discovered an unfortunate mistake in \cite{FG09}, which is not merely a typo but is easily fixable; namely, their intertwiner does not intertwine the representation they take in that paper for the non-commutative ring associated to each seed (Lem.\ref{lem:bf_K_prime_conjugation_on_bf_x}). One should replace it by their old representation which was used in a previous paper \cite{FG07}; one possible explanation for a reason why their old representation is more favorable than their new one is in terms of the `canonical' quantization of cotangent bundles, which may be discussed in the future work \cite{KS}. In addition, I replaced their `Heisenberg relations' by the corresponding `Weyl relations' which are more rigid, and also introduced the concept of \emph{special affine shift operators} on $L^2(\mathbb{R}^n)$ (\S\ref{subsec:special_affine_shift_operators}), which are induced by the `special affine' transformations $\mathrm{SL}_\pm(n,\mathbb{R}) \ltimes \mathbb{R}^n$ of $\mathbb{R}^n$, in order to facililate notations and computations for unitary operators, which may be adapted for future projects too. A rather easy yet crucial observation is that any special affine shift operators that is proportional to the identity operator is exactly the identity operator (Lem.\ref{lem:scalar_special_affine_shift_operator_is_identity}). The key point of the proof of the main result of the present paper is to show, using quantum dilogarithm identities, that the composition of the unitary intertwiners corresponding to the relevant mutation sequences for each `consistency' relation is a special affine shift operator. As Fock and Goncharov proved that these are scalar operators, it follows that these are the identity operators. Among the quantum dilogarithm identities is the famous pentagon identity of the quantum dilogarithm function. Certain generalizations of the pentagon identity are observed, namely the `hexagon' and the `octagon' identities (\S\ref{subsec:heuristic_proofs_of_hexagon_and_octagon_identities}), which I claim to follow from the pentagon; these identities will be treated more thoroughly in \cite{KY}.

\vs

The present `long version' of the paper contains some details and rigorous proofs regarding the functional analysis that is used in the construction of the quantization of cluster varieties, as well as friendly discussions on various aspects of the quantization. It is written in such a manner to help any reader, including myself, who is interested in understanding and eventually getting hands on the work of Fock and Goncharov, with the expense of being pretty long. In particular, it will serve as a basic framework for a joint work in progress with Carlos Scarinci on quantization of moduli spaces of $3$-dimensional globally hyperbolic spacetimes \cite{KS}. One can find a shortened version \cite{K16} prepared for journal submission, focusing on the original result, namely the computation of the phase constants appearing in the unitary representation of Fock-Goncharov's quantum cluster varieties. 

\vs

The triviality of these phases for example lets us remove the multiplicative constants in the quantum dilogarithm identities as stated in \cite[Thm.4.6]{KN}. But at the same time it unfortunately deprives of the representation theoretic meaning from the works \cite{FS10} \cite{X14} at the moment, leaving only the group theoretic meaning, which is further investigated in my work \cite{K12} in comparison with another quantization of Teichm\"uller spaces obtained by Kashaev. It was observed in \cite{K12} that these two quantizations of Teichm\"uller spaces are different in a very interesting sense, related to braid groups. In my opinion, this suggests that there exists a quantization of cluster varieties 1) with non-trivial phase constants, and furthermore 2) with these constants being (half-)integer powers of $\zeta = e^{-\frac{\pi\sqrt{-1}}{6}} c_\hbar$, where $c_\hbar=e^{-\frac{\pi\sqrt{-1}}{12}(\hbar+\hbar^{-1})}$ \eqref{eq:c_hbar} and $\hbar$ is the real quantization parameter. However, such a quantization has not been constructed yet, hence calls for a future investigation.

\vs

\noindent{\bf Acknowledgments.} I would like to thank Myungho Kim, Dylan Allegretti, Ivan Chi-ho Ip, Seung-Jo Jung, Woocheol Choi, Louis Funar, Vladimir V. Fock, and Alexander B. Goncharov for helpful discussions.

\section{Cluster $\mcal{A}$-, $\mcal{X}$- and $\mcal{D}$-varieties}

I shall describe Fock-Goncharov's cluster $\mcal{A}$-varieties, cluster $\mcal{X}$-varieties, and cluster $\mcal{D}$-varieties. Definitions and treatment in this section are from \cite{FG09} and references therein, with certain modifications. As in \cite{FG09}, we confine ourselves to simpler cases where there are no `frozen variables', although it is not hard to include them as well.

\subsection{Seed and seed tori}

In the present paper, we choose and fix one positive integer $n$, which can be regarded as the `rank' of the relevant cluster algebras/varieties.
\begin{definition}
\label{def:skew_symmetrizable_matrix}
Let $K$ be a field. We say that an $n\times n$ matrix $\varepsilon = (\varepsilon_{ij})_{i,j=1,\ldots,n}$ with entries in $K$ is \ul{\em skew-symmetrizable} if there exist $d_1,\ldots,d_n \in K^*$ such that the matrix $(\wh{\varepsilon}_{ij})_{i,j}$ defined as
\begin{align}
\label{eq:wh_varepsilon}
\wh{\varepsilon}_{ij} := \varepsilon_{ij} \, d_j^{-1}
\end{align}
is skew-symmetric, i.e. $\wh{\varepsilon}_{ij} = - \wh{\varepsilon}_{ji}$, $\forall i,j$. We call such a collection $d = (d_i)_{i=1,\ldots,n}$ a \ul{\em skew-symmetrizer} \footnote{not a standard terminology} of the matrix $(\varepsilon_{ij})_{i,j}$.
\end{definition}
For any skew-symmetrizable matrix over $K=\mathbb{Q}$, one can always find a skew-symmetrizer consisting only of integers.

\begin{definition}[seed]
\label{def:seed}
An \ul{\em $\mcal{A}$-seed} $\Gamma$ is a triple $(\varepsilon, d, \{A_i\}_{i=1}^n)$ of an $n\times n$ skew-symmetrizable integer matrix $\varepsilon$, a skew-symmetrizer $d$ consisting of positive integers, and an ordered set of formal variables $A_1,A_2,\ldots,A_n$ called the \ul{\em cluster $\mcal{A}$-variables}.

\vs

An \ul{\em $\mcal{X}$-seed} $\Gamma$ is a triple $(\varepsilon, d, \{X_i\}_{i=1}^n)$ of $\varepsilon$, $d$ as above, and an ordered set of variables $X_1,\ldots,X_n$ called the \ul{\em cluster $\mcal{X}$-variables}.
\vs

An \ul{\em $\mcal{D}$-seed} $\Gamma$ is a triple $(\varepsilon, d,\{B_i,X_i\}_{i=1}^n)$ of $\varepsilon,d$ as above, and an ordered set of variables $B_1,\ldots,B_n,X_1,\ldots,X_n$ called the \ul{\em cluster $\mcal{D}$-variables}.

\vs

We call any of these a \ul{\em seed}, $\varepsilon$ the \ul{\em exchange matrix} of the seed, the relevant formal variables the \ul{\em cluster variables}\,\footnote{Fock and Goncharov call these `cluster coordinates' in \cite{FG09}.} of the seed, and the relevant symbols $\mcal{A}$, $\mcal{X}$, $\mcal{D}$ the \ul{\em kind} of the seed. We write a seed by $\Gamma = (\varepsilon, d, *)$ if it is clear from the context what the cluster variables are. \footnote{This symbol $\Gamma$ has nothing to do with the cluster mapping class group (cluster modular group)  in \cite{FG09}.}
\end{definition}

Some words must be put in order. What Fock and Goncharov \cite{FG09} call a `feed' is essentially a seed without the cluster variables\,\footnote{Goncharov says that the uncommon terminology `feed' used in \cite{FG09} was supposed to be a joke.}, i.e. the exchange matrix $\varepsilon$ and the skew-symmetrizer $d$. I think it makes more sense to use seeds as defined in Def.\ref{def:seed} instead of feeds, in the construction of `cluster varieties', whence I do so in the present paper. An $\mcal{A}$-seed, with the data of $d$ left out, is what is called a `labeled seed' with `trivial coefficients' in the theory of cluster algebras initiated by Fomin and Zelevinsky \cite{FoZ02}.  A $\mcal{D}$-seed minus the data of $d$ is what is called a `labeled seed' with `coefficients' in the `universal semifield' in \cite{FoZ07}; $\{B_i\}_{i=1}^n$ is their `cluster' and $X_i$'s their coefficients. Readers will be able to verify these comparison when they get to the  formulas for the `mutations' in \S\ref{subsec:seed_mutations}. On the other hand, the cluster $\mcal{X}$-variables are not the usual `cluster variables' in the sense of Fomin-Zelevinsky. 
So, a `seed' as defined in the above Def.\ref{def:seed} is, say,  a `generalized seed together with the choice of a skew-symmetrizer'.  

\vs

For a $\mcal{D}$-seed there is an interesting set of redundant variables, which are included in what we will casually refer to as cluster ($\mcal{D}-$)variables:
\begin{definition}[tilde variables for $\mcal{D}$-seed]
For a $\mcal{D}$-seed $\Gamma = (\varepsilon, d, \{B_i,X_i\}_{i=1}^n)$, define new variables $\til{X}_i$, $i=1,\ldots,n$, by
$$
\til{X}_i := X_i \, \prod_{j = 1}^n B_j^{\varepsilon_{ij}}, \qquad \forall i = 1,\ldots,n.
$$
\end{definition}

For each seed, we consider a parametrized topological space whose coordinate functions are identified with the cluster variables.
\begin{definition}[seed tori]
\label{def:seed_tori}
To the seeds in Def.\ref{def:seed}, all denoted by $\Gamma$ by abuse of notations, we assign respectively the split algebraic tori
\begin{align*}
\mcal{A}_\Gamma := (\mathbb{G}_m)^n, \qquad
\mcal{X}_\Gamma = (\mathbb{G}_m)^n, \qquad
\mcal{D}_\Gamma = (\mathbb{G}_m)^{2n},
\end{align*}
called the \ul{\em seed $\mcal{A}$-torus}, the \ul{\em seed $\mcal{X}$-torus}, and the \ul{\em seed $\mcal{D}$-torus}, where $\mathbb{G}_m$ is the multiplicative algebraic group, defined as $\mathbb{G}_m (K) = K^*$ for any field $K$. We call these \ul{\em seed tori}, collectively. We identify the canonical coordinate functions of $\mcal{A}_\Gamma$ with the seed's cluster variables $A_i$'s, those of $\mcal{X}_\Gamma$ with $X_i$'s, and those of $\mcal{D}_\Gamma$ with $B_i$'s and $X_i$'s.
\end{definition}

So, one $\mcal{A}$-seed $\Gamma$ gives us one topological space $\mcal{A}_\Gamma$. We shall consider $\mcal{A}_\Gamma$ for different $\Gamma$'s, and glue them together by certain algebraic formulas, to construct a `cluster $\mcal{A}$-variety'. We do likewise for the spaces $\mcal{X}_\Gamma$ for different $\Gamma$'s to construct a `cluster $\mcal{X}$-variety', and similarly for $\mcal{D}_\Gamma$'s to construct a `cluster $\mcal{D}$-variety'. I will explain these gluings in the subsequent subsections. Note that we never glue seed tori of different kinds, whence the abuse of notations committed by labeling the different kinds of seeds by the same letter $\Gamma$ will not be too harmful. 

\vs

When constructing a scheme by gluing local `patches', why glue tori $(\mathbb{G}_m)^n$ (or $(\mathbb{G}_m)^{2n}$), instead of, say, affine spaces $\mathbb{A}^n$ (or $\mathbb{A}^{2n}$), where $\mathbb{A}(K)=K$ for any field $K$? This is because we would like the ring of regular functions on each local patch to be the ring of Laurent polynomials in the coordinate functions on the patch, for a reason I will speak about later (Rem.\ref{rem:Laurent_phenomenon}). At the moment, I ask the reader to just accept this preference for Laurent polynomial rings.

\vs

We construct a local patch as an irreducible affine variety defined as the ${\rm Spec}$ of the specified {\em ring of regular functions} on it. This ring of regular functions is what we care more than the topology of the local patch. We take this ring to be the ring of Laurent polynomials in $n$ (resp. $2n$) variables over $\mathbb{Z}$ (instead of, say, the ring of polynomials in $n$ variables). The corresponding affine variety, i.e. the ${\rm Spec}$ of this ring, is the split algebraic torus of rank $n$ (or $2n$), showing why it is sensible to have Def.\ref{def:seed_tori}.

\vs

What make the cluster varieties a lot richer and more interesting are the remarkable geometric structures they possess on top of their topology, which are defined on their local patches as follows.
\begin{definition}[geometric structures on seed tori]
\label{def:geometric_structures_on_seed_tori}
We equip the following $2$-form on $\mcal{A}_\Gamma$ 
\begin{align}
\label{eq:2-form_on_A_Gamma}
\Omega_\Gamma = \sum_{i,j \in \{1,\ldots,n\}} \, \til{\varepsilon}_{ij} \, d \log A_i \wedge d\log A_j, \quad \mbox{where} \quad \til{\varepsilon}_{ij} := d_i \varepsilon_{ij},
\end{align}
and the following Poisson structure on $\mcal{X}_\Gamma$
$$
\{X_i, X_j\} = \wh{\varepsilon}_{ij} \, X_i X_j, \qquad \forall i,j \in \{1,\ldots,n\}.
$$
On $\mcal{D}_\Gamma$ we consider the Poisson structure defined by
\begin{align}
\label{eq:Poisson_on_D}
\{B_i, B_j\} = 0, \quad
\{X_i, B_j\} = d_i^{-1} \delta_{ij} X_i B_j, \quad
\{X_i, X_j\} = \wh{\varepsilon}_{ij} X_i X_j,
\end{align}
where $\delta_{ij}$ is the Kronecker delta, as well as the following $2$-form
\begin{align}
\label{eq:2-form_on_D}
- \frac{1}{2} \sum_{i,j\in \{1,\ldots,n\}} \til{\varepsilon}_{ij} \, d\log B_i \wedge d \log B_j - \sum_{i \in \{1,\ldots,n\}} d_i \, d\log B_i \wedge d\log X_i.
\end{align}
\end{definition}

\begin{remark}
Fock and Goncharov say that $\Omega_\Gamma$ on $\mcal{A}_\Gamma$ is of `$K_2$ class'.
\end{remark}

\begin{lemma}[\cite{FG09}]
The $2$-form \eqref{eq:2-form_on_D} on $\mcal{D}_\Gamma$ is a symplectic form, and is compatible with the Poisson structure \eqref{eq:Poisson_on_D} on it.
\end{lemma}

\subsection{Seed mutations and seed automorphisms}
\label{subsec:seed_mutations}

In order to study the gluing of the seed tori, we first need to investigate the following transformation rules for seeds, which we can also understand as a way of recursively creating new seeds from previously constructed ones.
\begin{definition}[seed mutation]
\label{def:mutation}
For $k\in\{1,2,\ldots,n\}$, a seed $\Gamma' = (\varepsilon', d', *')$  is said to be obtained by applying the \ul{\em seed mutation $\mu_k$ in the direction $k$} to a seed $\Gamma = (\varepsilon, d, *)$, if all of the following hold: 

\begin{itemize}
\item the exchange matrices are related by
\begin{align}
\label{eq:varepsilon_prime_formula}
\varepsilon'_{i j} = \left\{
\begin{array}{ll}
-\varepsilon_{ij} & \mbox{if $i=k$ or $j=k$,} \\
\varepsilon_{ij} + \frac{1}{2}( |\varepsilon_{ik}| \varepsilon_{kj} + \varepsilon_{ik} | \varepsilon_{kj} | )& \mbox{otherwise}.
\end{array}
\right.
\end{align}

\item the skew symmetrizers are related by
\begin{align}
\label{eq:d_prime_formula}
d'_i = d_i, \qquad \forall i.
\end{align}

\item the cluster variables are related by:
\begin{enumerate}
\item[\rm 1)] In the case of $\mcal{A}$-seeds: 
\begin{align}
\label{eq:mu_k_on_A}
A'_i = \left\{
\begin{array}{ll}
A_i & \mbox{if $i\neq k$}, \\
\displaystyle A_k^{-1} \left( \prod_{j \, | \, \varepsilon_{kj}>0 } A_j^{\varepsilon_{kj}} + \prod_{j \, | \, \varepsilon_{kj}<0} A_j^{-\varepsilon_{kj}} \right) & \mbox{if $i=k$},
\end{array}
\right.
\end{align}
where $\Gamma = (\varepsilon, d, \{A_i\}_{i=1}^n)$ and $\Gamma' = (\varepsilon', d', \{A_i'\}_{i=1}^n)$.

\vs

\item[\rm 2)] In the case of $\mcal{X}$-seeds: 
\begin{align}
\label{eq:mu_k_on_X}
X'_i = \left\{
\begin{array}{ll}
X_k^{-1} & \mbox{if $i=k$}, \\
X_i \left( 1 + X_k^{{\rm sgn}(-\varepsilon_{ik})} \right)^{-\varepsilon_{ik}} & \mbox{if $i\neq k$,}
\end{array}
\right.
\end{align}
where $\Gamma = (\varepsilon, d, \{X_i\}_{i=1}^n)$ and $\Gamma' = (\varepsilon', d', \{X_i'\}_{i=1}^n)$, and
$$
{\rm sgn}(a)= \left\{
\begin{array}{ll}
1 & \mbox{if $a\ge 0$} \\
-1 & \mbox{if $a<0$}.
\end{array}
\right.
$$

\vs

\item[\rm 3)] In the case of $\mcal{D}$-seeds: the formulas in  \eqref{eq:mu_k_on_X}, together with
\begin{align}
\label{eq:mu_k_on_B}
B'_i = \left\{
\begin{array}{ll}
B_i & \mbox{if $i\neq k$}, \\
\displaystyle \frac{ \left(\prod_{j \, | \, \varepsilon_{kj}<0 } B_j^{-\varepsilon_{kj}}\right) + X_k \left(\prod_{j\, | \, \varepsilon_{kj}>0} B_j^{\varepsilon_{kj}}\right) }{B_k(1+X_k)} & \mbox{if $i=k$,} 
\end{array}
\right.
\end{align}
where $\Gamma = (\varepsilon, d, \{B_i,X_i\}_{i=1}^n)$ and $\Gamma' = (\varepsilon', d', \{B_i',X_i'\}_{i=1}^n)$.
\end{enumerate}

\end{itemize}

We denote such a situation by $\mu_k(\Gamma) = \Gamma'$, or $\Gamma \stackrel{k}{\to} \Gamma'$. We call this procedure a \ul{\em mutation}.
\end{definition}

The new cluster variables are viewed as elements of the `ambient field', the field of all rational functions in the previous cluster variables over $\mathbb{Q}$. Each of \eqref{eq:mu_k_on_A}--\eqref{eq:mu_k_on_B} should be thought of as an equality in the ambient field.

\begin{remark}
\cite{FG09} uses the notations $\mathbb{B}^+_k := \prod_{j \, | \, \varepsilon_{kj}>0} B_j^{\varepsilon_{kj}}$ and $\mathbb{B}^-_k := \prod_{j \, | \, \varepsilon_{kj}<0} B_j^{-\varepsilon_{kj}}$, which help simplifying the formula \eqref{eq:mu_k_on_B}.
\end{remark}

\begin{definition}[seed automorphism]
\label{def:seed_automorphism}
A seed $\Gamma' = (\varepsilon',d',*')$ is said to be obtained by applying the \ul{\em seed automorphism $P_\sigma$} to a seed $\Gamma = (\varepsilon,d,*)$ for a permutation $\sigma$ of $\{1,2,\ldots,n\}$, if
$$
d'_{\sigma(i)} = d_i, \qquad \varepsilon'_{\sigma(i) \, \sigma(j)} = \varepsilon_{ij}, \qquad \forall i,j,
$$
and
$$
A'_{\sigma(i)} = A_i, \quad(\mbox{for $\mcal{A}$-seeds}), \qquad
X'_{\sigma(i)} = X_i, \quad(\mbox{for $\mcal{X}$-, $\mcal{D}$-seeds}), \qquad
B'_{\sigma(i)} = B_i, \quad(\mbox{for $\mcal{D}$-seeds}),
$$
for all $i=1,2,\ldots,n$. We denote such a situation by $P_\sigma(\Gamma) = \Gamma'$.\footnote{$P_\sigma$ is denoted simply by $\sigma$ in \cite{FG09}.}
\end{definition}
I will often just say a `permutation' when I refer to a seed automorphism. We can consider the mutations and seed automorphisms as being applied to seeds from the left, and use the usual notation $\circ$ for the composition of them. So, when we apply a finite sequence of mutations and permutations to a seed, we read the sequence from right to left.

\begin{definition}[cluster transformations]
A \ul{\em cluster transformation} is a finite sequence of mutations and seed automorphisms, together with a seed $\Gamma$ on which the sequence is to be applied to. If the resulting seed is $\Gamma'$, we say that this cluster transformation \ul{\em connects} $\Gamma$ \ul{\em to} $\Gamma'$. We call mutations and seed automorphisms \ul{\em elementary cluster transformations}.
\end{definition}

Note that, each of the symbols $\mu_k$ and $P_\sigma$ stands for many different elementary cluster transformations, for it can be thought of as being applied to different seeds, which could be of any of the three kinds.

\subsection{Cluster modular groupoids}
\label{subsec:cluster_modular_groupoids}

In practice, we start from a single seed, which we often call an `initial seed', and produce new seeds by applying cluster transformations to it. Then the set of seeds thus created will be in bijection with the set of all finite sequences of mutations and permutations. However, I would like to identify two seeds whenever they are `essentially' the same, in the sense I explain now. 

\vs

Given any cluster transformation, say, connecting a seed $\Gamma$ to a seed $\Gamma'$, one obtains a well-defined identification of the cluster variables for $\Gamma'$ as rational functions in the cluster variables for $\Gamma$, by `composing' the formulas in Definitions \ref{def:mutation} and \ref{def:seed_automorphism}.
\begin{definition}[trivial cluster transformations]
\label{def:trivial_cluster_transformations}
A cluster transformation connecting a seed $\Gamma = (\varepsilon,d,*)$ to a seed $\Gamma'=(\varepsilon',d',*')$ is said to be \ul{\em weakly trivial} if $\varepsilon'=\varepsilon$ and $d'=d$. A cluster transformation is said to be \ul{\em $\mcal{A}$-trivial} if it is weakly trivial and induces the identity map between the cluster $\mcal{A}$-variables, i.e. the $i$-th cluster $\mcal{A}$-variable $A_i$ for $\Gamma$ equals the $i$-th cluster $\mcal{A}$-variable $A_i'$ for $\Gamma'$. The notions \ul{\em $\mcal{X}$-trivial} and \ul{\em $\mcal{D}$-trivial} are defined analogously. A cluster transformation is said to be \ul{\em trivial} if it is $\mcal{A}$-trivial, $\mcal{X}$-trivial, or $\mcal{D}$-trivial. 
\end{definition}

\begin{remark}
In \cite{FG09} a `(feed) cluster transformation' means a sequence of mutations and permutations applied to a feed $(\varepsilon,d)$, and it is said to be `trivial' if it is $\mcal{A}$-trivial and $\mcal{X}$-trivial at the same time, when applied to $\mcal{A}$-seeds and $\mcal{X}$-seeds whose underlying feeds are $(\varepsilon,d)$.
\end{remark}

\begin{definition}[identification of seeds]
\label{def:identification_of_seeds}
If a cluster transformation connecting an $\mcal{A}$-seed $\Gamma$ to an $\mcal{A}$-seed $\Gamma'$ is $\mcal{A}$-trivial, we identify $\Gamma$ and $\Gamma'$ as $\mcal{A}$-seeds, and write $\Gamma = \Gamma'$. Likewise for $\mcal{X}$-seeds and $\mcal{D}$-seeds.
\end{definition}
Keeping this identification of seeds in mind, we consider:
\begin{definition}[equivalence of seeds]
\label{def:equivalence_of_seeds}
The two seeds are \ul{\em equivalent} if they are connected by a cluster transformation. For a seed $\Gamma$, denote by $\mathscr{C} = |\Gamma|$ the equivalence class of $\Gamma$. \footnote{Identifying seeds connected by seed automorphisms amounts to considering Fomin-Zelevinsky's notion of `unlabeled' seeds.}
\end{definition}

Let us now investigate some examples of trivial cluster transformations. It is a standard and straightforward exercise to show that the following example is indeed a trivial cluster transformation.
\begin{lemma}[involution identity of a mutation]
\label{lem:involution_identity}
$\mu_k \circ \mu_k$ is a trivial cluster transformation on any seed of any kind. That is, if we write $\Gamma' = \mu_k(\Gamma)$ and $\Gamma'' = \mu_k(\Gamma')$, with $\Gamma'' = (\varepsilon'',d'',*'')$ and $\Gamma = (\varepsilon,d,*)$, then $\varepsilon''=\varepsilon$, $d''=d$, and $A_i'' = A_i$, $\forall i$ (for $\mcal{A}$-seeds), $X_i'' = X_i'$, $\forall i$ (for $\mcal{X}$- and $\mcal{D}$-seeds), and $B_i'' = B_i$, $\forall i$ (for $\mcal{D}$-seeds). \qed
\end{lemma}
In particular, we identify $\mu_k(\mu_k(\Gamma))$ and $\Gamma$. Some more easy-to-see trivial cluster transformations involving the seed automorphisms are as follows, which I also omit the proof of.
\begin{lemma}[permutation identities]
\label{lem:permutation_identities}
For any permutations $\sigma,\gamma$ of $\{1,\ldots,n\}$ and any $k\in \{1,\ldots,n\}$,
\begin{align*}
P_\sigma \circ P_\gamma = P_{\sigma \circ \gamma}, \qquad
P_\sigma \circ \mu_k \circ P_{\sigma^{-1}} = \mu_{\sigma(k)}, \qquad
P_{{\rm id}} = {\rm id},
\end{align*}
all of which hold as identities when applied to any seed of any kind, where $\mathrm{id}$ stands either for the identity permutation or the identity cluster transformation. \qed
\end{lemma}
The first two identities in this statement can be translated into saying that $P_\sigma \circ P_\gamma \circ P_{\sigma\circ \gamma}^{-1}$ and $P_\sigma \circ \mu_k \circ P_{\sigma^{-1}} \circ \mu_{\sigma(k)}^{-1}$ are trivial cluster transformations on any seed of any kind.

\vs

Besides these, there remain certain interesting trivial cluster transformations, pointed out e.g. in eq.(20) in of \cite[p.238]{FG09}:
\begin{lemma}[The $(h+2)$-gon relations]
\label{lem:h_plus_2-gon_relations}
Suppose that a seed $\Gamma = (\varepsilon,d,*)$ of any kind \footnote{\cite{FG09} restricts only to $\mcal{A}$-seeds and $\mcal{X}$-seeds, but I think it also holds for $\mcal{D}$-seeds.} satisfies
$$
\varepsilon_{ij} = - p \, \varepsilon_{ji} \quad \mbox{and} \quad |\varepsilon_{ij}| = p
$$
for some $i,j\in \{1,\ldots,n\}$ and $p\in \{0,1,2,3\}$.\footnote{\cite{FG09} uses $\varepsilon_{ij} = - p \, \varepsilon_{ji} = -p$, but it doesn't matter.} Let $h=2,3,4,6$ for $p=0,1,2,3$, respectively. Denote by $(i\, j)$ the transposition permutation of $\{1,\ldots,n\}$ interchanging $i$ and $j$. Then
\begin{align}
\mbox{$(P_{(i\, j)} \circ \mu_i)^{h+2}$ applied to $\Gamma$ is a trivial cluster transformation.}
\end{align}
\end{lemma}

Let us stop here to reflect on the trivial cluster transformations just presented. The permutation identities in Lem.\ref{lem:permutation_identities} do not really require a proof. The involution identity in Lem.\ref{lem:involution_identity} can be thought of as coming from `rank one' cluster algebras, namely, of Dynkin type $A_1$. The more serious identities in Lem.\ref{lem:h_plus_2-gon_relations} can be thought of as coming from `rank two' cluster algebras of Dynkin types $A_1\times A_1$, $A_2$, $B_2$, $G_2$, respectively, as pointed out in \cite{FG09}. Fock and Goncharov say in \cite{FG09} that they do not know how to find more trivial cluster transformations that are not consequences of the ones already discussed, let alone how to describe the complete list of trivial cluster transformations, which certainly exists. Thus they take only these known relations to formulate the problem and the result of their quantization.

\vs

Here let me share what I learned from a personal discussion with Myungho Kim. A `seed' in the usual theory of cluster algebras \cite{FoZ02}, maybe with `coefficients' \cite{FoZ07}, is said to be of `finite type' if its equivalence class is a finite set, i.e. the number of distinct seeds that can be produced by applying cluster transformations to the initial seed is finite, where we identify some seeds according to Def.\ref{def:identification_of_seeds}. In their remarkable work \cite{FoZ03}, Fomin and Zelevinsky classified all finite type seeds, where the classification is described miraculously by the famous Dynkin diagrams of finite type, which  appear in the Cartan-Killing classification of finite dimensional semi-simple Lie algebras. The rank of the Dynkin diagram coincides with the `rank' $n$ of the seed, i.e. the number of non-frozen (i.e. mutable, or exchangeable) variables. Keeping this result in mind, in the case of rank one or rank two seeds, one can readily deduce that some identities like Lem.\ref{lem:involution_identity} and Lem.\ref{lem:h_plus_2-gon_relations} hold, even without any serious computational check; namely, it is easy to see that {\em any} cluster transformation in the equivalence class of a finite type seed is of finite order. In fact, for any seed that `locally looks like' rank one\footnote{this is always so.} or rank two seeds\footnote{this amounts to the condition of Lem.\ref{lem:h_plus_2-gon_relations} being satisfied.}, these relations still hold, because one can regard this as a seed with only one or two non-frozen variables, where other cluster variables are regarded as frozen variables, or coefficients. Now, take any finite type seed. Think of any infinite sequence of mutations and permutations, and apply these one by one, starting from the initial finite type seed. We then keep producing seeds, and at some point, two of the seeds created thus far must coincide. Hence we obtain a trivial cluster transformation. One can do the same for any seed $\Gamma$ that looks locally like a finite type seed, to obtain trivial cluster transformations in the equivalence class $|\Gamma|$. However, we do not get any essentially new trivial cluster transformation this way, because it is possible to prove that any trivial cluster transformation on a finite type seed is a concatenation of finitely many trivial cluster transformations coming from rank 1 and rank 2. A proof of this statement can be extracted from \cite{FoZ03}. Next question is whether there exists a trivial cluster transformation that is not a consequence of the ones coming from the rank 1 and rank 2 cases. I remark that this question is answered negatively in \cite{FST}, in most cases of seeds of `surface type', coming from ideal triangulations of surfaces. More discussion on trivial cluster transformations will appear in a joint work in progress with Masahito Yamazaki \cite{KY}.

\begin{remark}
The `known' examples presented here and used in \cite{FG09} are all examples of sequences of mutations and permutations applied to seeds of any kind with designated underlying `feed' that are trivial in the sense of the last sentence in Def.\ref{def:trivial_cluster_transformations}. Is there an example of a sequence that is trivial in some kind but not on another?
\end{remark}

An idealistic structure that the equivalent seeds form is the following:
\begin{definition}
The \ul{\em $\mcal{A}$-cluster modular groupoid $\mcal{G}^\mcal{A} = \mcal{G}^\mcal{A}_\mathscr{C}$} associated to an equivalence class $\mathscr{C}$ of $\mcal{A}$-seeds is a category whose set of objects is $\mathscr{C}$, and whose set of morphisms ${\rm Hom}_{\mcal{G}^\mcal{A}} (\Gamma, \Gamma')$ from an object $\Gamma$ to $\Gamma'$ is the set of all cluster transformations from $\Gamma$ to $\Gamma'$ modulo $\mcal{A}$-trivial cluster transformations.

\vs

Likewise for \ul{\em the $\mcal{X}$-cluster modular groupoid $\mcal{G}^\mcal{X}$} and \ul{\em the $\mcal{D}$-cluster modular groupoid $\mcal{G}^\mcal{D}$}.
\end{definition}
Here `modulo' means the following. Consider cluster transformations ${\bf m}_1, {\bf m}_2, {\bf m}_3$, where some of them may be empty sequences of elementary cluster transformations. If ${\bf m}_2$ is a $\mcal{A}$-trivial cluster transformation, then ${\bf m}_1 \circ {\bf m}_2 \circ {\bf m}_3$ and ${\bf m}_1 \circ {\bf m}_3$ are viewed as the same morphism in the category $\mcal{G}^\mcal{A}$. Likewise for $\mcal{G}^\mcal{X}$ and $\mcal{G}^\mcal{D}$. I myself have always found the following description of these groupoids useful.
\begin{lemma}
\label{lem:only_one_morphism}
The set of morphisms ${\rm Hom}_{\mcal{G}^\mcal{A}}(\Gamma,\Gamma')$ between any two objects of $\mcal{G}^\mcal{A} = \mcal{G}^\mcal{A}_\mathscr{C}$ consists of exactly one element. Likewise for $\mcal{G}^\mcal{X}$ and $\mcal{G}^\mcal{D}$.
\end{lemma}

{\it Proof.} By definition of $\mathscr{C}$, the set ${\rm Hom}_{\mcal{G}^\mcal{A}}(\Gamma,\Gamma')$ has at least one element. Suppose ${\bf m}_1$ and ${\bf m}_2$, which are sequences of mutations and permutations being applied to $\Gamma$, are its elements. Then ${\bf m}_1(\Gamma) = \Gamma' = {\bf m}_2(\Gamma)$. Let ${\bf m}_2^{-1}$ be the sequence obtained by inverting each entry of the sequence ${\bf m}_2$ and putting in the reverse order, being applied to $\Gamma'$; here, `inverting' means replacing each $\mu_i$ and $P_\sigma$ by $\mu_i$ and $P_{\sigma^{-1}}$ respectively. Then, a repeated application of Lemmas \ref{lem:involution_identity} and \ref{lem:permutation_identities} tells us $({\bf m}_2^{-1} \circ {\bf m}_2)(\Gamma) = \Gamma$, hence ${\bf m}_2^{-1}(\Gamma') = \Gamma$. Define cluster transformations ${\bf m}_3 := {\bf m}_2^{-1} \circ {\bf m}_1$, being applied to $\Gamma$, and ${\bf m}_4 := {\bf m}_2 \circ {\bf m}_2^{-1}$, being applied to $\Gamma'$. Then ${\bf m}_3 (\Gamma) = {\bf m}_2^{-1}(\Gamma) = \Gamma$ and ${\bf m}_4(\Gamma') = {\bf m}_2(\Gamma) = \Gamma'$, hence ${\bf m}_3$ and ${\bf m}_4$ are trivial cluster transformations. Thus ${\bf m}_2$ and ${\bf m}_2 \circ {\bf m}_3$ are the same morphisms in $\mcal{G}^\mcal{A}$ from $\Gamma$ to $\Gamma'$. Note ${\bf m}_2 \circ {\bf m}_3 = {\bf m}_2 \circ {\bf m}_2^{-1} \circ {\bf m}_1 = {\bf m}_4 \circ {\bf m}_1$, which is the same morphism as ${\bf m}_1$. So ${\bf m}_2$ is the same morphism as ${\bf m}_1$, which is the desired result for $\mcal{G}^\mcal{A}$. Similarly for $\mcal{G}^\mcal{X}$ and $\mcal{G}^\mcal{D}$. \qed

\vs

A slightly less idealistic version is the following, taking into account only the known relations.

\begin{definition}
\label{eq:saturated-modular_groupoid}
The \ul{\em saturated $\mcal{A}$-cluster modular groupoid $\wh{\mcal{G}}^\mcal{A} = \wh{\mcal{G}}^\mcal{A}_\mathscr{C}$} associated to an equivalence class $\mathscr{C}$ of $\mcal{A}$-seeds is a category whose set of objects is $\mathscr{C}$, and whose set of morphisms ${\rm Hom}_{\wh{\mcal{G}}^\mcal{A}}(\Gamma, \Gamma')$ from an object $\Gamma$ to $\Gamma'$ is the set of all cluster transformations from $\Gamma$ to $\Gamma'$ modulo only the trivial cluster transformations that are described in Lemmas \ref{lem:involution_identity}, \ref{lem:permutation_identities} and \ref{lem:h_plus_2-gon_relations}. Likewise for $\wh{\mcal{G}}^\mcal{X}$ and $\wh{\mcal{G}}^\mcal{D}$.
\end{definition}

\begin{definition}
In these groupoids, an \ul{\em elementary morphism} is a morphism representing an elementary cluster transformation.\footnote{probably not a completely standard terminology (used by Funar, Sergiescu and collaborators)}
\end{definition}

As we shall soon see, these groupoids provide a handy way to formulate the construction of `cluster varieties' and the quantization. However, Fock and Goncharov \cite{FG09} define groupoids based on `feeds' $(\varepsilon,d)$, rather than seeds. Using theirs has an advantage as it allows to define the notion of (saturated) cluster modular group, or (saturated) cluster mapping class group\footnote{in the case of cluster variety `coming from a surface', this group coincides with the `mapping class group' of the surface.}, as the group of automorphisms of one feed; in particular, upon quantization, one gets projective representations of the cluster mapping class group. However, in view of the properties of the scheme that we will soon construct, feeds do not seem to be the best way to go.

\subsection{Cluster varieties}

Each elementary cluster transformation, say from a seed $\Gamma$ to a seed $\Gamma'$, induces a rational map from the seed torus for $\Gamma$ to that for $\Gamma'$, denoted by the same symbol as the elementary cluster transformation itself, defined in the level of functions by the identification formulas for the cluster variables, obtained by `composing' the formulas in Definitions \ref{def:mutation} and \ref{def:seed_automorphism} as mentioned earlier. For example, we define the rational map $\mu_k : \mcal{X}_\Gamma \to \mcal{X}_{\Gamma'}$ between the seed $\mcal{X}$-tori by describing what the pullbacks of the coordinate functions of the torus $\mcal{X}_{\Gamma'}$ are, in terms of coordinate functions of $\mcal{X}_\Gamma$:
\begin{align*}
\mu_k^*(X'_i) = \left\{
\begin{array}{ll}
X_k^{-1} & \mbox{if $i=k$}, \\
X_i \left( 1 + X_k^{{\rm sgn}(-\varepsilon_{ik})} \right)^{-\varepsilon_{ik}} & \mbox{if $i\neq k$,}
\end{array}
\right.
\end{align*}
and $\mu_k^*( {X_i'}^{-1} ) = {\mu_k^*(X_i')}^{-1}$, $\forall i$.

\vs

It is easy to see that the map $P_\sigma$ between two seed tori is an isomorphism of varieties. 
The map $\mu_k$ between two seed tori is indeed a rational map, for it is defined on the subsets of the respective seed tori on which the denominators of the transformation formulas are not zero, which are open subsets whose complements in the tori are of lower dimension. As each $\mu_k$ has a rational inverse map, namely $\mu_k$ (Lem.\ref{lem:involution_identity}), we see that $\mu_k$ is in fact a birational map.

\vs

So we can now say that any cluster transformation from a seed $\Gamma$ to a seed $\Gamma'$ induces a birational map between the two seed tori. We finally construct cluster varieties by gluing the seed tori along these birational maps.
\begin{definition}[cluster varieties]
\label{def:cluster_varieties}
The \ul{\em cluster $\mcal{A}$-variety} for an equivalence class $\mathscr{C}$ of $\mcal{A}$-seeds, is a scheme obtained by gluing all the seed tori $\mcal{A}_{\Gamma}$ for $\Gamma$ in this equivalence class $\mathscr{C}$ using the birational maps associated to cluster transformations. We denote it by $\mcal{A}_{\mathscr{C}}$, or just $\mcal{A}$ if the equivalence class is apparent from the context.

\vs

Define the \ul{\em cluster $\mcal{X}$-variety} $\mcal{X} = \mcal{X}_{\mathscr{C}}$ and the \ul{\em cluster $\mcal{D}$-variety} $\mcal{D} = \mcal{D}_{\mathscr{C}}$ in a similar fashion.
\end{definition}
In order for $\mcal{A}_\mathscr{C}$ to be well-defined, for each pair of seeds $(\Gamma, \Gamma')$ in $\mathscr{C}$ we must have a unique birational map $\mcal{A}_\Gamma \to \mcal{A}_{\Gamma'}$ along which we glue the two tori. There may be many different sequences of mutations and permutations that connect $\Gamma$ to $\Gamma'$, and we should make sure that they induce the same birational map on the tori. Readers can easily deduce this by using Lem.\ref{lem:only_one_morphism}, using the fact that trivial cluster transformations induce identity maps between the tori. Similarly for $\mcal{X}_\mathscr{C}$ and $\mcal{D}_\mathscr{C}$ too.

\vs

This definition is different from that of Fock and Goncharov \cite{FG09}. Their cluster variety is defined as the ${\rm Spec}$ of the ring of global regular functions on the scheme constructed in a similar way as in Def.\ref{def:cluster_varieties} using `feed tori' instead of seed tori. In our seed-version, this ring could be defined as follows.
\begin{definition}
For a seed $\Gamma$ of any kind, denote by $R_\Gamma$ the ring of regular functions on the seed torus for $\Gamma$, and let $\mathbb{Q}_\Gamma$ be the field of fractions of $R_\Gamma$.

For equivalent seeds $\Gamma$ and $\Gamma'$, denote by $\mu_{\Gamma,\Gamma'}^*$ the unique isomorphism $\mathbb{Q}_{\Gamma'} \to \mathbb{Q}_\Gamma$ of fields as described above on the generators.

For an equivalence class $\mathscr{C}$ of seeds, the \emph{\ul{ring of global regular functions}} $\mathbb{L}_\mathscr{C}$ is defined as the following subring of $\mathbb{Q}_\Gamma$
\begin{align}
\label{eq:ring_of_global_regular_functions}
\mathbb{L}_\Gamma ~ := ~ \bigcap_{\Gamma' \in \mathscr{C}} \mu^*_{\Gamma,\Gamma'}(R_{\Gamma'}) \quad \subset \quad \mathbb{Q}_\Gamma,
\end{align}
for any chosen $\Gamma \in \mathscr{C}$.
\end{definition}
To be a bit more explicit about the above notation, the ring $R_\Gamma$ for an $\mathcal{A}$-seed $\Gamma = (\varepsilon,d, \{A_i\}_{i=1}^n)$, for example, is the ring $\mathbb{Z}[A_1,A_1^{-1}, \ldots,A_n,A_n^{-1}]$.
Note that \eqref{eq:ring_of_global_regular_functions} `does not depend on the choice of $\Gamma \in \mathscr{C}$' in the sense that
$$
\mathbb{L}_\Gamma = \mu_{\Gamma,\Gamma'}^*( \mathbb{L}_{\Gamma'} ), \qquad \forall \Gamma, \Gamma' \in \mathscr{C}.
$$
We collectively denote all $\mathbb{L}_\Gamma$ for $\Gamma\in \mathscr{C}$ by the symbol $\mathbb{L}_\mathscr{C}$.
\begin{remark}
As $R_\Gamma$ is the ring of Laurent polynomials (in $n$ or $2n$ variables), Fock and Goncharov calls $\mathbb{L}_\mathscr{C}$ the ring of `universally Laurent polynomials'; $\mathbb{L}$ stands for the `L'aurent. In fact, they defined such a ring only in the quantum case, as we shall encounter later. In the case of a genuine cluster algebra setting, e.g. for $\mathcal{A}$-seeds, $\mathbb{L}_\mathscr{C}$ is called the \emph{upper cluster algebra}.
\end{remark}
Fock and Goncharov say that thus finally obtained variety $\mathrm{Spec}(\mathbb{L}_\mathscr{C})$ is the `affine closure'\footnote{This may not be a standard usage of the term `affine closure'.} of the one constructed just by gluing the tori. Although taking an affine closure may guarantee somewhat nicer topological properties, I am not sure whether it preserves the scheme structures, which I would like to keep for the sake of upcoming quantization, in the cases of $\mcal{X}$-variety and $\mcal{D}$-variety. So, in the present article, I stick to Def.\ref{def:cluster_varieties} without taking the `affine closure'.

\begin{remark}
If there is only one torus for one feed, which possibly underlies many seeds, then this feed torus may be glued to itself in weird ways, creating potential analytical or even topological problems. So we better associate a torus to a seed.
\end{remark}

\begin{remark}
\label{rem:Laurent_phenomenon}
The `Laurent phenomenon' of cluster algebras \cite{FoZ02} says that when $\Gamma$ is an $\mcal{A}$-seed, then for each $\mcal{A}$-seed $\Gamma'$ equivalent to $\Gamma$, the image under $\mu^*_{\Gamma,\Gamma'}$ of each cluster $\mcal{A}$-variable of $\Gamma'$ is in $\mathbb{L}_\Gamma$. The subalgebra of $\mathbb{L}_\Gamma$ generated by all these images is the \ul{\em cluster algebra} for the seed $\Gamma$. So, there are `many' global regular functions the cluster $\mcal{A}$-variety. If we had declared that the ring of regular functions on the local patch associated to a seed should be the ring of polynomials for example, instead of Laurent polynomials, then there would not be so many globally regular functions.
\end{remark}

The cluster varieties defined in Def.\ref{def:cluster_varieties} can be formulated in the language of functors on the (saturated) cluster modular groupoids considered in \S\ref{subsec:cluster_modular_groupoids}. Little more generally, a scheme with an atlas whose charts (or, the corresponding `local patches') are enumerated by $\mathscr{C}$ can be thought of as a contravariant functor
$$
\eta : \mbox{the cluster modular groupoid} ~ \mcal{G}_\mathscr{C} ~ \longrightarrow ~ \mbox{certain category of commutative rings},
$$
where $\mcal{G}_\mathscr{C}$ stands for one of $\mcal{G}^\mcal{A}_\mathscr{C}$, $\mcal{G}^\mcal{X}_\mathscr{C}$, and $\mcal{G}^\mcal{D}_\mathscr{C}$, while morphisms from an object to another in the category in the RHS are homomorphisms from the field of fraction of the latter object to that of the former. For our situation, one may further require that the objects of the RHS category must be isomorphic to the ring of Laurent polynomials in $n$ (or $2n$) variables over $\mathbb{Z}$, and that we only allow homomorphisms whose image on the generators can be written as rational expressions using only additions, multiplications, and divisions (i.e. `subtraction-free'). Then the local patch corresponding to an element of $\mathscr{C}$ is set to be the ${\rm Spec}$ of the image of that element under $\eta$, and we glue the patches corresponding to each two elements of $\mathscr{C}$ by the birational map given by the image under $\eta$ of the unique morphism in $\mcal{G}_\mathscr{C}$ from the first to the latter; recall from Lem.\ref{lem:only_one_morphism} that in $\mcal{G}_\mathscr{C}$ there is exactly one morphism from any object to an object.  Such functor $\eta$ is a slight reformulation of Fock-Goncharov's `positive space' \cite[Def.1.3]{FG09b}, in which the objects in the target category (denoted by ${\rm Pos}$) are split algebraic tori and the morphisms are positive rational maps between the tori. 

\vs

This groupoid formulation of the cluster varieties is not just an unnecessary luxury to have, for such a formulation in terms of groupoids and rings is the only known way to describe the quantum versions of cluster varieties; there is no actual topological space in the quantum world.

\vs

The following geometric structures, which are extra data on the cluster varieties, are crucial in the story of quantization.
\begin{lemma}[geometric structures on cluster varieties]
\label{lem:geometric_structures_on_cluster_varieties}
The geometric structures on the seed tori defined in Def.\ref{def:geometric_structures_on_seed_tori} induce well-defined corresponding geometric structures on the respective cluster varieties.
\end{lemma}
One checks this lemma by verifying that $\mu_k$ and $P_\sigma$ preserve the geometric structures on the seed tori; I omit the computation.

\subsection{Maps among three varieties}

Only in the present subsection, the ambiguity of the symbol $\Gamma$ comes in action. Here, let it denote a seed of any of the three kinds, with a fixed underlying data $(\varepsilon,d)$, which is called a `feed' in \cite{FG09}. For each kind, consider the equivalence class $\mathscr{C}=|\Gamma|$ of seeds for $\Gamma$, and construct the corresponding cluster variety; denote them by $\mathcal{A}$, $\mathcal{X}$, $\mathcal{D}$, respectively. We first describe the map
$$
p : \mathcal{A} \to \mathcal{X},
$$
given by gluing the regular maps
$$
p_\Gamma : \mathcal{A}_\Gamma \to \mathcal{X}_\Gamma
$$
for each $\Gamma \in \mathscr{C}$, by the pullback formula (see e.g. \cite[Thm.2.3.(b)]{FG09})
\begin{align}
  \label{eq:p_Gamma_formula}
  p^*_\Gamma X_i = \prod_{j=1}^n A_j^{\varepsilon_{ij}}, \qquad \forall i=1,\ldots,n.
\end{align}
First, one can check that different $p_\Gamma$'s are compatible with respect to the mutations and seed automorphisms, hence they indeed glue to give a map $p: \mathcal{A} \to \mathcal{X}$. In general, this map $p$ is neither surjective nor injective. However, it nicely `respects' the geometric structures discussed in Def.\ref{def:geometric_structures_on_seed_tori} and Lem.\ref{lem:geometric_structures_on_cluster_varieties}. Namely, the part of $\mathcal{A}$ that is killed by the map $p$ is exactly the null foliation of the degenerate $2$-form on $\mathcal{A}$; more precisely, the restriction onto $p(\mathcal{A})$ of the Poisson structure on $\mathcal{X}$ is symplectic and the pullpack of this symplectic $2$-form under $p$ coincides with the degenerate $2$-form on $\mathcal{A}$. The triple $(p,\mathcal{A},\mathcal{X})$ is what is called a `cluster ensemble' by Fock and Goncharov in \cite{FG09b}. Moreover, there are commutative diagrams (\cite[Thm.2.3.(d)]{FG09})
\begin{align*}
\xymatrix@R-4mm{
\mathcal{A} \times \mathcal{A} \ar[rd]^-{\varphi} \ar[dd]^{p\times p} & \\
& \mathcal{D} \ar[ld]^-{\pi} \\
\mathcal{X} \times \mathcal{X}^{\rm op} &
}
  \qquad\qquad\qquad
\xymatrix@R+5mm@C-3mm{
\mathcal{X} \ar@{^(->}[r]^j \ar[d] & \mathcal{D} \ar[d]^{\pi} \\
\Delta_\mathcal{X} \ar[r] & \mathcal{X} \times \mathcal{X}^{\rm op}
},
\end{align*}
with appropriate maps satisfying certain desirable properties, where $\mcal{X}^{\rm op}$ means the `opposite' Possion variety of $\mcal{X}$. I refrain from saying any more detail here, and just refer the readers to \cite{FG09}.

\section{Algebraic quantization}

In \cite{FG09} and their previous works, Fock and Goncharov obtained certain `deformation quantization' of $\mcal{X}$-variety `along' its Poisson structure and that of $\mcal{D}$-variety `along' its symplectic structure; some version of the former follows from the latter, which I think is more canonical. Instead of giving a comprehensive explanation of the meaning of `deformation quantization' in general, I only present here what were obtained in \cite{FG09}.

\subsection{Basic formulas for algebraic quantization}

To a $\mcal{D}$-seed $\Gamma$ we associate the following family of non-commutative associative $*$-algebras that deforms the coordinate ring of the seed $\mcal{D}$-torus $\mcal{D}_\Gamma$ `in the direction of' its symplectic or Poisson structure.
\begin{definition}[quantum torus algebra for a seed $\mcal{D}$-torus]
\label{def:quantum_D-torus_algebra}
Let $\Gamma=(\varepsilon,d,\{B_i,X_i\}_{i=1}^n)$ be a $\mcal{D}$-seed, and $N$ be the smallest positive integer such that $\wh{\varepsilon}_{ij} = \varepsilon_{ij}/d_j \in \frac{1}{N}\mathbb{Z}$, $\forall i,j$ \footnote{One could just use $N=$ the least common multiple of all $d_i$'s.}. The \ul{\em seed quantum $\mcal{D}$-torus algebra ${\bf D}_\Gamma^q$} associated to $\Gamma$ for a quantum parameter $q$, which can be regarded as a formal parameter at the moment, is the free associate $*$-algebra over $\mathbb{Z}[q^{1/N}, q^{-1/N}]$ generated by ${\bf B}_1,\ldots,{\bf B}_n,{\bf X}_1,\ldots,{\bf X}_n$ and their inverses\footnote{I use bold faced letters in order to distinguish from the `classical' cluster variables.}, mod out by the relations
\begin{align*}
\begin{array}{rcll}
q_j^{-\varepsilon_{ij}} {\bf X}_i {\bf X}_j & = & q_i^{-\varepsilon_{ji}} {\bf X}_j {\bf X}_i, & \quad \forall i,j\in \{1,\ldots,n\}, \\
q_i^{-1} {\bf X}_i {\bf B}_i & = & q_i \, {\bf B}_i {\bf X}_i, & \quad \forall i \in \{1,\ldots,n\}, \\ 
{\bf B}_i {\bf X}_j & = & {\bf X}_j {\bf B}_i, & \quad \mbox{whenever $i\neq j$}, \\
{\bf B}_i {\bf B}_j & = & {\bf B}_j {\bf B}_i, & \quad \forall i,j \in \{1,\ldots,n\},
\end{array}
\end{align*}
where
\begin{align}
\label{eq:q_i}
q_i := (q^{1/N})^{N/d_i}, \qquad \forall i \in \{1,\ldots,n\},
\end{align}
with the $*$-structure defined as the unique ring anti-homomorphism satisfying
$$
*{\bf X}_i = {\bf X}_i, \qquad *{\bf B}_i = {\bf B}_i, ~~~ \forall i \in \{1,\ldots,n\}, \qquad *q^{1/N} = q^{-1/N}.
$$
Denote by $\mathbb{D}^q_\Gamma$ the skew field of fractions of ${\bf D}^q_\Gamma$. \footnote{`skew field' is a synonym of `division ring'}
\end{definition}
\begin{remark}
Eq.\eqref{eq:q_i} can be viewed as $q_i = q^{1/d_i}$.
\end{remark}
I think that it is a mistake that the relations ${\bf B}_i {\bf B}_j = {\bf B}_j {\bf B}_i$, which do not follow from others, are omitted in eq.(58) in p.252 \cite[\S3]{FG09}.

\begin{remark}
As in p232 of \cite{FG09}, to any lattice $\Lambda$, i.e. a free abelian group, equipped with a skew-symmetric bilinear form $\Lambda \times \Lambda \to \frac{1}{N} \mathbb{Z}$, one can associate a `quantum torus algebra', which is an associative $*$-algebra over $\mathbb{Z}[q^{1/N}, q^{-1/N}]$, with generators enumerated by $\Lambda$. Def.\ref{def:quantum_D-torus_algebra} is the case when $\Lambda$ is of rank $2n$, whose basis vectors corresponds to the generators $\mathbf{B}_1,\ldots,\mathbf{B}_n,\mathbf{X}_1,\ldots,\mathbf{X}_n$, where the skew-symmetric form comes from the Poisson bracket among the classical variables $B_1,\ldots,B_n,X_1,\ldots,X_n$ as we see in \eqref{eq:Poisson_on_D}.
\end{remark}

\begin{definition}[tilde variables for the seed quantum $\mcal{D}$-torus algebra]
For each $i\in \{1,\ldots,n\}$, define an element $\til{{\bf X}}_i$ of ${\bf D}^q_\Gamma$ by
$$
\til{{\bf X}}_i := {\bf X}_i \, \prod_{j =1}^n {\bf B}_j^{\varepsilon_{ij}}.
$$
\end{definition}

Then one can check
$$
q_j^{\varepsilon_{ij}} \til{\mathbf{X}}_i \til{\mathbf{X}}_j = q_i^{\varepsilon_{ji}} \til{\mathbf{X}}_j \til{\mathbf{X}}_i, \qquad
q_i^{-\delta_{i,j}} \til{\mathbf{X}}_i \, \mathbf{B}_j = q_i^{\delta_{i,j}} \mathbf{B}_j \til{\mathbf{X}}_i, \qquad
\mathbf{X}_i \til{\mathbf{X}}_j = \til{\mathbf{X}}_j \mathbf{X}_i,
$$
for all $i,j \in \{1,\ldots,n\}$.

\vs

To a mutation $\Gamma \stackrel{k}{\to} \Gamma'$ we associate the following $*$-isomorphism $\mathbb{D}^q_{\Gamma'} \to \mathbb{D}^q_\Gamma$ of skew fields.
\begin{definition}[quantum mutation map]
\label{def:quantum_mutation_map}
Let the two $\mcal{D}$-seeds $\Gamma=(\varepsilon,d,*)$ and $\Gamma'=(\varepsilon',d',*')$ be related by the mutation along $k$, that is, $\Gamma' = \mu_k(\Gamma)$. Define the map $\mu^q_k : \mathbb{D}^q_{\Gamma'} \to \mathbb{D}^q_\Gamma$ by
$$
\mu^q_k := \mu^{\sharp q}_k \circ \mu'_k,
$$
where $\mu^{\sharp q}_k$ is the automorphism of $\mathbb{D}^q_\Gamma$ given by the following formulas on generators
\begin{align}
\label{eq:quantum_mutation_sharp_on_B_i}
\mu^{\sharp q}_k ({\bf B}_i) & = \left\{
\begin{array}{ll}
{\bf B}_i & \mbox{if $i\neq k$}, \\
{\bf B}_k (1+q_k {\bf X}_k) (1+q_k \til{{\bf X}}_k)^{-1} & \mbox{if $i=k$},
\end{array}
\right. \\
\label{eq:quantum_mutation_sharp_of_X_i}
\mu^{\sharp q}_k ({\bf X}_i) & = {\bf X}_i \prod_{r=1}^{|\varepsilon_{ik}|} ( 1 + (q_k^{{\rm sgn}(-\varepsilon_{ik})} )^{2r-1} {\bf X}_k )^{ {\rm sgn}(-\varepsilon_{ik}) }, \qquad \forall i \in \{1,\ldots,n\},
\end{align}
and $\mu'_k$ is induced by the map $\mu'_k : {\bf D}^q_{\Gamma'} \to {\bf D}^q_\Gamma$ given by
\begin{align}
\nonumber
\mu'_k({\bf B}_i') = \left\{
\begin{array}{ll}
{\bf B}_i & \mbox{if $i\neq k$,} \\
{\bf B}_k^{-1} \prod_{j=1}^n {\bf B}_j^{[-\varepsilon_{kj}]_+}  & \mbox{if $i = k$,}
\end{array}
\right.
\qquad
\mu'_k({\bf X}_i') = \left\{
\begin{array}{ll}
q_k^{-\varepsilon_{ik}[\varepsilon_{ik}]_+} {\bf X}_i \, ({\bf X}_k)^{[\varepsilon_{ik}]_+} & \mbox{if $i \neq k$}, \\
{\bf X}_k^{-1} & \mbox{if $i=k$,}
\end{array}
\right.
\end{align}
on the generators of ${\bf D}^q_{\Gamma'}$, where $[~~]_+$ denotes the `positive part' of a real number:
$$
[a]_+ := \frac{a+  |a|  }{2} = \left\{ \begin{array}{ll} a, & \mbox{if $a\ge 0$}, \\ 0, & \mbox{otherwise} \end{array} \right., \qquad \forall a\in \mathbb{R}.
$$
\end{definition}

\begin{definition}[quantum permutation map]
\label{def:quantum_permutation_map}
Let $\Gamma' = P_\sigma(\Gamma)$ for a permutation $\sigma$ of $\{1,\ldots,n\}$. Define the map $\mathbb{P}_\sigma : \mathbb{D}^q_{\Gamma'} \to \mathbb{D}^q_\Gamma$ to be the one induced by the map $\mathbb{P}_\sigma : {\bf D}^q_{\Gamma'} \to {\bf D}^q_\Gamma$ given by
$$
\mathbb{P}_\sigma ({\bf B}'_{\sigma(i)})= {\bf B}_i, \qquad
\mathbb{P}_\sigma ({\bf X}'_{\sigma(i)})= {\bf X}_i, \qquad
\forall i \in \{1,\ldots,n\}.
$$
\end{definition}

For a quantization, we would like to assign to each cluster transformation a well-defined quantum map between the skew fields of fractions of the respective seed quantum $\mcal{D}$-torus algebras, such that 1) it is a $*$-isomorphism of skew fields, and 2) it recovers the classical map between the seed $\mcal{D}$-tori in the classical limit $q\to 1$. The above defined maps $\mu_k^q$ and $\mathbb{P}_\sigma$ are what we assign to elementary cluster transformations, and it is a straightforward task to check that these satisfy the two conditions 1) and 2) just mentioned. 

\vs

A general cluster transformation is the composition of a sequence of elementary ones, and we would assign to it the composition of the sequence of corresponding quantum maps $\mu_k^q$ and $\mathbb{P}_\sigma$ for the elementary cluster transformations, although in a reverse order, for the maps $\mu_k^q$ and $\mathbb{P}_\sigma$ are constructed in a `contravariant' manner. Problem is, there can be several sequences of elementary cluster transformations expressing the same classical cluster transformation, and we have to make sure that these sequences yield the same quantum map. Equivalently, we must make sure that each trivial cluster transformation is assigned the identity map of the relevant quantum $\mcal{D}$-torus algebra. This `consistency' is the key aspect in the algebraic part of the program of quantization of cluster varieties. A standard way of formulating this problem is via the (saturated) cluster modular groupoids considered in \S\ref{subsec:cluster_modular_groupoids}, and we do this in the following subsection.

\subsection{Groupoid formulation of algebraic quantization}

\begin{definition}
The \ul{\em category of seed quantum $\mcal{D}$-torus algebras ${\rm QDTor}^q={\rm QDTor}^q_{\mathscr{C}}$} associated to an equivalence class $\mathscr{C}$ of $\mcal{D}$-seeds with quantum parameter $q$ is a category whose objects are ${\bf D}^q_\Gamma$ with $\Gamma$ in the equivalence class $\mathscr{C}$, and the set morphisms ${\rm Hom}_{{\rm QDTor}^q}({\bf D}^q_\Gamma, {\bf D}^q_{\Gamma'})$ is the set of all $*$-isomomorphisms from the skew field $\mathbb{D}^q_{\Gamma'}$ to the skew field $\mathbb{D}^q_\Gamma$ whose images of the generators of $\mathbb{D}^q_{\Gamma'}$ are subtraction-free elements of $\mathbb{D}^q_\Gamma$. \quad (see Def.\ref{def:quantum_D-torus_algebra} for ${\bf D}^q_\Gamma$ and $\mathbb{D}^q_\Gamma$.)
\end{definition}
Here, a `subtraction-free' element is one that can be expressed in terms of generators and $q^{\pm 1/N}$ using only additions, multiplications, and divisions.
\begin{definition}
\label{def:quantum_D-space}
By a \ul{\em quantum $\mcal{D}$-space $\mcal{D}^q = \mcal{D}^q_{\mathscr{C}}$}, or \ul{\em quantum double}, for a cluster $\mcal{D}$-variety $\mcal{D} = \mcal{D}_{\mathscr{C}}$ with quantum parameter $q$ we mean a contravariant functor
\begin{align}
\label{eq:algebraic_quantization_functor}
\eta^q ~:~ \mbox{{\em the saturated cluster modular groupoid}} ~~
\wh{\mcal{G}}^\mcal{D}_{\mathscr{C}} ~ \longrightarrow ~ {\rm QDTor}^q_{\mathscr{C}}
\end{align}
whose images of elementary morphisms recover the formulas in Definitions \ref{def:mutation} and \ref{def:seed_automorphism} in the `classical limit' $q\to 1$. 

\vs

Let us call such a functor $\eta^q$ an \ul{\em algebraic quantum cluster ($\mcal{D}$-)variety}.
\end{definition}
\begin{remark}
\cite{FG09} views a functor $\eta^q$ in \eqref{eq:algebraic_quantization_functor} as a `quantum scheme'.
\end{remark}
Of course, an ultimate goal to attain is a functor \eqref{eq:algebraic_quantization_functor} from the cluster modular groupoid $\mcal{G}^\mcal{D}_\mathscr{C}$, not just from the saturated cluster modular groupoid $\wh{\mcal{G}}^\mcal{D}_\mathscr{C}$. What has been obtained so far is only \eqref{eq:algebraic_quantization_functor}. We know the images of the objects (Def.\ref{def:quantum_D-torus_algebra}), as well as images of the elementary morphisms (Def.\ref{def:quantum_mutation_map}, \ref{def:quantum_permutation_map}). Hence we also know the images of general morphisms, for they are sequences of elementary ones. The only problem is then the well-definedness of this construction.
\begin{lemma}[Lem.3.4 of \cite{FG09}]
\label{lem:eta_q_is_well-defined}
The functor $\eta^q$ described above is well-defined.
\end{lemma}
Fock-Goncharov \cite{FG09} says that this can be obtained either by direct computation, or by using the results of \cite[Thm.6.1]{BZ05}. In the present paper we shall prove a different statement that imply these relations.

\section{Review on functional analysis}

As is often the case in a `physical' quantization, one does not stop at the non-commutative algebras, but considers representations of them on Hilbert spaces, respecting the $*$-structure. In the present section we review some functional analytic background knowledge, which is necessary for understanding of the formulation of representations of quantum cluster varieties, as well as for the proof of the main result of the present paper. However, this section is not to be taken as a comprehensive overview of the theory of operators on Hilbert spaces, and for this I refer the readers to \cite{RS70} and \cite{Y}, which are great sources to learn the subject, and from which I collected most of the definitions and statements for this section.

\subsection{Densely defined operators on a Hilbert space}
\label{subsec:densely_defined_operators}

I first recall some basic notions on operators on Hilbert spaces. Sometimes I elaborated on some terms or modified some notations from those used in \cite{RS70} and \cite{Y}, trying to conform with and not to deviate too much from standard usage in modern days. All Hilbert spaces here are over $\cplx$. For a Hilbert space $V$ we denote by
$$
\langle \cdot , \cdot \rangle_V : V \times V \to \cplx
$$
its inner product, which is complex linear in the first argument and conjugate linear in the second argument. The norm
$$
||v||_V := \sqrt{\langle v,v\rangle_V}, \quad \forall v\in V
$$
on $V$ defines a topology on $V$, which we refer to as the \emph{standard Hilbert space topology} on $V$.

\begin{definition}[densely defined operators on a Hilbert space]
Let $V$ be a Hilbert space, and let $T$ be a $\cplx$-linear map $T : D(T) \to V$ for a $\cplx$-vector subspace $D(T)$ of $V$ that is dense in $V$ with respect to the standard Hilbert space topology on $V$. We say that $D(T)$ is the \ul{\em domain} of $T$ and that $T$ is a \ul{\em densely defined operator on $V$}. 
\end{definition}

\begin{remark}
The domain is also a part of the data of a `densely defined operator'.
\end{remark}

\begin{definition}
Let $D$ be a $\cplx$-linear subspace of a Hilbert space $V$, and let $T : D\to V$ be a $\cplx$-linear map. We say that \ul{\em $D$ is invariant under $T$} or \ul{\em $T$ leaves $D$ invariant} if $Tv \in D$, $\forall v \in D$.
\end{definition}

\begin{definition}
Let $T_1, T_2$ be densely defined operators on a Hilbert space $V$. We say that $T_2$ \ul{\em extends} $T_1$ (or $T_2$ is an \ul{\emph{extension}} of $T_1$) if $D(T_1) \subset D(T_2)$ and $T_1=T_2$ on $D(T_1)$. We write this as
$$
T_1 \subset T_2.
$$
\end{definition}

\begin{definition}
\label{def:symmetric_operator}
A densely defined operator $T$ on a Hilbert space $V$ is said to be \ul{\em symmetric} if
$$
\langle Tv, w\rangle_V = \langle v, Tw\rangle_V, \qquad \forall v,w\in D(T).
$$
\end{definition}

Let us now review the definition of the {\em adjoint} of a densely defined operator. Let $T$ be a densely defined operator on a Hilbert space $V$, with its domain $D(T)$. Consider $v\in V$ such that the following holds:
\begin{align}
\label{eq:adjoint_domain}
w \mapsto \langle Tw, v\rangle_V \mbox{ is a bounded linear functional on $D(T)$.}
\end{align}
A {\em bounded linear functional} on $D(T)$ means a $\cplx$-linear map $\rho : D(T) \to \cplx$ such that ${\displaystyle \sup_{w\in D(T)}} \frac{ \rho(w) }{ ||w||_V }$ is finite. Then the Riesz Representation Theorem (or the `Riesz Lemma') tells us that there exists a unique vector $u$ in $V$ such that $\langle Tw, v \rangle_V = \langle w, u\rangle_V$ for all $w\in D(T)$, for $D(T)$ is dense in $V$. Denote such $u$ by $u = T^* v$. Thus we obtain a certain $\cplx$-linear map $T^* : D(T^*) \to V$ where
$$
D(T^*) := \{ v\in V \, | \, \eqref{eq:adjoint_domain} \mbox{ holds}\},
$$
which is not necessarily dense in $V$.

\begin{definition}[adjoint of a densely defined operator]
\label{def:adjoint_of_an_operator}
Let $T$ be a densely defined operator on a Hilbert space $V$. The unique $\cplx$-linear map $D(T^*) \to V$ constructed above is called the \ul{\em adjoint} of $T$.
\end{definition}
Obviously we have $T\subset T^*$ for each symmetric operator $T$. In fact this is a characterization of symmetric operators among densely defined operators. 
\begin{definition}[self-adjoint and essentially self-adjoint operators]
We say that a densely operator $T$ on a Hilbert space $V$ is \ul{\em self-adjoint} if and only if $T = T^*$. 

\vs

A densely defined symmetric operator $T$ is called \ul{\em essentially self-adjoint} if and only if it has a unique self-adjoint extension.
\end{definition}

Densely defined self-adjoint operators on Hilbert spaces\footnote{The notion of self-adjoint operators applies only to densely defined operators, so we might just say `self-adjoint operators' for densely defined self-adjoint operators.} are one of the most basic objects of study in quantum theory. However, as it is tricky to deal with their domains and ranges, we often turn the situation into the story of another class of operators, which are far more nicely behaved.
\begin{definition}
A \ul{\emph{unitary operator}} on a Hilbert space $V$ is a bijective $\cplx$-linear operator $T : V\to V$ such that $\langle Tv, Tw \rangle_V = \langle v,w\rangle_V$, $\forall v,w\in V$.
\end{definition}

In particular, a unitary operator on $V$ is defined on the whole $V$. If $T$ is a densely defined operator on $V$ satisfying $\langle Tv, Tw \rangle_V = \langle v,w\rangle_V$ for all $v,w$ in its domain, we can uniquely extend $T$ to the whole $V$, for it is continuous, i.e. bounded (use the B.L.T. theorem; see e.g. \cite[Thm.I.7]{RS70}); we will then identify such $T$ with its unique extension.

\vs

For a densely defined self-adjoint operator $T$, there are two ways to obtain unitary operators. One is to consider the family $e^{\sqrt{-1} tT}$, $t\in \mathbb{R}$, of unitary operators (Stone's Theorem, discussed at the end of \S\ref{subsec:spectral_theorem}), and the other is to consider the `Cayley transform' $(T-\sqrt{-1}\cdot\mathrm{id})(T+\sqrt{-1}\cdot\mathrm{id})^{-1}$, which is also a unitary operator. For more detalied treatments, see \cite{RS70} and \cite{Y}.

\subsection{Spectral theorem, continuous functional calculus of self-adjoint operators, and Stone's theorem}
\label{subsec:spectral_theorem}

If $T$ is a bounded operator on $V$, then for any one-variable polynomial $f$ over $\mathbb{C}$, the operator $f(T)$ is well-defined as a bounded operator on $T$. One can also easily make sense of $f(T)$ for any complex-valued real-analytic function $f$ on $\mathbb{R}$ whose radius of convergence at $0$ is greater than the operator norm $||T|| := \displaystyle \sup_{||v||_V=1} ||Tv||_V$. On the other hand, a na\"ive approach does not work when $T$ is an unbounded operator, i.e. when $||T|| = \infty$. However, the `spectral theorem' allows us to consider $f(T)$ when $T$ is a densely defined (unbounded) self-adjoint operator and $f$ is any continuous $\cplx$-valued function on $\mathbb{R}$, or even when $f$ is merely a measurable function. In the present subsection we present, without proofs, one version of the theorem among its several different guises, which may not be in its most general form but is sufficient for our purposes.

\begin{definition}
An \ul{\em orthogonal projection} on a Hilbert space $V$ is a bounded $\cplx$-linear operator $P : V \to V$ such that $P^2=P$ and $P = P^*$.
\end{definition}
Orthogonal projections on $V$ are in one-to-one correspondence with closed $\cplx$-linear subspaces of $V$. Namely, the range $W:=\mathrm{Ran}(P) = \{Pv : v\in V\}$ of $P$ is a closed $\cplx$-linear subspace of $V$, and we have the vector space direct sum decomposition $V = W \oplus W^\perp$ which is also orthogonal with respect to $\langle \,\, , \,\rangle_V$, where $W^\perp = \{ u \in V : \langle u,w\rangle_V=0, ~\forall w\in W\}$; the projection $P$ then sends each $w+u$ to $w$, $\forall w\in W$, $u\in W^\perp$.

\begin{definition}[see {\cite[Chap.XI.5]{Y}}]
A family of orthogonal projections $E(\lambda)$, $\lambda\in\mathbb{R}$, on a Hilbert space $V$, is called a (real) \emph{resolution of the identity} if it satisfies
\begin{align}
  & E(\lambda) \, E(\mu) = E(\min(\lambda,\mu)), \quad\forall \lambda,\mu \in \mathbb{R}, \\
  & E(-\infty)=0, \quad E(+\infty)=\mathrm{id}_V, \\
  & E(\lambda+0) = E(\lambda), \quad \forall \lambda\in \mathbb{R},
\end{align}
where
\begin{align}
  E(-\infty) v := \lim_{\lambda\to-\infty} E(\lambda) v, \quad
  E(\infty) v := \lim_{\lambda\rightarrow\infty} E(\lambda) v, \quad
  E(\lambda+0) v := \lim_{\epsilon\searrow 0} \, E(\lambda+\epsilon)v,
\end{align}
for all $v\in V$ and $\lambda\in\mathbb{R}$, where the limits are taken with respect to the standard Hilbert space topology on $V$.
\end{definition}
For each two real numbers $\alpha,\beta$ with $\alpha<\beta$ we define
$$
E(\alpha,\beta] := E(\beta) - E(\alpha),
$$
which is an orthogonal projection.

\begin{defprop}[see {\cite[Chap.XI.5]{Y}}]
Let $E(\lambda)$ be a resolution of the identity on $V$. For any complex valued continuous function $f$ on $\mathbb{R}$ and any two real numbers $\alpha,\beta$ with $\alpha<\beta$, for each $v\in V$ we can define
\begin{align}
  \int_\alpha^\beta f(\lambda) \, dE(\lambda) \, v
\end{align}
as the limit of the Riemann sums
\begin{align}
  \sum_j f(\lambda_j') \, E(\lambda_j, \lambda_{j+1}] \, v, \quad \mbox{where} \quad \alpha = \lambda_1 < \lambda_2 < \cdots < \lambda_n = \beta, \qquad \lambda'_j \in (\lambda_j,\lambda_{j+1}]
\end{align}
as $\max_j |\lambda_{j+1}-\lambda_j|$ tends to zero.

\vs

Define
\begin{align}
  \int_{-\infty}^\infty f(\lambda) \, dE(\lambda) \, v
= \lim_{\alpha\rightarrow -\infty, \, \beta\rightarrow\infty} \,\int_\alpha^\beta f(\lambda) \, dE(\lambda) \, v
\end{align}
when the limit exists, where the limit is taken with respect to the standard Hilbert space topology on $V$.
\end{defprop}

\begin{theorem}[see {\cite[Chap.XI.5]{Y}}]
\label{thm:spectral_theorem1}
Let $f(\lambda)$ be a real-valued continuous function on $\mathbb{R}$. Let $E(\lambda)$ be a resolution of the identity on a Hilbert space $V$. Then the formula
\begin{align}
\label{eq:T_E_f}
  T_{E,f} (v) = \int_{-\infty}^\infty f(\lambda) \, d E(\lambda) \, v, \qquad \forall v\in D(T_{E,f}),
\end{align}
for 
\begin{align}
  D(T_{E,f}) := \left\{ \, v \in V \,\, \left| \,\,\, \int_{-\infty}^\infty \right. |f(\lambda)|^2 \, d||E(\lambda) \, v||_V^2 < \infty \, \right\}
\end{align}
defines a $\cplx$-linear map $T_{E,f} : D(T_{E,f}) \rightarrow V$ which is a densely defined self-adjoint operator on $V$ with the domain $D(T_{E,f})$. Moreover,
\begin{align}
\label{eq:T_E_f_norm}
||T_{E,f}(v)||_V = \int_{-\infty}^\infty |f(\lambda)|^2 \, d||E(\lambda) \, v||_V^2, \qquad \forall v\in D(T_{E,f}).
\end{align}
\end{theorem}
Often in the literature (e.g. also in \cite[Chap.XI.5]{Y}), \eqref{eq:T_E_f} is written as
\begin{align}
  \label{eq:T_E_f2}
  \langle T_{E,f}(v), \,w \rangle_V = \int_{-\infty}^\infty f(\lambda) \, d\langle E(\lambda) v, w \rangle_V, \qquad \forall v \in D(T_{E,f}), \quad \forall w\in V.
\end{align}
As a particular case of Thm.\ref{thm:spectral_theorem1}, for the identity function $f(\lambda) \equiv \lambda$, we have a densely defined operator $T := T_{E,\mathrm{id}}$ on $V$ defined by
\begin{align}
\label{eq:spectral_resolution}
T= T_{E,\mathrm{id}}(v) = \int_{-\infty}^\infty \lambda \, d E(\lambda) \, v, \qquad \forall v\in D(T_{E,\mathrm{id}}),
\end{align}
on its domain
$$
D(T_{E,\mathrm{id}}) = \left\{ \, v \in V \,\, \left| \,\,\, \int_{-\infty}^\infty \right. |\lambda|^2 \, d||E(\lambda) \, v||^2_V < \infty \, \right\}.
$$
We write \eqref{eq:spectral_resolution} symbolically as
\begin{align}
\label{eq:spectral_resolution2}
  T = \int_{-\infty}^\infty \, \lambda \, dE(\lambda),
\end{align}
and call \eqref{eq:spectral_resolution2} as the \emph{spectral resolution} of the self-adjoint operator $T$.

\vs

\begin{theorem}[see {\cite[Chap.XI.5]{Y}}]
\label{thm:spectral_theorem2}
Each (densely defined) self-adjoint operator $T$ on a Hilbert space $V$ has a unique spectral resolution.
\end{theorem}
Theorems \ref{thm:spectral_theorem1} and \ref{thm:spectral_theorem2} together can be thought of as a version of the `spectral theorem of self-adjoint operators'. One can then establish the `functional calculus' of self-adjoint operators as follows. We only discuss the version for continuous functions here, but a version for bounded Borel (measurable) functions can also be considered (see \cite[Chap.VIII]{RS70}, \cite[Chap.XI.12]{Y})

\vs

Let $T$ be a densely defined self-adjoint operator $T$ on a Hilbert space $V$. Then there exists a unique resolution of the identity $E(\lambda)$ such that $T = \int_{-\infty}^\infty \lambda \, dE(\lambda)$. Now, for any real-valued continuous function $f$ on $\mathbb{R}$, one defines the operator $f(T)$ on the domain
\begin{align}
\label{eq:D_f_T}
  D(f(T)) := \left\{ \, v \in V \,\, \left| \,\,\, \int_{-\infty}^\infty \right. |f(\lambda)|^2 \, d||E(\lambda) \, v||_V^2 < \infty \, \right\}
\end{align}
by the formula
\begin{align}
\label{eq:real_f_T}
  f(T) \, v = \int_{-\infty}^\infty f(\lambda) \, d E(\lambda) \, v, \qquad \forall v\in D(f(T)).
\end{align}
This formula \eqref{eq:real_f_T} indeed gives a well-defined $f(T) \, v \in V$ for each $v\in D(T)$ due to Thm.\ref{thm:spectral_theorem1}, which also says that $f(T)$ is a densely defined self-adjoint operator on $V$ with its domain being $D(f(T))$ in \eqref{eq:D_f_T}.

\vs

From Thm.\ref{thm:spectral_theorem1} one can also deduce that for a densely defined self-adjoint operator $T = \int_{-\infty}^\infty \lambda \, dE(\lambda)$ on a Hilbert space $V$ and a complex-valued continuous function $f$ on $\mathbb{R}$, one can define the operator $f(T)$ as $f(T) v = \int_{-\infty}^\infty f(\lambda) \, dE(\lambda) v$ on the domain $D(f(T))$ defined as in \eqref{eq:D_f_T}, which is a dense subspace of $V$. In case when $|f(\lambda)|=1$ for all $\lambda\in \mathbb{R}$, eq.\eqref{eq:T_E_f_norm} of Thm.\ref{thm:spectral_theorem1} implies $||f(T)(v)||_V^2 = \int_{-\infty}^\infty |f(\lambda)|^2 \, d||E(\lambda)v||_V^2 = \int_{-\infty}^\infty d||E(\lambda) v||_V^2 = \displaystyle \lim_{\alpha\to-\infty, \, \beta\rightarrow\infty} (||E(\beta)v||_V^2 - ||E(\alpha)v||_V^2) = ||E(\infty)v||_V^2 - ||E(-\infty)v||_V^2 = ||v||^2_V$ for all $v \in D(f(T))$. Hence $D(f(T)) = V$ and $f(T)$ is a unitary operator on $V$. Let us state this as a lemma:
\begin{lemma}
[unitary operators from functional calculus]
\label{lem:unitary_operators_from_functional_calculus}
A complex-valued function on a measure space or a topological space is called a \ul{\emph{unitary function}}\footnote{not a standard terminology} if its values lie in $\mathrm{U}(1) \subset \mathbb{C}$, i.e. are of modulus $1$. Suppose that $f$ is a unitary function on $\mathbb{R}$ and $T$ is a self-adjoint operator on a Hilbert space $V$. Then the operator $f(T)$, defined by the functional calculus described so far, is unitary. \qed
\end{lemma}

So, for a self-adjoint operator $T$ on a Hilbert space $V$ and each fixed real number $t$, one can apply the functional calculus for the complex function $f(\lambda) = e^{\sqrt{-1} \, t \lambda}$ and $T$, to get a unitary operator $f(T) = e^{\sqrt{-1} \, tT}$ on $V$.
\begin{theorem}[see e.g. {\cite[Chap.VIII.4]{RS70}}]
\label{thm:unitary_family_for_self-adjoint_operator}
If $T$ is a self-adjoint operator on a Hilbert space $V$, then the family $\{U_T^{(t)}:= e^{\sqrt{-1} \, tT}\}_{t\in \mathbb{R}}$ of unitary operators on $V$ satisfies
\begin{align}
\label{eq:Stone_condition1}
U_T^{(t)} \, U_T^{(s)} = U_T^{(t+s)}, \qquad \forall s,t\in \mathbb{R},
\end{align}
and
\begin{align}
\label{eq:Stone_condition2}
\lim_{t\rightarrow t_0} U_T^{(t)} v = U_T^{(t_0)} v, \qquad \forall t_0\in \mathbb{R} \quad\mbox{and}\quad \forall v\in V,
\end{align}
where the limit is taken with respect to the standard Hilbert space topology.
\end{theorem}

\begin{definition}[{\cite[Chap.VIII.4]{RS70}}]
\label{def:strongly_continuous_one-parameter_unitary_group}
If $\{U^{(t)}\}_{t\in \mathbb{R}}$ is a family of unitary operators on a Hilbert space $V$ satisfying \eqref{eq:Stone_condition1} and \eqref{eq:Stone_condition2}, i.e. if $U^{(t)} U^{(s)} = U^{(t+s)}$, $\forall s,t$, and $\lim_{t\rightarrow t_0} U^{(t)}v = U^{(t_0)}v$, $\forall t_0\in \mathbb{R}$, $\forall v\in V$, where the limit is taken with respect to the standard Hilbert space topology, we call this family a \ul{\em strongly continuous one-parameter unitary group}.
\end{definition}

\begin{theorem}[`Stone's theorem' \cite{St32}; see e.g. {\cite[Chap.VIII.4]{RS70}}]
\label{thm:Stones_theorem}
If $\{U^{(t)}\}_{t\in \mathbb{R}}$ is a strongly continuous one-parameter unitary group on a Hilbert space $V$, there exists a unique densely defined self-adoint operator $T$ on $V$ such that $U^{(t)} = e^{\sqrt{-1} \, tT}$, $\forall t\in \mathbb{R}$. We call $T$ the \ul{\em infinitesimal generator} of $\{U^{(t)}\}_{t\in \mathbb{R}}$.
\end{theorem}

\begin{theorem}[{\cite[Chap.VIII.4]{RS70}}]
\label{thm:strong_derivative}
Let $\{U^{(t)}\}_{t\in \mathbb{R}}$ be a strongly continuous one-parameter unitary group on a Hilbert space $V$. If $D$ is a dense $\cplx$-linear subspace of $V$ such that $U^{(t)}v \in D$ for all $v\in D$ and 
\begin{align}
\label{eq:Av}
Av := \lim_{t\to 0} \frac{ U^{(t)} v - v }{t}
\end{align}
exists for all $v\in D$ where the limit is taken with respect to the standard Hilbert space topology, then $\frac{1}{\sqrt{-1}} A$ is essentially self-adjoint on $D$, and its unique self-adjoint extension is the infinitesimal generator of $\{U^{(t)}\}_{t\in \mathbb{R}}$.
\end{theorem}

Let me add the following lemmas which are easy to prove:
\begin{lemma}[unitarily equivalent essentially self-adjoint operators]
\label{lem:unitarily_equivalent_essentially_self-adjoint_operators}
Let $T$ be an essentially self-adjoint operator on a Hilbert space $V$ with its dense domain $D(T)$, and let $U$ be a unitary operator on $V$. Then,
\begin{enumerate}
\item[\rm 1)] the operator $S := UTU^{-1}$ defined on the dense subspace $D(S) := U(D(T))$ is also essentially self-adjoint,

\item[\rm 2)] if $\mathbf{T}$ and $\mathbf{S}$ denote the unique self-adjoint extensions of $T$ and $S$ respectively, then $D(\mathbf{S}) = U(D(\mathbf{T}))$ and $\mathbf{S} = U \mathbf{T} U^{-1}$ holds on $D(\mathbf{S})$. \qed
\end{enumerate}

\end{lemma}

\begin{remark}
In fact we have $\mathbf{T}=T^*$ and $\mathbf{S}=S^*$, but let's not bother.
\end{remark}

\begin{lemma}[unitary conjugation commutes with functional calculus]
\label{lem:unitary_conjugation_commutes_with_functional_calculus}
If $T$ is a densely defined self-adjoint operator on a Hilbert space $V$ with its domain $D(T)$, and if $U$ is a unitary operator on $V$, then $U T U^{-1}$ is a densely defined self-adjoint operator on $V$ with domain $D(U T U^{-1}) = U(D(T)) := \{ Uv \, | \, v\in D(T) \}$. For any complex-valued continuous function $f$ on $\mathbb{R}$, the densely defined operators $f(UTU^{-1})$ and $f(T)$, obtained by the functional calculus applied to $UTU^{-1}$ and $T$ respectaively with the same function $f$, are related by
\begin{align*}
f(UTU^{-1}) = Uf(T)U^{-1},
\end{align*}
where $D(f(UTU^{-1})) = U(D(f(T)))$.
\end{lemma}

\textit{Proof.} First, we define the operator $UTU^{-1}$ on $D(UTU^{-1}) := U(D(T))$ as $(UTU^{-1})(Uv) = U(Tv)$ for any element $Uv$ of $U(D(T))$, which makes sense. As $U$ is a homeomorphism from $V$ to itself with respect to the standard Hilbert space topology, $U(D(T))$ is dense in $V$. Pick any two elements $Uv$ and $Uw$ of $U(D(T))$, with $v,w\in D(T)$. Then
$  \langle Uv, (UTU^{-1})(Uw) \rangle_V = \langle Uv, \, U(Tw) \rangle_V = \langle v,Tw \rangle_V = \langle Tv, w \rangle_V = \langle UTv, Uw \rangle_V = \langle (UTU^{-1})(Uv), Uw \rangle_V
$, therefore $UTU^{-1}$ is symmetric; hence $D(UTU^{-1}) \subset D((UTU^{-1})^*)$. Now, suppose $v\in D((UTU^{-1})^*)$, that is, $f\mapsto \langle (UTU^{-1})f, v\rangle_V$ is a bounded linear functional on $D(UTU^{-1}) = U(D(T))$. Consider the map $D(T) \to \cplx$ given by $g\mapsto \langle Tg, U^{-1} v \rangle_V$, which is clearly $\cplx$-linear. As the map $D(T) \to U(D(T))$, $g\mapsto Ug = f$, is a bijective $\cplx$-linear map preserving the norm, one has
$$
\sup_{g\in D(T)} \frac{ |\langle T g, U^{-1} v\rangle_V| }{||g||_V} = \sup_{f \in U(D(T))} \frac{ |\langle TU^{-1}f, U^{-1}v\rangle_V|}{||f||_V} = \sup_{f\in U(D(T))} \frac{ |\langle UTU^{-1} f, v \rangle_V | }{||f||_V},
$$
which we know to be finite from the fact that $v\in D((UTU^{-1})^*)$. Thus $g\mapsto \langle Tg, U^{-1} v\rangle_V$ is a bounded linear functional on $D(T)$, hence $U^{-1} v \in D(T^*)$ by the definition of $D(T^*)$. We thus have $U^{-1} v \in D(T)$, since $D(T)=D(T^*)$, for $T$ is self-adjoint. Therefore $v\in U(D(T)) = D(UTU^{-1})$, hence $D((UTU^{-1})^*) = D(UTU^{-1})$, showing that $UTU^{-1}$ is self-adjoint.

\vs

Write $T = \int_{-\infty}^\infty \lambda \, dE(\lambda)$, the unique spectral resolution of $T$. So $E(\lambda)$ is a resolution of the identity and we have
$$
D(T) = \left\{ v \in V \, \left| \, \int_{-\infty}^\infty \right. |\lambda|^2 \, d||E(\lambda) \, v||^2_V < \infty \right\}.
$$
It is easy to check that $F(\lambda) := U E(\lambda) U^{-1}$ is also a resolution of the identity. Thus we may form a self-adjoint operator $T' := \int_{-\infty}^\infty \lambda \, dF(\lambda)$, on its domain 
$$
D(T') = \left\{ v \in V \, \left| \, \int_{-\infty}^\infty \right. |\lambda|^2 \, d||F(\lambda) \, v||^2_V < \infty \right\}.
$$
As $||F(\lambda) (Uv) ||^2_V = ||UE(\lambda) v||^2_V = ||E(\lambda) v||^2_V$ for any $v\in V$, one can easily check that $D(T') = U(D(T))$. Following the construction in the present subsection, using the fact that $U$ is a norm-preserving homeomorphism of $V$ to itself with respect to the standard Hilbert space topology, one can show that $T' = UTU^{-1}$. Likewise, one can show that the densely defined operator $f(T') = \int_{-\infty}^\infty f(\lambda) \, dF(\lambda)$ obtained from the functional calculus of $T'$ for the function $f$ exactly coincides with $U f(T) U^{-1}$, where $f(T) = \int_{-\infty}^\infty f(\lambda) \, dE(\lambda)$ is obtained from the functional calculus of $T$ for $f$, and also that $D(f(T')) = \{ v\in V | \int_{-\infty}^\infty |f(\lambda)|^2 ||dF(\lambda)v||^2_V<\infty \}$ coincides $U(D(f(T)))$, where $D(f(T)) = \{v \in V | \int_{-\infty}^\infty |f(\lambda)|^2 ||dE(\lambda)||^2_V<\infty\}$. \qed

\begin{lemma}
\label{lem:unitary_conjugation_corollary}
Let $T,T'$ be densely defined self-adjoint operators on a Hilbert space $V$ and $U$ be a unitary operator on $V$. Suppose $U e^{\sqrt{-1} t T} U^{-1} = e^{\sqrt{-1} t T'}$ for all $t\in \mathbb{R}$. Then $T' = UTU^{-1}$. Furthermore, if $D$ is a dense subspace of $V$ invariant under $U, U^{-1}, T, T', e^{\sqrt{-1}tT}, e^{\sqrt{-1} tT'}$, $\forall t\in \mathbb{R}$, then $(T'\restriction D) = U(T\restriction D)U^{-1}$ on $D$ \quad (where the symbol $\restriction$ means the restriction).
\end{lemma}

\textit{Proof.} Let $T$ be densely defined self-adjoint and $U$ be unitaryon the Hilbert space $V$. By Lem.\ref{lem:unitary_conjugation_commutes_with_functional_calculus}, $T'' := U T U^{-1}$ is a densely defined self-adjoint operator on $V$, and $f_t(T'') = U f_t(T) U^{-1}$ holds on the whole $V$, where the unitary operators $f_t(T'')$ and $f_t(T)$ on $V$ are obtained by the functional calculus applied respectively to $T''$ and $T$ for the function $f_t(\lambda) = e^{\sqrt{-1} t \lambda}$. By assumption we have a densely defined self-adjoint operator $T'$ on $V$ satisfying $Ue^{\sqrt{-1} t T} U^{-1} = e^{\sqrt{-1} t T'}$ for all $t\in \mathbb{R}$. As $U e^{\sqrt{-1} t T} U^{-1} = U f_t(T) U^{-1} = f_t(T'') = e^{\sqrt{-1}tT''}$, we have $e^{\sqrt{-1}tT'} = e^{\sqrt{-1}t T''}$. By the uniqueness of an infinitesimal generator of the strongly continuous unitary group $\{e^{\sqrt{-1}tT'}\}_{t\in\mathbb{R}}$ (Theorems \ref{thm:unitary_family_for_self-adjoint_operator}, \ref{thm:Stones_theorem}), we have $T'=T''$, that is, $T' = UTU^{-1}$.

\vs

So $T'v = UTU^{-1}v$ for all $v\in D(T') = D(UTU^{-1}) = U(D(T))$ (Lem.\ref{lem:unitary_conjugation_commutes_with_functional_calculus}). Let $v$ be any element of the space $D$ in the statement of the Lemma. By assumption, $D\subseteq D(T')$, so $v\in D(T')$, so we have $T'v = UTU^{-1}v$, and $T'v = (T'\restriction D)v$. By assumption $U^{-1}v \in D \subseteq D(T)$, so $TU^{-1}v = (T\restriction D)(U^{-1} v)$. Hence indeed $(T'\restriction D)v = U(T\restriction D)U^{-1} v$. \qed

\subsection{The Stone-von Neumann Theorem}

In the present subsection I introduce basic ingredients of the construction of representations, namely the operators forming the `Schr\"odinger representation' on the Hilbert space $L^2(\mathbb{R},dx)$. We first define a special $\mathbb{C}$-vector subspace $D$ of $L^2(\mathbb{R},dx)$.
\begin{align}
\label{eq:D}
D := \mathrm{span}_\cplx \{e^{-\alpha x^2 - \beta x} P(x) \, | \, \alpha, \, \beta \in \cplx, ~ \mathrm{Re}(\alpha)>0, ~ \mbox{$P$ is a polynomial in $x$} \},
\end{align}
which is dense in $L^2(\mathbb{R}, dx)$ with respect to the standard Hilbert space topology; one way of seeing that this subspace is dense is to observe that the `Hermite functions', which form a topological countable basis of $L^2(\mathbb{R},dx)$, are in $D$. One nice property of $D$ is:
\begin{lemma}
\label{lem:D_is_preserved_by_Fourier_transforms}
Let $\mathcal{F} : L^2(\mathbb{R},dx) \to L^2(\mathbb{R},dx)$ be the Fourier transform; more precisely, on $L^1(\mathbb{R}, dx) \cap L^2(\mathbb{R}, dx)$ it is given by the formula
\begin{align}
\label{eq:mcal_F}
(\mcal{F} f)(x) = \int_\mathbb{R} e^{-2\pi \sqrt{-1} xy} f(y) \, dy.
\end{align}
Then, $\mcal{F}(D) = D$.
\end{lemma}
{\it Proof.} A key fact to use is that $\mcal{F}$ sends $e^{-\alpha x^2}$ to $\frac{\sqrt{\pi}}{\sqrt{\alpha}} \, e^{-\pi^2 x^2/\alpha}$, for any complex $\alpha$ with $\mathrm{Re}(\alpha)>0$, where $\sqrt{\alpha}$ is taken to be such that $\mathrm{Re}(\sqrt{\alpha})>0$; one can prove this using $\int_\mathbb{R} e^{-x^2} dx = \sqrt{\pi}$ and a basic contour integral technique. I omit the details. \qed

\vs

We shall use the well-known fact that $\mcal{F}$ is unitary, and that its inverse is given by the formula $(\mcal{F}^{-1}f)(x) = \int_\mathbb{R} e^{2\pi\sqrt{-1} xy} f(y)dy$ for $f\in L^1\cap L^2$.

\begin{definition}[the Schr\"odinger representation]
\label{def:Schrodinger_representation}
Let $\hbar \in \mathbb{R}_{>0}\setminus\mathbb{Q}$. The operators $\wh{\mathbf{q}}^\hbar$ and $\wh{\mathbf{p}}^\hbar$ are defined as the densely defined operators on $L^2(\mathbb{R},dx)$ on the domain $D$ \eqref{eq:D} given by \footnote{the superscripts $\hbar$ in $\wh{\mathbf{q}}_i^\hbar$ and $\wh{\mathbf{p}}_i^\hbar$ do not mean powers.}
\begin{align}
\label{eq:Schrodinger_representation}
  \wh{\mathbf{q}}^\hbar := x \qquad\mbox{and}\qquad \wh{\mathbf{p}}^\hbar := 2\pi\sqrt{-1} \, \hbar \, \frac{d}{d x} \qquad \mbox{on}\quad D,
\end{align}
i.e.
$$
(\wh{\mathbf{q}}^\hbar f)(x) = x \, f(x), \qquad
(\wh{\mathbf{p}}^\hbar f)(x) = 2\pi\sqrt{-1} \, \hbar \, \frac{d f}{d x} (x) \qquad \forall f(x) \in D.
$$
\end{definition}
It is straightforward to check:
\begin{lemma}
\label{lem:Schrodinger_operators_lemma1}
Each of $\wh{\mathbf{q}}^\hbar$ and $\wh{\mathbf{p}}^\hbar$ leaves $D$ invariant, and is symmetric in the sense of Def.\ref{def:symmetric_operator}. \qed
\end{lemma}

\begin{lemma}
\label{lem:Schrodinger_operators_lemma2}
Each of $\wh{\mathbf{q}}^\hbar$ and $\wh{\mathbf{p}}^\hbar$ is essentially self-adjoint on $D$, in the Hilbert space $L^2(\mathbb{R},dx)$. On $D$ they satisfy the `Heisenberg relation'
\begin{align}
\label{eq:Heisenberg_relation_for_p_and_q}
[\wh{\mathbf{p}}^\hbar, \wh{\mathbf{q}}^\hbar] = 2\pi \sqrt{-1} \, \hbar \, \cdot \mathrm{id},
\end{align}
i.e. $\wh{\mathbf{p}}^\hbar (\wh{\mathbf{q}}^\hbar f) - \wh{\mathbf{q}}^\hbar (\wh{\mathbf{p}}^\hbar f) = 2\pi\sqrt{-1}\hbar \cdot f$ for all $f\in D$. Moreover, if we let $P$ and $Q$ be the unique self-adjoint extensions of $\wh{\mathbf{p}}^\hbar$ and $\wh{\mathbf{q}}^\hbar$ respectively, the corresponding strongly continuous unitary groups $\{ e^{\sqrt{-1} \alpha P}\}_{\alpha\in\mathbb{R}}$ and $\{ e^{\sqrt{-1} \beta Q} \}_{\beta\in\mathbb{R}}$, defined by the functional calculus in \S\ref{subsec:spectral_theorem} for $P$ and $Q$ respectively with the functions $\lambda \mapsto e^{\sqrt{-1}\alpha \lambda}$ and $\lambda \mapsto e^{\sqrt{-1} \beta \lambda}$ which are unitary functions in the sense of Lem.\ref{lem:unitary_operators_from_functional_calculus}, satisfy the `Weyl relations'
\begin{align}
\label{eq:Weyl_relations_for_P_and_Q}
  e^{\sqrt{-1} \alpha P} \, e^{\sqrt{-1} \beta Q} = e^{- 2\pi \sqrt{-1} \, \hbar \, \alpha \beta} \, e^{\sqrt{-1} \beta Q} \, e^{\sqrt{-1} \alpha P}.
\end{align}
\end{lemma}
{\it Proof.} The Heisenberg relation is easy to check directly. For each real numbers $\alpha,\beta\in\mathbb{R}$ define operators $U^{(\alpha)}$ and $V^{(\beta)}$ on $L^2(\mathbb{R},dx)$ to itself as
$$
(U^{(\alpha)}f)(x) := f(x - 2\pi \hbar \, \alpha), \qquad (V^{(\beta)} f)(x) := e^{\sqrt{-1} \beta x} \, f(x), \qquad \forall f(x) \in L^2(\mathbb{R},dx).
$$
One can indeed check that these are well-defined operators on the whole $L^2(\mathbb{R},dx)$, that they are unitary, and that they are one-parameter groups, i.e. $U^{(\alpha_1)} U^{(\alpha_2)} = U^{(\alpha_1+\alpha_2)}$ and $V^{(\beta_1)} V^{(\beta_2)} = V^{(\beta_1+\beta_2)}$. One can moreover check that $\{U^{(\alpha)}\}_{\alpha\in\mathbb{R}}$ and $\{V^{(\beta)}\}_{\beta\in\mathbb{R}}$ are strongly continuous one-parameter unitary groups in the sense of Def.\ref{def:strongly_continuous_one-parameter_unitary_group}. We now use Thm.\ref{thm:strong_derivative} to express their infinitesimal generators. First, one can easily observe that each of $U^{(\alpha)}$ and $V^{(\beta)}$ preserves $D$. In fact, we shall find later that a more general version of what we want to prove for $U^{(\alpha)}$ must be proved; in order to save time and paper I refer to this `future' proof. Namely, as a special case $n=1$ of Lem.\ref{lem:shift_operator_as_exponential}, we see that the operator $\wh{\mathbf{p}}^\hbar$ on $D$ is essentially self-adjoint and its unique self-adjoint extension $P$ coincides with the infinitesimal generator of $\{U^{(\alpha)}\}_{\alpha\in\mathbb{R}}$; so we can write $U^{(\alpha)} = e^{\sqrt{-1} \alpha P}$. One can directly prove that $\{V^{(\beta)}\}_{\beta\in\mathbb{R}}$ is a strongly continuous one-parameter unitary group, but we take a shortcut, whose method will become useful later. 

\vs

Recall the Fourier transform $\mcal{F}$ \eqref{eq:mcal_F}; we will see how it `exchanges' $U^{(\alpha)}$ and $V^{(\beta)}$, and also $\wh{\mathbf{p}}^\hbar$ and $\wh{\mathbf{q}}^\hbar$. For $f\in D$, one has $U^{(\alpha)} f \in D \subset L^1(\mathbb{R}) \cap L^2(\mathbb{R})$, hence
\begin{align*}
(\mcal{F}^{\pm 1} (U^{(\alpha)}f))(x) & = \int_\mathbb{R} e^{\mp 2\pi\sqrt{-1} xy} \, (U^{(\alpha)}f)(y) \, dy  = \int_\mathbb{R} e^{\mp 2\pi\sqrt{-1} xy} \, f(y-2\pi\hbar\, \alpha) dy \\
& = \int_\mathbb{R} e^{\mp 2\pi\sqrt{-1} x (Y + 2\pi\hbar\, \alpha)} f(Y) dY \qquad (\because Y = y - 2\pi\hbar \, \alpha) \\
& = e^{\sqrt{-1} ( \mp (2\pi)^2 \hbar \, \alpha) x} \cdot (\mcal{F}^{\pm 1} f)(x)
= (V^{(\mp (2\pi)^2 \hbar \, \alpha)} (\mcal{F}^{\pm 1} f))(x),
\end{align*}
where all the symbols $\pm$ and $\mp$ are to be coherent. As $\mcal{F}$ is unitary and $\mcal{F}(D)=D$, one thus has
$$
\mcal{F} \, U^{(\alpha)} \, \mcal{F}^{-1} = V^{(-(2\pi)^2\hbar \, \alpha)} \qquad\mbox{and}\qquad
\mcal{F} \, V^{((2\pi)^2 \hbar \, \alpha)} \, \mcal{F}^{-1} = U^{(\alpha)},
$$
on $D$, and therefore on the whole $L^2(\mathbb{R})$. For $f\in D$, we have $\wh{\mathbf{p}}^\hbar f \in D \subset L^1(\mathbb{R})\cap L^2(\mathbb{R})$, so
\begin{align*}
(\mcal{F}^{\pm 1} \, (\wh{\mathbf{p}}^\hbar f))(x) & = 2\pi\sqrt{-1} \hbar \, \int_\mathbb{R} e^{\mp 2\pi\sqrt{-1} xy} \cdot \frac{df(y)}{dy} \, dy \\
& = - 2\pi \sqrt{-1} \hbar \, \int_\mathbb{R} \left( \frac{d}{dy} e^{\mp 2\pi\sqrt{-1} xy} \right) \cdot f(y) \, dy \quad (\because \mbox{integration by parts}) \\
& = \mp (2\pi)^2 \hbar \, x \cdot \int_\mathbb{R} e^{\mp 2\pi\sqrt{-1} xy} \, f(y) \, dy = (\mp (2\pi)^2 \hbar \, \wh{\mathbf{q}}^\hbar (\mcal{F}^{\pm 1} f))(x).
\end{align*}
Thus, one deduces the following equalities of operators $D\to D$:
\begin{align}
\label{eq:mcal_F_conjugation_on_p_and_q}
\mcal{F} \, \wh{\mathbf{p}}^\hbar \, \mcal{F}^{-1} = - (2\pi)^2 \, \hbar \, \wh{\mathbf{q}}^\hbar \qquad\mbox{and}\qquad
\mcal{F} \, ( (2\pi)^2 \hbar \, \wh{\mathbf{q}}^\hbar ) \, \mcal{F}^{-1} = \wh{\mathbf{p}}^\hbar.
\end{align}
In particular, $\wh{\mathbf{p}}^\hbar$ on $D$ being essentially self-adjoint, together with $\mcal{F}$ being unitary and $\mcal{F}(D)=D$, implies that $\wh{\mathbf{q}}^\hbar$ on $D$ is also essentially self-adjoint, thanks to Lem.\ref{lem:unitarily_equivalent_essentially_self-adjoint_operators}.

\vs

Lem.\ref{lem:unitary_conjugation_commutes_with_functional_calculus} tells us that the densely defined operator $Q := \mcal{F} (-\frac{1}{(2\pi)^2\hbar }\, P) \mcal{F}^{-1}$ on $D(Q) := \mcal{F}(D(-\frac{1}{(2\pi)^2\hbar} \, P)) = \mcal{F}(D(P))$ is self-adjoint, and that $\mcal{F} \, e^{\sqrt{-1}\beta(-\frac{1}{(2\pi)^2\hbar}) P} \, \mcal{F}^{-1} = e^{\sqrt{-1}\beta Q}$ for each $\beta\in \mathbb{R}$. So, from $U^{(\alpha)} = e^{\sqrt{-1}\alpha P}$ and $\mcal{F} \, U^{(\alpha)} \, \mcal{F}^{-1} = V^{(-(2\pi)^2\, \hbar\, \alpha)}$ one deduces $V^{(\beta)} = e^{\sqrt{-1}\beta Q}$. As the restrictions of $-\frac{1}{(2\pi)^2 \hbar} \, P$ and $Q$ on $D$ are $-\frac{1}{(2\pi)^2\hbar} \, \wh{\mathbf{p}}^\hbar$ and $\wh{\mathbf{q}}^\hbar$ respectively, 
Lem.\ref{lem:unitary_conjugation_corollary} applied to this situation with the space $D$ tells us that $\mcal{F} ( - \frac{1}{(2\pi)^2 \hbar} \, P \restriction D) \, \mcal{F}^{-1} = Q\restriction D$; since the restriction of $- \frac{1}{(2\pi)^2 \hbar} \, P$ is $- \frac{1}{(2\pi)^2 \hbar} \, \wh{\mathbf{p}}^\hbar$ and we know $\mcal{F} \, ( - \frac{1}{(2\pi)^2 \hbar} \, \wh{\mathbf{p}}^\hbar ) \, \mcal{F}^{-1} = \wh{\mathbf{q}}^\hbar$, we find $Q\restriction D = \wh{\mathbf{q}}^\hbar$, hence $Q$ is the unique self-adjoint extension of $\wh{\mathbf{q}}^\hbar$.

\vs

Lastly, the Weyl relations are easily checked as follows:
\begin{align*}
( U^{(\alpha)} (V^{(\beta)} f) )(x) & = (V^{(\beta)} f)(x-2\pi\hbar \, \alpha) = e^{\sqrt{-1} \beta(x-2\pi\hbar\, \alpha)} \, f(x-2\pi\hbar\, \alpha) \\
& = e^{-2\pi \sqrt{-1} \hbar\, \alpha \beta} e^{\sqrt{-1}\beta x} \cdot (U^{(\alpha)} f)(x)
=  e^{-2\pi \sqrt{-1} \hbar\, \alpha \beta} \, (V^{(\beta)} (U^{(\alpha)} f))(x),
\end{align*}
for each $f\in L^2(\mathbb{R})$ and each real numbers $\alpha,\beta$.\qed

\vs

What is remarkable is that the pair $(P,Q)$ of self-adjoint operators given by the formula \eqref{eq:Schrodinger_representation} on a dense subspace is the unique `irreducible' pair satisfying the Weyl relations \eqref{eq:Weyl_relations_for_P_and_Q}. However, note that the Heisenberg relation \eqref{eq:Heisenberg_relation_for_p_and_q} on a Hilbert space may not imply the Weyl relations \eqref{eq:Weyl_relations_for_P_and_Q}; see `Nelson's example' and its perturbation in \cite[Chap.VIII.5]{RS70}. So the Weyl relations are what make the example \eqref{eq:Schrodinger_representation} \emph{the standard representation}, not the Heisenberg relation.

\begin{theorem}
[Stone von-Neumann Theorem \cite{vN31}; {\cite[Thm.VIII.14]{RS70}}]
\label{thm:SvN}
Let $\mathscr{H}$ be a separable Hilbert space. Let $\{U^{(\alpha)}\}_{\alpha\in \mathbb{R}}$ and $\{V^{(\beta)}\}_{\beta\in\mathbb{R}}$ be strongly continuous one-parameter unitary groups on $\mathscr{H}$, satisfying the Weyl relations
$$
U^{(\alpha)} \, V^{(\beta)} = e^{-\sqrt{-1} \alpha \beta} \, V^{(\beta)} U^{(\alpha)}, \qquad \forall \alpha,\beta\in\mathbb{R}.
$$
Then there are closed subspaces $\mathscr{H}_\ell$, $\ell=1,\ldots,N$, where $N$ is a positive integer or $\infty$, such that
\begin{itemize}
\item[\rm 1)] $\mathscr{H} = \bigoplus_{\ell=1}^N \, \mathscr{H}_\ell$,

\item[\rm 2)] Each of $U^{(\alpha)}$ and $V^{(\beta)}$ preserves each $\mathscr{H}_\ell$,

\item[\rm 3)] For each $\ell$, there exists a unitary operator $T_\ell : \mathscr{H}_\ell \to L^2(\mathbb{R},dx_\ell)$ such that
$$
( (T_\ell \,  U^{(\alpha)} \, T_\ell^{-1}) \, f)(x_\ell) = f(x_\ell \, -\, \alpha), \qquad
( (T_\ell \,  V^{(\beta)} \, T_\ell^{-1}) \, f)(x_\ell) = e^{\sqrt{-1}\beta x_\ell} \cdot f(x_\ell),
$$
for all $f(x_\ell) \in L^2(\mathbb{R}, dx_\ell)$.
\end{itemize}
\end{theorem}

\section{Representations of algebraic quantum cluster varieties}
\label{sec:representations_of_algebraic_quantum_cluster_varieties}

\subsection{Positive integrable representations of quantum torus algebras}
\label{subsec:positive_representations}

For a representation of a seed quantum $\mcal{D}$-torus algebra ${\bf D}^q_\Gamma$, we would like to look for operators $\wh{\bf B}_i = \pi_\Gamma({\bf B}_i)$, $\wh{\bf X}_i = \pi_\Gamma({\bf X}_i)$, $i=1,\ldots,n$, on a Hilbert space $\mathscr{H}_\Gamma$, that satisfy the algebraic relations of ${\bf B}_i$ and ${\bf X}_i$. Since ${\bf B}_i$, ${\bf X}_i$, and their corresponding operators, are supposed to be quantum counterparts of the positive real valued functions $B_i$, $X_i$ on the space $\mcal{D}_\mathscr{C}(\mathbb{R}_{>0})$, i.e. the positive real points of our cluster $\mcal{D}$-variety, it is natural to require that the operators be (densely defined) self-adjoint operators that are `positive-definite', by which I mean:
\begin{definition}
A densely defined self-adjoint operator $T$ on a Hilbert space $V$ with its domain $D(T)$ is said to be \ul{\em positive-definite} if $\langle Tv, v\rangle_V >0$ for all nonzero $v\in D(T)$.
\end{definition}
\begin{remark}
In the theory of bounded operators, the condition $\langle Tv, v\rangle_V\ge 0$, $\forall v \in V$, implies self-adjointness.
\end{remark}
As in usual quantization stories, $\mathscr{H}_\Gamma$ is the space of all square-integrable functions on some Euclidean (measure) space $\mathbb{R}^n$, and each of these positive-definite self-adjoint operators are defined only on a dense subspace of $\mathscr{H}_\Gamma$. So it is a tricky matter to deal with the commutation relations of such operators, for we have to be careful about the domains and ranges of these operators. I suggest the following definition as a way to define and characterize certain `nicely-behaved' class of representations of the quantum $\mcal{D}$-torus algebra $\mathbf{D}^q_\Gamma$, using the `Weyl relations', as employed already by von Neumann.

\begin{definition}
\label{def:positive_integrable_representation}
Let $\hbar \in \mathbb{R}_{>0}\setminus\mathbb{Q}$, and let
$$
q = e^{\pi \sqrt{-1} \hbar}. 
$$
A \ul{\em positive integrable $*$-representation} of the seed quantum $\mcal{D}$-torus algebra ${\bf D}^q_\Gamma$ is a triple $(\mathscr{H}_\Gamma, D_\Gamma, \pi^q_\Gamma)$, where $\mathscr{H}_\Gamma$ is a Hilbert space, $D_\Gamma$ is a dense $\cplx$-linear subspace of $\mathscr{H}_\Gamma$, and $\pi^q_\Gamma$ consists of essentially self-adjoint operators $\wh{\mathbf{b}}^\hbar_i$, $\wh{\mathbf{x}}^\hbar_i$, $i=1,\ldots,n$, on $D_\Gamma$, such that\footnote{the superscripts $\hbar$ in $\wh{\mathbf{b}}_i^\hbar$ and $\wh{\mathbf{x}}_i^\hbar$ do not mean powers, but just indicate the associated parameter $\hbar$} 
\begin{itemize}
\item[\rm 1)] for any
\begin{align}
\label{eq:alpha_and_beta}
(\alpha) = (\alpha_1,\ldots,\alpha_n) \in \mathbb{R}^n \quad\mbox{and}\quad
(\beta) = (\beta_1,\ldots,\beta_n) \in \mathbb{R}^n
\end{align}
the following unitary operators on $\mathscr{H}_\Gamma$
$$
\wh{\mathbf{B}}^{(\alpha)}_i := e^{\sqrt{-1} \alpha_i \wh{\mathbf{b}}^\hbar_i}, \qquad \wh{\mathbf{X}}^{(\beta)}_i := e^{\sqrt{-1} \beta_i \wh{\mathbf{x}}^\hbar_i}, \qquad i=1,\ldots,n,
$$
defined by the functional calculus in \S\ref{subsec:spectral_theorem} for the unique self-adjoint extensions of the operators $\wh{\mathbf{b}}^\hbar_i$'s and $\wh{\mathbf{x}}^\hbar_i$'s with the (unitary) functions $\lambda \mapsto e^{\sqrt{-1} \alpha_i \lambda}$ and $\lambda \mapsto e^{\sqrt{-1} \beta_i \lambda}$ respectively, leave $D_\Gamma$ invariant and satisfy the relations
\begin{align}
\label{eq:Weyl_relations_for_D_q_Gamma}
\left\{
{\renewcommand{\arraystretch}{1.3} \begin{array}{rcll}
e^{\pi \sqrt{-1} \hbar_j \varepsilon_{ij} \beta_i \beta_j} \, \wh{\bf X}_i^{(\beta)} \wh{\bf X}_j^{(\beta)} & = & e^{\pi \sqrt{-1} \hbar_i \varepsilon_{ji} \beta_j \beta_i} \, \wh{\bf X}_j^{(\beta)} \wh{\bf X}_i^{(\beta)}, & \quad \forall i,j\in \{1,\ldots,n\}, \\
e^{\pi\sqrt{-1} \hbar_i \beta_i \alpha_i} \, \wh{\bf X}_i^{(\beta)} \wh{\bf B}_i^{(\alpha)} & = & e^{-\pi\sqrt{-1} \hbar_i \alpha_i \beta_i} \, \wh{\bf B}_i^{(\alpha)} \wh{\bf X}_i^{(\beta)}, & \quad \forall i \in \{1,\ldots,n\},  \\
\wh{\bf B}_i^{(\alpha)} \wh{\bf X}_j^{(\beta)} & = & \wh{\bf X}_j^{(\beta)} \wh{\bf B}_i^{(\alpha)}, & \quad \mbox{whenever $i\neq j$}, \\
\wh{\bf B}_i^{(\alpha)} \wh{\bf B}_j^{(\alpha)} & = & \wh{\bf B}_j^{(\alpha)} \wh{\bf B}_i^{(\alpha)}, & \quad \forall i,j \in \{1,\ldots,n\},
\end{array} }
\right.
\end{align}
where
\begin{align}
\label{eq:hbar_i}
\hbar_i := \hbar/d_i, \qquad \forall i=1,\ldots,n, \qquad \mbox{and}
\end{align}

\item[\rm 2)] the following densely defined self-adjoint operators on $\mathscr{H}_\Gamma$
\begin{align}
\label{eq:pi_Gamma_of_generators}
\pi^q_\Gamma(\mathbf{B}_i^{\pm 1}) := e^{\pm \wh{\mathbf{b}}^\hbar_i}, \qquad \pi^q_\Gamma(\mathbf{X}_i^{\pm 1}) := e^{\pm \wh{\mathbf{x}}^\hbar_i}, \qquad i=1,\ldots,n,
\end{align}
defined by the functional calculus in \S\ref{subsec:spectral_theorem} for the unique self-adjoint extensions of the operators $\wh{\mathbf{b}}^\hbar_i$'s and $\wh{\mathbf{x}}^\hbar_i$'s with the functions $e^{\pm \lambda}$ respectively, leave $D_\Gamma$ invariant.
\end{itemize}
\end{definition}
Few words on the terminology `positive integrable $*$-representation'. From \eqref{eq:T_E_f2} and the fact that $e^{\pm \lambda}>0$, one can easily see that the densely defined operators in \eqref{eq:pi_Gamma_of_generators} which represent the generators of $\mathbf{D}^q_\Gamma$ which are $*$-invariant are positive-definite self-adjoint, hence justifying the words `positive' and `$*$-'. The word `integrable' hints to the `Weyl relations' \eqref{eq:Weyl_relations_for_D_q_Gamma}. In this regard, we note that the `Weyl relations' \eqref{eq:Weyl_relations_for_D_q_Gamma} imply the `Heisenberg relations'
\begin{align}
\label{eq:Heisenberg_relations_for_D_q_Gamma}
[\wh{\mathbf{x}}^\hbar_i,\wh{\mathbf{x}}^\hbar_j] = 2\pi \sqrt{-1} \, \hbar_j \, \varepsilon_{ij} \cdot {\rm id}, \quad
[\wh{\mathbf{x}}^\hbar_i, \wh{\mathbf{b}}^\hbar_j]=2\pi \sqrt{-1} \, \hbar_i \, \delta_{i,j} \cdot {\rm id}, \quad
[\wh{\mathbf{b}}^\hbar_i,\wh{\mathbf{b}}^\hbar_j]=0, \quad \forall i,j,
\end{align}
on $D_\Gamma$ (see e.g. Corollary of Theorem VIII.14 of \cite{RS70}), where $\delta_{i,j}$ is the Kronecker delta, but the `Heisenberg relations' \eqref{eq:Heisenberg_relations_for_D_q_Gamma} may not imply the `Weyl relations' \eqref{eq:Weyl_relations_for_D_q_Gamma} in a general Hilbert space; see `Nelson's example' and its perturbation in \cite[Chap.VIII.5]{RS70}. Finally, one can prove that the operators \eqref{eq:pi_Gamma_of_generators} satisfy the algebraic relations satisfied by the corresponding generators of $\mathbf{D}^q_\Gamma$, when applied to elements in $D_\Gamma$. One way of showing this is to use the `Stone-von Neumann Theorem' (see e.g. Thm.VIII.14 of \cite{RS70}) for each $i=1,\ldots, n$ (one at a time), to write down $\wh{\mathbf{b}}_i$ and $\wh{\mathbf{x}}_i$ as explicit operators on an explicit Hilbert space, and check directly; for instance, one can check for the concrete examples of positive integrable $*$-representations that I will soon present. So, $\pi_\Gamma$ can indeed be thought of as a `representation' of the algebra $\mathbf{D}^q_\Gamma$.

\vs

The `tilde' operators are defined by the formula
\begin{align}
\label{eq:tilde_operators}
\wh{\til{\mathbf{x}}}^\hbar_i := \wh{\mathbf{x}}^\hbar_i + \sum_{j=1}^n \varepsilon_{ij} \, \wh{\mathbf{b}}^\hbar_j, \qquad \pi^q_\Gamma(\til{\mathbf{X}}^{\pm 1}_i) := e^{\pm \wh{\til{\mathbf{x}}}^\hbar_i}, \qquad \forall i=1,\ldots,n,
\end{align}
as essentially self-adjoint operators on $D_\Gamma$.

\vs

Def.\ref{def:positive_integrable_representation} is designed in such a way to deal with the analytical issues for the algebraic relations that must be satisfied by the positive-definite self-adjoint operators corresponding to the generators of ${\bf D}^q_\Gamma$. At the moment, let us not bother about the operators corresponding to other elements of the skew field of frations $\mathbb{D}^q_\Gamma$ of $\mathbf{D}^q_\Gamma$; we will later discuss how Fock and Goncharov \cite{FG09} dealt with the analytical issues for the operators for sufficiently many elements of $\mathbb{D}^q_\Gamma$ simultaneously.

\begin{definition}
\label{def:weakly_irreducible}
A positive integrable $*$-representation $(\mathscr{H}_\Gamma, D_\Gamma, \pi^q_\Gamma)$ of the seed quantum $\mcal{D}$-torus algebra $\mathbf{D}^q_\Gamma$ is said to be \ul{\em weakly irreducible} if $D_\Gamma$ has no nonzero $\cplx$-vector subspace $D$ that is invariant under all $\wh{\mathbf{B}}^{(\alpha)}_i$, $\wh{\mathbf{X}}^{(\beta)}_i$, $\pi^q_\Gamma(\mathbf{B}_i^{\pm 1})$, $\pi^q_\Gamma(\mathbf{X}_i^{\pm 1})$, $i=1,\ldots,n$, such that the closure of $D$ in $\mathscr{H}_\Gamma$ is a proper Hilbert subspace of $\mathscr{H}_\Gamma$.
\end{definition}

\begin{question}
Is a weakly irreducible positive integrable $*$-representation of a seed quantum $\mathcal{D}$-torus algebra $\mathbf{D}^q_\Gamma$ unique in some sense, up to unitary equivalence (e.g. via the Stone-von Neumann theorem)?
\end{question}

Let us now present two examples of weakly irreducible positive integrable $*$-representations of $\mathbf{D}^q_\Gamma$, due to Fock-Goncharov \cite{FG09}. 

\begin{definition}
\label{def:H_Gamma_and_D_Gamma}
For a $\mcal{D}$-seed $\Gamma$, define the Hilbert space
\begin{align}
\label{eq:H_Gamma}
\mathscr{H}_\Gamma := L^2(\mathbb{R}^n, \, da_1\, da_2 \, \ldots \, da_n),
\end{align}
where the measure on $\mathbb{R}^n$ is the product of Lebesgue measures $da_i$ on $\mathbb{R}$, and the inner product $\langle \,\, , \, \rangle_{\mathscr{H}_\Gamma}$ is the usual one
$$
\langle f, g \rangle_{\mathscr{H}_\Gamma} = \int_{\mathbb{R}^n} f \, \ol{g} \, da_1 \ldots da_n.
$$
Define $D_\Gamma$ as the following $\cplx$-vector subspace of $\mathscr{H}_\Gamma$:
\begin{align}
\label{eq:D_Gamma}
\hspace{-3mm}
D_\Gamma := {\rm span}_\cplx\left\{ ~ e^{ - \sum_{i=1}^n (\alpha_i a_i^2 + \beta_i a_i) } \, P(a_1,\ldots,a_n) ~ \, \left| \, ~ \begin{array}{l} \alpha_i,  \beta_i \in \cplx, ~ \mathrm{Re}(\alpha_i)>0, ~\mbox{and $P$ is} \\ \mbox{a polynomial in $a_1,\ldots,a_n$ over $\mathbb{C}$} \end{array} \right. \right\}.
\end{align}
\end{definition}
\begin{remark}
The definitions of $\mathscr{H}_\Gamma$ and $D_\Gamma$ do not really use the data of a $\mathcal{D}$-seed $\Gamma$, except for the number $n$. The subscript $\Gamma$ is mainly for the sake of keeping track of notations. Readers will find this useful later when we deal with multiple $\Gamma$'s.
\end{remark}
Recall that an element of $\mathscr{H}_\Gamma$ is the equivalence class of a square-integrable measurable function on $\mathbb{R}^n$, where two functions are set to be equivalent if they coincide except at a set of measure zero. Usually we will be working with one representative for each equivalence class.

\begin{lemma}
In Def.\ref{def:H_Gamma_and_D_Gamma}, $D_\Gamma$ is dense in $\mathscr{H}_\Gamma$, with respect to the standard Hilbert space topology.
\end{lemma}
\textit{Proof.} We observe that $\mathscr{H}_\Gamma$ is canonically isomorphic, as Hilbert spaces, to the Hilbert space tensor product of all $L^2(\mathbb{R}, da_i)$, $i=1,\ldots,n$ (see e.g. \cite[Thm.II.10]{RS70}). We can then note that, under this canonical identification, the algebraic tensor product of the subspaces
\begin{align}
\label{eq:D_i}
D_i := \mathrm{span}_\cplx \{e^{-\alpha_i a_i^2 - \beta_i a_i} P_i(a_i) \, | \, \alpha_i, \, \beta_i \in \cplx, ~ \mathrm{Re}(\alpha_i)>0, ~ \mbox{$P_i$ is a polynomial in $a_i$} \}
\end{align}
of $L^2(\mathbb{R},da_i)$ for $i=1,\ldots,n$, coincides with $D_\Gamma$. It is easy to show that the algebraic tensor product of dense subspaces of (finitely many) Hilbert spaces is dense in the Hilbert space tensor product of the Hilbert spaces. \qed

\vs

This $D_\Gamma$ contains Fock-Goncharov's dense subspace in \cite{FG09}, which is denoted by $W_\mathbf{i}$ there; they require $\alpha_1 = \alpha_2 = \cdots = \alpha_n>0$. One nice property of $D_\Gamma$ (or $W_{\bf i}$) is:

\begin{lemma}
\label{lem:D_is_preserved_by_Fourier_transforms}
Let $\mathcal{F}_i : \mathscr{H}_\Gamma \to \mathscr{H}_\Gamma$ be the Fourier transform in the $a_i$ variable; more precisely, on $L^1(\mathbb{R}^n, da_1\, \cdots \, da_n) \cap \mathscr{H}_\Gamma$ it is given by the formula
\begin{align}
\label{eq:mcal_F_i}
(\mathcal{F}_i f)(a_1,\ldots,a_i,\ldots,a_n) = \int_\mathbb{R} e^{-2\pi \sqrt{-1} a_i b_i} f(a_1,\ldots,b_i,\ldots,a_n) \, db_i.
\end{align}
Then, $\mathcal{F}_i(D_\Gamma) = D_\Gamma$. \qed
\end{lemma}

\begin{definition}[the Schr\"odinger representation]
\label{def:Schrodinger_representation}
Let $\hbar \in \mathbb{R}_{>0}\setminus\mathbb{Q}$ and let $D_\Gamma$, $\mathscr{H}_\Gamma$ be as in Def.\ref{def:H_Gamma_and_D_Gamma}. Define the operators $\wh{\mathbf{q}}^\hbar_i$, $\wh{\mathbf{p}}^\hbar_i$ : $D_\Gamma \to \mathscr{H}_\Gamma$, $i=1,\ldots,n$, by the formulas\footnote{the superscripts $\hbar$ in $\wh{\mathbf{q}}_i^\hbar$ and $\wh{\mathbf{p}}_i^\hbar$ do not mean powers.}
\begin{align}
\label{eq:Schrodinger_representation}
  \wh{\mathbf{q}}^\hbar_i := a_i \qquad\mbox{and}\qquad \wh{\mathbf{p}}^\hbar_i := 2\pi\sqrt{-1} \, \hbar \, \frac{\partial}{\partial a_i} \qquad \mbox{on}\quad D_\Gamma,
\end{align}
for each $i=1,\ldots,n$, i.e.
$$
(\wh{\mathbf{q}}^\hbar_i f)(a_1,\ldots,a_n) = a_i \, f(a_1,\ldots,a_n), \qquad
(\wh{\mathbf{p}}^\hbar_i f)(a_1,\ldots,a_n) = 2\pi\sqrt{-1} \, \hbar \, \frac{\partial f}{\partial a_i} (a_1,\ldots,a_n)
$$
for all $f(a_1,\ldots,a_n) \in D_\Gamma$.
\end{definition}
It is straightforward to check:
\begin{lemma}
\label{lem:Schrodinger_operators_lemma1}
Each of $\wh{\mathbf{q}}^\hbar_i$'s and $\wh{\mathbf{p}}^\hbar_i$'s leaves $D_\Gamma$ invariant, and is symmetric in the sense of Def.\ref{def:symmetric_operator}. \qed
\end{lemma}

\begin{lemma}
\label{lem:Schrodinger_operators_lemma2}
Each of $\wh{\mathbf{q}}^\hbar_i$'s and $\wh{\mathbf{p}}^\hbar_i$'s is essentially self-adjoint on $D_\Gamma$, in the Hilbert space $\mathscr{H}_\Gamma$. On $D_\Gamma$ they satisfy the `Heisenberg relations'
\begin{align}
\label{eq:Heisenberg_relations_of_p_and_q}
[\wh{\mathbf{p}}^\hbar_i, \wh{\mathbf{q}}^\hbar_j] = 2\pi \sqrt{-1} \, \hbar \, \delta_{i,j} \cdot \mathrm{id}.
\end{align}
Moreover, the corresponding strongly continuous unitary groups $e^{\sqrt{-1} \alpha_i \wh{\mathbf{q}}^\hbar_i}$, $e^{\sqrt{-1} \beta_i \wh{\mathbf{p}}^\hbar_i}$, defined by the functional calculus in \S\ref{subsec:spectral_theorem} for the unique self-adjoint extensions of the operators $\wh{\mathbf{q}}^\hbar_i$'s and $\wh{\mathbf{p}}^\hbar_i$'s with the functions $e^{\sqrt{-1}\alpha_i \lambda}$ and $e^{\sqrt{-1} \beta_i \lambda}$ respectively, for any real numbers $\alpha_i$'s and $\beta_i$'s, satisfy the `Weyl relations'
\begin{align}
\nonumber
  e^{\sqrt{-1} \alpha_i \wh{\mathbf{p}}^\hbar_i} \, e^{\sqrt{-1} \beta_i \wh{\mathbf{q}}^\hbar_i} = e^{- 2\pi \sqrt{-1} \, \hbar \, \alpha_i \beta_i} \, e^{\sqrt{-1} \beta_i \wh{\mathbf{q}}^\hbar_i} \, e^{\sqrt{-1} \alpha_i \wh{\mathbf{p}}^\hbar_i}.
\end{align}
\end{lemma}

\begin{lemma}[new $*$-representation of ${\bf D}^q_\Gamma$; {\cite[\S5]{FG09}}]
\label{lem:new_representation}
Let $\hbar \in\mathbb{R}_{>0}\setminus \mathbb{Q}$, and $q = e^{\pi \sqrt{-1}  \hbar}$. The following operators defined on the dense subspace $D_\Gamma$ of the Hilbert space $\mathscr{H}_\Gamma$ given by
\begin{align}
\label{eq:new_representation}
\displaystyle \wh{\mathbf{b}}^\hbar_i = (\wh{\mathbf{b}}^\hbar_i)^{\rm new} := \wh{\mathbf{q}}^\hbar_i, \quad\qquad \wh{\mathbf{x}}^\hbar_i = (\wh{\mathbf{x}}^\hbar_i)^{\rm new} := d_i^{-1} \, \wh{\mathbf{p}}^\hbar_i - \sum_{j=1}^n [\varepsilon_{ij}]_+ \, \wh{\mathbf{q}}^\hbar_j, \quad\qquad \forall i =1,\ldots,n,
\end{align}
where $\wh{\mathbf{q}}^\hbar_i$'s and $\wh{\mathbf{p}}^\hbar_i$'s are as in Def.\ref{def:Schrodinger_representation}, define a weakly irreducible positive integrable $*$-representation of the seed quantum $\mathcal{D}$-torus algebra $\mathbf{D}^q_\Gamma$, denoted by $(\mathscr{H}_\Gamma, D_\Gamma, \pi^q_\Gamma)$.
\end{lemma}
This representation is what is used in \cite{FG09} by Fock and Goncharov; they changed from another choice of representation from their previous paper \cite{FG07}. We will soon see that opting to use this new representation was not a good decision. Anyways, the `tilde' counterpart \eqref{eq:tilde_operators} of the new representation is given by
$$
\displaystyle \wh{\til{\mathbf{x}}}^\hbar_i = (\wh{\til{\mathbf{x}}}^\hbar_i)^{\rm new} := (\wh{\mathbf{x}}^\hbar_i)^{\rm new} + \sum_{j=1}^n \varepsilon_{ij} \, (\wh{\mathbf{b}}^\hbar_j)^{\rm new}
= d_i^{-1} \, \wh{\mathbf{p}}^\hbar_i - \sum_{j=1}^n [-\varepsilon_{ij}]_+ \, \wh{\mathbf{q}}^\hbar_j \qquad (\mbox{on $D_\Gamma$}).
$$
Their `old' representation is as follows (eq.(78) of \cite{FG09}, or \cite[\S4.1]{FG07}):
\begin{lemma}[old $*$-representation of ${\bf D}^q_\Gamma$; \cite{FG07}]
\label{lem:old_representation}
Let $\hbar \in\mathbb{R}_{>0}\setminus \mathbb{Q}$, and $q = e^{\pi \sqrt{-1} \hbar}$. The following operators defined on the dense subspace $D_\Gamma$ of the Hilbert space $\mathscr{H}_\Gamma$ given by
\begin{align}
\label{eq:old_representation}
(\wh{\mathbf{b}}^\hbar_i)^{\rm old} := 2\wh{\mathbf{q}}_i, \qquad
\displaystyle (\wh{\mathbf{x}}^\hbar_i)^{\rm old} := \frac{1}{2} d_i^{-1} \, \wh{\mathbf{p}}^\hbar_i - \sum_{j=1}^n \varepsilon_{ij} \, \wh{\mathbf{q}}^\hbar_j, \qquad \forall i=1,\ldots,n,
\end{align}
where $\hbar_i$'s are as in \eqref{eq:hbar_i}, define a weakly irreducible positive integrable $*$-representation of the seed quantum $\mcal{D}$-torus algebra $\mathbf{D}^q_\Gamma$, denoted by $(\mathscr{H}_\Gamma, D_\Gamma, (\pi^q_\Gamma)^{\rm old})$.
\end{lemma}
Again, the operators for the `tilde' counterpart \eqref{eq:tilde_operators} are given by
\begin{align}
\label{eq:old_representation_tilde}
\displaystyle (\wh{\til{\mathbf{x}}}^\hbar_i)^{\rm old} := (\wh{\mathbf{x}}^\hbar_i)^{\rm old} + \sum_{j=1}^n \varepsilon_{ij} \, (\wh{\mathbf{b}}^\hbar_j)^{\rm old}
= \frac{1}{2} d_i^{-1} \, \wh{\mathbf{p}}^\hbar_i + \sum_{j=1}^n \varepsilon_{ij} \, \wh{\mathbf{q}}^\hbar_j \qquad (\mbox{on $D_\Gamma$}).
\end{align}

In Lemmas \ref{lem:Schrodinger_operators_lemma2}, \ref{lem:new_representation}, and \ref{lem:old_representation}, the essential self-adjointness of the operators is something to be proven, although it may be taken to be a well-known fact, at least in the case of Lem.\ref{lem:Schrodinger_operators_lemma2}; I shall show this for Lem.\ref{lem:old_representation}, from which readers can easily figure out proofs for the other two lemmas. For convenience of presentation, we postpone the proof until the end of \S\ref{subsec:special_affine_shift_operators}.

\subsection{The Langlands modular dual}

For later purposes, it is necessary to define the notion of the `Langlands dual' of each seed. For a motivation, consult e.g. the `cluster linear algebra' in \cite[\S2.3]{FG09}. Here I just formulate the definitions that we will use, without explaining the reasoning behind them.

\begin{definition}
\label{def:Gamma_vee}
For a seed $\Gamma = (\varepsilon,d,*)$ of any kind, define the \ul{\em Langlands dual seed $\Gamma^\vee = (\varepsilon^\vee, d^\vee, *^\vee)$} as
\begin{align}
\label{eq:varepsilon_vee}
\varepsilon^\vee_{ij}:= d_i d_j^{-1} \varepsilon_{ij}, \qquad
d_i^\vee := d_i^{-1}, \qquad\forall i,j=1,\ldots,n.
\end{align}
\end{definition}
In particular, note that $\varepsilon^\vee_{ij} (d_j^\vee)^{-1} = d_i \varepsilon_{ij} \stackrel{\eqref{eq:2-form_on_A_Gamma}}{=} \til{\varepsilon}_{ij}$ is skew-symmetric. Here I point out an unfortunate typo in \cite[eq.(5)]{FG07} and \cite[end of \S2.1]{FG09}, which says $\varepsilon^\vee_{ij} = - \varepsilon^\vee_{ji}$, as if the matrix $\varepsilon^\vee$ is skew-symmetric. It should be corrected to
$$
\varepsilon^\vee_{ij} = - \varepsilon_{ji},
$$
which one can easily deduce from $\varepsilon^\vee_{ij} = d_i d_j^{-1} \varepsilon_{ij}$ and the skew-symmetry of $\wh{\varepsilon}_{ij} = d_j^{-1} \varepsilon_{ij}$; in particular, $\varepsilon^\vee$ is an integer matrix. The duality assignement $\varepsilon \leadsto \varepsilon^\vee$ of exchange matrices commutes with matrix mutations, i.e.
$$
(\varepsilon_{ij}')^\vee = (\varepsilon^\vee_{ij})'
$$
where the prime ${}'$ in the LHS means mutation $\mu_k$ applied to the matrix $\varepsilon$ of $\Gamma$ and the one in the RHS means mutation $\mu^\vee_k$  applied to the Langlands dual matrix $\varepsilon^\vee$ of $\Gamma^\vee$, which follows the same formula as the one of $\mu_k$ for the usual exchange matrix:
$$
(\varepsilon^\vee_{ij})' := \left\{
\begin{array}{ll}
-\varepsilon_{ij}^\vee & \mbox{if $i=k$ or $j=k$,} \\
\varepsilon_{ij}^\vee + \frac{1}{2}( |\varepsilon_{ik}^\vee| \varepsilon^\vee_{kj} + \varepsilon^\vee_{ik} | \varepsilon^\vee_{kj} | )& \mbox{otherwise}.
\end{array}
\right.
$$
Meanwhile, $(d_i^\vee)' = d_i^\vee$ for all $i$, under the mutation $\mu^\vee_k$. We declare that the cluster variables mutate as if $(\varepsilon^\vee,d^\vee,*^\vee)$; that is, we put $\vee$ at each cluster variable and $\varepsilon_{\cdot \, \cdot}$ appearing in the formulas \eqref{eq:mu_k_on_A}, \eqref{eq:mu_k_on_X}, \eqref{eq:mu_k_on_B}. We use the notation $P^\vee_\sigma$ for a seed automorphism of a Langlands dual seed, whose definition is same as $P_\sigma$ with $\vee$ put everywhere appropriate. So, a Langlands dual seed of a (genuine) seed can be viewed a generalized seed; it is not exactly a seed in the sense of Def.\ref{def:seed} because the skew-symmetrizer $(d^\vee)$ does not consist of positive integers.
\begin{definition}
\label{def:bf_m_vee}
For any cluster transformation $\mathbf{m}$ of (genuine) seeds, denote by $\mathbf{m}^\vee$ the corresponding cluster transformation of Langlands dual seeds.
\end{definition}

\begin{definition}
Let $\Gamma=(\varepsilon,d,\{B_i,X_i\}_{i=1}^n)$ be a $\mathcal{D}$-seed, and $\Gamma^\vee = (\varepsilon^\vee, d^\vee, \{B_i^\vee,X_i^\vee\}_{i=1}^n)$ be its Langlands dual $\mathcal{D}$-seed. The \ul{\em Langlands dual quantum $\mathcal{D}$-torus algebra $\mathbf{D}^{q^\vee}_{\Gamma^\vee}$}, for a formal quantum parameter $q^\vee$, is the free associative $*$-algebra over $\mathbb{Z}[q^\vee,(q^\vee)^{-1}]$ generated by $\mathbf{B}_i^\vee$, $\mathbf{X}_i^\vee$, $i=1,\ldots,n$, and their inverses, mod out by the relations
\begin{align*}
\begin{array}{rcll}
(q_j^\vee)^{- \varepsilon^\vee_{ij} } \mathbf{X}_i^\vee \mathbf{X}_j^\vee & = & (q^\vee_i)^{- \varepsilon^\vee_{ji}} \mathbf{X}_j^\vee \mathbf{X}_i^\vee, & \quad \forall i,j \in \{1,\ldots,n\}, \\
(q^\vee_i)^{-1} \mathbf{X}_i^\vee \mathbf{B}_i^\vee & = & (q^\vee_i) \mathbf{B}_i^\vee \mathbf{X}_i^\vee, & \quad \forall i=1,\ldots,n, \\
\mathbf{B}_i^\vee \mathbf{X}_j^\vee & = & \mathbf{X}_j^\vee \mathbf{B}_i^\vee, & \quad \mbox{whenever $i\neq j$}, \\
\mathbf{B}^\vee_i \mathbf{B}^\vee_j & = & \mathbf{B}^\vee_j \mathbf{B}^\vee_i, & \quad \forall i,j \in \{1,\ldots,n\},
\end{array}
\end{align*}
where
\begin{align}
\label{eq:q_i_vee}
q_i^\vee := (q^\vee)^{1/d_i^\vee} = (q^\vee)^{d_i} \in \mathbb{Z}[q^\vee,(q^\vee)^{-1}], \qquad \forall i=1,\ldots,n,
\end{align}
with the $*$-structure defined as the unique ring anti-homomorphism satisfying
$$
*\mathbf{X}^\vee_i = \mathbf{X}^\vee_i, \qquad *\mathbf{B}^\vee_i = \mathbf{B}^\vee_i, ~~~ \forall i \in \{1,\ldots,n\}, \qquad *q^\vee = (q^\vee)^{-1}.
$$
Denote by $\mathbb{D}^{q^\vee}_{\Gamma^\vee}$ the skew field of fractions of $\mathbf{D}^{q^\vee}_{\Gamma^\vee}$. For each $i\in \{1,\ldots,n\}$ define
$$
\til{\mathbf{X}}_i^\vee := \mathbf{X}_i^\vee \, \prod_{j=1}^n (\mathbf{B}_j^\vee)^{\varepsilon^\vee_{ij}} ~ \in ~ \mathbf{D}^{q^\vee}_{\Gamma^\vee}.
$$
\end{definition}

\begin{definition}[quantum mutation map for Langlands dual]
\label{def:quantum_mutation_map_vee}
Let the two $\mcal{D}$-seeds $\Gamma=(\varepsilon,d,*)$ and $\Gamma'=(\varepsilon',d',*')$ be related by the mutation along $k$, that is, $\Gamma' = \mu_k(\Gamma)$, so that the Langlands dual seeds $\Gamma^\vee = (\varepsilon^\vee, d^\vee, *^\vee)$ and $(\Gamma')^\vee = ( (\varepsilon')^\vee, (d')^\vee, (*')^\vee)$ are related by $(\Gamma')^\vee = \mu_k(\Gamma^\vee)$. Define the map $\mu^{q^\vee}_k : \mathbb{D}^{q^\vee}_{(\Gamma')^\vee} \to \mathbb{D}^{q^\vee}_{\Gamma^\vee}$ by
$$
\mu^{q^\vee}_k := \mu^{\sharp q^\vee}_k \circ (\mu'_k)^\vee,
$$
where $\mu^{\sharp q^\vee}_k$ is the automorphism of $\mathbb{D}^{q^\vee}_{\Gamma^\vee}$ given by the following formulas on generators
\begin{align}
\label{eq:quantum_mutation_sharp_on_B_i_vee}
\mu^{\sharp q^\vee}_k ({\bf B}^\vee_i) & = \left\{
\begin{array}{ll}
{\bf B}^\vee_i & \mbox{if $i\neq k$}, \\
{\bf B}^\vee_k (1+q^\vee_k {\bf X}^\vee_k) (1+q^\vee_k \til{{\bf X}}^\vee_k)^{-1} & \mbox{if $i=k$},
\end{array}
\right. \\
\label{eq:quantum_mutation_sharp_of_X_i_vee}
\mu^{\sharp q^\vee}_k ({\bf X}^\vee_i) & = {\bf X}^\vee_i \prod_{r=1}^{|\varepsilon^\vee_{ik}|} ( 1 + ((q_k^\vee)^{{\rm sgn}(-\varepsilon^\vee_{ik})} )^{2r-1} {\bf X}^\vee_k )^{ {\rm sgn}(-\varepsilon^\vee_{ik}) }, \qquad \forall i \in \{1,\ldots,n\},
\end{align}
and $(\mu'_k)^\vee$ is induced by the map $(\mu'_k)^\vee : {\bf D}^{q^\vee}_{(\Gamma')^\vee} \to {\bf D}^{q^\vee}_{\Gamma^\vee}$ given by
\begin{align}
\nonumber
& (\mu'_k)^\vee( ({\bf B}_i')^\vee ) = \left\{
\begin{array}{ll}
{\bf B}^\vee_i & \mbox{if $i\neq k$,} \\
({\bf B}_k^\vee)^{-1} \prod_{j=1}^n ({\bf B}^\vee_j)^{[-\varepsilon^\vee_{kj}]_+}  & \mbox{if $i = k$,}
\end{array}
\right.
\\
\nonumber
& (\mu'_k)^\vee(({\bf X}_i')^\vee) = \left\{
\begin{array}{ll}
(q_k^\vee)^{\varepsilon^\vee_{ik}[\varepsilon^\vee_{ik}]_+} {\bf X}^\vee_i \, ({\bf X}^\vee_k)^{[\varepsilon^\vee_{ik}]_+} & \mbox{if $i \neq k$}, \\
({\bf X}^\vee_k)^{-1} & \mbox{if $i=k$,}
\end{array}
\right.
\end{align}
on the generators of ${\bf D}^{q^\vee}_{(\Gamma')^\vee}$.
\end{definition}

\begin{lemma}
\label{lem:eta_q_vee}
Each quantum $\mcal{D}$-space functor $\eta^q$ \eqref{eq:algebraic_quantization_functor} induces a corresponding well-defined contravariant functor for the Langlands dual, denoted by $\eta^{q^\vee}$. \qed
\end{lemma}

\subsection{The Schwartz space for universally Laurent elements, and the unitary intertwiners as quantum mutation operators}
\label{subsec:Schwartz_space}

Is the notion of a (weakly irreducible) positive integrable $*$-representation $(\mathscr{H}_\Gamma, D_\Gamma, \pi^q_\Gamma)$ of the seed quantum $\mathcal{D}$-torus algebra $\mathbf{D}^q_\Gamma$ defined and discussed in \S\ref{subsec:positive_representations} satisfactory? Let's think about two issues. First, some operators were essentially self-adjoint on $D_\Gamma$, but the maximal extension of each operator may have a different domain from that of another. So it would be nice if it is possible to consider the \emph{common} maximal domain of `sufficiently many' operators. The more serious question arises when we try to build a relationship between positive integrable $*$-representations of the seed quantum $\mathcal{D}$-torus algebras for different seeds. 
Namely, for a representation $(\mathscr{H}_\Gamma, D_\Gamma,\pi^q_\Gamma)$ of $\mathbf{D}^q_\Gamma$ and a representation $(\mathscr{H}_{\Gamma'}, D_{\Gamma'}, \pi^q_{\Gamma'})$ of $\mathbf{D}^q_{\Gamma'}$, one would like to find a unitary map $\mathscr{H}_\Gamma \to \mathscr{H}_{\Gamma'}$ that `intertwines' the two representations $\pi^q_\Gamma$ and $\pi^q_{\Gamma'}$ of the two quantum torus algebras which are related by a certain composition of quantum mutation maps and quantum permutation maps. Then, would we want $D_\Gamma$ to be sent to $D_{\Gamma'}$ by this unitary map? In fact, what we want to be preserved by this unitary intertwining map is the common maximal domain mentioned above.

\vs

Let us be more precise. Suppose $(\mathscr{H}_\Gamma,D_\Gamma,\pi^q_\Gamma)$ is a positive integrable $*$-representation of $\mathbf{D}^q_\Gamma$. We then have operators $\pi^q_\Gamma(\mathbf{B}_i^{\pm 1})$ and $\pi^q_\Gamma(\mathbf{X}_i^{\pm 1})$ on $D_\Gamma$, leaving $D_\Gamma$ invariant and representing the generators of the algebra $\mathbf{D}^q_\Gamma$. As these operators satisfy the defining algebraic relations of the corresponding generators of $\mathbf{D}^q_\Gamma$, for any $\mathbf{u}\in \mathbf{D}^q_\Gamma$ they induce the well-defined operator $\pi^q_\Gamma(\mathbf{u})$ on $D_\Gamma$, leaving $D_\Gamma$ invariant. For any $\Gamma' \in |\Gamma|$ (Def.\ref{def:equivalence_of_seeds}) and a cluster transformation $\mathbf{m} \in \mathrm{Hom}_{\wh{\mathcal{G}}^\mcal{D}_{|\Gamma|}}(\Gamma,\Gamma')$ from $\Gamma$ to $\Gamma'$ in the saturated cluster modular groupoid $\wh{\mcal{G}}^\mcal{D}_{|\Gamma|}$, recall that $\eta^q(\mathbf{m}) : \mathbb{D}^q_{\Gamma'} \to \mathbb{D}^q_\Gamma$ (Def.\ref{def:quantum_D-space}) is the quantum map corresponding to $\mathbf{m}$ relating the two quantum torus algebras, which is an isomorphism of skew fields $\mathbb{D}^q_{\Gamma'}$ and $\mathbb{D}^q_\Gamma$. So, for as many as possible $\mathbf{u} \in \mathbb{D}^q_\Gamma$ we would like the densely defined operator $\pi^q_\Gamma(\mathbf{u})$ on $\mathscr{H}_\Gamma$ to correspond to the densely defined operator $\pi^q_{\Gamma'} ( (\eta^q(\mathbf{m}))^{-1} (\mathbf{u}))$ on $\mathscr{H}_{\Gamma'}$, via the sought-for unitary map $\mathscr{H}_\Gamma \to \mathscr{H}_{\Gamma'}$. However, at the moment, the operator $\pi^q_\Gamma(\mathbf{u})$ is known only when $\mathbf{u} \in \mathbf{D}^q_\Gamma$, while the operator $\pi^q_{\Gamma'} ( (\eta^q(\mathbf{m}))^{-1} (\mathbf{u}))$ is known only when $(\eta^q(\mathbf{m}))^{-1} (\mathbf{u}) \in \mathbf{D}^q_{\Gamma'}$. So we require $\mathbf{u} \in \mathbb{D}^q_\Gamma$ to be in the intersection of $\mathbf{D}^q_\Gamma$ and $\eta^q(\mathbf{m})(\mathbf{D}^q_{\Gamma'})$. So we first collect such elements.
\begin{definition}[quantum ring of universally Laurent polynomials]
Let $\mathscr{C}$ be an equivalence class of $\mathcal{D}$-seeds. For any $\Gamma \in \mathscr{C}$, the \ul{\emph{$*$-algebra of universally Laurent elements $\mathbb{L}^q_\mathscr{C}$}} is defined as the following subring of the skew field of fractions $\mathbb{D}^q_\Gamma$ of the seed quantum $\mathcal{D}$-torus algebra $\mathbf{D}^q_\Gamma$:
\begin{align}
\label{eq:bf_L_q_Gamma}
  \mathbb{L}^q_\Gamma := \bigcap_{\Gamma' \in \mathscr{C}} \,\, \bigcap_{\mathbf{m}\in \mathrm{Hom}_{\wh{\mathcal{G}}^\mathcal{D}_{\mathscr{C}}}(\Gamma,\Gamma')} \eta^q(\mathbf{m}) (\mathbf{D}^q_{\Gamma'}) \quad \subset \quad \mathbb{D}^q_\Gamma.
\end{align}
\end{definition}
As each $\eta^q(\mathbf{m}) : \mathbb{D}^q_{\Gamma'} \to \mathbb{D}^q_\Gamma$ is a $*$-algebra homomorphism, one observes that \eqref{eq:bf_L_q_Gamma} indeed defines a $*$-subring of $\mathbb{D}^q_\Gamma$. Also, this definition of $\mathbb{L}^q_\Gamma$ is `independent' of the choice of a $\mathcal{D}$-seed $\Gamma$ in the equivalence class $\mathscr{C}$. Namely, pick any $\Gamma' \in \mathscr{C}$ and $\mathbf{m} \in \mathrm{Hom}_{\wh{\mcal{G}}^\mcal{D}_\mathscr{C}}(\Gamma,\Gamma')$, then $\eta^q(\mathbf{m})$ provides a $*$-algebra isomorphism $\mathbb{L}^q_{\Gamma'} \to \mathbb{L}^q_{\Gamma}$; one may verify this, as well as the `consistency' of this identification, from $\eta^q(\mathbf{m}_1) \circ \eta^q(\mathbf{m}_2) = \eta^q(\mathbf{m}_2 \circ \mathbf{m}_1)$, which holds because $\eta^q$ is a contravariant functor. We collectively denote all $\mathbb{L}^q_\Gamma$ for $\Gamma \in \mathscr{C}$ by $\mathbb{L}^q_\mathscr{C}$.

\vs

So we have the following primitive version of a wish-list
: 
\begin{itemize}
\item[\rm (wish 1)] for each $\mathcal{D}$-seed $\Gamma \in \mathscr{C}$, a dense subspace $\mathscr{S}^q_\Gamma$ of $\mathscr{H}_\Gamma$ that is a `common maximal domain' on which $\pi^q_\Gamma(\mathbf{u})$ is defined for all $\mathbf{u}\in \mathbb{L}^q_\Gamma$,

\item[\rm (wish 2)] for each $\Gamma,\Gamma'\in \mathscr{C}$ and each cluster transformation $\mathbf{m} \in \mathrm{Hom}_{\wh{\mcal{G}}^\mcal{D}_\mathscr{C}}(\Gamma,\Gamma')$, a unitary map $\mathbf{K}^q(\mathbf{m}) : \mathscr{H}_{\Gamma'} \to \mathscr{H}_\Gamma$ such that
\begin{align}
\nonumber
\mathbf{K}^q(\mathbf{m}) (\mathscr{S}^q_{\Gamma'}) = \mathscr{S}^q_\Gamma,
\end{align}
intertwining the actions, in the sense that
\begin{align}
\label{eq:intertwining_equation1}
&  \pi^q_{\Gamma}(\mathbf{u}) \, \mathbf{K}^q(\mathbf{m}) \, v = \mathbf{K}^q(\mathbf{m}) \, \pi^q_{\Gamma'}((\eta^q(\mathbf{m}))^{-1}(\mathbf{u})) \, v, \qquad \forall v\in \mathscr{S}^q_{\Gamma'}, \quad \forall \mathbf{u}\in \mathbb{L}^q_{\Gamma},
\end{align}

\item[\rm (wish 3)] for each $\Gamma,\Gamma',\Gamma''\in \mathscr{C}$, $\mathbf{m}_1 \in \mathrm{Hom}_{\wh{\mcal{G}}^\mcal{D}_\mathscr{C}}(\Gamma,\Gamma')$, and  $\mathbf{m}_2 \in \mathrm{Hom}_{\wh{\mcal{G}}^\mcal{D}_\mathscr{C}}(\Gamma',\Gamma'')$ , the `consistency'
$$
\mathbf{K}^q(\mathbf{m}_1) \circ \mathbf{K}^q(\mathbf{m}_2) = c_{\mathbf{m}_1,\mathbf{m}_2} \, \mathbf{K}^q(\mathbf{m}_2\circ\mathbf{m}_1)
$$
holds for some complex constant $c_{\mathbf{m}_1,\mathbf{m}_2}$.
\end{itemize}
One interprets this situation as saying that $\mathbf{K}^q(\mathbf{m})$ intertwines the actions $\pi^q_\Gamma$ and $\pi^q_{\Gamma'}$ of $\mathbb{L}^q_\Gamma$ and $\mathbb{L}^q_{\Gamma'}$ on $\mathscr{S}^q_\Gamma$ and $\mathscr{S}^q_{\Gamma'}$ related by the quantum algebra map $\eta^q(\mathbf{m})$ in a consistent manner.

\vs

It turns out that the known solution $\mathbf{K}^q(\mathbf{m})$ to this problem admits the `modular double' phenomenon\footnote{`modular double' was first discovered in \cite{F95}.}. Namely, it also intertwines the actions of the Langlands modular dual counterparts (Def.\ref{def:Gamma_vee}), namely the actions $\pi^{q^\vee}_{\Gamma^\vee}$ and $\pi^{q^\vee}_{(\Gamma')^\vee}$ of $\mathbb{L}^{q^\vee}_{\Gamma^\vee}$ and $\mathbb{L}^{q^\vee}_{(\Gamma')^\vee}$ related by $\eta^{q^\vee}(\mathbf{m}^\vee)$ (Def.\ref{def:bf_m_vee}, Lem.\ref{lem:eta_q_vee}), where $\Gamma^\vee$ and $(\Gamma')^\vee$ are the Langlands dual $\mathcal{D}$-seeds of $\Gamma$ and $\Gamma'$ (Def.\ref{def:Gamma_vee}), while the quantum parameter for the Langlands dual
$$
q^\vee := e^{\pi\sqrt{-1}(1/\hbar)}
$$
is obtained by applying the `modular transformation' $\hbar \mapsto 1/\hbar$ to $q=e^{\pi \sqrt{-1} \hbar}$ \footnote{a genuine modular transformation would be $\hbar\to -1/\hbar$; so, in fact, it might be better to put $q^\vee = e^{\pi \sqrt{-1}(-1/\hbar)}$ instead, but $q^\vee=e^{\pi \sqrt{-1}(1/\hbar)}$ is what has been used in the literature.}. Moreover, the action of $\mathbb{L}^q_{\Gamma}$ via $\pi^q_{\Gamma}$ commutes with that of $\mathbb{L}^{q^\vee}_{\Gamma^\vee}$ via $\pi^{q^\vee}_{\Gamma^\vee}$ on a dense subspace of $\mathscr{H}_\Gamma \equiv \mathscr{H}_{\Gamma^\vee}$. This allows to put together these two actions to form a representation of the bigger algebra
\begin{align}
\label{eq:bf_L}
\mathbf{L}^\hbar_\Gamma := \mathbb{L}^q_\Gamma \otimes_\mathbb{Z} \mathbb{L}^{q^\vee}_{\Gamma^\vee}
\end{align}
called the `(Langlands) modular double' of $\mathbb{L}^q_\Gamma$, making the representation `strongly' irreducible, in the sense that there is no more operator that commutes with all the operators representing the elements of the algebra via the representation. As a matter of fact, \eqref{eq:bf_L} is not the full description of the Langlands modular double; we also impose the transcendental relations, saying that the $\frac{1}{\hbar_i}$-th power of each generator $\mathbf{B}_i$ (resp. $\mathbf{X}_i$) of $\mathbb{L}^q_\Gamma$ must `coincide with' the generator $\mathbf{B}_i^\vee$ (resp. $\mathbf{X}_i^\vee$) of $\mathbb{L}^{q^\vee}_{\Gamma^\vee}$, which makes sense when considering representations, as we shall see (see the definition of `modular double' e.g. in \cite{FK12} \cite{F95}), and which determines the precise relationship between $\pi^q_\Gamma$ and $\pi^{q^\vee}_{\Gamma^\vee}$.

\vs

Define the following operators on the dense subspace $D_\Gamma$ of the Hilbert space $\mathscr{H}_\Gamma$:
\begin{align}
\label{eq:checked_operators}
  (\wh{\mathbf{x}}^{\hbar}_i)^\vee := \frac{1}{\hbar_i} \wh{\mathbf{x}}^\hbar_i, \qquad\quad
  (\wh{\mathbf{b}}^\hbar_i)^\vee := \frac{1}{\hbar_i} \wh{\mathbf{b}}^\hbar_i, \qquad \forall i =1,\ldots,n,
\end{align}
where $\wh{\mathbf{x}}^\hbar_i$ and $\wh{\mathbf{b}}^\hbar_i$ are in \eqref{eq:new_representation} or in \eqref{eq:old_representation} and $\hbar_i = \hbar / d_i$. Then one immediately observes from \eqref{eq:Weyl_relations_for_D_q_Gamma} that the Weyl relations hold: namely, for $(\alpha)$ and $(\beta)$ as in \eqref{eq:alpha_and_beta}, the unitary operators
$$
(\wh{\mathbf{B}}^\vee_i)^{(\alpha)} := e^{\sqrt{-1} \alpha_i (\wh{\mathbf{b}}^\hbar_i)^\vee}, \qquad \wh{\mathbf{X}}^{(\beta)}_i := e^{\sqrt{-1} \beta_i (\wh{\mathbf{x}}^\hbar_i)^\vee}, \qquad i=1,\ldots,n,
$$
satisfy
\begin{align}
\nonumber
\left\{
{\renewcommand{\arraystretch}{1.3} \begin{array}{rcll}
e^{\pi \sqrt{-1} \hbar_j^\vee \varepsilon^\vee_{ij} \beta_i \beta_j} \, (\wh{\bf X}^\vee_i)^{(\beta)} (\wh{\bf X}^\vee_j)^{(\beta)} & = & e^{\pi \sqrt{-1} \hbar^\vee_i \varepsilon^\vee_{ji} \beta_j \beta_i} \, (\wh{\bf X}_j^\vee)^{(\beta)} (\wh{\bf X}^\vee_i)^{(\beta)}, & \quad \forall i,j\in \{1,\ldots,n\}, \\
e^{\pi\sqrt{-1} \hbar^\vee_i \beta_i \alpha_i} \, (\wh{\bf X}_i^\vee)^{(\beta)} (\wh{\bf B}_i^\vee)^{(\alpha)} & = & e^{-\pi\sqrt{-1} \hbar^\vee_i \alpha_i \beta_i} \, (\wh{\bf B}_i^\vee)^{(\alpha)} (\wh{\bf X}_i^\vee)^{(\beta)}, & \quad \forall i \in \{1,\ldots,n\},  \\
(\wh{\bf B}_i^\vee)^{(\alpha)} (\wh{\bf X}^\vee_j)^{(\beta)} & = & (\wh{\bf X}^\vee_j)^{(\beta)} (\wh{\bf B}^\vee_i)^{(\alpha)}, & \quad \mbox{whenever $i\neq j$}, \\
(\wh{\bf B}^\vee_i)^{(\alpha)} (\wh{\bf B}^\vee_j)^{(\alpha)} & = & (\wh{\bf B}^\vee_j)^{(\alpha)} (\wh{\bf B}^\vee_i)^{(\alpha)}, & \quad \forall i,j \in \{1,\ldots,n\},
\end{array} }
\right.
\end{align}
because, for example, the $e$-power in the LHS of the first equation is $\pi\sqrt{-1} \beta_i \beta_j$ times
$$
\hbar_i^{-1}\hbar_j^{-1} (\hbar_j \varepsilon_{ij}) = \hbar_i^{-1} \varepsilon_{ij} = \hbar^{-1} d_i \, \varepsilon_{ij} = \hbar^{-1} d_j (d_i d_j^{-1} \varepsilon_{ij}) = \hbar^\vee_j \, \varepsilon^\vee_{ij},
$$
where we denote
$$
\hbar^\vee_i := \hbar^\vee/d_i^\vee = 1/\hbar_i \quad\mbox{for each $i=1,\ldots,n$.}
$$
These relations imply the Heisenberg counterparts as an analog of \eqref{eq:Heisenberg_relations_for_D_q_Gamma}, which are easier to grasp than the Weyl relations: on $D_\Gamma$,
\begin{align*}
  [(\wh{\mathbf{x}}^\hbar_i)^\vee,(\wh{\mathbf{x}}^\hbar_j)^\vee] = 2\pi \sqrt{-1} \hbar_j^\vee \varepsilon_{ij}^\vee \cdot \mathrm{id}, \quad
  [(\wh{\mathbf{x}}^\hbar_i)^\vee,(\wh{\mathbf{b}}^\hbar_j)^\vee] = 2\pi \sqrt{-1} \hbar_i^\vee \delta_{i,j} \cdot \mathrm{id}, \quad
  [(\wh{\mathbf{b}}^\hbar_i)^\vee,(\wh{\mathbf{b}}^\hbar_j)^\vee] = 0,
\end{align*}
for all $i,j \in \{1,\ldots,n\}$. Moreover, again using \eqref{eq:Heisenberg_relations_for_D_q_Gamma} one finds
\begin{align}
\label{eq:commutation_with_vee_parts}
  [\wh{\mathbf{x}}^\hbar_i, (\wh{\mathbf{x}}^\hbar_j)^\vee] = 2\pi \sqrt{-1} \varepsilon_{ij} \cdot \mathrm{id}, \quad
  [\wh{\mathbf{x}}^\hbar_i, (\wh{\mathbf{b}}^\hbar_j)^\vee] = 2\pi \sqrt{-1} \delta_{i,j} \cdot \mathrm{id}, \quad
  [\wh{\mathbf{b}}^\hbar_i, (\wh{\mathbf{b}}^\hbar_j)^\vee] = 0,
\end{align}
for all $i,j$; note that all the commutators in \eqref{eq:commutation_with_vee_parts} are identity times a scalar in $2\pi \sqrt{-1} \, \mathbb{Z}$.

\vs

One can then check, e.g. by using the Stone-von Neumann theorem or by dealing with the explicit example we presented, that
\begin{lemma}
\begin{enumerate}
\item[\rm 1)] the operators $\pi^q_\Gamma(\mathbf{B}_i^{\pm 1}) := e^{\pm \wh{\mathbf{b}}^\hbar_i}$, $\pi^q_\Gamma(\mathbf{X}_i^{\pm 1}) := e^{\pm \wh{\mathbf{x}}^\hbar_i}$, $i=1,\ldots,n$, satisfy the defining relations of the generators $\mathbf{B}_i^{\pm 1}$, $\mathbf{X}_i^{\pm 1}$, $i=1,\ldots,n$, of the algebra $\mathbf{D}^q_\Gamma$ on $D_\Gamma$,

\item[\rm 2)] the operators $\pi^{q^\vee}_{\Gamma^\vee}( (\mathbf{B}_i^\vee)^{\pm 1}) := e^{\pm (\wh{\mathbf{b}}^\hbar_i)^\vee}$, $\pi^{q^\vee}_{\Gamma^\vee}((\mathbf{X}_i^\vee)^{\pm 1}) := e^{\pm (\wh{\mathbf{x}}^\hbar_i)^\vee}$, $i=1,\ldots,n$, preserve $D_\Gamma$, are essentially self-adjoint on $D_\Gamma$, and satisfy the defining relations of the generators $(\mathbf{B}_i^\vee)^{\pm 1}$, $(\mathbf{X}_i^\vee)^{\pm 1}$, $i=1,\ldots,n$, of the algebra $\mathbf{D}^{q^\vee}_{\Gamma^\vee}$ on $D_\Gamma$,

\item[\rm 3)] the actions of $\mathbf{D}^q_\Gamma$ and $\mathbf{D}^{q^\vee}_{\Gamma^\vee}$ on $D_\Gamma$ via $\pi^q_\Gamma$ and $\pi^{q^\vee}_{\Gamma^\vee}$ commute, in the sense that
\begin{align*}
  \pi^q_\Gamma(\mathbf{u}) \, \pi^{q^\vee}_{\Gamma^\vee}(\mathbf{v}) =   \pi^{q^\vee}_{\Gamma^\vee}(\mathbf{v})  \, \pi^q_\Gamma(\mathbf{u}) \quad \mbox{on} ~ D_\Gamma,
\end{align*}
for all generators $\mathbf{u}$ and $\mathbf{v}$ (hence for all elements) of $\mathbf{D}^q_\Gamma$ and $\mathbf{D}^{q^\vee}_{\Gamma^\vee}$.

\item[\rm 4)] One has
\begin{align}
\label{eq:pi_vee}
\pi^{q^\vee}_{\Gamma^\vee}( \mathbf{B}_i^\vee ) = ( \pi^q_\Gamma(\mathbf{B}_i) )^{1/\hbar_i}, \qquad \pi^{q^\vee}_{\Gamma^\vee}( \mathbf{X}_i^\vee ) = ( \pi^q_\Gamma(\mathbf{X}_i) )^{1/\hbar_i}, \qquad \forall i=1,\ldots,n,
\end{align}
where the right hand sides are defined using the functional calculus of the unique self-adjoint extensions of $\pi^q_\Gamma(\mathbf{B}_i)$ and $\pi^q_\Gamma(\mathbf{X}_i)$.
\end{enumerate}
\end{lemma}
In particular, 3) can be formally proved by using the Baker-Campbell-Hausdorff (BCH) formula
$$
e^X e^Y = e^{X+Y+\frac{1}{2}[X,Y] + \frac{1}{12}( [X,[X,Y]] + [Y,[Y,X]]) + \cdots}
$$
and \eqref{eq:commutation_with_vee_parts}. Note that 4) is the `transcendental' relation mentioned earlier. We could interpret the above 1)--4) as having a `positive integrable $*$-representation of the Langlands modular double algebra $\mathbf{L}^\hbar_\Gamma$'. Let us refer to this representation as
\begin{align}
\label{eq:pi_hbar}
\pi^\hbar_{\Gamma} := \pi^q_\Gamma \otimes \pi^{q^\vee}_{\Gamma^\vee}, 
\end{align}
for convenience, by a slight abuse of notation.

\vs

Finally, we shall define the sought-for Schwartz space $\mathscr{S}^\hbar_\Gamma$ (I use this notation instead of $\mathscr{S}^q_\Gamma$, hinting the presence of the Langlands modular double) as the common domain for the operators representing the elements of $\mathbf{L}^\hbar_\Gamma$ via $\pi^q_\Gamma \otimes \pi^{q^\vee}_{\Gamma^\vee}$, which form the promised collection of `sufficiently many' operators to be intertwined.
\begin{definition}[the Schwartz space of Fock and Goncharov; see e.g. {\cite[Def.5.2]{FG09}}]
For a $\mathcal{D}$-seed $\Gamma$, define the \\
\ul{\em seed $\mathcal{D}$-Schwartz space $\mathscr{S}_\Gamma$} as
\begin{align}
\label{eq:FG_Schwartz_space}
  \mathscr{S}^\hbar_\Gamma := \left\{ f \in \mathscr{H}_\Gamma \, \left| \, \begin{array}{l} \mbox{for each $\mathbf{u}\in \mathbf{L}^\hbar_\Gamma$ \eqref{eq:bf_L}, the linear functional $v \mapsto \langle \pi^\hbar_\Gamma(\mathbf{u}) v, \, f \rangle_{\mathscr{H}_\Gamma}$} \\ \mbox{is a continuous functional on $D_\Gamma$ \eqref{eq:D_Gamma}} \end{array} \right. \right\}
\end{align}
\end{definition}
Let $f \in \mathscr{S}^\hbar_\Gamma$. Then, for each $\mathbf{u} \in \mathbf{L}^\hbar_\Gamma$, how do we define $\pi^\hbar_\Gamma(\mathbf{u})f$? From the definition, the linear functional $v\mapsto \langle \pi^\hbar_\Gamma(*\mathbf{u})v, \, f \rangle_{\mathscr{H}_\Gamma}$ on $D_\Gamma$ is continuous, because $*\mathbf{u} \in \mathbf{L}^\hbar_\Gamma$. Then, since $D_\Gamma$ is dense in $\mathscr{H}_\Gamma$, by Riesz Representation Theorem (or Riesz Lemma), there exists a unique vector $g \in \mathscr{H}_\Gamma$ `representing' this functional, i.e. $\langle \pi^\hbar_\Gamma(*\mathbf{u})v, \, f\rangle_{\mathscr{H}_\Gamma} = \langle v, g \rangle_{\mathscr{H}_\Gamma}$ for all $v\in D_\Gamma$. We now declare $\pi^\hbar_\Gamma(\mathbf{u})f$ to be this $g$. 

\vs

Let us see why this is natural. Now, suppose that for each $\mathbf{u}\in \mathbf{L}^\hbar_\Gamma$, the symbol $\pi^\hbar_\Gamma(\mathbf{u})$ stands for a densely defined operator on some chosen domain, say $D_\mathbf{u}$, which contains $D_\Gamma$. For each $\mathbf{u}\in \mathbf{L}^\hbar_\Gamma$, the adjoint $(\pi^\hbar_\Gamma(\mathbf{u}))^*$ is uniquely determined, together with its well-defined domain (Def.\ref{def:adjoint_of_an_operator}); we do not know yet whether its domain is dense. For all $v,w\in D_\Gamma$, $\pi^\hbar_\Gamma(\mathbf{u})v$ and $\pi^\hbar_\Gamma(*\mathbf{u})w$ are well-defined elements of $D_\Gamma$, and we can verify that $\langle \pi^\hbar_\Gamma(\mathbf{u})v, w\rangle_{\mathscr{H}_\Gamma} = \langle v, \pi^\hbar_\Gamma(*\mathbf{u})w\rangle_{\mathscr{H}_\Gamma}$. Thus, for each $w \in D_\Gamma$, the linear functional $v\mapsto \langle \pi^\hbar_\Gamma(\mathbf{u})v, w\rangle_{\mathscr{H}_\Gamma}$ on $D_\Gamma$ is continuous; namely, the norm of this functional is bounded by $||\pi^\hbar_\Gamma(*\mathbf{u})w||_{\mathscr{H}_\Gamma} < \infty$. This means that each $w\in D_\Gamma$ is in the domain of $(\pi^\hbar_\Gamma(\mathbf{u}))^*$; in particular, $(\pi^\hbar_\Gamma(\mathbf{u}))^*$ is a densely defined operator.

\vs

We just saw that both the domain of $\pi^\hbar_\Gamma(\mathbf{u})$ and that of $(\pi^\hbar_\Gamma(\mathbf{u}))^*$ contain $D_\Gamma$. If we are to require $\pi^\hbar_\Gamma$ to be a $*$-representation, for each $\mathbf{u}\in\mathbf{L}^\hbar_\Gamma$ we must have
$$
(\pi^\hbar_\Gamma(\mathbf{u}))^* = \pi^\hbar_\Gamma(*\mathbf{u})
$$
preferably with the domains matching exactly. Let us weaken this a little bit, to require just
\begin{align}
\label{eq:pi_to_be_star-representation}
\pi^\hbar_\Gamma(*\mathbf{u}) \subset (\pi^\hbar_\Gamma(\mathbf{u}))^*.
\end{align}
Suppose that, some $f\in \mathscr{H}_\Gamma$ belongs to the domain of $\pi^\hbar_\Gamma(\mathbf{u})$, for some $\mathbf{u}\in\mathbf{L}^\hbar_\Gamma$. Then \eqref{eq:pi_to_be_star-representation}, together with the observation $*(*\mathbf{u})=\mathbf{u}$, says that $f$ is in the domain of $(\pi^\hbar_\Gamma(*\mathbf{u}))^*$. Hence $v\mapsto \langle \pi^\hbar_\Gamma(*\mathbf{u})v,f \rangle_{\mathscr{H}_\Gamma}$ is a continuous linear functional on the domain of $\pi^\hbar_\Gamma(*\mathbf{u})$, and therefore also on $D_\Gamma$. 

\vs

So, if $f \in \mathscr{H}_\Gamma$ belongs to the domain of $\pi^\hbar_\Gamma(\mathbf{u})$ for every $\mathbf{u}\in\mathbf{L}^\hbar_\Gamma$, then it must belong to $\mathscr{S}^\hbar_\Gamma$. In that case, for each $\mathbf{u}\in\mathbf{L}^\hbar_\Gamma$, the Riesz Representation Theorem for the continuous linear functional $v\mapsto \langle \pi^\hbar_\Gamma(*\mathbf{u})v,f\rangle_{\mathscr{H}_\Gamma}$ on $D_\Gamma$ says that $(\pi^\hbar_\Gamma(*\mathbf{u}))^*f$ is the unique element of $\mathscr{H}_\Gamma$ such that $\langle \pi^\hbar_\Gamma(*\mathbf{u})v, f\rangle_{\mathscr{H}_\Gamma} = \langle v, (\pi^\hbar_\Gamma(*\mathbf{u}))^* f\rangle_{\mathscr{H}_\Gamma}$ for all $v \in D_\Gamma$. We then observe from \eqref{eq:pi_to_be_star-representation} that $\pi^\hbar_\Gamma(\mathbf{u})f = (\pi^\hbar_\Gamma(*\mathbf{u}))^*f$. This is consistent with our earlier description of $\pi^\hbar_\Gamma(\mathbf{u})f$ for each $f\in \mathscr{S}^\hbar_\Gamma$.

\vs

Another crucial aspect of the Schwartz space $\mathscr{S}^\hbar_\Gamma$ is the Fr\'echet topology on it \cite{G08} \cite{FG09}, given by the countable family of seminorms $\rho_\mathbf{u}$ for $\mathbf{u}$ running through a basis in $\mathbf{L}^\hbar_\Gamma$, where
$$
\rho_\mathbf{u}(v) := ||\pi^\hbar_\Gamma(\mathbf{u})v||_{\mathscr{H}_\Gamma}, \quad \forall v\in \mathscr{S}^\hbar_\Gamma.
$$
Although their subspaces $\mathscr{S}_{\bf i}$ and $W_{\bf i}$ of the Hilbert space $L^2(\mathscr{A}^+_{\bf i}) \equiv L^2(\mathbb{R}^n)$ are different from our subspaces $\mathscr{S}^\hbar_\Gamma$ and $D_\Gamma$ of $\mathscr{H}_\Gamma \equiv L^2(\mathbb{R}^n)$, I believe a same proof as the one for \cite[Prop.5.8]{FG09} would give the following, which can be interpreted as saying that the algebra of operators $\pi^\hbar_\Gamma(\mathbf{L}^\hbar_\Gamma)$ acting on the dense subspace $\mathscr{S}^\hbar_\Gamma$ of the Hilbert space $\mathscr{H}_\Gamma$ is essentially self-adjoint, as mentioned by Fock-Goncharov \cite{FG09}.
\begin{proposition}
The space $D_\Gamma$ is dense in the Schwartz space $\mathscr{S}^\hbar_\Gamma$ in the Fr\'echet topology. \qed
\end{proposition}

\subsection{Groupoid formulation of representations}

Here we present a groupoid formulation of representations of (algebraic) quantum cluster $\mathcal{D}$-varieties, which summarizes and refines the discussion in the previous subsection. First, recall from Definitions \ref{def:positive_integrable_representation} and \ref{def:weakly_irreducible} the notion of a weakly irreducible positive $*$-representation $(\mathscr{H}_\Gamma,D_\Gamma,\pi^q_\Gamma)$ of a $\mathcal{D}$-seed $\Gamma$, where $q=e^{\pi\sqrt{-1}\hbar}$ with $\hbar \in \mathbb{R}_{>0}\setminus\mathbb{Q}$; the datum $\pi^q_\Gamma$ consists of essentially self-adjoint operators $\wh{\mathbf{b}}^\hbar_{\Gamma;i}$, $\wh{\mathbf{x}}^\hbar_{\Gamma;i}$, $i=1,\ldots,n$, on a dense subspace $D_\Gamma$ of a complex Hilbert space $\mathscr{H}_\Gamma$, satisfying the `Weyl relations' as in Def.\ref{def:positive_integrable_representation}.1). This entails some subsequent notions: the `tilde' operators $\wh{\til{\mathbf{x}}}^\hbar_{\Gamma;i}$ \eqref{eq:tilde_operators}, $i=1,\ldots,n$, the `checked' operators $(\wh{\mathbf{b}}^\hbar_{\Gamma;i})^\vee$, $(\wh{\mathbf{x}}^\hbar_{\Gamma;i})^\vee$ \eqref{eq:checked_operators}, $i=1,\ldots,n$, the representation $\pi^\hbar_\Gamma$ of the algebra $\mathbf{L}^\hbar_\Gamma$ \eqref{eq:bf_L} of universally Laurent elements of the Langlands modular double quantum $\mathcal{D}$-torus algebra, and the Schwartz space $\mathscr{S}_\Gamma^\hbar$ \eqref{eq:FG_Schwartz_space}.

\begin{definition}
\label{def:category_of_representations}
Let $\mathscr{C}$ be an equivalence class of $\mathcal{D}$-seeds. Let $\hbar \in \mathbb{R}_{>0}\setminus\mathbb{Q}$ be a quantum parameter, and let $q = e^{\pi \sqrt{-1}\hbar}$. The \ul{\em category of representations of quantum cluster $\mathcal{D}$-variety $\mathrm{DRep}_\mathscr{C}^\hbar$} associated to $\mathscr{C}$ is a category whose objects are weakly irreducible positive $*$-representation $(\mathscr{H}_\Gamma, D_\Gamma, \pi^q_\Gamma)$ with $\Gamma$ in $\mathscr{C}$, and whose set of morphisms from $(\mathscr{H}_\Gamma, D_\Gamma, \pi^q_\Gamma)$ to $(\mathscr{H}_{\Gamma'}, D_{\Gamma'}, \pi^q_{\Gamma'})$ is the set of all intertwiners between these two objects, where an \ul{\em intertwiner} for a cluster transformation $\mathbf{m} \in \mathrm{Hom}_{\wh{\mcal{G}}^\mcal{D}_\mathscr{C}}(\Gamma,\Gamma')$ is a unitary map
$
\mathbf{K}^\hbar_{\Gamma,\Gamma'} : \mathscr{H}_{\Gamma'} \to \mathscr{H}_\Gamma
$
satisfying the conditions
\begin{itemize}
\item[\rm (Sch 1)] $\mathbf{K}^\hbar_{\Gamma,\Gamma'}(\mathscr{S}^\hbar_{\Gamma'}) = \mathscr{S}^\hbar_\Gamma$,

\item[\rm (Sch 2)] $\mathbf{K}^\hbar_{\Gamma,\Gamma'} \, \pi^\hbar_{\Gamma'}(\mathbf{u}') \, v = \pi^\hbar_\Gamma(\eta^\hbar(\mathbf{m})(\mathbf{u}')) \, \mathbf{K}^\hbar_{\Gamma,\Gamma'} \, v$ for all $v \in \mathscr{S}_{\Gamma'}$ and $\mathbf{u}' \in \mathbf{L}^\hbar_{\Gamma'}$,
\end{itemize}
where $\eta^\hbar(\mathbf{m}) : \mathbf{L}^\hbar_{\Gamma'} \to \mathbf{L}^\hbar_\Gamma$ is the restriction of
\begin{align}
\label{eq:eta_hbar}
\eta^\hbar(\mathbf{m}) := \eta^q(\mathbf{m}) \otimes \eta^{q^\vee}(\mathbf{m}^\vee) : \mathbb{D}^q_{\Gamma'} \otimes \mathbb{D}^{q^\vee}_{(\Gamma')^\vee} \to \mathbb{D}^q_\Gamma \otimes \mathbb{D}^{q^\vee}_{\Gamma^\vee}.
\end{align}
\end{definition}

\begin{definition}
Let $\mathscr{C}$, $\hbar$ be as in Def.\ref{def:category_of_representations}. By a \ul{\em representation of quantum cluster $\mathcal{D}$-variety} associated to $\mathscr{C}$ we mean a contravariant projective functor
\begin{align}
\label{eq:representation_functor}
\mathbf{K}^\hbar ~:~ \mbox{{\em the saturated cluster modular groupoid}} ~~
\wh{\mcal{G}}^\mcal{D}_{\mathscr{C}} ~~ ({\rm Def}.\ref{eq:saturated-modular_groupoid}) ~ \longrightarrow ~ {\rm DRep}^\hbar_{\mathscr{C}}
\end{align}
such that $\mathbf{K}^\hbar(\mathbf{m})$ is an intertwiner for $\mathbf{m}$ (Def.\ref{def:category_of_representations}) for each morphism $\mathbf{m}$ of $\wh{\mcal{G}}^\mcal{D}_\mathscr{C}$, where a \ul{\em projective} functor means that the composition of morphisms are preserved only up to constants: for morphisms $\mathbf{m}_1,\mathbf{m}_2$ such that $\mathbf{m}_1 \circ \mathbf{m}_2$ is well-defined, there is a constant $c_{\mathbf{m}_1,\mathbf{m}_2} \in \mathrm{U}(1) \subset \mathbb{C}^\times$ such that
$$
\mathbf{K}^\hbar(\mathbf{m}_1\circ \mathbf{m}_2) = \, c_{\mathbf{m}_1,\mathbf{m}_2} \cdot \mathbf{K}^\hbar(\mathbf{m}_2) \, \mathbf{K}^\hbar(\mathbf{m}_1).
$$
\end{definition}

Recall from \S\ref{subsec:cluster_modular_groupoids}  that the groupoid $\wh{\mcal{G}}^\mcal{D}_\mathscr{C}$ can be presented as `generators and relations', with generators being mutations $\mu_k$ (for $k\in \{1,\ldots,n\}$), together with seed automorphisms $P_\sigma$ (for $\sigma$ a permutation of the set $\{1,\ldots,n\}$), and relations being the ones in Lemmas \ref{lem:involution_identity}, \ref{lem:permutation_identities}, and \ref{lem:h_plus_2-gon_relations}. So, in order to give a construction of a representation $\mathbf{K}^\hbar$ of quantum cluster $\mathcal{D}$-variety associated to $\mathscr{C}$, it suffices to describe what $\mathbf{K}^\hbar(\mu_k)$ and $\mathbf{K}^\hbar(P_\sigma)$ are, and verify that they satisfy the relations (written in the reverse orders) up to constants. For $\mathbf{K}^\hbar(P_\sigma)$ we take
\begin{align}
\label{eq:bf_K_P_sigma}
\mathbf{K}^\hbar(P_\sigma) := \mathbf{P}_\sigma ~ : ~ \mathscr{H}_{P_\sigma(\Gamma)} \to \mathscr{H}_\Gamma,
\end{align}
given by
$$
(\mathbf{P}_\sigma f)(a_1,\ldots,a_n) := f(a_1,\ldots,a_n), \qquad \forall f \in \mathscr{H}_{P_\sigma(\Gamma)} = L^2(\mathbb{R}^n, \, da_{\sigma(1)} \, \cdots \, da_{\sigma(n)}).
$$
This description of $\mathbf{P}_\sigma$ calls for a few words. Here, elements of $\mathscr{H}_{P_\sigma(\Gamma)}$ are functions in the variables $a_{\sigma(1)},\ldots,a_{\sigma(n)}$, where the notation is such that the $i$-th argument is $a_{\sigma(i)}$. To see what is happening, it is wise to extend the definition of $\mathbf{P}_\sigma$ to the space of `all' functions on $\mathbb{R}^n$ in the variables $a_{\sigma(1)},\ldots,a_{\sigma(n)}$, by the same formula. For example, the symbol $f:= a_{\sigma(i)}$, regarded as a function on $\mathbb{R}^n$ in the variables $a_{\sigma(1)}, \ldots, a_{\sigma(n)}$, is the function that spits out the $i$-th coordinate as its value (in the $P_\sigma(\Gamma)$ world); one can write $f(a_{\sigma(1)},\ldots,a_{\sigma(n)}) = a_{\sigma(i)}$ explicitly. The above formula then says $(\mathbf{P}_\sigma f)(a_1,\ldots,a_n) = a_i$; so, $\mathbf{P}_\sigma f$, regarded as a function on $\mathbb{R}^n$ in the variables $a_1,\ldots,a_n$, is what spits out the $i$-th coordinate, hence can be conveniently wrote as $a_i$ (in the $\Gamma$ world). Hence
$$
\mathbf{P}_\sigma \, a_{\sigma(i)} = a_i, \qquad \forall i=1,\ldots,n.
$$
The following two subsections \S\ref{subsec:special_affine_shift_operators} and \S\ref{subsec:non-compact_QD} are devoted to the two main ingredients for the construction in \S\ref{subsec:formulas_for_intertwiners} of the intertwiner $\mathbf{K}^\hbar(\mu_k)$ for a mutation $\mu_k$. The relations which must be satisfied by $\mathbf{K}^\hbar(\mu_k)$ and $\mathbf{K}^\hbar(P_\sigma)$ are the subject of the next section \S\ref{sec:operator_identities}.

\subsection{Special affine shift operators}
\label{subsec:special_affine_shift_operators}

Let $(M,\mu)$ and $(N,\nu)$ be some measure spaces, and let $\phi : M \to N$ be an invertible measure-preserving map whose inverse is also measure-preserving. It induces a natural unitary isomorphism of Hilbert spaces $L^2(M,d\mu) \to L^2(N,d\nu)$ given by $f \mapsto (\phi^{-1})^* f = f \circ \phi^{-1}$; this simple observation is generally referred to as `Koopman's Lemma' \cite{K31}. In the present subsection we define and investigate some special cases, when both measure spaces are the usual Euclidean space $(\mathbb{R}^n, da_1 \, \ldots \, da_n)$ and $\phi: \mathbb{R}^n \to \mathbb{R}^n$ is of the following type.
\begin{definition}
Write each element of $\mathbb{R}^n$ as a row vector $\mathbf{a} = (a_1 \, a_2 \, \cdots \, a_n)$ with real entries. 

A bijective map $\phi : \mathbb{R}^n \to \mathbb{R}^n$ is called an \ul{\emph{special affine transformation}} if there exist some matrix $\mathbf{c} = (c_{ij})_{i,j\in\{1,\ldots,n\}} \in \mathrm{SL}_\pm(n,\mathbb{R})$ and some $\mathbf{t} \in \mathbb{R}^n$, where
\begin{align}
\mathrm{SL}_\pm(n,\mathbb{R}) := \{ \, \mathbf{c} \in \mathrm{GL}(n,\mathbb{R}) \, : \, |\det \mathbf{c}|=1 \, \},
\end{align}
such that
\begin{align}
\label{eq:phi_inverse}
\phi^{-1}(\mathbf{a}) = \mathbf{a} \, \mathbf{c} + \mathbf{t}, \quad \forall \mathbf{a} \in \mathbb{R}^n.
\end{align}
\end{definition}
The notions and results in the present subsection are nothing fundamentally new, but as I have never seen them formulated as a special class of operators, I find it convenient to investigate them as I shall do here now.

\begin{lemma}
The above correspondence $\phi \leadsto (\mathbf{c},\mathbf{t})$ is a group isomorphism between the group of all special affine transformations of $\mathbb{R}^n$ and the following semi-direct product group
\begin{align}
\mathrm{SL}_\pm(n,\mathbb{R}) \ltimes \mathbb{R}^n = \{ (\mathbf{c}, \mathbf{t}) \, | \, \mathbf{c}\in \mathrm{SL}_\pm(n,\mathbb{R}), \,\, \mathbf{t} \in \mathbb{R}^n \},
\end{align}
whose multiplication is given by
\begin{align}
\label{eq:multiplication_of_semidirect_product}
(\mathbf{c}, \mathbf{t}) \, (\mathbf{c}', \mathbf{t}') = (\mathbf{c} \, \mathbf{c}', \mathbf{t} \, \mathbf{c}' + \mathbf{t}').
\end{align}
\end{lemma}
\textit{Proof.} Suppose two special affine transformations $\phi$ and $\psi$ are given by $\phi^{-1} (\mathbf{a}) = \mathbf{a} \, \mathbf{c} + \mathbf{t}$ and $\psi^{-1} (\mathbf{a}) = \mathbf{a} \, \mathbf{c}' + \mathbf{t}'$. Then $(\phi \circ \psi)^{-1}(\mathbf{a}) = \psi^{-1}(\phi^{-1}(\mathbf{a})) = \phi^{-1} (\mathbf{a}) \, \mathbf{c}' + \mathbf{t}' = (\mathbf{a} \, \mathbf{c} + \mathbf{t}) \, \mathbf{c}' + \mathbf{t}' = \mathbf{a} \, \mathbf{c} \, \mathbf{c}' + \mathbf{t} \, \mathbf{c}' + \mathbf{t}'$. And the identity special affine transformation corresponds to $(\mathbf{id}, \mathbf{0}) \in \mathrm{SL}_\pm(n,\mathbb{R}) \ltimes \mathbb{R}^n$, where $\mathbf{id}$ denotes the identity matrix and $\mathbf{0} = (0 \, \cdots \, 0)$. \qed

\vs

Let us adapt the notation for the Hilbert space $\mathscr{H}_\Gamma$ and its dense subspace $D_\Gamma$, defined in \eqref{eq:H_Gamma} and \eqref{eq:D_Gamma} in Def.\ref{def:H_Gamma_and_D_Gamma}. Note that the $\mathcal{D}$-seed data $\Gamma$, except for its `rank' $n$, is not playing any role in the definitions of $\mathscr{H}_\Gamma$ and $D_\Gamma$; it does only when we consider multiple $\mathcal{D}$-seeds at the same time, which is not the case in the present subsection. So, in principle we can drop the subscript $\Gamma$ for now, but I keep it in order to avoid any unnecessary confusion.

\begin{lemma}
For any special affine transformation $\phi$, the map $f \mapsto (\phi^{-1})^* f = f \circ \phi^{-1}$ is a unitary isomorphism from the Hilbert space $\mathscr{H}_\Gamma$ \eqref{eq:H_Gamma} to itself.
\end{lemma}

\textit{Proof.} Denote the volume form on $\mathbb{R}^n$ by $\mathrm{vol} := da_1 \wedge da_2 \wedge \cdots \wedge da_n$. Then for any $f \in \mathscr{H}_\Gamma$ we have
$$
||f||_{\mathscr{H}_\Gamma}^2 = \int_{\mathbb{R}^n} |f|^2 \, da_1\cdots da_n
= \int_{\mathbb{R}^n} |f|^2 \, \mathrm{vol},
$$
where the last integral is viewed as an integral of an $n$-form on the $n$-dimensional parametrized oriented manifold $\mathbb{R}^n$. Note now that
$$
\int_{\mathbb{R}^n} |f|^2 \, \, \mathrm{vol}
= \int_{\phi(\mathbb{R}^n)} (\phi^{-1})^*( |f|^2 \,\, \mathrm{vol})
= \int_{\phi(\mathbb{R}^n)} |(\phi^{-1})^*f|^2 \,\, (\phi^{-1})^* \mathrm{vol}.
$$
Observe that $(\phi^{-1})^* \mathrm{vol} = \det(D(\phi^{-1})) \cdot \mathrm{vol}$, where $D(\phi^{-1})$ is the Jacobian matrix of $\phi^{-1}$: the $(i,j)$-th entry is the first order partial derivative of the $i$-th component function of $\phi^{-1} :\mathbb{R}^n \to \mathbb{R}^n$ with respect to the $j$-th variable $a_j$. One can easily verify from \eqref{eq:phi_inverse} that $\det(D(\phi^{-1}))=\pm 1$; it being $1$ means $\phi$ being orientation-preserving, while $-1$ means orientation-reversing. In either case we have
$$
\int_{\phi(\mathbb{R}^n)} |(\phi^{-1})^*f|^2 \,\, \det(D(\phi^{-1})) \cdot \mathrm{vol} = \int_{\mathbb{R}^n} |(\phi^{-1})^*f|^2 \,\, \mathrm{vol} = ||(\phi^{-1})^* f||_{\mathscr{H}_\Gamma}^2,
$$
as desired. \qed

\vs

\begin{definition}
\label{def:special_affine_shift_operator}
For each special affine transformation $\phi$, given by $\phi^{-1}(\mathbf{a}) = \mathbf{a} \, \mathbf{c} + \mathbf{t}$, denote the unitary automorphism $f\mapsto (\phi^{-1})^*f$ of the Hilbert space $\mathscr{H}_\Gamma$ \eqref{eq:H_Gamma} by $\mathbf{S}_{(\mathbf{c},\mathbf{t})}$, which we call a \ul{\em special affine shift operator}.
\end{definition}
More explicitly, $\mathbf{S}_{(\mathbf{c},\mathbf{t})}$ is given by
\begin{align}
\label{eq:def_of_S_c_t}
  (\mathbf{S}_{(\mathbf{c},\mathbf{t})}f)(\mathbf{a}) = f(\mathbf{a} \, \mathbf{c} + \mathbf{t}), \qquad \forall f\in \mathscr{H}_\Gamma=L^2(\mathbb{R}^n, da_1\cdots da_n), \quad \forall \mathbf{a}\in \mathbb{R}^n.
\end{align}

\begin{lemma}
The correspondence $(\mathbf{c}, \mathbf{t}) \leadsto \mathbf{S}_{(\mathbf{c}, \mathbf{t})}$ is a group homomorphism, i.e.
\begin{align}
  \label{eq:bf_S_is_group_homomorphism}
  \mathbf{S}_{(\mathbf{c}, \mathbf{t}) (\mathbf{c}',\mathbf{t}')} = \mathbf{S}_{(\mathbf{c}, \mathbf{t})} \, \mathbf{S}_{(\mathbf{c}', \mathbf{t}')}, \qquad\quad \mathbf{S}_{(\mathbf{id}, \mathbf{0})} = \mathrm{id}_{\mathscr{H}_\Gamma},
\end{align}
and is injective.
\end{lemma}
\textit{Proof.} Straighforward exercise. \qed

\begin{corollary}
\label{cor:special_affine_shift_operators_form_a_group}
The special affine shift operators on $\mathscr{H}_\Gamma$ \eqref{eq:H_Gamma} form a group, isomorphic to $\mathrm{SL}_\pm(n,\mathbb{R})\ltimes \mathbb{R}^n$. \qed
\end{corollary}

From the definiton \eqref{eq:def_of_S_c_t} we immediately obtain the following observation, which becomes useful in the proof of the main result of the present paper:
\begin{lemma}
\label{lem:scalar_special_affine_shift_operator_is_identity}
If $\mathbf{S}_{(\mathbf{c},\mathbf{t})} = c \cdot \mathrm{id}_{\mathscr{H}_\Gamma}$ for some complex scalar $c$, then $c=1$. \qed
\end{lemma}
The following lemma is also immediate:
\begin{lemma}
\label{lem:bf_S_preserves_D}
$\mathbf{S}_{(\mathbf{c},\mathbf{t})} (D_\Gamma) = D_\Gamma$. 
\end{lemma}
{\it Proof.}
It is easy to see that each $\mathbf{S}_{(\mathbf{c},\mathbf{t})}$ leaves $D_\Gamma$ invariant. The result then follows from the invertibility of $\mathbf{S}_{(\mathbf{c},\mathbf{t})}$.
\qed

\vs

The permutation operator $\wh{\mathbf{P}}_\sigma$ associated to a permutation $\sigma$ of $\{1,\ldots,n\}$, which is given by $(\wh{\mathbf{P}}_\sigma f)(a_1,\ldots,a_n) = f(a_{\sigma(1)},\ldots,a_{\sigma(n)})$, is an example of a special affine shift operator. Another class of examples is given by the usual `shift operators $\mathbf{S}_{(\mathbf{id},\mathbf{t})}$: for any $\mathbf{t} \in \mathbb{R}^n$ observe
\begin{align}
\label{eq:bf_S_as_shift_operator}
(\mathbf{S}_{(\mathbf{id}, \mathbf{t})}f)(\mathbf{a}) = f(\mathbf{a}+\mathbf{t}). 
\end{align}
In view of the Taylor series expansion, one may expect that we could write $\mathbf{S}_{(\mathbf{id},\mathbf{t})} = e^{\sum_{i=1}^n t_i \frac{\partial}{\partial a_i}}$, which a priori makes sense only formally. The following lemma allows us to do this rigorously.
\begin{lemma}[shift operators form a strongly continuous unitary group]
\label{lem:shift_operator_as_exponential}
For any $\mathbf{t} = (t_1 \, \cdots \, t_n) \in \mathbb{R}^n$, the family $\{\mathbf{S}_{(\mathbf{id},\alpha \mathbf{t})}\}_{\alpha \in \mathbb{R}}$ is a strongly continuous one-parameter unitary group on the Hilbert space $\mathscr{H}_\Gamma$ \eqref{eq:H_Gamma} in the sense of Def.\ref{def:strongly_continuous_one-parameter_unitary_group}, and it preserves $D_\Gamma$ \eqref{eq:D_Gamma}. Its infinitesimal generator, in the sense of Thm.\ref{thm:Stones_theorem}, is the unique self-adjoint extension of the essentially self-adjoint operator on the dense subspace $D_\Gamma$ of $\mathscr{H}_\Gamma$ given as $- \sum_{i=1}^n \frac{t_i}{2\pi \hbar} \wh{\mathbf{p}}^\hbar_i $, where $\wh{\mathbf{p}}^\hbar_i$'s are as defined in Def.\ref{def:Schrodinger_representation}.
\end{lemma}

\textit{Proof.} The condition \eqref{eq:Stone_condition1} for this family comes easily from \eqref{eq:bf_S_is_group_homomorphism} and \eqref{eq:multiplication_of_semidirect_product}. Note now that
\begin{align*}
  ||\mathbf{S}_{(\mathbf{id},\alpha \mathbf{t})}f - \mathbf{S}_{(\mathbf{id},\alpha_0 \mathbf{t})}f ||_{\mathscr{H}_\Gamma}^2
= \int_{\mathbb{R}^n} | f(\mathbf{a} + \alpha \mathbf{t}) - f(\mathbf{a} + \alpha_0 \mathbf{t}) |^2 \, d\mathbf{a}.
\end{align*}
For now, view $f$ as an everywhere-defined square-integrable measurable function on $\mathbb{R}^n$; that is, pick one representative for the equivalence class of functions, which is an element of $\mathscr{H}_\Gamma$. Further, assume that $f$ is continuous. Then for each $\mathbf{a}$ we have $f(\mathbf{a}+\alpha \mathbf{t}) - f(\mathbf{a}+\alpha_0 \mathbf{t}) \to 0$ as $\alpha \to \alpha_0$. On the other hand, $| f(\mathbf{a} + \alpha \mathbf{t}) - f(\mathbf{a} + \alpha_0 \mathbf{t}) |^2 \le (|f(\mathbf{a}+\alpha\mathbf{t})| + |f(\mathbf{a}+\alpha_0 \mathbf{t})|)^2 \le 2|f(\mathbf{a}+\alpha\mathbf{t})|^2 + 2|f(\mathbf{a}+\alpha_0\mathbf{t})|^2 = 2|(\mathbf{S}_{(\mathbf{id},\alpha\mathbf{t})}f)(\mathbf{a})|^2 + 2|(\mathbf{S}_{(\mathbf{id},\alpha_0\mathbf{t})}f)(\mathbf{a})|^2$, which is integrable because $\mathbf{S}_{(\mathbf{id},\alpha\mathbf{t})}f$ and $\mathbf{S}_{(\mathbf{id},\alpha_0\mathbf{t})}f$ are square-integrable. So by the Dominated Convergence Theorem we get $||\mathbf{S}_{(\mathbf{id},\alpha \mathbf{t})}f - \mathbf{S}_{(\mathbf{id},\alpha_0 \mathbf{t})}f ||_{\mathscr{H}_\Gamma}^2 \to \int_{\mathbb{R}^n} 0 \, d\mathbf{a} = 0$ as $\alpha\rightarrow\alpha_0$. Using the facts that continuous square-integrable functions is dense in $\mathscr{H}_\Gamma$ and that $\mathbf{S}_{(\mathbf{id},\alpha \mathbf{t})}$ and $\mathbf{S}_{(\mathbf{id},\alpha \mathbf{t})}$ are unitary, one can deduce $||\mathbf{S}_{(\mathbf{id},\alpha \mathbf{t})}f - \mathbf{S}_{(\mathbf{id},\alpha_0 \mathbf{t})}f ||_{\mathscr{H}_\Gamma}^2 \to0$ as $\alpha\rightarrow\alpha_0$ for any $f\in \mathscr{H}_\Gamma$, establishing the condition \eqref{eq:Stone_condition2}.

\vs

We shall now use Thm.\ref{thm:strong_derivative} to obtain the infinitesimal generator. Lem.\ref{lem:bf_S_preserves_D} says that each $\mathbf{S}_{(\mathbf{id},\alpha\mathbf{t})}$ leaves $D_\Gamma$ invariant. Let $f\in D_\Gamma\subset \mathscr{H}_\Gamma$, and pick its representative so that it literally is an element of the space described in the RHS of \eqref{eq:D_Gamma}. Note that for each $\mathbf{a}$ we have ${\displaystyle \lim_{\alpha\rightarrow 0}} \frac{ (\mathbf{S}_{(\mathbf{id},\alpha\mathbf{t})}f)(\mathbf{a}) - f(\mathbf{a})}{\alpha} = {\displaystyle \lim_{\alpha\rightarrow 0}} \frac{f(\mathbf{a}+\alpha\mathbf{t}) - f(\mathbf{a})}{\alpha}=\,[\,$the directional derivative of $f$ at $\mathbf{a}$ in the direction of the vector $\mathbf{t} \, ] \, = \sum_{i=1}^n t_i \left( \frac{\partial f}{\partial a_i}(\mathbf{a}) \right)=\left( \sum_{i=1}^n t_i \frac{\partial f}{\partial a_i}\right)(\mathbf{a})$, and that $\sum_{i=1}^n t_i \frac{\partial f}{\partial a_i} = (\sum_{i=1}^n t_i \, \frac{\partial}{\partial a_i})f \in D_\Gamma$ (see e.g. Lem.\ref{lem:Schrodinger_operators_lemma1}). As $f$ is $C^\infty$ smooth on $\mathbb{R}^n$, from basic multivariable calculus (the Taylor expansion) we get
\begin{align}
\label{eq:Taylor_expansion}
  f(\mathbf{a}+\alpha\mathbf{t}) - f(\mathbf{a}) 
& = 
\alpha \sum_{i=1}^n t_i \left( \frac{\partial f}{\partial a_i} (\mathbf{a}) \right) +
\int_0^\alpha \int_0^\beta \sum_{i,j =1}^n t_i t_j \left( \frac{\partial^2 f}{\partial a_j \partial a_i} (\mathbf{a} + \gamma \mathbf{t}) \right) d\gamma \, d\beta
\end{align}
for any $\mathbf{a},\mathbf{t}\in\mathbb{R}^n$ and any $\alpha\in\mathbb{R}$, where $\beta,\gamma$ are real dummy variables. To check the validity of \eqref{eq:Taylor_expansion}, regard both sides as functions in $\alpha$, with fixed $\mathbf{a}$ and $\mathbf{t}$; put $\alpha=0$ to both sides, put $\alpha=0$ to both sides after applying $\frac{d}{d\alpha}$, and put $\alpha=0$ to both sides after applying $\frac{d^2}{d\alpha^2}$. From \eqref{eq:Taylor_expansion}, for $\alpha\neq 0$ we get
\begin{align*}
  \left| \frac{f(\mathbf{a}+\alpha\mathbf{t}) - f(\mathbf{a}) }{\alpha} - \left(\sum_{i=1}^n t_i \, \frac{\partial f}{\partial a_i}\right)(\mathbf{a}) \right| \le \frac{1}{|\alpha|} \, \left| \int_0^\alpha \int_0^\beta \left| \left( \sum_{i,j=1}^n t_i t_j \frac{\partial^2 f}{\partial a_j \partial a_i} \right)(\mathbf{a}+\gamma\mathbf{t}) \right| \, d\gamma \, d\beta \, \right|.
\end{align*}
For a fixed $\mathbf{t}=(t_1 \, \cdots \, t_n) \in \mathbb{R}^n$ and $f\in D_\Gamma$, let
$$
F := \sum_{i,j=1}^n t_i t_j \frac{\partial^2 f}{\partial a_j \partial a_i}.
$$
As each $\frac{\partial^2 f}{\partial a_j \partial a_i} = \frac{\partial}{\partial a_j}(\frac{\partial}{\partial a_i})f \in D_\Gamma$ (e.g. from Lem.\ref{lem:Schrodinger_operators_lemma1}), we have $F\in D_\Gamma$. Denote $|\mathbf{a}| := \sqrt{a_1^2+\cdots+a_n^2}$ for any $\mathbf{a}\in \mathbb{R}^n$. Considering the basis of $D_\Gamma$ as written in \eqref{eq:D}, it is a straightforward exercise to show that here exists positive real numbers $R,k$ such that $|F(\mathbf{a})|\le e^{-k |\mathbf{a}|^2}$ whenever $|\mathbf{a}| \ge R$. Now let $L := {\displaystyle \max_{|\mathbf{a}|\le 3R}} \, |F(\mathbf{a})|$, which is a positive real number.

\vs

Suppose from now on that $|\alpha|$ is small enough such that $|\alpha| \, |\mathbf{t}| < R$. Our dummy variable $\gamma$ runs in between $0$ and $\beta$, while $\beta$ in between $0$ and $\alpha$; hence $|\gamma| \, |\mathbf{t}| \le |\alpha|\, |\mathbf{t}| <R$. If $|\mathbf{a}| \le 2R$, then $|\mathbf{a}+\gamma \mathbf{t}| \le |\mathbf{a}| + |\gamma| \, |\mathbf{t}| < 3R$, therefore $|F(\mathbf{a}+\gamma\mathbf{t})|\le L$. Now suppose $|\mathbf{a}| \ge 2R$. Then $|\mathbf{a}+\gamma \mathbf{t}| \ge |\mathbf{a}| - |\gamma| \, |\mathbf{t}| > 2R - R = R$, therefore $|F(\mathbf{a}+\gamma \mathbf{t})| \le e^{-k|\mathbf{a}+\gamma\mathbf{t}|^2}$. On the other hand $|\mathbf{a}+\gamma\mathbf{t}| \ge |\mathbf{a}| - |\gamma|\, |\mathbf{t}| > |\mathbf{a}| - R \ge |\mathbf{a}| - \frac{|\mathbf{a}|}{2} = \frac{|\mathbf{a}|}{2}$, therefore $e^{-k|\mathbf{a}+\gamma\mathbf{t}|^2} \le e^{-k|\mathbf{a}|^2/4}$. Now let $L' := {\displaystyle \min_{|\mathbf{a}|\le 3R}} \, |e^{-k|\mathbf{a}|^2/4}|$, which is a positive real number. Put $L'' := \max\{ \frac{L}{L'}, 1\}$, which is a positive real number. Then, for $|\mathbf{a}|\le 2R$ we have $|F(\mathbf{a}+\gamma\mathbf{t})| \le L \le L'' L' \le L'' e^{-k|\mathbf{a}|^2/4}$, while for $|\mathbf{a}|\ge 2R$ we have $|F(\mathbf{a}+\gamma\mathbf{t})| \le e^{-k|\mathbf{a}|^2/4} \le L'' e^{-k|\mathbf{a}|^2/4}$. 

\vs

What have we proved? We fixed $f$ and $\mathbf{t}$, hence $F$ is fixed; then we can find some $R,k$ from $F$. Then the positive real constants $L,L,L''$ can be obtained from $F$ and $R$. Require that $|\alpha|$ is small enough such that $|\alpha| \, |\mathbf{t}| < R$. We proved that $|F(\mathbf{a}+\gamma\mathbf{t})|\le L'' e^{-k|\mathbf{a}|^2/4}$ holds for all $\mathbf{a} \in \mathbb{R}^n$ and for all values of the real dummy variable $\gamma$ in between $0$ and $\alpha$. So we have
\begin{align*}
  \left| \frac{f(\mathbf{a}+\alpha\mathbf{t}) - f(\mathbf{a}) }{\alpha} - \left(\sum_{i=1}^n t_i \, \frac{\partial f}{\partial a_i}\right)(\mathbf{a}) \right| & \le \frac{1}{|\alpha|} \, \left| \int_0^\alpha \int_0^\beta \left| F(\mathbf{a}+\gamma\mathbf{t}) \right| \, d\gamma \, d\beta \, \right| \\
& \le \frac{1}{|\alpha|} \, L'' e^{-k|\mathbf{a}|^2/4} \, \left| \int_0^\alpha \int_0^\beta d\gamma \, d\beta \, \right| \\
& = \frac{|\alpha|}{2} \, L'' \, e^{-k|\mathbf{a}|^2/4}.
\end{align*}
For our fixed $f\in D_\Gamma$ and $\mathbf{t}\in \mathbb{R}^n$, for each real number $\alpha$ satisfying $|\alpha| \, |\mathbf{t}| < R$ and $\frac{|\alpha|}{2} \cdot L'' \le 1$ define a function $F_\alpha : \mathbb{R}^n \to \mathbb{R}$ as
\begin{align*}
  F_\alpha(\mathbf{a}) := \frac{f(\mathbf{a}+\alpha\mathbf{t}) - f(\mathbf{a}) }{\alpha} - \left(\sum_{i=1}^n t_i \, \frac{\partial f}{\partial a_i}\right)(\mathbf{a}), \qquad \forall\mathbf{a}\in\mathbb{R}^n.
\end{align*}
In particular, $F_\alpha \in D_\Gamma$. We just proved that $|F_\alpha(\mathbf{a})| \le \frac{|\alpha|}{2} L'' e^{-k|\mathbf{a}|^2/4} \le e^{-k|\mathbf{a}|^2/4}$ holds for all $\mathbf{a}\in \mathbb{R}^n$, where $k$ is a positive real constant obtained from $f$ and $\mathbf{t}$. Recall that the pointwise limit of the funciton $F_\alpha$ is the zero function, i.e. for each $\mathbf{a}\in \mathbb{R}^n$ we have ${\displaystyle \lim_{\alpha\to 0}} F_\alpha(\mathbf{a})=0$, hence for each $\mathbf{a}\in\mathbb{R}^n$ we have ${\displaystyle \lim_{\alpha\rightarrow 0}} |F_\alpha(\mathbf{a})|^2=0$. On the other hand, the function $|e^{-k|\mathbf{a}|^2/4}|^2$ is integrable on $\mathbb{R}^n$ (note $e^{-k|\mathbf{a}|^2/4} \in D_\Gamma$). Therefore by the Dominated Convergence Theorem, we have ${\displaystyle \lim_{\alpha\to 0}} \int_{\mathbb{R}^n} |F_\alpha(\mathbf{a})|^2 \, d\mathbf{a} = \int_{\mathbb{R}^n} \lim_{\alpha\to0} |F_\alpha(\mathbf{a})|^2 \, d\mathbf{a} = \int_{\mathbb{R}^n} 0 \, d\mathbf{a} =0$. Hence ${\displaystyle \lim_{\alpha\to 0}} ||F_\alpha(\mathbf{a})-0||_{\mathscr{H}_\Gamma}=0$, so ${\displaystyle \lim_{\alpha\to 0}} F_\alpha = 0$ with respect to the standard Hilbert space topology of $\mathscr{H}_\Gamma$.

\vs

In conclusion, ${\displaystyle \lim_{\alpha\rightarrow 0}} \frac{ (\mathbf{S}_{(\mathbf{id},\alpha\mathbf{t})}f)(\mathbf{a}) - f(\mathbf{a})}{\alpha} = {\displaystyle \lim_{\alpha\rightarrow 0}} \frac{f(\mathbf{a}+\alpha\mathbf{t}) - f(\mathbf{a})}{\alpha} = \left( \sum_{i=1}^n t_i \frac{\partial f}{\partial a_i}\right)(\mathbf{a})$, where the two limits are not just pointwise limits, but are limits of elements of $\mathscr{H}_\Gamma$ with respect to the standard Hilbert space topology of $\mathscr{H}_\Gamma$. So we can write
$$
\lim_{\alpha\to 0} \frac{ \mathbf{S}_{(\mathbf{id},\alpha\mathbf{t})} f - f }{\alpha} =\left( \sum_{i=1}^n \frac{\partial}{\partial a_i} \right)f, \qquad \forall f\in D_\Gamma,
$$
where the limit is with respect to the standard Hilbert space topology of $\mathscr{H}_\Gamma$. Then, from Thm.\ref{thm:strong_derivative} we obtain the desired result. \qed

\vs

It is handy to have the following application of Lem.\ref{lem:shift_operator_as_exponential}:
\begin{lemma}
\label{lem:bf_S_conjugation_on_bf_p_and_bf_q}
Let $\mathbf{c} \in \mathrm{SL}_\pm(n,\mathbb{R})$ and consider $\mathbf{S}_{(\mathbf{c},\mathbf{0})} : \mathscr{H}_\Gamma \to \mathscr{H}_\Gamma$. Then, for each $i=1,\ldots,n$,
\begin{align}
\label{eq:conjugation_of_bf_S_on_bf_p}
  \mathbf{S}_{(\mathbf{c},\mathbf{0})} \, \wh{\mathbf{p}}^\hbar_i \,\, \mathbf{S}_{(\mathbf{c},\mathbf{0})}^{-1} = \sum_{j=1}^n c^{ij} \, \wh{\mathbf{p}}^\hbar_j \qquad \mbox{on}~ D_\Gamma, \quad \mbox{where} \quad \mathbf{c}^{-1} = (c^{ij})_{i,j\in\{1,\ldots,n\}}.
\end{align}
This equality also holds when $\wh{\mathbf{p}}^\hbar_i$ and $\sum_{j=1}^n c^{ij}\, \wh{\mathbf{p}}^\hbar_j$ are replaced by their unique self-adjoint extensions, on the domain of the self-adjoint extension of $\sum_{j=1}^n c^{ij} \, \wh{\mathbf{p}}^\hbar_j$.

\vs

One also has for each $i=1,\ldots,n$,
\begin{align}
\label{eq:conjugation_of_bf_S_on_bf_q}
  \mathbf{S}_{(\mathbf{c},\mathbf{0})} \,\, \wh{\mathbf{q}}^\hbar_i \,\, \mathbf{S}_{(\mathbf{c},\mathbf{0})}^{-1} = \sum_{j=1}^n c_{ji} \, \wh{\mathbf{q}}^\hbar_j \qquad \mbox{on}~ D_\Gamma, \quad \mbox{where} \quad \mathbf{c} = (c_{ij})_{i,j\in\{1,\ldots,n\}}.
\end{align}
\end{lemma}

{\it Proof.} Observe for any $\mathbf{c}\in \mathrm{SL}_\pm(n,\mathbb{R})$, $\mathbf{t}\in\mathbb{R}^n$, $\alpha\in\mathbb{R}$ that
\begin{align*}
  \mathbf{S}_{(\mathbf{c},\mathbf{0})} \mathbf{S}_{(\mathbf{id},\alpha \mathbf{t})} \mathbf{S}_{(\mathbf{c},\mathbf{0})}^{-1} = \mathbf{S}_{(\mathbf{c},\mathbf{0})(\mathbf{id},\alpha\mathbf{t})(\mathbf{c},\mathbf{0})^{-1}} = \mathbf{S}_{(\mathbf{c}, \alpha\mathbf{t})(\mathbf{c}^{-1},\mathbf{0})} = \mathbf{S}_{(\mathbf{id},\alpha \mathbf{t} \mathbf{c}^{-1})}
\end{align*}
Put $\mathbf{t} = - 2\pi \hbar \, \mathbf{e}_i$, where $\mathbf{e}_i := (0 \, \cdots \, 1\, \cdots \, 0)$, the standard $i$-th basis row vector of $\mathbb{R}^n$. Lem.\ref{lem:shift_operator_as_exponential} tells us that the unique self-adjoint extension of the essentially self-adjoint operator $\frac{2\pi \hbar}{2\pi \hbar} \, \wh{\mathbf{p}}^\hbar_i = \wh{\mathbf{p}}^\hbar_i$ (on $D_\Gamma$) is the infinitesimal generator of the strongly continuous unitary group $\{\mathbf{S}_{(\mathbf{id},\alpha\mathbf{t})}\}_{\alpha\in\mathbb{R}}$, while the self-adjoint extension of $ \sum_{j=1}^n \frac{ 2\pi \hbar \, c^{ij} }{2\pi \hbar } \wh{\mathbf{p}}^\hbar_j = \sum_{j=1}^n c^{ij} \, \wh{\mathbf{p}}^\hbar_j$, where $c^{ij}$ is the $(i,j)$-th entry of the matrix $\mathbf{c}^{-1}$, is the infinitesimal generator of $\{\mathbf{S}_{(\mathbf{id},\alpha\mathbf{t}\mathbf{c}^{-1})}\}_{\alpha\in\mathbb{R}}$. As every special linear shift operator is a unitary operator leaving $D_\Gamma$ invariant and since each $\wh{\mathbf{p}}^\hbar_i$ leaves $D_\Gamma$ invariant, from Lem.\ref{lem:unitary_conjugation_corollary} we get \eqref{eq:conjugation_of_bf_S_on_bf_p}. Meanwhile, \eqref{eq:conjugation_of_bf_S_on_bf_q} is an easy check. \qed

\vs

In particular, we can now get to a promised proof of the essential self-adjointness of the operators appearing in Lem.\ref{lem:old_representation}:

\vs

{\it Proof of the essential self-adjointness in Lem.\ref{lem:old_representation}.} \quad Recall the Fourier transforms $\mcal{F}_i : D_\Gamma \to D_\Gamma$ from \eqref{eq:mcal_F_i} of Lem.\ref{lem:D_is_preserved_by_Fourier_transforms}; we know these are (more precisely, extend to) unitary operators. We can deduce from \eqref{eq:mcal_F_conjugation_on_p_and_q} how $\mcal{F}_i$ `exchanges' the operators $\wh{\mathbf{q}}^\hbar_i = a_i$ and $\wh{\mathbf{p}}^\hbar_i = 2\pi \sqrt{-1} \hbar \frac{\partial}{\partial a_i}$ on $D_\Gamma$, which preserve $D_\Gamma$: 
\begin{align}
\label{eq:mcal_F_i_conjugation}
\mcal{F}_i^{-1} \,\, \wh{\mathbf{q}}^\hbar_i \, \mcal{F}_i = -\frac{1}{(2\pi)^2 \hbar} \, \wh{\mathbf{p}}^\hbar_i \qquad \mbox{as operators $D_\Gamma \to D_\Gamma$}.
\end{align}
It is easy to verify $\mcal{F}_i^{-1} \, \wh{\mathbf{p}}^\hbar_j \, \mcal{F}_i = \wh{\mathbf{p}}^\hbar_j$ as operators $D_\Gamma \to D_\Gamma$, whenever $i\neq j$. Consider
$$
\til{\mcal{F}}_i := \prod_{j\neq i} \mcal{F}_j.
$$  
The product order does not matter, for the factors are mutually commuting unitary operators; we note that $\til{\mcal{F}}_i$ is unitary, and preserves $D_\Gamma$. Observe that
\begin{align*}
  \til{\mcal{F}}_i^{-1} \,\, \wh{\mathbf{x}}^\hbar_i  \, \til{\mcal{F}}_i \,\, \stackrel{\eqref{eq:old_representation}, ~ \varepsilon_{ii}=0}{=} \,\, \prod_{j\neq i} \mcal{F}_j^{-1} \left( \frac{d_i^{-1}}{2} \wh{\mathbf{p}}^\hbar_i - \sum_{j\neq i} \varepsilon_{ij} \, \wh{\mathbf{q}}^\hbar_j \right) \, \prod_{j\neq i} \mcal{F}_j
= \frac{d_i^{-1}}{2}\wh{\mathbf{p}}^\hbar_i + \sum_{j\neq i} \frac{\varepsilon_{ij}}{(2\pi)^2 \hbar} \, \wh{\mathbf{p}}^\hbar_j =: P_i,
\end{align*}
as operators $D_\Gamma \to D_\Gamma$. In Lem.\ref{lem:shift_operator_as_exponential} we saw that this operator $P_i$ on $D_\Gamma$ appears as the essentially self-adjoint operator whose unique self-adjoint extension is the infinitesimal generator of a strongly continuous one-parameter unitary group. All we need here is the fact that $P_i$ on $D_\Gamma$ is essentially self-adjoint. As we have $\til{\mcal{F}}_i^{-1} \, \, \wh{\mathbf{x}}^\hbar_i \, \til{\mcal{F}}_i = P_i$ as operators $D_\Gamma\to D_\Gamma$, where $\til{\mcal{F}}_i$ is unitary with $\til{F}_i(D_\Gamma) = D_\Gamma$, one observes that $\wh{\mathbf{x}}^\hbar_i$ on $D_\Gamma$ is essentially self-adjoint as well, thanks to Lem.\ref{lem:unitarily_equivalent_essentially_self-adjoint_operators}. Almost the same argument proves the essential self-adjointness of $\wh{\til{\mathbf{x}}}^\hbar_i$. As for the operator $\wh{\mathbf{b}}^\hbar_i = 2\wh{\mathbf{q}}^\hbar_i$ on $D_\Gamma$, from the equality $\mcal{F}_i^{-1} \wh{\mathbf{b}}^\hbar_i \mcal{F}_i = - \frac{2}{(2\pi)^2 \hbar} \wh{\mathbf{p}}^\hbar_i$ of operators $D_\Gamma \to D_\Gamma$ whose RHS appears as an essentially self-adjoint operator on $D_\Gamma$ in Lem.\ref{lem:shift_operator_as_exponential}, we see that $\wh{\mathbf{b}}^\hbar_i$ is essentially self-adjoint on $D_\Gamma$, again thanks to Lem.\ref{lem:unitarily_equivalent_essentially_self-adjoint_operators}. \qed

\subsection{Quantum dilogarithm functions}
\label{subsec:non-compact_QD}

First recall the `compact' quantum dilogarithm of Faddeev-Kashaev \cite{FK94} \cite{F95}.
\begin{definition}[the compact quantum dilogarithm]
For any complex number $\mathbf{q}$ with $|\mathbf{q}|<1$ define a merophorphic function $\Psi^\mathbf{q}(z)$ on $\mathbb{C}$ as
\begin{align}
\label{eq:Psi_q}
\Psi^\mathbf{q}(z) := \prod_{i=1}^\infty (1 + \mathbf{q}^{2i-1} z)^{-1} = \frac{1}{(1+\mathbf{q}z)(1+\mathbf{q}^3 z)(1+\mathbf{q}^5 z) \cdots}
\end{align}
\end{definition}
\begin{lemma}
The above infinite product absolutely converges for any $z$ except at the poles $z = \mathbf{q}^{-(2i-1)}$, $i=1,2,3,\ldots$. \qed
\end{lemma}
\begin{lemma}
One has the functional equation
\begin{align}
\label{eq:functional_equation_of_compact_QD}
  \Psi^\mathbf{q}(\mathbf{q}^2 z) = (1+\mathbf{q}z) \, \Psi^\mathbf{q}(z). \qed
\end{align}
\end{lemma}
For the moment, regard $\mathbf{q}$ as a formal parameter; then one can make sense of $\Psi^\mathbf{q}(z)^{-1}$ as a formal power series in $z$ with coefficients in $\mathbb{Z}_{\ge 0}[\mathbf{q}] \subset \mathbb{Z}[\mathbf{q}]$. Let $N$ be as in Def.\ref{def:quantum_D-torus_algebra}, and assume that $\mathbf{q}^{1/N}$ is a well-defined formal variable, so that $\mathbf{q}_k := \mathbf{q}^{1/k} = (\mathbf{q}^{1/N})^{N/k}$ is well-defined. Then, one can verify that the automorphism $\mu_k^{\sharp \mathbf{q}}$ of $\mathbb{D}_\Gamma^\mathbf{q}$ in Def.\ref{def:quantum_mutation_map} can be realized as conjugation by a single expression $\Psi^{\mathbf{q}_k}(\mathbf{X}_k) \Psi^{\mathbf{q}_k}(\til{\mathbf{X}}_k)^{-1}$;
\begin{align}
\label{eq:mu_k_sharp_q_as_conjugation}
\mu_k^{\sharp \mathbf{q}}(u) = \Psi^{\mathbf{q}_k} (\mathbf{X}_k) \, \Psi^{\mathbf{q}_k} (\til{\mathbf{X}}_k)^{-1} \,\, u \,\, \Psi^{\mathbf{q}_k}(\til{\mathbf{X}}_k) \, \Psi^{\mathbf{q}_k}(\mathbf{X}_k)^{-1}, \qquad \forall u \in \mathbb{D}_\Gamma^\mathbf{q}.
\end{align}
Fock and Goncharov even took this formal conjugation \eqref{eq:mu_k_sharp_q_as_conjugation} as the definition of $\mu_k^{\sharp \mathbf{q}}$; see \cite[Def.3.1]{FG09}. The RHS of \eqref{eq:mu_k_sharp_q_as_conjugation} is not a priori an element of $\mathbb{D}_\Gamma^\mathbf{q}$, but one can check that it is, for each generator $u$ of $\mathbb{D}^\mathbf{q}_\Gamma$. For example,
\begin{align*}
\Psi^{\mathbf{q}_k} (\mathbf{X}_k) \, \Psi^{\mathbf{q}_k} (\til{\mathbf{X}}_k)^{-1} \,\, \mathbf{X}_i \,\, \Psi^{\mathbf{q}_k}(\til{\mathbf{X}}_k) \, \Psi^{\mathbf{q}_k}(\mathbf{X}_k)^{-1} & = \Psi^{\mathbf{q}_k} (\mathbf{X}_k) \,\, \mathbf{X}_i \,\, \Psi^{\mathbf{q}_k}(\mathbf{X}_k)^{-1} \\
 & = \mathbf{X}_i \Psi^{\mathbf{q}_k}(\mathbf{q}_k^{-2\varepsilon_{ik}} \mathbf{X}_k) \, \Psi^{\mathbf{q}_k}(\mathbf{X}_k)^{-1},
\end{align*}
which becomes \eqref{eq:quantum_mutation_sharp_of_X_i} after a successive application of \eqref{eq:functional_equation_of_compact_QD} on the part $\Psi^{\mathbf{q}_k}(\mathbf{q}_k^{-2\varepsilon_{ik}} \mathbf{X}_k) \, \Psi^{\mathbf{q}_k}(\mathbf{X}_k)^{-1}$.

\vs

However, for (physical) quantization of cluster varieties one needs the case when $|\mathbf{q}|=1$. Thus I recall the `non-compact' quantum dilogarithm; among its many guises, here I introduce a version as used in \cite{FG09}. The non-compact quantum dilogarithm can be thought of as a suitable limit of the compact quantum dilogarithm as $|\mathbf{q}|\to 1$. The na\"ive limit does not exist; only certain ratio of compact quantum dilogarithm does. Note first that if we write $\mathbf{q}=e^{\pi\sqrt{-1}h}$ for $h\in \mathbb{C}$, then $|\mathbf{q}|<1$ if and only if $\mathrm{Im}(h)>0$, and $|\mathbf{q}|=1$ means $\mathrm{Im}(h)=0$. As $z\mapsto -1/z$ is an example of a $\mathrm{PSL}(2,\mathbb{R})$ M\"obius transformation of the upper half plane, we see that $\mathrm{Im}(h)>0$ implies $\mathrm{Im}(-1/h)>0$. So, if we define $\mathbf{q}^\vee := e^{\pi \sqrt{-1}/h}$ for $\mathrm{Im}(h)>0$, we have $| 1/\mathbf{q}^\vee |<1$. When we form the ratio $\Psi^\mathbf{q}(e^z)/\Psi^{1/\mathbf{q}^\vee}(e^{z/h})$, the poles and zeroes make appropriate cancellation, so that the limit as $\mathrm{Im}(h) \to 0$ exists; this would then coincide with the non-compact quantum dilogarithm $\Phi^h(z)$ with $h\in \mathbb{R}$. Here I write a little bit stronger statement, included in the following summary of some of the basic properties of $\Phi^h$:
\begin{lemma}
For any $h\in \mathbb{C}$ with $\mathrm{Im}(h)\ge 0$ and $\mathrm{Re}(h)>0$, let
$$
\mathbf{q}:=e^{\pi\sqrt{-1} h}, \qquad \mathbf{q}^\vee := e^{\pi\sqrt{-1}/h}.
$$
\begin{enumerate}
\item[\rm a)] 
The integral in the expression
$$
\Phi^h(z) = \exp\left( - \frac{1}{4} \int_\Omega \frac{e^{-ipz}}{\sinh(\pi p) \sinh(\pi h p)} \right),
$$
where $\Omega$ is a contour on the real line that avoids the origin via a small half circle above the origin, absolutely converges for $z$ in the strip $\mathrm{Im}(z) < \pi(1 + \mathrm{Re}(h))$, yielding a non-vanishing complex analytic function $\Phi^h(z)$ on this strip. Each of the functional equations
\begin{align}
\label{eq:Phi_h_difference_equations}
\left\{ {\renewcommand{\arraystretch}{1.3} \begin{array}{rcl}
\Phi^h(z+2\pi\sqrt{-1} h) & = & (1+\mathbf{q}e^z) \Phi^h(z) \\
\Phi^h(z+2\pi\sqrt{-1}) & = & (1+\mathbf{q}^\vee e^{z/h}) \Phi^h(z)
\end{array}} \right.
\end{align}
holds when the two arguments of $\Phi^h$ are in the strip. These functional equations let us to analytically continue $\Phi^h$ to a meromorphic function on the whole plane $\mathbb{C}$, with
\begin{align*}
\mbox{the set of zeros} & = \{ (2\ell+1) \pi\sqrt{-1} +  (2m+1) \pi\sqrt{-1} h \, | \, \ell,m\in \mathbb{Z}_{\ge 0} \}, \quad \mbox{and} \\
\mbox{the set of poles} & =  \{ - (2\ell+1) \pi\sqrt{-1} - (2m+1) \pi\sqrt{-1} h \, | \, \ell,m\in \mathbb{Z}_{\ge 0} \}.
\end{align*}

\item[\rm b)] \emph{(relationship between compact and non-compact)} When $\mathrm{Im}(h)>0$, one has the equality
\begin{align}
\label{eq:Phi_h_as_ratio}
\frac{\Psi^\mathbf{q}(e^z)}{\Psi^{1/\mathbf{q}^\vee}(e^{z/h})} = \Phi^h(z),
\end{align}
for any $z$ that is not a pole of $\Phi^h$. 

\item[\rm c)] For any sequence $\{h_n\}_{n=1,2,\ldots}$ in $\mathbb{C}$ with $\mathrm{Im}(h_n)>0$ and $\mathrm{Re}(h_n)>0$ that converges to $\hbar \in \mathbb{R}_{>0}$, one has $\lim_{n\to \infty} \Phi^{h_n}(z) = \Phi^\hbar(z)$ for almost every $z$.

\item[\rm d)] When $h = \hbar \in \mathbb{R}_{>0}$, every pole and zero of $\Phi^h$ is simple if and only if $h\notin \mathbb{Q}$.

\item[\rm e)] \emph{($\hbar \leftrightarrow 1/\hbar$ duality)} When $h = \hbar \in \mathbb{R}_{>0}$, one has
\begin{align}
\label{eq:QD_self-duality}
  \Phi^{1/\hbar}(z/\hbar) = \Phi^\hbar(z).
\end{align}

\item[\rm f)] \emph{(unitarity)} When $h = \hbar \in \mathbb{R}_{>0}$, one has
$$
|\Phi^\hbar(z)|=1, \qquad \forall z\in \mathbb{R}.
$$

\item[\rm g)] \emph{(involutivity)} One has
\begin{align}
\label{eq:QD_identity_quadratic}
\Phi^h(z) \Phi^h(-z) = c_h\, \exp\left( \frac{z^2}{4\pi \sqrt{-1} \hbar} \right),
\end{align}
where
\begin{align}
\label{eq:c_hbar}
  c_h := e^{- \frac{\pi \sqrt{-1}}{12} (h + h^{-1})} \in \mathbb{C}^\times.
\end{align}
In particular, $|c_\hbar|=1$ for $\hbar = h \in \mathbb{R}_{>0}$.
\end{enumerate}
\qed
\end{lemma}

\begin{remark}
In the definition of $\Phi^h$, one may use a contour that avoids the origin by a small circle below the origin. The resulting value $\Phi^h(z)$ is same as before.
\end{remark}

Besides the `difference equations' \eqref{eq:functional_equation_of_compact_QD} and \eqref{eq:Phi_h_difference_equations}, the quantum dilogarithm functions satisfy a crucial equation called the `pentagon equation'.

\begin{proposition}[the pentagon equation for the compact quantum dilogarithm; \cite{FK94}]
Let $X,Y$ be formal variables satisfying $XY = \mathbf{q}^2YX$. Then one has
\begin{align}
\label{eq:compact_QD_pentagon}
  \Psi^\mathbf{q}(Y)^{-1} \, \Psi^\mathbf{q}(X)^{-1} = \Psi^\mathbf{q}(X)^{-1} \, \Psi^\mathbf{q}(\mathbf{q}XY)^{-1} \, \Psi^\mathbf{q}(Y)^{-1},
\end{align}
where each factor and each side is regarded as a formal power series in the variables $X$ and $Y$ with coefficients in $\mathbb{Z}_{\ge 0}[\mathbf{q}] \subset \mathbb{Z}[\mathbf{q}]$. \qed
\end{proposition}

When $X$ and $Y$ are formal variables with $XY = \mathbf{q}^2 YX$, we declare $X^{1/h}$ and $Y^{1/h}$ to be the new formal variables that satisfy $X^{1/h} \, Y^{1/h} = (\mathbf{q}^{1/\vee})^2 \, Y^{1/h} \, X^{1/h}$ and commute with $X$ and $Y$. Then from \eqref{eq:compact_QD_pentagon} for the compact quantum dilogarithms $\Psi^\mathbf{q}$ and $\Psi^{1/\mathbf{q}^\vee}$ one can formally deduce the corresponding pentagon equation for the non-compact quantum dilogarithm $\Phi^h$, which is given as ratio of these \eqref{eq:Phi_h_as_ratio}. It is best to formulate the pentagon equation of non-compact quantum dilogarithm in terms of an identity of unitary operators on a Hilbert space in the case when $h = \hbar \in \mathbb{R}_{>0} \setminus \mathbb{Q}$, so I postpone doing this until \S\ref{subsec:known_operator_identities}.

\subsection{Formulas for intertwiners}
\label{subsec:formulas_for_intertwiners}

I finally present Fock-Goncharov's formula for the intertwiner $\mathbf{K}^\hbar(\mu_k)$ for the mutation $\mu_k$ from a $\mcal{D}$-seed $\Gamma$ to $\Gamma'$. As a way of remembering where $\mu_k$ is being applied to and ends up with, we write
$$
\mathbf{K}^\hbar(\mu_k) = \mathbf{K}^\hbar_{\Gamma\mut{k}\Gamma'} ~ : ~ \mathscr{H}_{\Gamma'} \to \mathscr{H}_\Gamma, \quad \mbox{when} \quad \Gamma' = \mu_k(\Gamma).
$$
\begin{definition}[the intertwiner for a mutation; {\cite[Def.5.1]{FG09}}]
\label{def:bf_K}
Let $\Gamma = (\varepsilon,d,\{B_i,X_i\})_{i=1}^n$ be a $\mathcal{D}$-seed, and let $\mu_k(\Gamma) = \Gamma' = (\varepsilon',d',\{B_i',X_i'\}_{i=1}^n)$ for some $k\in \{1,\ldots,n\}$. Let $\hbar \in \mathbb{R}_{>0}\setminus\mathbb{Q}$. We define the unitary map $\mathbf{K}^\hbar_{\Gamma\mut{k}\Gamma'} : \mathscr{H}_{\Gamma'} \to \mathscr{H}_\Gamma$ as
$$
\mathbf{K}^\hbar_{\Gamma\mut{k}\Gamma'} := \mathbf{K}^{\sharp \hbar}_{\Gamma\mut{k}\Gamma'} \circ \mathbf{K}'_{\Gamma\mut{k}\Gamma'},
$$
where 
$$
\mathbf{K}^{\sharp \hbar}_{\Gamma\mut{k}\Gamma'} := \Phi^{\hbar_k}(\wh{\mathbf{x}}^\hbar_{\Gamma;k}) \left( \Phi^{\hbar_k}(\wh{\til{\mathbf{x}}}^\hbar_{\Gamma;k}) \right)^{-1}  ~ : ~ \mathscr{H}_\Gamma \to \mathscr{H}_\Gamma
$$
is the composition of two unitary operators obtained by applying the functional calculus in \S\ref{subsec:spectral_theorem} to the unique self-adjoint extensions of $\wh{\mathbf{x}}^\hbar_{\Gamma;k}$ and $\wh{\til{\mathbf{x}}}^\hbar_{\Gamma;k}$ respectively which are defined either by \eqref{eq:new_representation} or by \eqref{eq:old_representation} for the unitary functions $z\mapsto \Phi^{\hbar_k}(z)$ and $z\mapsto (\Phi^{\hbar_k}(z))^{-1}$, while
\begin{align}
\label{eq:bf_K_prime}
  \mathbf{K}'_{\Gamma\mut{k}\Gamma'} := \mathbf{S}_{(\mathbf{c}_{\Gamma\mut{k}\Gamma'}\,,\,\mathbf{0})} \circ \mathbf{I}_{\Gamma\mut{k}\Gamma'}: \mathscr{H}_{\Gamma'} \to \mathscr{H}_\Gamma,
\end{align}
where $\mathbf{I}_{\Gamma\mut{k}\Gamma'} : \mathscr{H}_{\Gamma'}\to \mathscr{H}_\Gamma$ is induced by identifying each $a_i'$ with $a_i$, i.e.
\begin{align}
\label{eq:bf_I}
(\mathbf{I}_{\Gamma\mut{k}\Gamma'} \, f)(a_1,\ldots,a_n) := f(a_1',\ldots,a_n'), \qquad \forall f \in \mathscr{H}_{\Gamma'} = L^2(\mathbb{R}^n, da_1'\, \cdots da_n'),
\end{align}
and $\mathbf{S}_{(\mathbf{c}_{\Gamma\mut{k}\Gamma'}\,, \,\mathbf{0})} : \mathscr{H}_\Gamma \to \mathscr{H}_\Gamma$ is defined as in \eqref{eq:def_of_S_c_t} of \S\ref{subsec:special_affine_shift_operators}, where $\mathbf{c}_{\Gamma\mut{k}\Gamma'} = (c_{ij})_{i,j\in \{1,\ldots,n\}} \in \mathrm{SL}_\pm(n,\mathbb{R})$ is given by
\begin{align}
\label{eq:c_Gamma_Gamma_prime_entries}
c_{ii} = 1, \quad \forall i \neq k, \qquad c_{kk} = -1, \qquad c_{ik} = [-\varepsilon_{ki}]_+, \quad \forall i\neq k, \qquad c_{ij}=0 \quad\mbox{otherwise}.
\end{align}
\end{definition}
So, as in eq.(73) of \cite{FG09}, one can view $\mathbf{K}'_{\Gamma\mut{k}\Gamma'} : \mathscr{H}_{\Gamma'} \to \mathscr{H}_\Gamma$ as induced by the map
\begin{align}
\label{eq:what_bf_K_prime_does}
a_i' \mapsto \left\{
\begin{array}{ll}
a_i & \mbox{if $i\neq k$}, \\
-a_k + \sum_{j=1}^n [-\varepsilon_{kj}]_+ \, a_j  & \mbox{if $i=k$.}
\end{array}
\right.
\end{align}

From Lem.\ref{lem:bf_S_preserves_D} it follows that:
\begin{lemma}
$\mathbf{K}'_{\Gamma\mut{k}\Gamma'}(D_{\Gamma'}) = D_\Gamma$. \qed
\end{lemma}

In order to verify that $\mathbf{K}^\hbar_{\Gamma\mut{k}\Gamma'}$ indeed intertwines the relevant actions, we begin with the conjugation action of $\mathbf{K}'_{\Gamma\mut{k}\Gamma'}$ on the operators $\wh{\mathbf{b}}^\hbar_{\Gamma';i}$, $\wh{\mathbf{x}}^\hbar_{\Gamma';i}$, $i=1,\ldots,n$. For this, we apply Lem.\ref{lem:bf_S_conjugation_on_bf_p_and_bf_q} to $\mathbf{c}=\mathbf{c}_{\Gamma\mut{k}\Gamma'}$. From \eqref{eq:c_Gamma_Gamma_prime_entries} one finds that $\mathbf{c}_{\Gamma\mut{k}\Gamma'}$ is its own inverse:
\begin{align}
\label{eq:bf_c_is_self-inverse}
(\mathbf{c}_{\Gamma\mut{k}\Gamma'})^{-1} = \mathbf{c}_{\Gamma\mut{k}\Gamma'}.  
\end{align}
\begin{lemma}
\label{lem:bf_K_prime_conjugation_on_bf_b}
Let $\Gamma,\Gamma',k,\hbar$ be as in Def.\ref{def:bf_K}. One has
\begin{align}
\label{eq:bf_K_prime_conjugation_on_bf_b_i}
\mathbf{K}'_{\Gamma\mut{k}\Gamma'} \,\, \wh{\mathbf{b}}^\hbar_{\Gamma';i} \,\, (\mathbf{K}'_{\Gamma\mut{k}\Gamma'})^{-1} & = \wh{\mathbf{b}}^\hbar_{\Gamma;i}, \qquad \forall i\neq k, \\
\label{eq:bf_K_prime_conjugation_on_bf_b_k}
\mathbf{K}'_{\Gamma\mut{k}\Gamma'} \,\, \wh{\mathbf{b}}^\hbar_{\Gamma';k} \,\, (\mathbf{K}'_{\Gamma\mut{k}\Gamma'})^{-1} & = - \wh{\mathbf{b}}^\hbar_{\Gamma;k} + \sum_{j=1}^n [-\varepsilon_{kj}]_+ \, \wh{\mathbf{b}}^\hbar_j,
\end{align}
as equalities on $D_\Gamma$, when the operators $\wh{\mathbf{b}}^\hbar_{\Gamma';i}$ and $\wh{\mathbf{b}}^\hbar_{\Gamma;i}$ on $D_\Gamma$ are taken to be given by the `old' formulas \eqref{eq:old_representation} or by the the `new' formulas \eqref{eq:new_representation}. If $\wh{\mathbf{b}}^\hbar_{\Gamma';i}$ and $\wh{\mathbf{b}}^\hbar_{\Gamma;i}$ are replaced by their unique self-adjoint extensions, each equality still holds on the domain of the self-adjoint extension of either side of the equation.
\end{lemma}

{\it Proof.} For each $i=1,\ldots,n$ observe
\begin{align*}
& \mathbf{K}'_{\Gamma\mut{k}\Gamma'} \, \wh{\mathbf{q}}^\hbar_{\Gamma';i} \, (\mathbf{K}_{\Gamma\mut{k}\Gamma'})^{-1}
= (\mathbf{S}_{(\mathbf{c}_{\Gamma\mut{k}\Gamma'}\,,\,\mathbf{0})} \circ \mathbf{I}_{\Gamma\mut{k}\Gamma'}) \left( \wh{\mathbf{q}}^\hbar_{\Gamma';i} \right) (\mathbf{I}_{\Gamma\mut{k}\Gamma'}^{-1} \circ \mathbf{S}_{(\mathbf{c}_{\Gamma\mut{k}\Gamma'}\,,\,\mathbf{0})}^{-1} ) \\
& = \mathbf{S}_{(\mathbf{c}_{\Gamma\mut{k}\Gamma'}\,,\,\mathbf{0})} \left( \wh{\mathbf{q}}^\hbar_{\Gamma;i} \right) \mathbf{S}_{(\mathbf{c}_{\Gamma\mut{k}\Gamma'}\,,\,\mathbf{0})}^{-1} \stackrel{\vee}{=} \left\{
  \begin{array}{ll}
    \wh{\mathbf{q}}^\hbar_{\Gamma;i} & i\neq k, \\
    - \wh{\mathbf{q}}^\hbar_{\Gamma;k} + \sum_{j=1}^n [-\varepsilon_{kj}]_+ \, \wh{\mathbf{q}}^\hbar_{\Gamma;j}, & i=k,
  \end{array}
\right.
\end{align*}
where the checked equality is from \eqref{eq:conjugation_of_bf_S_on_bf_q}, \eqref{eq:c_Gamma_Gamma_prime_entries}, and $\varepsilon_{kk}=0$. From \eqref{eq:new_representation} and \eqref{eq:old_representation} we get the desired result. \qed

\vs

We go on.
\begin{lemma}
\label{lem:bf_K_prime_conjugation_on_bf_x}
Let $\Gamma,\Gamma',k,\hbar$ be as in Def.\ref{def:bf_K}. One has
\begin{align}
\label{eq:bf_K_prime_conjugation_on_bf_x_i}
\mathbf{K}'_{\Gamma\mut{k}\Gamma'} \,\, \wh{\mathbf{x}}^\hbar_{\Gamma';i} \,\, (\mathbf{K}'_{\Gamma\mut{k}\Gamma'})^{-1} & = \wh{\mathbf{x}}^\hbar_{\Gamma;i} + [\varepsilon_{ik}]_+ \, \wh{\mathbf{x}}^\hbar_{\Gamma;k}, \qquad \forall i\neq k, \\
\label{eq:bf_K_prime_conjugation_on_bf_til_x_i}
\mathbf{K}'_{\Gamma\mut{k}\Gamma'} \,\, \wh{\til{\mathbf{x}}}^\hbar_{\Gamma';i} \,\, (\mathbf{K}'_{\Gamma\mut{k}\Gamma'})^{-1} & = \wh{\til{\mathbf{x}}}^\hbar_{\Gamma;i} + [\varepsilon_{ik}]_+ \, \wh{\til{\mathbf{x}}}^\hbar_{\Gamma;k}, \qquad \forall i\neq k, \\
\label{eq:bf_K_prime_conjugation_on_bf_x_k}
\mathbf{K}'_{\Gamma\mut{k}\Gamma'} \,\, \wh{\mathbf{x}}^\hbar_{\Gamma';k} \,\, (\mathbf{K}'_{\Gamma\mut{k}\Gamma'})^{-1} & = - \wh{\mathbf{x}}^\hbar_{\Gamma;k}, \qquad\quad
\mathbf{K}'_{\Gamma\mut{k}\Gamma'} \,\, \wh{\til{\mathbf{x}}}^\hbar_{\Gamma';k} \,\, (\mathbf{K}'_{\Gamma\mut{k}\Gamma'})^{-1} = - \wh{\til{\mathbf{x}}}^\hbar_{\Gamma;k},
\end{align}
as equalities on $D_\Gamma$, when the operators $\wh{\mathbf{x}}^\hbar_{\Gamma';i}$ and $\wh{\mathbf{x}}^\hbar_{\Gamma;i}$ are taken to be given by the `old' formulas \eqref{eq:old_representation}. These equalities do NOT hold when the operators are taken from the `new' formulas \eqref{eq:new_representation}. If $\wh{\mathbf{x}}^\hbar_{\Gamma';i}$ and $\wh{\mathbf{x}}^\hbar_{\Gamma;i}$ are replaced by their unique self-adjoint extensions, each equality still holds on the domain of the self-adjoint extension of either side of the equation.
\end{lemma} 

{\it Proof.} First, take the `old' representation \eqref{eq:old_representation}; then one writes for each $i=1,\ldots,n$,
\begin{align}
\label{eq:old_representation_revisited}
\left\{
{\renewcommand{\arraystretch}{1.4}
\begin{array}{ll}
(\wh{\mathbf{x}}^\hbar_{\Gamma;i})^{\rm old} = \frac{d_i^{-1}}{2} \wh{\mathbf{p}}^\hbar_{\Gamma;i} -\sum_{j=1}^n \varepsilon_{ij} \wh{\mathbf{q}}^\hbar_{\Gamma;j},
& (\wh{\til{\mathbf{x}}}^\hbar_{\Gamma;i})^{\rm old} = \frac{d_i^{-1}}{2} \wh{\mathbf{p}}^\hbar_{\Gamma;i} + \sum_{j=1}^n \varepsilon_{ij} \wh{\mathbf{q}}^\hbar_{\Gamma;j}, \\
(\wh{\mathbf{x}}^\hbar_{\Gamma';i})^{\rm old} = \frac{d_i^{-1}}{2} \wh{\mathbf{p}}^\hbar_{\Gamma';i} -\sum_{j=1}^n \varepsilon_{ij} \wh{\mathbf{q}}^\hbar_{\Gamma';j}, 
& (\wh{\til{\mathbf{x}}}^\hbar_{\Gamma';i})^{\rm old} = \frac{d_i^{-1}}{2} \wh{\mathbf{p}}^\hbar_{\Gamma';i} + \sum_{j=1}^n \varepsilon'_{ij} \wh{\mathbf{q}}^\hbar_{\Gamma';j},
\end{array}}
\right.
\end{align}
where $\wh{\mathbf{p}}^\hbar_{\Gamma;i}$ and $\wh{\mathbf{q}}^\hbar_{\Gamma;j}$ are defined on $D_{\Gamma} \subset \mathscr{H}_{\Gamma}$ by  \eqref{eq:Schrodinger_representation}, while $\wh{\mathbf{p}}^\hbar_{\Gamma';i}$ and $\wh{\mathbf{q}}^\hbar_{\Gamma';j}$ are defined on $D_{\Gamma'} \subset \mathscr{H}_{\Gamma'}$ by \eqref{eq:Schrodinger_representation} with $a_i$, $\frac{\partial}{\partial a_i}$, $d_i$ replaced by the primed versions $a_i'$, $\frac{\partial}{\partial a_i'}$, $d_i' = d_i$. First, observe
\begin{align*}
& \mathbf{K}'_{\Gamma\mut{k}\Gamma'} \,\, \left( \frac{d_k^{-1}}{2} \wh{\mathbf{p}}^\hbar_{\Gamma';k} \mp \sum_{j=1}^n \varepsilon'_{kj} \wh{\mathbf{q}}^\hbar_{\Gamma';j}\right) (\mathbf{K}'_{\Gamma\mut{k}\Gamma'})^{-1} 
= \mathbf{S}_{(\mathbf{c}_{\Gamma\mut{k}\Gamma'}\,,\,\mathbf{0})} \left( \frac{d_k^{-1}}{2} \wh{\mathbf{p}}^\hbar_{\Gamma;k} \mp \sum_{j=1}^n (\varepsilon'_{kj}) \wh{\mathbf{q}}^\hbar_{\Gamma;j}\right) \mathbf{S}_{(\mathbf{c}_{\Gamma\mut{k}\Gamma'}\,,\,\mathbf{0})}^{-1} \\ 
& \stackrel{\vee}{=} \frac{d_k^{-1}}{2} ( - \wh{\mathbf{p}}^\hbar_{\Gamma;k}) \mp \sum_{j\neq k} (-\varepsilon_{kj}) \wh{\mathbf{q}}^\hbar_{\Gamma;j} 
= - \left( \frac{d_k^{-1}}{2} \wh{\mathbf{p}}^\hbar_{\Gamma;k} \mp \sum_{j=1}^n \varepsilon_{kj} \wh{\mathbf{q}}^\hbar_{\Gamma;j} \right),
\end{align*}
where in the checked equality we used \eqref{eq:conjugation_of_bf_S_on_bf_p},  \eqref{eq:conjugation_of_bf_S_on_bf_q}, \eqref{eq:bf_c_is_self-inverse}, \eqref{eq:c_Gamma_Gamma_prime_entries}, $\varepsilon_{kk}=0$, and $\varepsilon'_{kj} = -\varepsilon_{kj}$. In view of \eqref{eq:old_representation_revisited}, we can translate this result into the following:
$$
\mathbf{K}'_{\Gamma\mut{k}\Gamma'} \,\, (\wh{\mathbf{x}}^\hbar_{\Gamma';k})^{\rm old} \,\, (\mathbf{K}'_{\Gamma\mut{k}\Gamma'})^{-1} = - (\wh{\mathbf{x}}^\hbar_{\Gamma;k})^{\rm old}, \qquad\quad
\mathbf{K}'_{\Gamma\mut{k}\Gamma'} \,\, (\wh{\til{\mathbf{x}}}^\hbar_{\Gamma';k})^{\rm old} \,\, (\mathbf{K}'_{\Gamma\mut{k}\Gamma'})^{-1} = - (\wh{\til{\mathbf{x}}}^\hbar_{\Gamma;k})^{\rm old},
$$
which is the desired eq.\eqref{eq:bf_K_prime_conjugation_on_bf_x_k} for the old representation.

\vs

For each $i\neq k$  observe
\begin{align}
\nonumber
& \mathbf{K}'_{\Gamma\mut{k}\Gamma'} \,\, \left( \frac{d_i^{-1}}{2} \wh{\mathbf{p}}^\hbar_{\Gamma';i} \mp \sum_{j=1}^n \varepsilon'_{ij} \wh{\mathbf{q}}^\hbar_{\Gamma';j} \right) \,\, (\mathbf{K}'_{\Gamma\mut{k}\Gamma'})^{-1}  = \mathbf{S}_{(\mathbf{c}_{\Gamma\mut{k}\Gamma'}\,,\,\mathbf{0})} \left( \frac{d_i^{-1}}{2} \wh{\mathbf{p}}^\hbar_{\Gamma;i} \mp \sum_{j=1}^n \varepsilon'_{ij} \wh{\mathbf{q}}^\hbar_{\Gamma;j} \right) \mathbf{S}_{(\mathbf{c}_{\Gamma\mut{k}\Gamma'}\,,\,\mathbf{0})}^{-1} \\
\label{eq:bf_K_prime_conjugation_eq1}
&
= \frac{d_i^{-1}}{2}\left(\wh{\mathbf{p}}^\hbar_{\Gamma;i} + [-\varepsilon_{ki}]_+ \, \wh{\mathbf{p}}^\hbar_{\Gamma;k}\right) \mp  \sum_{j\neq k} \varepsilon'_{ij} \, \wh{\mathbf{q}}^\hbar_{\Gamma;j}
\mp \varepsilon'_{ik} \left( - \wh{\mathbf{q}}^\hbar_{\Gamma;k} + \sum_{\ell\neq k} [-\varepsilon_{k\ell}]_+ \wh{\mathbf{q}}^\hbar_{\Gamma;\ell} \right),
\end{align}
using \eqref{eq:conjugation_of_bf_S_on_bf_p},  \eqref{eq:conjugation_of_bf_S_on_bf_q}, \eqref{eq:bf_c_is_self-inverse}, \eqref{eq:c_Gamma_Gamma_prime_entries}. Note that $[-\varepsilon_{ki}]_+ \frac{d_k}{d_i} \stackrel{\eqref{eq:wh_varepsilon}}{=} [-\wh{\varepsilon}_{ki} d_k]_+ = [\wh{\varepsilon}_{ik} d_k]_+ \stackrel{\eqref{eq:wh_varepsilon}}{=} [\varepsilon_{ik}]_+$. We use the matrix mutation formula \eqref{eq:varepsilon_prime_formula} $\varepsilon'_{ij} = \varepsilon_{ij} + \frac{1}{2}(|\varepsilon_{ik}|\varepsilon_{kj} + \varepsilon_{ik} |\varepsilon_{kj}|)$ and $\varepsilon'_{ik}=-\varepsilon_{ik}$. Hence
\begin{align*}
  \eqref{eq:bf_K_prime_conjugation_eq1}
= \left(\frac{d_i^{-1}}{2}\wh{\mathbf{p}}^\hbar_{\Gamma;i} \mp \sum_{j=1}^n  \varepsilon_{ij} \wh{\mathbf{q}}^\hbar_{\Gamma;j}\right)
+ \left(\frac{d_k^{-1}}{2} [\varepsilon_{ik}]_+ \wh{\mathbf{p}}^\hbar_{\Gamma;k} \mp \sum_{j\neq k} (\underbrace{\frac{|\varepsilon_{ik}|\varepsilon_{kj} + \varepsilon_{ik}|\varepsilon_{kj}|}{2} - \varepsilon_{ik} [-\varepsilon_{kj}]_+ }) \wh{\mathbf{q}}^\hbar_{\Gamma;j} \right)
\end{align*}
Putting $|\varepsilon_{ik}| = 2 [\varepsilon_{ik}]_+ - \varepsilon_{ik}$ and $|\varepsilon_{kj}| = 2[-\varepsilon_{kj}]_+ + \varepsilon_{kj}$ makes the underbraced part to
$$
\frac{ (2[\varepsilon_{ik}]_+-\varepsilon_{ik})\varepsilon_{kj} + \varepsilon_{ik}(2[-\varepsilon_{kj}]_+ + \varepsilon_{kj}) }{2} - \varepsilon_{ik} [-\varepsilon_{kj}]_+ = [\varepsilon_{ik}]_+ \varepsilon_{kj}.
$$
Thus, keeping in mind $\varepsilon_{kk}=0$, we get
$$
\eqref{eq:bf_K_prime_conjugation_eq1} = \left(\frac{d_i^{-1}}{2}\wh{\mathbf{p}}^\hbar_{\Gamma;i} \mp \sum_{j=1}^n  \varepsilon_{ij} \wh{\mathbf{q}}^\hbar_{\Gamma;j}\right)
+ [\varepsilon_{ik}]_+ \left( \frac{d_k^{-1}}{2} \wh{\mathbf{p}}^\hbar_{\Gamma;k} \mp \sum_{j=1}^n  \varepsilon_{kj} \wh{\mathbf{q}}^\hbar_{\Gamma;j}\right),
$$
which is the desired eq.\eqref{eq:bf_K_prime_conjugation_on_bf_x_i}, \eqref{eq:bf_K_prime_conjugation_on_bf_til_x_i} for the old representation.

\vs

On the other hand, let us try to use the new representation: in this case we have
$$
(\wh{\mathbf{x}}^\hbar_{\Gamma;i})^{\rm new} = d_i^{-1} \wh{\mathbf{p}}^\hbar_{\Gamma;i} - \sum_{j=1}^n  [\varepsilon_{ij}]_+ \wh{\mathbf{q}}^\hbar_{\Gamma;j}, \qquad (\wh{\mathbf{x}}^\hbar_{\Gamma';i})^{\rm new} = d_i^{-1} \wh{\mathbf{p}}^\hbar_{\Gamma';i} - \sum_{j=1}^n  [\varepsilon'_{kj}]_+ \wh{\mathbf{q}}^\hbar_{\Gamma';j}, \qquad \forall i=1,\ldots,n.
$$
We just check one case. Observe
\begin{align*}
\hspace{-6mm}
\mathbf{K}'_{\Gamma\mut{k}\Gamma'} (\wh{\mathbf{x}}^\hbar_{\Gamma';k})^{\rm new} (\mathbf{K}'_{\Gamma\mut{k}\Gamma'})^{-1} = \mathbf{S}_{(\mathbf{c}_{\Gamma\mut{k}\Gamma'}\,,\,\mathbf{0})} \left( d_k^{-1} \wh{\mathbf{p}}^\hbar_{\Gamma;k} - \sum_{j=1}^n [\varepsilon'_{kj}]_+ \wh{\mathbf{q}}^\hbar_{\Gamma;j} \right) \mathbf{S}_{(\mathbf{c}_{\Gamma\mut{k}\Gamma'}\,,\,\mathbf{0})}^{-1}
= - d_k^{-1} \wh{\mathbf{p}}^\hbar_{\Gamma;k} - \sum_{j\neq k} [-\varepsilon_{kj}]_+ \, \wh{\mathbf{q}}^\hbar_{\Gamma;j},
\end{align*}
where we used \eqref{eq:conjugation_of_bf_S_on_bf_p},  \eqref{eq:conjugation_of_bf_S_on_bf_q}, \eqref{eq:bf_c_is_self-inverse}, \eqref{eq:c_Gamma_Gamma_prime_entries}, together with $\varepsilon'_{kk}=0$. By inspection, one immediately finds that this does NOT equal $-(\wh{\mathbf{x}}^\hbar_{\Gamma;k})^{\rm new}$, unless $[-\varepsilon_{kj}]_+ = - [\varepsilon_{kj}]_+$, $\forall j\neq k$, which implies $\varepsilon_{kj}=0$, $\forall j\neq k$, which is not always the case. In particular, if we require $\mathbf{K}'_{\Gamma\mut{k}\Gamma'} (\wh{\mathbf{x}}^\hbar_{\Gamma';k})^{\rm new} (\mathbf{K}'_{\Gamma\mut{k}\Gamma'})^{-1}$ to equal $-(\wh{\mathbf{x}}^\hbar_{\Gamma;k})^{\rm new}$ for each $k\in \{1,\ldots,n\}$, then we would have to have $\varepsilon_{kj}=0$ for all $j,k \in \{1,\ldots,n\}$, which is undesirable. Eq.\eqref{eq:bf_K_prime_conjugation_on_bf_x_i} and \eqref{eq:bf_K_prime_conjugation_on_bf_til_x_i} for the new representation lead to more complicated equations for $\varepsilon_{ij}$'s, which do not always hold. The last statement of the lemma follows from Lem.\ref{lem:unitarily_equivalent_essentially_self-adjoint_operators}. \qed

\vs

So, here I point out one mistake in \cite{FG09}. The right identity in eq.(82) of Thm.5.6 in \cite[\S5.3]{FG09} says $\mathbf{K}' \wh{A} w = \wh{\mu'_k(A)} \mathbf{K}^\sharp w$ holds for all $w\in W_{\mathbf{i}'}$ and all $A\in \mathbf{L}$, where $\wh{A}$ and $\wh{\mu'_k(A)}$ stand for the operatos for $A$ and $\mu_k'(A)$ corresponding to their `new' representation, and $\mu_k'$ is as in Def.\ref{def:quantum_mutation_map}; in the proof they say that this is ``straightforward". First, the right hand side of this identity must read $\wh{\mu'_k(A)} \mathbf{K}' w$, which is an obvious minor typo. As mentioned earlier, their $W_{\mathbf{i}'}$ is a subspace of our $D_{\Gamma'}$. Notice that this asserted equation can be obtained from the exponentials of our (logarithmic) equations \eqref{eq:bf_K_prime_conjugation_on_bf_x_i}, \eqref{eq:bf_K_prime_conjugation_on_bf_til_x_i}, and \eqref{eq:bf_K_prime_conjugation_on_bf_x_k}. For example, the exponential of \eqref{eq:bf_K_prime_conjugation_on_bf_x_k} would yield
$$
\mathbf{K}'_{\Gamma\mut{k}\Gamma'}\, \pi^\hbar_{\Gamma'}(\mathbf{X}_{\Gamma';k}) (\mathbf{K}'_{\Gamma\mut{k}\Gamma'})^{-1} = \mathbf{K}'_{\Gamma\mut{k}\Gamma'} \, e^{\wh{\mathbf{x}}^\hbar_{\Gamma';k}} (\mathbf{K}'_{\Gamma\mut{k}\Gamma'})^{-1} \stackrel{\eqref{eq:bf_K_prime_conjugation_on_bf_x_k}}{=} e^{-\wh{\mathbf{x}}^\hbar_{\Gamma;k}} = \pi^\hbar_\Gamma(\mathbf{X}_{\Gamma;k}^{-1}) \stackrel{{\rm Def.}\ref{def:quantum_mutation_map}}{=} \pi^\hbar_\Gamma(\mu'_k(\mathbf{X}_{\Gamma';k})).
$$
However, recall from Lem.\ref{lem:bf_K_prime_conjugation_on_bf_x} that \eqref{eq:bf_K_prime_conjugation_on_bf_x_k} holds only for the `old' representation, but not for the `new' representation which is what is used in \cite{FG09}. Likewise, the exponential of \eqref{eq:bf_K_prime_conjugation_on_bf_x_i} would yield, for each $i\neq k$,
\begin{align*}
& \mathbf{K}'_{\Gamma\mut{k}\Gamma'}\, \pi^\hbar_{\Gamma'}(\mathbf{X}_{\Gamma';i}) (\mathbf{K}'_{\Gamma\mut{k}\Gamma'})^{-1} = \mathbf{K}'_{\Gamma\mut{k}\Gamma'} \, e^{\wh{\mathbf{x}}^\hbar_{\Gamma';i}} (\mathbf{K}'_{\Gamma\mut{k}\Gamma'})^{-1} \stackrel{\eqref{eq:bf_K_prime_conjugation_on_bf_x_i}}{=} e^{\wh{\mathbf{x}}^\hbar_{\Gamma;i} + [\varepsilon_{ik}]_+ \, \wh{\mathbf{x}}^\hbar_{\Gamma;k}} \\
& \stackrel{\vee}{=} q_k^{-\varepsilon_{ik}[\varepsilon_{ik}]_+} e^{\wh{\mathbf{x}}^\hbar_{\Gamma;i}} (e^{\wh{\mathbf{x}}^\hbar_{\Gamma;k}})^{[\varepsilon_{ik}]_+} 
= \pi^\hbar_\Gamma(q_k^{-\varepsilon_{ik}[\varepsilon_{ik}]_+} \mathbf{X}_i \mathbf{X}_k^{[\varepsilon_{ik}]_+}) \stackrel{{\rm Def.}\ref{def:quantum_mutation_map}}{=} \pi^\hbar_\Gamma (\mu'_k(\mathbf{X}_{\Gamma';i})),
\end{align*}
where $q_k = e^{\pi \sqrt{-1} \hbar_k} = e^{\pi\sqrt{-1}\hbar/d_k}$; the checked equality is from the BCH formula and $[\wh{\mathbf{x}}^\hbar_{\Gamma;i}, \wh{\mathbf{x}}^\hbar_{\Gamma;k}] \stackrel{\eqref{eq:Heisenberg_relations_for_D_q_Gamma}}{=} 2\pi \sqrt{-1} \hbar_k \varepsilon_{ik}\cdot\mathrm{id}$. Similarly for \eqref{eq:bf_K_prime_conjugation_on_bf_til_x_i}. But again, these hold only for the `old' representation, not the `new' one.

\vs

Notice that in a previous paper \cite{FG07} Fock and Goncharov used the `old' representation along with the operator $\mathbf{K}'_{\Gamma\mut{k}\Gamma'}$ given as in our \eqref{eq:bf_K_prime}. Then in \cite{FG09} they decided to use the `new' representation, mainly for aesthetic purpose in my opinion, but using the exactly same operator $\mathbf{K}'_{\Gamma\mut{k}\Gamma'}$. We just noted that $\mathbf{K}'_{\Gamma\mut{k}\Gamma'}$ intertwine the `old' operators correctly, but not the `new' ones. Even though the `old' operators and the `new' operators both satisfy the same Heisenberg algebra relations respectively, the `old' ones are kind of special with respect to the choice of the operator $\mathbf{K}'_{\Gamma\mut{k}\Gamma'}$ \eqref{eq:bf_K_prime}, which looks pretty canonical at first glance. What is special about the `old' operators? Why did Fock and Goncharov expect that their `new' representation would be intertwined correctly by the same operator $\mathbf{K}'_{\Gamma\mut{k}\Gamma'}$ as well? As for the first question, I think it might not have been a very good idea to consider the Hilbert space $\mathscr{H}_\Gamma$ as the space of functions in the logarithmic $\mcal{A}$-variables $a_1,\ldots,a_n$, as done by Fock and Goncharov \cite{FG09}. I mean, it seems that they are viewing these variables as the actual logarithms of the $\mcal{A}$-variables $A_1,\ldots,A_n$ which in turn are related to the $\mcal{X}$-variables by the formula \eqref{eq:p_Gamma_formula}. This may cause no problem in the classical setting, but in the quantum setting we should be careful. In this regard, I give one remark. Although in general there is no natural quantization of the $\mcal{A}$-coordinate functions $A_1,\ldots,A_n$, suppose that they are quantized to be some operators; it is natural to expect that these operators are related via a formula like \eqref{eq:p_Gamma_formula} to the operators quantizing the $\mcal{X}$-coordinate variables, which do exist. As the quantum $\mcal{X}$ operators do not mutually commute, these `imaginary' quantum $\mcal{A}$ operators must not mutually commute. However, $a_1,\ldots,a_n$ regarded as operators on $L^2(\mathbb{R}^n; da_1\wedge \cdots da_n) = \mathscr{H}_\Gamma$ \emph{do} mutually commute; therefore these can NOT be the `log' version of the above mentioned `imaginary' quantum $\mcal{A}$ operators. So we should NOT expect that the fact the operator $\mathbf{K}'_{\Gamma\mut{k}\Gamma'}$ is induced by the transformation formula \eqref{eq:what_bf_K_prime_does} would automatically guarantee that the conjugation action of $\mathbf{K}'_{\Gamma\mut{k}\Gamma'}$ on the quantum $\mcal{X}$ operators realizes the transformation $\mu'_k$ defined in Def.\ref{def:quantum_mutation_map} which looks like being induced from \eqref{eq:what_bf_K_prime_does} and \eqref{eq:p_Gamma_formula}. A possible answer to the second question is that the fact that $\mathbf{K}'_{\Gamma\mut{k}\Gamma'}$ is induced by \eqref{eq:what_bf_K_prime_does} \emph{does}, somewhat misleadingly as we just saw, guarantee the correct conjugation action on the $\wh{\mathbf{b}}^\hbar_{\Gamma';i}$ operators for both the old and the new representations, as we showed in Lem.\ref{lem:bf_K_prime_conjugation_on_bf_x}.

\vs

So the `new' representation used in \cite{FG09} must be retracted unless a corresponding correct operator $\mathbf{K}'_{\Gamma\mut{k}\Gamma'}$ is found. Hence, from now on, we stick to their `old' representation \eqref{eq:old_representation}; the operators $\wh{\mathbf{b}}^\hbar_{\Gamma;i}$ and $\wh{\mathbf{x}}^\hbar_{\Gamma;i}$ for each $\mathcal{D}$-seed $\Gamma$ and $i\in\{1,\ldots,n\}$ would mean $(\wh{\mathbf{b}}^\hbar_{\Gamma;i})^{\rm old}$ and $(\wh{\mathbf{x}}^\hbar_{\Gamma;i})^{\rm old}$ respectively.

\vs

First, one can notice that the generators $\mathbf{B}_{\Gamma';i}$, $i=1,\ldots,n$, as well as the `checked' generators $\mathbf{B}_{\Gamma';i}^\vee$, $\mathbf{X}_{\Gamma';i}^\vee$, $i=1,\ldots,n$, are intertwined by $\mathbf{K}'_{\Gamma\mut{k}\Gamma'}$ in the way we want them to be, with respect to $\mu_k' : \mathbb{D}^q_{\Gamma'} \to \mathbb{D}^q_\Gamma$ and $(\mu_k')^\vee : \mathbb{D}^{q^\vee}_{(\Gamma')^\vee} \to \mathbb{D}^{q^\vee}_{\Gamma^\vee}$; these two $\mu_k'$ are collectively denoted by $\mu_k'$ when applied on elements of $\mathbf{L}^\hbar_{\Gamma'}$, by abuse of notation. So we have
$$
\mathbf{K}'_{\Gamma\mut{k}\Gamma'} \, \pi^\hbar_{\Gamma'}(\mathbf{u}') = \pi^\hbar_{\Gamma}(\mu_k'(\mathbf{u}')) \, \mathbf{K}'_{\Gamma\mut{k}\Gamma'}, \qquad \forall \mathbf{u}' \in \mathbf{L}^\hbar_{\Gamma'}.
$$
One may take this equality as an equality of operators $D_{\Gamma'} \to D_\Gamma$, and then take unique self-adjoint extensions of $\pi^\hbar_{\Gamma'}(\mathbf{u}')$ and $\pi^\hbar_{\Gamma}(\mu_k'(\mathbf{u}'))$.

\vs

Let us now check the conjugation action by $\mathbf{K}^{\sharp \hbar}_{\Gamma\mut{k}\Gamma'}$; we argue this here only at a formal level, for a heuristic purpose. See \cite[Thm.5.6]{FG09} for a rigorous treatment.

\begin{lemma}
Let $\Gamma,\Gamma',k,\hbar$ be as in Def.\ref{def:bf_K}. For each $\mathbf{u} \in \mathbf{L}^\hbar_\Gamma$ one has
\begin{align*}
  \mathbf{K}^{\sharp \hbar}_{\Gamma\mut{k}\Gamma'} \, \pi^\hbar_\Gamma (\mathbf{u}) = \pi^\hbar_\Gamma(\mu^{\sharp \hbar}_k(\mathbf{u})) \, \mathbf{K}^{\sharp \hbar}_{\Gamma\mut{k}\Gamma'} \quad\mbox{on $D_\Gamma$},
\end{align*}
where $\mu^{\sharp \hbar}_k = \mu^{\sharp q}_k \otimes \mu^{\sharp q^\vee}_k$.
\end{lemma}

{\it Formal heuristic `proof'.} Let $f$ be a meromorphic function on the complex plane, and $A,B$ be self-adjoint operators, with
$$
[A,B] = \sqrt{-1} c \cdot \mathrm{id}
$$
for some real number $c$, where the equality holds on on some dense subspace, say; we assume that $f$ does not have a pole on $\mathbb{R}$ nor on $\mathbb{R}+\sqrt{-1} c$. To be more rigorous one can assume the Welyl-relations-version of this equality, namely $e^{\sqrt{-1} \alpha A} e^{\sqrt{-1} \beta B} = e^{- \alpha \beta \sqrt{-1} c} e^{\sqrt{-1} \beta B} e^{\sqrt{-1} \alpha A}$ for each $\alpha,\beta \in \mathbb{R}$, which is an equality of unitary operators. We then assert
$$
e^A f(B) e^{-A} = f(B + \sqrt{-1} c \cdot \mathrm{id}),
$$
as densely defined operators, at least formally.

\vs

Note
$$
[\wh{\mathbf{x}}^\hbar_{\Gamma;i}, \wh{\mathbf{x}}^\hbar_{\Gamma;k}] 
\stackrel{\eqref{eq:Heisenberg_relations_for_D_q_Gamma}}{=} 2\pi \sqrt{-1} \, \hbar_k \, \varepsilon_{ik} \cdot {\rm id}
$$

Observe
\begin{align}
\nonumber
  \mathbf{K}^{\sharp \hbar}_{\Gamma\mut{k}\Gamma'} \, \pi^\hbar_\Gamma(\mathbf{X}_{\Gamma;i}) \, (\mathbf{K}^{\sharp \hbar}_{\Gamma\mut{k}\Gamma'})^{-1}
& = \Phi^{\hbar_k}(\wh{\mathbf{x}}^\hbar_{\Gamma;k}) \, \ul{ \left( \Phi^{\hbar_k}(\wh{\til{\mathbf{x}}}^\hbar_{\Gamma;k} ) \right)^{-1} \, e^{\wh{\mathbf{x}}^\hbar_{\Gamma;i}} \, \Phi^{\hbar_k}(\wh{\til{\mathbf{x}}}^\hbar_{\Gamma;k} ) } \, \left( \Phi^{\hbar_k}(\wh{\mathbf{x}}^\hbar_{\Gamma;k}) \right)^{-1} \\
\nonumber
& = \Phi^{\hbar_k}(\wh{\mathbf{x}}^\hbar_{\Gamma;k}) \, e^{\wh{\mathbf{x}}^\hbar_{\Gamma;i}} \, \left( \Phi^{\hbar_k}(\wh{\mathbf{x}}^\hbar_{\Gamma;k}) \right)^{-1} \\
\label{eq:bf_K_sharp_conjugation_eq1}
& = e^{\wh{\mathbf{x}}^\hbar_{\Gamma;i}} \, \Phi^{\hbar_k}(\wh{\mathbf{x}}^\hbar_{\Gamma;k} - 2\pi\sqrt{-1} \hbar_k \, \varepsilon_{ik} \cdot \mathrm{id}) \, \left( \Phi^{\hbar_k}(\wh{\mathbf{x}}^\hbar_{\Gamma;k}) \right)^{-1}
\end{align}
From one of the two equations \eqref{eq:Phi_h_difference_equations} note that
\begin{align*}
\Phi^\hbar(z \pm 2\pi\sqrt{-1}\hbar) \, \left( \Phi^\hbar(z) \right)^{-1} = ( 1 + q^{\pm 1} \, e^z)^{\pm 1}
\end{align*}
holds as meromorphic functions on the complex plane, where the three $\pm$ symbols are either all $+$ at the same time or all $-$ at the same time. By an easy induction one obtains 
\begin{align}
\label{eq:QD_identity_M}
  \Phi^\hbar(z + 2\pi \sqrt{-1} \hbar M) \, \left( \Phi^\hbar(z) \right)^{-1} = \prod_{m=1}^{|M|} ( 1 + q^{(2m-1) {\rm sgn}(M)} \, e^z)^{{\rm sgn}(M)},
\end{align}
as an equality of meromorphic functions, for each integer $M$. Either side of \eqref{eq:QD_identity_M} does not have a pole on the real line, hence restricts to a real analytic function on $z\in \mathbb{R}$. Functional calculus in \S\ref{subsec:spectral_theorem} for the self-adjoint operator $\wh{\mathbf{x}}^\hbar_{\Gamma;k}$ applied to these two functions with $\hbar_k$ in place of $\hbar$ (thus $q_k$ in place of $q$) and $M = -\varepsilon_{ik}$ yield two identical densely defined operators on $\mathscr{H}_\Gamma$:
\begin{align}
\label{eq:QD_identity_M_operator}
\Phi^{\hbar_k} (\wh{\mathbf{x}}^\hbar_{\Gamma;k} - 2\pi \sqrt{-1} \hbar_k \, \varepsilon_{ik} \cdot \mathrm{id}) \, \left( \Phi^\hbar(\wh{\mathbf{x}}^\hbar_{\Gamma;k}) \right)^{-1} = \prod_{m=1}^{|\varepsilon_{ik}|} ( 1 + q_k^{(2m-1) {\rm sgn}(-\varepsilon_{ik})} \, e^{\wh{\mathbf{x}}^\hbar_{\Gamma;k}})^{{\rm sgn}(-\varepsilon_{ik})},
\end{align}
with their domains matching; note that the product order in the RHS nor that in the LHS does not matter, as the factors commute with one another in each side. So, from \eqref{eq:bf_K_sharp_conjugation_eq1} we get
\begin{align*}
  \mathbf{K}^{\sharp \hbar}_{\Gamma\mut{k}\Gamma'} \, e^{\wh{\mathbf{x}}^\hbar_{\Gamma;i}} \, (\mathbf{K}^{\sharp \hbar}_{\Gamma\mut{k}\Gamma'})^{-1}
= e^{\wh{\mathbf{x}}^\hbar_{\Gamma;i}} \, \prod_{m=1}^{|\varepsilon_{ik}|} ( 1 + q_k^{(2m-1) {\rm sgn}(-\varepsilon_{ik})} \, e^{\wh{\mathbf{x}}^\hbar_{\Gamma;k}})^{{\rm sgn}(-\varepsilon_{ik})}
\stackrel{\eqref{eq:quantum_mutation_sharp_of_X_i}}{=} \pi^\hbar_\Gamma(\mu^{\sharp q}_k(\mathbf{X}_{\Gamma;i})).
\end{align*}
One can likewise prove $\mathbf{K}^{\sharp \hbar}_{\Gamma\mut{k}\Gamma'} \, \pi^\hbar_\Gamma(\mathbf{X}_i^{-1}) \, (\mathbf{K}^{\sharp \hbar}_{\Gamma\mut{k}\Gamma'})^{-1} =\pi^\hbar_\Gamma(\mu^{\sharp q}_k(\mathbf{X}_{\Gamma;i}^{-1}))$.

\vs

Let us now investigate the conjugation action on $\pi^\hbar_\Gamma(\mathbf{B}_i^{\pm 1}) = e^{\pm\wh{\mathbf{b}}^\hbar_i}$. From \eqref{eq:Heisenberg_relations_for_D_q_Gamma} and \eqref{eq:tilde_operators}
\begin{align}
\label{eq:x_and_b_commutation}
  [\wh{\til{\mathbf{x}}}^\hbar_{\Gamma;k}, \, \wh{\mathbf{b}}^\hbar_{\Gamma;i}] = [\wh{\mathbf{x}}^\hbar_{\Gamma;k}, \, \wh{\mathbf{b}}^\hbar_{\Gamma;i}] = 2\pi \sqrt{-1} \hbar_k\, \delta_{k,i} \cdot \mathrm{id},
\end{align}
hence for each $i\neq k$, one observes
\begin{align}
\nonumber
  \mathbf{K}^{\sharp \hbar}_{\Gamma\mut{k}\Gamma'} \, \pi^\hbar_\Gamma(\mathbf{B}_{\Gamma;i}^{\pm 1}) \, (\mathbf{K}^{\sharp \hbar}_{\Gamma\mut{k}\Gamma'})^{-1}
& = \Phi^{\hbar_k}(\wh{\mathbf{x}}^\hbar_{\Gamma;k}) \, \left( \Phi^{\hbar_k}(\wh{\til{\mathbf{x}}}^\hbar_{\Gamma;k} ) \right)^{-1} \, e^{\pm \wh{\mathbf{b}}^\hbar_{\Gamma;i}} \, \Phi^{\hbar_k}(\wh{\til{\mathbf{x}}}^\hbar_{\Gamma;k} )  \, \left( \Phi^{\hbar_k}(\wh{\mathbf{x}}^\hbar_{\Gamma;k}) \right)^{-1} \\
\nonumber
& = e^{\pm \wh{\mathbf{b}}^\hbar_{\Gamma;i}} = \pi^\hbar_\Gamma(\mathbf{B}_{\Gamma;i}^{\pm 1}),
\end{align}
while
\begin{align*}
& \mathbf{K}^{\sharp \hbar}_{\Gamma\mut{k}\Gamma'} \, \pi^\hbar_\Gamma(\mathbf{B}_{\Gamma;k}) \, (\mathbf{K}^{\sharp \hbar}_{\Gamma\mut{k}\Gamma'})^{-1} \\
& = \Phi^{\hbar_k}(\wh{\mathbf{x}}^\hbar_{\Gamma;k}) \, \ul{ \left( \Phi^{\hbar_k}(\wh{\til{\mathbf{x}}}^\hbar_{\Gamma;k} ) \right)^{-1} \, e^{\wh{\mathbf{b}}^\hbar_{\Gamma;k}} } \, \Phi^{\hbar_k}(\wh{\til{\mathbf{x}}}^\hbar_{\Gamma;k} )  \, \left( \Phi^{\hbar_k}(\wh{\mathbf{x}}^\hbar_{\Gamma;k}) \right)^{-1} \\
& = \ul{ \Phi^{\hbar_k}(\wh{\mathbf{x}}^\hbar_{\Gamma;k}) \, e^{\wh{\mathbf{b}}^\hbar_{\Gamma;k}}  } \, \left( \Phi^{\hbar_k}(\wh{\til{\mathbf{x}}}^\hbar_{\Gamma;k} + 2\pi \sqrt{-1} \hbar_k \cdot \mathrm{id}) \right)^{-1} \Phi^{\hbar_k}(\wh{\til{\mathbf{x}}}^\hbar_{\Gamma;k}) \left( \Phi^{\hbar_k}(\wh{\mathbf{x}}^\hbar_{\Gamma;k}) \right)^{-1} \\
& = e^{\wh{\mathbf{b}}^\hbar_{\Gamma;k}} \, \Phi^{\hbar_k}(\wh{\mathbf{x}}^\hbar_{\Gamma;k} + 2\pi\sqrt{-1}\hbar_k \cdot \mathrm{id}) \, \ul{ \left( \Phi^{\hbar_k}(\wh{\til{\mathbf{x}}}^\hbar_{\Gamma;k} + 2\pi \sqrt{-1} \hbar_k \cdot \mathrm{id}) \right)^{-1} \Phi^{\hbar_k}(\wh{\til{\mathbf{x}}}^\hbar_{\Gamma;k}) } \,\, \ul{ \left( \Phi^{\hbar_k}(\wh{\mathbf{x}}^\hbar_{\Gamma;k}) \right)^{-1} } \\
& = e^{\wh{\mathbf{b}}^\hbar_{\Gamma;k}}\, \ul{ \Phi^{\hbar_k}(\wh{\mathbf{x}}^\hbar_{\Gamma;k} + 2\pi\sqrt{-1}\hbar_k \cdot \mathrm{id}) \left( \Phi^{\hbar_k}(\wh{\mathbf{x}}^\hbar_{\Gamma;k}) \right)^{-1} } \,\, \ul{ \left( \Phi^{\hbar_k}(\wh{\til{\mathbf{x}}}^\hbar_{\Gamma;k} + 2\pi \sqrt{-1} \hbar_k \cdot \mathrm{id}) \right)^{-1} \Phi^{\hbar_k}(\wh{\til{\mathbf{x}}}^\hbar_{\Gamma;k}) } \\
& = e^{\wh{\mathbf{b}}^\hbar_{\Gamma;k}} (1 + q_k \, e^{\wh{\mathbf{x}}^\hbar_{\Gamma;k}}) \, ( 1 + q_k \, e^{\wh{\til{\mathbf{x}}}^\hbar_{\Gamma;k}})^{-1}
= \pi^\hbar_\Gamma( \mu^{\sharp q}_k(\mathbf{B}_{\Gamma;k})),
\end{align*}
and similarly when $\mathbf{B}_{\Gamma;k}$ is replaced by $\mathbf{B}_{\Gamma;k}^{-1}$. 

\vs

What about the `checked' generators? Note that
$$
[\frac{\wh{\mathbf{x}}^\hbar_{\Gamma;i}}{\hbar_i}, \wh{\mathbf{x}}^\hbar_{\Gamma;k}] = 2\pi\sqrt{-1} \frac{\hbar_k}{\hbar_i} \varepsilon_{ik} \cdot \mathrm{id}
= 2\pi \sqrt{-1} \frac{d_i}{d_k} \varepsilon_{ik} \cdot \mathrm{id}
\stackrel{\eqref{eq:varepsilon_vee}}{=} 2\pi \sqrt{-1} \varepsilon_{ik}^\vee \cdot \mathrm{id}, \qquad
[\frac{\wh{\mathbf{x}}^\hbar_{\Gamma;i}}{\hbar_i}, \wh{\til{\mathbf{x}}}^\hbar_{\Gamma;k}]=0,
$$
hold on $D_\Gamma$, as well as their corresponding Weyl-relations.
\begin{align*}
\mathbf{K}^{\sharp \hbar}_{\Gamma\mut{k}\Gamma'} \, \pi^\hbar_\Gamma(\mathbf{X}_{\Gamma;i}^\vee) \, (\mathbf{K}^{\sharp \hbar}_{\Gamma\mut{k}\Gamma'})^{-1}
& \stackrel{\eqref{eq:pi_vee},\eqref{eq:pi_hbar}}{=} \Phi^{\hbar_k}(\wh{\mathbf{x}}^\hbar_{\Gamma;k}) \, \ul{ \left( \Phi^{\hbar_k}(\wh{\til{\mathbf{x}}}^\hbar_{\Gamma;k} ) \right)^{-1} \, e^{\wh{\mathbf{x}}^\hbar_{\Gamma;i}/\hbar_i} \, \Phi^{\hbar_k}(\wh{\til{\mathbf{x}}}^\hbar_{\Gamma;k} ) } \, \left( \Phi^{\hbar_k}(\wh{\mathbf{x}}^\hbar_{\Gamma;k}) \right)^{-1} \\
& = \Phi^{\hbar_k}(\wh{\mathbf{x}}^\hbar_{\Gamma;k}) \, e^{\wh{\mathbf{x}}^\hbar_{\Gamma;i}/\hbar_i} \, \left( \Phi^{\hbar_k}(\wh{\mathbf{x}}^\hbar_{\Gamma;k}) \right)^{-1} \\
& = e^{\wh{\mathbf{x}}^\hbar_{\Gamma;i}/\hbar_i} \, \underbrace{ \Phi^{\hbar_k}(\wh{\mathbf{x}}^\hbar_{\Gamma;k} - 2\pi\sqrt{-1} \varepsilon^\vee_{ik} \cdot \mathrm{id}) \, \left( \Phi^{\hbar_k}(\wh{\mathbf{x}}^\hbar_{\Gamma;k}) \right)^{-1} }.
\end{align*}
Recall that $\varepsilon^\vee_{ik}$ is an integer. From one of the two equations \eqref{eq:Phi_h_difference_equations} we have
\begin{align*}
\Phi^\hbar(z \pm 2\pi\sqrt{-1}) \, \left( \Phi^\hbar(z) \right)^{-1} = ( 1 + (q^\vee)^{\pm 1} \, e^{z/\hbar})^{\pm 1},
\end{align*}
where the three $\pm$ symbols are either all $+$ at the same time or all $-$ at the same time. By induction one obtains, for each integer $M$,
\begin{align}
\label{eq:QD_identity_M_2}
  \Phi^\hbar(z + 2\pi \sqrt{-1} M) \, \left( \Phi^\hbar(z) \right)^{-1} = \prod_{m=1}^{|M|} ( 1 + (q^\vee)^{(2m-1) {\rm sgn}(M)} \, e^{z/\hbar})^{{\rm sgn}(M)},
\end{align}
which is an analog of \eqref{eq:QD_identity_M}. By formally applying the functional calculus for the self-adjoint extension of $\wh{\mathbf{x}}^\hbar_{\Gamma;k}$ on both sides of \eqref{eq:QD_identity_M_2}, while keeping in mind $(q_k)^\vee := e^{\pi\sqrt{-1}/\hbar_k} = (e^{\pi\sqrt{-1}/\hbar})^{d_k} = (q^\vee)^{d_k} \stackrel{\eqref{eq:q_i_vee}}{=} q^\vee_k$, one can transform the underbraced part above, so that we have
$$
\mathbf{K}^{\sharp \hbar}_{\Gamma\mut{k}\Gamma'} \, \pi^\hbar_\Gamma(\mathbf{X}_{\Gamma;i}^\vee) \, (\mathbf{K}^{\sharp \hbar}_{\Gamma\mut{k}\Gamma'})^{-1}
= e^{\wh{\mathbf{x}}^\hbar_{\Gamma;i}/\hbar_i} \, \prod_{m=1}^{|\varepsilon_{ik}^\vee|} ( 1+(q_k^\vee)^{(2m-1){\rm sgn}(-\varepsilon_{ik}^\vee)} e^{\wh{\mathbf{x}}^\hbar_{\Gamma;k}/\hbar_k})^{{\rm sgn}(-\varepsilon_{ik}^\vee)},
$$
which equals $\pi^\hbar_\Gamma(\mu^{\sharp q^\vee}_k(\mathbf{X}_{\Gamma;i}^\vee))$, in view of \eqref{eq:quantum_mutation_sharp_of_X_i_vee}, \eqref{eq:pi_vee}, and \eqref{eq:pi_hbar}. 

\vs

Also, from \eqref{eq:x_and_b_commutation} one has$$
[\wh{\til{\mathbf{x}}}^\hbar_{\Gamma;k}, \, \frac{\wh{\mathbf{b}}^\hbar_{\Gamma;i}}{\hbar_i}] = [\wh{\mathbf{x}}^\hbar_{\Gamma;k}, \, \frac{\wh{\mathbf{b}}^\hbar_{\Gamma;i}}{\hbar_i}] = 2\pi \sqrt{-1} \frac{\hbar_k}{\hbar_i} \, \delta_{k,i} \cdot \mathrm{id} = 2\pi\sqrt{-1} \delta_{k,i} \cdot \mathrm{id}
$$
on $D_\Gamma$, as well as its Weyl-relations version. Hence, for each $i\neq k$ one observes  \\
$\mathbf{K}^{\sharp \hbar}_{\Gamma\mut{k}\Gamma'} \, \pi^\hbar_\Gamma((\mathbf{B}_{\Gamma;i}^\vee)^{\pm 1}) \, (\mathbf{K}^{\sharp \hbar}_{\Gamma\mut{k}\Gamma'})^{-1} = \pi^\hbar_\Gamma((\mathbf{B}_{\Gamma;i}^\vee)^{\pm 1})$, while
\begin{align*}
& \mathbf{K}^{\sharp \hbar}_{\Gamma\mut{k}\Gamma'} \, \pi^\hbar_\Gamma(\mathbf{B}^\vee_{\Gamma;k}) \, (\mathbf{K}^{\sharp \hbar}_{\Gamma\mut{k}\Gamma'})^{-1} \\
& = \Phi^{\hbar_k}(\wh{\mathbf{x}}^\hbar_{\Gamma;k}) \, \ul{ \left( \Phi^{\hbar_k}(\wh{\til{\mathbf{x}}}^\hbar_{\Gamma;k} ) \right)^{-1} \, e^{\wh{\mathbf{b}}^\hbar_{\Gamma;k}/\hbar_k} } \, \Phi^{\hbar_k}(\wh{\til{\mathbf{x}}}^\hbar_{\Gamma;k} )  \, \left( \Phi^{\hbar_k}(\wh{\mathbf{x}}^\hbar_{\Gamma;k}) \right)^{-1} \\
& = \ul{ \Phi^{\hbar_k}(\wh{\mathbf{x}}^\hbar_{\Gamma;k}) \, e^{\wh{\mathbf{b}}^\hbar_{\Gamma;k}/\hbar_k}  } \, \left( \Phi^{\hbar_k}(\wh{\til{\mathbf{x}}}^\hbar_{\Gamma;k} + 2\pi \sqrt{-1} \cdot \mathrm{id}) \right)^{-1} \Phi^{\hbar_k}(\wh{\til{\mathbf{x}}}^\hbar_{\Gamma;k}) \left( \Phi^{\hbar_k}(\wh{\mathbf{x}}^\hbar_{\Gamma;k}) \right)^{-1} \\
& = e^{\wh{\mathbf{b}}^\hbar_{\Gamma;k}/\hbar_k} \, \Phi^{\hbar_k}(\wh{\mathbf{x}}^\hbar_{\Gamma;k} + 2\pi\sqrt{-1} \cdot \mathrm{id}) \, \ul{ \left( \Phi^{\hbar_k}(\wh{\til{\mathbf{x}}}^\hbar_{\Gamma;k} + 2\pi \sqrt{-1} \cdot \mathrm{id}) \right)^{-1} \Phi^{\hbar_k}(\wh{\til{\mathbf{x}}}^\hbar_{\Gamma;k}) } \,\, \ul{ \left( \Phi^{\hbar_k}(\wh{\mathbf{x}}^\hbar_{\Gamma;k}) \right)^{-1} } \\
& = e^{\wh{\mathbf{b}}^\hbar_{\Gamma;k}/\hbar_k}\, \ul{ \Phi^{\hbar_k}(\wh{\mathbf{x}}^\hbar_{\Gamma;k} + 2\pi\sqrt{-1} \cdot \mathrm{id}) \left( \Phi^{\hbar_k}(\wh{\mathbf{x}}^\hbar_{\Gamma;k}) \right)^{-1} } \,\, \ul{ \left( \Phi^{\hbar_k}(\wh{\til{\mathbf{x}}}^\hbar_{\Gamma;k} + 2\pi \sqrt{-1} \cdot \mathrm{id}) \right)^{-1} \Phi^{\hbar_k}(\wh{\til{\mathbf{x}}}^\hbar_{\Gamma;k}) } \\
& = e^{\wh{\mathbf{b}}^\hbar_{\Gamma;k}/\hbar_k} (1 + q_k^\vee \, e^{\wh{\mathbf{x}}^\hbar_{\Gamma;k}/\hbar_k}) \, ( 1 + q_k^\vee \, e^{\wh{\til{\mathbf{x}}}^\hbar_{\Gamma;k}/\hbar_k})^{-1}
= \pi^\hbar_\Gamma( \mu^{\sharp q^\vee}_k(\mathbf{B}^\vee_{\Gamma;k})).
\end{align*}
\qed

\vs

Thus for each $\mathbf{u}'\in\mathbb{L}^q_{\Gamma'} \cup \mathbb{L}^{q^\vee}_{(\Gamma')^\vee} \subset \mathbf{L}^\hbar_{\Gamma'}$  one has
\begin{align*}
\mathbf{K}^\hbar_{\Gamma\mut{k}\Gamma'} \, \pi^\hbar_{\Gamma'}(\mathbf{u}') \, (\mathbf{K}^\hbar_{\Gamma\mut{k}\Gamma'})^{-1}
& = \mathbf{K}^{\sharp \hbar}_{\Gamma\mut{k}\Gamma'} \, \mathbf{K}'_{\Gamma\mut{k}\Gamma'} \, \pi^\hbar_{\Gamma'}(\mathbf{u}') \, (\mathbf{K}'_{\Gamma\mut{k}\Gamma'})^{-1} \, (\mathbf{K}^{\sharp \hbar}_{\Gamma\mut{k}\Gamma'})^{-1} \\
& = \mathbf{K}^{\sharp \hbar}_{\Gamma\mut{k}\Gamma'} \, \pi^\hbar_\Gamma( \mu_k'(\mathbf{u}')) \, (\mathbf{K}^{\sharp \hbar}_{\Gamma\mut{k}\Gamma'})^{-1} \\
& = \left\{ \begin{array}{ll}
\pi^\hbar_\Gamma(\mu_k^{\sharp q} (\mu_k'(\mathbf{u}')))  & \mbox{if $\mathbf{u'}\in \mathbb{L}^q_{\Gamma'}$} \\
\pi^\hbar_\Gamma(\mu^{\sharp q^\vee}_k(\mu'_k(\mathbf{u}'))) & \mbox{if $\mathbf{u'}\in \mathbb{L}^{q^\vee}_{(\Gamma')^\vee}$}
\end{array} \right. \\
& = \left\{ \begin{array}{ll}
\pi^\hbar_\Gamma(\mu_k^q (\mathbf{u}'))  & \mbox{if $\mathbf{u'}\in \mathbb{L}^q_{\Gamma'}$} \\
\pi^\hbar_\Gamma(\mu^{q^\vee}_k(\mathbf{u}')) & \mbox{if $\mathbf{u'}\in \mathbb{L}^{q^\vee}_{(\Gamma')^\vee}$}
\end{array} \right. \\
& = \pi^\hbar_\Gamma(\mu^\hbar_k(\mathbf{u}')),
\end{align*}
where $\mu^\hbar_k$ denotes the restriction of
$$
\mu^\hbar_k := \mu^q_k \otimes \mu^{q^\vee}_k : \mathbb{D}^q_{\Gamma'} \otimes \mathbb{D}^{q^\vee}_{(\Gamma')^\vee} \to \mathbb{D}^q_\Gamma \otimes \mathbb{D}^{q^\vee}_{\Gamma^\vee}
$$
to a map $\mathbf{L}^\hbar_{\Gamma'} \to \mathbf{L}^\hbar_\Gamma$, just like \eqref{eq:eta_hbar}. Thus we have
$$
\mathbf{K}^\hbar_{\Gamma\mut{k}\Gamma'} \, \pi^\hbar_{\Gamma'}(\mathbf{u}') \, (\mathbf{K}^\hbar_{\Gamma\mut{k}\Gamma'})^{-1} = \pi^\hbar_\Gamma(\mu^\hbar_k(\mathbf{u}')), \qquad \forall \mathbf{u}' \in \mathbf{L}^\hbar_{\Gamma'},
$$
as desired. 

\vs

Notice that I provided only a heuristic argument of the above result. Although this is a crucial anayltic result to prove rigorously, I do not give a full proof here, and just refer the readers to \cite[\S5]{FG09}. There, the intertwining property is checked for certain `generating' elemnts of $\mathbf{L}^\hbar_{\Gamma'}$, applied to the elements of their space $W_{{\bf i}'}$ which is dense in their (feed) Schwartz space $\mathscr{S}_{{\bf i}'}$ in the Fr\'echet topology. In the next section, we shall concentrate on the `consistency' relations of the unitary intertwiners $\mathbf{K}^\hbar_{\Gamma\mut{k}\Gamma'}$; for those I will provide rigorous proofs.

\section{Computation of operator identities}
\label{sec:operator_identities}

In the present section, we always assume that
$$
\hbar\in \mathbb{R}_{>0}\setminus\mathbb{Q}.
$$

\subsection{Operator identities as consequences of irreducibility}

In \eqref{eq:bf_K_P_sigma} and Def.\ref{def:bf_K} of \S\ref{sec:representations_of_algebraic_quantum_cluster_varieties} we assigned unitary intertwiner maps $\mathbf{K}^\hbar(P_\sigma) = \mathbf{P}_\sigma$ and $\mathbf{K}^\hbar(\mu_k)$ to elementary morphisms $P_\sigma$ and $\mu_k$ of the saturated cluster modular groupoid $\wh{\mcal{G}}^\mcal{D}_\mathscr{C}$ (Def.\ref{eq:saturated-modular_groupoid}). In the present subsection we discuss the relations that they must satisfy: namely, the ones in Lemmas \ref{lem:involution_identity}, \ref{lem:permutation_identities}, and \ref{lem:h_plus_2-gon_relations}, each written in the reverse order, up to a constant. The relations in Lem.\ref{lem:permutation_identities}, i.e. the permutation identities, are easily checked to hold exactly, without multiplicative constants:
$$
\mathbf{P}_\gamma \, \mathbf{P}_\sigma = \mathbf{P}_{\sigma \circ \gamma}, \qquad
\mathbf{P}_{\sigma^{-1}} \, \mathbf{K}^\hbar(\mu_k) \, \mathbf{P}_\sigma = \mathbf{K}^\hbar(\mu_{\sigma(k)}), \qquad
\mathbf{P}_{\mathrm{id}} = \mathrm{id}.
$$
The way how Fock and Goncharov proves other relations in \cite{FG09} is somewhat indirect. Suppose $\mathbf{m}$ is a trivial morphism in $\wh{\mcal{G}}^\mcal{D}_\mathscr{C}$, i.e. a morphism from an object $\Gamma$ to itself, written as a particular sequence of mutations and seed automorphisms. Denote by $\mathbf{K}^\hbar(\mathbf{m}) : \mathscr{H}_\Gamma \to \mathscr{H}_\Gamma$ the composition of the reversed sequence of intertwiners corresponding to these mutations and seed automorphisms. 
From the intertwining properties of $\mathbf{K}^\hbar(\mu_k)$ and $\mathbf{K}^\hbar(P_\sigma)$ we have the following intertwining equality
$$
\mathbf{K}^\hbar(\mathbf{m}) \, \pi^\hbar_\Gamma(\mathbf{u})\, v = \pi^\hbar_\Gamma((\eta^\hbar(\mathbf{m}))(\mathbf{u}) ) \, \mathbf{K}^\hbar(\mathbf{m}) \, v, \qquad \forall \mathbf{u} \in \mathbf{L}^\hbar_\Gamma, \quad \forall v \in \mathscr{S}^\hbar_\Gamma,
$$
where $\eta^\hbar(\mathbf{m}) = \eta^q(\mathbf{m}) \otimes \eta^{q^\vee}(\mathbf{m}^\vee)$. As $\mathbf{m}$ is a trivial morphism, Lem.\ref{lem:eta_q_is_well-defined} tells us that both $\eta^q(\mathbf{m})$ and $\eta^{q^\vee}(\mathbf{m})$ are identity maps of algebras. Hence we have
$$
\mathbf{K}^\hbar(\mathbf{m}) \, \pi^\hbar_\Gamma(\mathbf{u}) = \pi^\hbar_\Gamma(\mathbf{u})\, \mathbf{K}^\hbar(\mathbf{m}), \qquad \forall \mathbf{u}\in \mathbf{L}^\hbar_\Gamma,
$$
when applied to elements of $v\in \mathscr{S}_\Gamma^\hbar$. Fock and Goncharov \cite{FG09} established that the representation $\pi^\hbar_\Gamma$ of $\mathbf{L}^\hbar_\Gamma$ (on $\mathscr{S}_\Gamma^\hbar$) is `strongly irreducible', in the sense that any bounded operator $\mathscr{S}^\hbar_\Gamma\to \mathscr{S}^\hbar_\Gamma$ commuting with $\pi^\hbar_\Gamma(\mathbf{u})$ for all $\mathbf{u}\in \mathbf{L}^\hbar_\Gamma$ is a scalar operator; they described such a situation as saying that the algebra $\mathbf{L}^\hbar_\Gamma$ is `big enough'. Anyhow, they proved:
\begin{theorem}[\cite{FG09}]
\label{thm:FG_constant}
For any trivial morphism $\mathbf{m}$ in $\wh{\mcal{G}}^\mcal{D}_\mathscr{C}$ from $\Gamma$ to itself written as a sequence of mutations and seed automorphisms, if $\mathbf{K}^\hbar(\mathbf{m})$ denotes the composition of the reversed sequence of intertwiners corresponding to the mutations and seed automorphisms constituting the sequence $\mathbf{m}$, then
$$
\mathbf{K}^\hbar(\mathbf{m}) = c_{\mathbf{m}} \cdot \mathrm{id}_{\mathscr{H}_\Gamma},
$$
for some constant $c_{\mathbf{m}} \in \mathrm{U}(1) \subset \mathbb{C}^\times$. \qed
\end{theorem}
\begin{corollary}
The constructed intertwiners $\mathbf{K}^\hbar(\mu_k)$ and $\mathbf{K}^\hbar(P_\sigma)$ for the elementary cluster transformations induce a well-defined projective functor $\mathbf{K}^\hbar$ \eqref{eq:representation_functor}. \qed
\end{corollary}
The constant $c_{\mathbf{m}}$ is denoted by $\lambda$ in \cite[Thm.5.5]{FG09}, and is not precisely determined as an explicit number depending on $\hbar$. One of the purposes of the present paper is to show by computation that this constant $c_{\mathbf{m}}$ is always $1$. For this, we write the result of the above theorem more explicitly in terms of the relations in Lemmas \ref{lem:involution_identity} and \ref{lem:h_plus_2-gon_relations}, like in \cite[Thm.5.5]{FG09} but even more concretely. I separate this into several cases.

\begin{proposition}[rank $1$ identity; `twice-flip is identity'; \cite{FG09}]
\label{prop:rank_1_constant}
Let $\Gamma=(\varepsilon,d,*)$ be a $\mathcal{D}$-seed, $k \in \{1,\ldots,n\}$, $\Gamma' = \mu_k(\Gamma) = (\varepsilon',d',*')$. Then $\Gamma = \mu_k(\Gamma')$, and there exists a constant $c_{A_1} \in \mathrm{U}(1)\subset \mathbb{C}^\times$ such that
\begin{align}
\label{eq:twice_flip_to_prove}
\mathbf{K}^\hbar_{\Gamma\mut{k}\Gamma'} \, \mathbf{K}^\hbar_{\Gamma'\mut{k}\Gamma}  = c_{A_1} \cdot \mathrm{id}_{\mathscr{H}_\Gamma}. \qquad \qed
\end{align}
\end{proposition}

Using the above Proposition, the $(h+2)$-gon relations can be made into the following forms.

\begin{proposition}[$A_1 \times A_1$ identity; `commuting identity'; \cite{FG09}]
\label{prop:FG_A1_times_A1_identity}
Let $\Gamma^{(0)} = (\varepsilon^{(0)}, d^{(0)}, *^{(0)})$ be a $\mathcal{D}$-seed, and assume that some two distinct $i,j\in \{1,\ldots,n\}$ satisfy
\begin{align}
\label{eq:A1_times_A1_condition}
\varepsilon^{(0)}_{ij} = 0 = \varepsilon^{(0)}_{ji}.
\end{align}
Let
$$
\Gamma^{(1)} := \mu_i (\Gamma^{(0)}), \quad \Gamma^{(2)} := \mu_j(\Gamma^{(1)}), \quad\mbox{and}\quad \Gamma^{(3)}:=\mu_j(\Gamma^{(0)}).
$$
Then $\Gamma^{(2)} = \mu_i(\Gamma^{(3)})$, and there exists a constant $c_{A_1\times A_1} \in \mathrm{U}(1) \subset\mathbb{C}^\times$ such that
\begin{align}
\label{eq:A1_times_A1_to_prove}
\mathbf{K}^\hbar_{\Gamma^{(0)}\mut{i}\Gamma^{(1)}} \, \mathbf{K}^\hbar_{\Gamma^{(1)}\mut{j}\Gamma^{(2)}} = c_{A_1 \times A_1} \cdot \mathbf{K}^\hbar_{\Gamma^{(0)}\mut{j}\Gamma^{(3)}}, \, \mathbf{K}^\hbar_{\Gamma^{(3)}\mut{i}\Gamma^{(2)}}. \qquad \qed
\end{align}
\end{proposition}

\begin{proposition}[$A_2$ identity; `pentagon identity'; \cite{FG09}]
\label{prop:FG_A2_identity}
Let $\Gamma^{(0)} = (\varepsilon^{(0)}, d^{(0)}, *^{(0)})$ be a $\mathcal{D}$-seed, and assume that some $i,j\in \{1,\ldots,n\}$ satisfy
$$
\left\{ {\renewcommand{\arraystretch}{1.3} \begin{array}{l}
\varepsilon^{(0)}_{ij} = 1, \\
\varepsilon^{(0)}_{ji} = -1,
\end{array}}\right. \qquad\mbox{or}\qquad
\left\{ {\renewcommand{\arraystretch}{1.3} \begin{array}{l}
\varepsilon^{(0)}_{ij} = -1, \\
\varepsilon^{(0)}_{ji} = 1.
\end{array}} \right.
$$
Let
$$
\Gamma^{(1)} := \mu_i (\Gamma^{(0)}), \quad \Gamma^{(2)} := \mu_j(\Gamma^{(1)}),\quad \Gamma^{(3)}:=\mu_i(\Gamma^{(2)}), \quad\mbox{and}\quad \Gamma^{(4)}:=\mu_j(\Gamma^{(0)}), \quad \Gamma^{(5)}:=\mu_i(\Gamma^{(4)}).
$$
Then $\Gamma^{(3)} = P_{(i\, j)}(\Gamma^{(5)})$, and there exists a constant $c_{A_2} \in \mathrm{U}(1) \subset \mathbb{C}^\times$ such that 
\begin{align}
\label{eq:A2_to_prove}
\mathbf{K}^\hbar_{\Gamma^{(0)}\mut{i}\Gamma^{(1)}} \, \mathbf{K}^\hbar_{\Gamma^{(1)}\mut{j}\Gamma^{(2)}} \, \mathbf{K}^\hbar_{\Gamma^{(2)}\mut{i}\Gamma^{(3)}} = c_{A_2} \cdot \mathbf{K}^\hbar_{\Gamma^{(0)}\mut{j}\Gamma^{(4)}} \, \mathbf{K}^\hbar_{\Gamma^{(4)}\mut{i}\Gamma^{(5)}} \, \mathbf{P}_{(i\, j)},
\end{align}
as operators from $\mathscr{H}_{\Gamma^{(3)}}$ to $\mathscr{H}_{\Gamma^{(0)}}$; here $\mathbf{P}_{(i\, j)}$ is regarded as an operator from $\mathscr{H}_{P_{(i\, j)}(\Gamma^{(5)})} \equiv \mathscr{H}_{\Gamma^{(3)}}$ to $\mathscr{H}_{\Gamma^{(5)}}$. \qed
\end{proposition}

\begin{proposition}[$B_2$ identity; `hexagon identity'; \cite{FG09}]
\label{prop:FG_B2_identity}
Let $\Gamma^{(0)} = (\varepsilon^{(0)}, d^{(0)}, *^{(0)})$ be a $\mathcal{D}$-seed, and assume that some $i,j\in \{1,\ldots,n\}$ satisfy
$$
\varepsilon^{(0)}_{ij} = - 2 \, \varepsilon^{(0)}_{ji}, \qquad |\varepsilon^{(0)}_{ij}| = 2.
$$
Let
\begin{align*}
& \Gamma^{(1)} := \mu_i (\Gamma^{(0)}), \quad \Gamma^{(2)} := \mu_j(\Gamma^{(1)}),\quad \Gamma^{(3)}:=\mu_i(\Gamma^{(2)}), \quad\mbox{and}\quad\\
& \Gamma^{(4)}:=\mu_j(\Gamma^{(0)}), \quad \Gamma^{(5)}:=\mu_i(\Gamma^{(4)}).
\end{align*}
Then $\Gamma^{(3)} = \mu_j(\Gamma^{(5)})$, and there exists a constant $c_{B_2} \in \mathrm{U}(1) \subset \mathbb{C}^\times$ such that 
\begin{align}
\label{eq:B2_to_prove}
\mathbf{K}^\hbar_{\Gamma^{(0)}\mut{i}\Gamma^{(1)}} \, \mathbf{K}^\hbar_{\Gamma^{(1)}\mut{j}\Gamma^{(2)}} \, \mathbf{K}^\hbar_{\Gamma^{(2)}\mut{i}\Gamma^{(3)}} = c_{B_2} \cdot \mathbf{K}^\hbar_{\Gamma^{(0)}\mut{j}\Gamma^{(4)}} \, \mathbf{K}^\hbar_{\Gamma^{(4)}\mut{i}\Gamma^{(5)}} \, \mathbf{K}^\hbar_{\Gamma^{(5)}\mut{j}\Gamma^{(3)}}. \qquad \qed
\end{align}
\end{proposition}

\begin{proposition}[$G_2$ identity; `octagon identity'; \cite{FG09}]
\label{prop:FG_G2_identity}
Let $\Gamma^{(0)} = (\varepsilon^{(0)}, d^{(0)}, *^{(0)})$ be a $\mathcal{D}$-seed, and assume that some $i,j\in \{1,\ldots,n\}$ satisfy
$$
\varepsilon^{(0)}_{ij} = - 3 \, \varepsilon^{(0)}_{ji}, \qquad |\varepsilon^{(0)}_{ij}| = 3.
$$
Let
\begin{align*}
& \Gamma^{(1)} := \mu_i (\Gamma^{(0)}), \quad \Gamma^{(2)} := \mu_j(\Gamma^{(1)}),\quad \Gamma^{(3)}:=\mu_i(\Gamma^{(2)}), \quad \Gamma^{(4)} := \mu_j(\Gamma^{(3)}), \quad\mbox{and}\quad\\
& \Gamma^{(5)}:=\mu_j(\Gamma^{(0)}), \quad \Gamma^{(6)}:=\mu_i(\Gamma^{(5)}), \quad \Gamma^{(7)} := \mu_j(\Gamma^{(6)}).
\end{align*}
Then $\Gamma^{(4)} = \mu_i(\Gamma^{(7)})$, and there exists a constant $c_{G_2} \in \mathrm{U}(1) \subset \mathbb{C}^\times$ such that 
\begin{align}
\label{eq:G2_to_prove}
{\renewcommand{\arraystretch}{1.4}\begin{array}{l}
\mathbf{K}^\hbar_{\Gamma^{(0)}\mut{i}\Gamma^{(1)}} \, \mathbf{K}^\hbar_{\Gamma^{(1)}\mut{j}\Gamma^{(2)}} \, \mathbf{K}^\hbar_{\Gamma^{(2)}\mut{i}\Gamma^{(3)}} \, \mathbf{K}^\hbar_{\Gamma^{(3)}\mut{j}\Gamma^{(4)}} \\
\quad = c_{G_2} \cdot \mathbf{K}^\hbar_{\Gamma^{(0)}\mut{j}\Gamma^{(5)}} \, \mathbf{K}^\hbar_{\Gamma^{(5)}\mut{i}\Gamma^{(6)}} \, \mathbf{K}^\hbar_{\Gamma^{(6)}\mut{j}\Gamma^{(7)}} \, \mathbf{K}^\hbar_{\Gamma^{(7)}\mut{i}\Gamma^{(4)}}. 
\end{array}} \qquad \qed
\end{align}
\end{proposition}

\begin{remark}
The terminologies `hexagon' and `octagon' which I used for the last two equalities as analogs of `pentagon' are probably not standard.
\end{remark}
\begin{remark}
These constants depend only on the underlying `feed' data $(\varepsilon,d)$, not on specific seeds $\Gamma$. So Prop.\ref{prop:rank_1_constant}--\ref{prop:FG_G2_identity} are slightly stronger than Thm.\ref{thm:FG_constant}. 
\end{remark}
One cannot determine these constants $c_{A_1}, c_{A_1 \times A_1}, c_{A_2}, c_{B_2}, c_{G_2}$ by means of the indirect argument of Fock and Goncharov. So we will make explicit computations about the unitary intertwining operators, and use some known operator identities involving the quantum dilogarithm function $\Phi^\hbar$.

\subsection{Known operator identities}
\label{subsec:known_operator_identities}

In the current subsection, we establish some operator identities, including the pentagon equation for the non-compact quantum dilogarithm $\Phi^\hbar$ promised in \S\ref{subsec:non-compact_QD}.

\vs

Recall the identity $\Phi^\hbar(z) \Phi^\hbar(-z) = c_\hbar \cdot \exp( \frac{z^2}{4\pi\sqrt{-1}\hbar})$ \eqref{eq:QD_identity_quadratic}, where $|c_\hbar|=1$. As $z\mapsto \Phi^\hbar(z)$, $z\mapsto \Phi^\hbar(-z)$, and $z\mapsto \exp\left( \frac{z^2}{4\pi\sqrt{-1}\hbar} \right)$ are all unitary functions in the sense of Lem.\ref{lem:unitary_operators_from_functional_calculus}, the application of the functional calculus in \S\ref{subsec:spectral_theorem} for a self-adjoint operator on these functions yield unitary operators.
\begin{lemma}
Let $T$ be a (densely-defined) self-adjoint operator on a Hilbert space $V$. Its functional calculus applied to the three unitary functions just mentioned yield three unitary operators $\Phi^\hbar(T)$, $\Phi^\hbar(-T)$, and $\exp(\frac{T^2}{4\pi\sqrt{-1}\hbar})$, satisfying
\begin{align}
\label{eq:QD_identity_quadratic}
\Phi^\hbar(T) \, \Phi^\hbar(-T) = c_\hbar \, \exp\left( \frac{T^2}{4\pi \sqrt{-1} \hbar} \right),
\end{align}
where $c_\hbar$ is as in \eqref{eq:c_hbar}. \qed
\end{lemma}
Proof of \eqref{eq:QD_identity_quadratic} is immediate from the construction of the functional calculus in \S\ref{subsec:spectral_theorem} using the unique spectral resolution of $T$.

\begin{lemma}
[linear combination of Heisenberg operators]
\label{lem:linear_combination_of_Heisenberg_operators}
Let $P$ and $Q$ be densely defined self-adjoint operators on a separable Hilbert space $V$ satisfying the Weyl relations \eqref{eq:Weyl_relations_for_P_and_Q}
$$
e^{\sqrt{-1} \alpha P} \, e^{\sqrt{-1} \beta Q} = e^{-2\pi\sqrt{-1} \hbar \, \alpha\beta} \, e^{\sqrt{-1}\beta Q} \, e^{\sqrt{-1} \alpha P}, \qquad \forall \alpha,\beta\in\mathbb{R},
$$
corresponding to the Heisenberg relation
$$
[P,Q] = 2\pi\sqrt{-1} \, \hbar \cdot \mathrm{id}
$$
which makes sense on a dense subspace of $V$. Then there exists a dense subspace $D$ of $V$ such that
\begin{itemize}
\item[\rm 1)] each of the restrictions $\wh{\mathbf{p}}^\hbar := P \restriction D$ and $\wh{\mathbf{q}}^\hbar := Q \restriction D$ preserves $D$ and is essentially self-adjoint,

\item[\rm 2)] for each real numbers $\alpha_0, \beta_0$, the operator $\alpha_0 \wh{\mathbf{p}}^\hbar + \beta_0 \wh{\mathbf{q}}^\hbar$ on $D$ is essentially self-adjoint.
\end{itemize}
Denote by $\alpha_0 P+\beta_0 Q$ the unique self-adjoint extension of $\alpha_0 \wh{\mathbf{p}}^\hbar + \beta_0 \wh{\mathbf{q}}^\hbar$.
\end{lemma}

{\it Proof.} As a corollary of Thm.\ref{thm:SvN} we can assume that $V = L^2(\mathbb{R},dx)$, and $(e^{\sqrt{-1}\alpha P}f)(x) = f(x-2\pi\hbar \,\alpha)$ and $(e^{\sqrt{-} \beta Q}f)(x) = e^{\sqrt{-1}\beta x} \, f(x)$ for all $f\in V = L^2(\mathbb{R})$. One could also use the dense subspace $D$ \eqref{eq:D} which is preserved by $P$ and $Q$, whose restrictions on $D$ are the operators $\wh{\mathbf{p}}^\hbar = P\restriction D = 2\pi\sqrt{-1} \hbar \, \frac{d}{dx}$ and $\wh{\mathbf{q}}^\hbar = Q\restriction D = x$ which are essentially self-adjoint, establishing the assertion 1).

\vs

If $\alpha_0=0$ then there is nothing to prove for 2), so let $\alpha_0\neq 0$; one may assume $\alpha_0 = 1$ by scaling. For each $\gamma\in\mathbb{R}$ consider the operator $U^{(\gamma)} : V \to V$ given by
\begin{align}
\label{eq:U_gamma}
(U^{(\gamma)}f)(x) = e^{\sqrt{-1}\gamma x^2} f(x).
\end{align}
Then $U^{(\gamma)}$ is unitary and preserves $D$. It is easy to see
$
U^{(\gamma)} \, \wh{\mathbf{q}}^\hbar \, (U^{(\gamma)})^{-1} = \wh{\mathbf{q}}^\hbar
$
as operators $D\to D$, as each factor is multiplication operator (hence commutes with one another). Notice that $(U^{(\gamma)})^{-1} = U^{(-\gamma)}$. For each $f\in D$ observe
\begin{align*}
( U^{(\gamma)} ( \wh{\mathbf{p}}^\hbar ( (U^{(\gamma)})^{-1} f ) ) )(x) & = e^{\sqrt{-1} \gamma x^2} \cdot ( \wh{\mathbf{p}}^\hbar ( U^{(-\gamma)} f ) )(x)
= e^{\sqrt{-1} \gamma x^2} \cdot (2\pi\sqrt{-1} \hbar) \cdot \frac{d}{dx} \left( e^{-\sqrt{-1}\gamma x^2} \, f(x) \right) \\
& = \cancel{ e^{\sqrt{-1} \gamma x^2} } (2\pi\sqrt{-1} \hbar) \, \cancel{ e^{-\sqrt{-1}\gamma x^2} } \cdot \left( - 2 \sqrt{-1} \gamma x \, f(x) + \frac{df(x)}{dx}  \right) \\
& = ( 4\pi \gamma \hbar \, \wh{\mathbf{q}}^\hbar f)(x) + (\wh{\mathbf{p}}^\hbar f)(x),
\end{align*}
hence $U^{(\gamma)} \, \wh{\mathbf{p}}^\hbar \, (U^{(\gamma)})^{-1} = \wh{\mathbf{p}}^\hbar + 4\pi \gamma \hbar \, \wh{\mathbf{q}}^\hbar$ as operators $D\to D$. Thus, putting $\gamma = - \frac{\beta_0}{4\pi \hbar}$ one obtains
\begin{align}
\label{eq:U_gamma_conjugation}
U^{(-\frac{\beta_0}{4\pi\hbar})} \, (\wh{\mathbf{p}}^\hbar + \beta_0  \wh{\mathbf{q}}^\hbar) \, (U^{(-\frac{\beta_0}{4\pi \hbar})})^{-1} = \wh{\mathbf{p}}^\hbar
\end{align}
as operators $D\to D$. As the operator $\wh{\mathbf{p}}^\hbar$ on $D$ is essentially self-adjoint and $U^{(-\frac{\beta_0}{4\pi\gamma \hbar})}$ is a unitary operator preserving $D$, one can easily deduce e.g. using Lem.\ref{lem:unitary_conjugation_commutes_with_functional_calculus} that $\wh{\mathbf{p}}^\hbar + \beta_0  \wh{\mathbf{q}}^\hbar$ is also essentially self-adjoint, i.e. has a unique self-adjoint extension. \qed

\begin{remark}
The reason why we allow $\alpha$ appearing in the definition \eqref{eq:D} of the space $D \subset L^2(\mathbb{R},dx)$ is to make $D$ to be invariant under the operators $U^{(\gamma)} = e^{\sqrt{-1}\gamma x^2}$. 
\end{remark}

\begin{theorem}[the pentagon equation of the non-compact quantum dilogarithm]
Let $P$, $Q$, and $V$ be as in Lem.\ref{lem:linear_combination_of_Heisenberg_operators}. Then
$$
\Phi^\hbar(P) \, \Phi^\hbar(Q) = \Phi^\hbar(Q) \, \Phi^\hbar(Q+P) \Phi^\hbar(P).
$$
\end{theorem}
As pointed out in \cite{FG09}, this pentagon equation was suggested in \cite{F95} and proved in \cite{FKV01}, \cite{W}, and \cite{G08}.

\begin{corollary}
\label{cor:QD_identity_pentagon2}
Let $P$, $Q$, and $V$ be as in Lem.\ref{lem:linear_combination_of_Heisenberg_operators}. Then
\begin{align}
\label{eq:QD_identity_pentagon2}
\Phi^\hbar(P) \Phi^\hbar(Q) \Phi^\hbar(-P) = c_\hbar \, \Phi^\hbar(Q) \, \Phi^\hbar(Q+P) \, \exp\left(\frac{P^2}{4\pi \sqrt{-1} \hbar}\right). \qed
\end{align}
\end{corollary}

It is sensible to expect that the following lemma holds.
\begin{lemma}
Let $P$ and $Q$ be densely defined self-adjoint operators on a Hilbert space $V$ that `strongly commute', i.e. their corresponding strongly continuous one-parameter unitary groups commute:
$$
e^{\sqrt{-1}\alpha P} \, e^{\sqrt{-1} \beta Q} = e^{\sqrt{-1} \beta Q} \, e^{\sqrt{-1} \alpha P}, \qquad \forall \alpha,\beta \in \mathbb{R}.
$$
Let $f_1, f_2: \mathbb{R} \to \mathbb{C}$ be any unitary functions in the sense of Lem.\ref{lem:unitary_operators_from_functional_calculus}. Then the unitary operators $f_1(P)$ and $f_2(Q)$ commute:
$$
f_1(P) f_2(Q) = f_2(Q) f_1(P).
$$
\end{lemma}
Instead of trying to prove the above lemma in full generality, we show only the following, which suffices for our purpose.
\begin{lemma}
\label{lem:commuting}
Let $\Gamma$ be a $\mathcal{D}$-seed with $n\ge 2$, and $D_\Gamma \subset \mathscr{H}_\Gamma = L^2(\mathbb{R}^n, \, da_1 \, \cdots da_n)$ be as in \eqref{eq:H_Gamma} and \eqref{eq:D_Gamma}. Suppose $\mathbf{o}_1$ and $\mathbf{o}_2$ are operators on $D_\Gamma$, preserving $D_\Gamma$, given by $\mathbb{R}$-linear combinations of $\wh{\mathbf{p}}^\hbar_i$, $\wh{\mathbf{q}}^\hbar_i$ \eqref{eq:Schrodinger_representation}, $i=1,\ldots,n$. Assume that none of the two is the zero operator, and that
$$
[\mathbf{o}_1, \mathbf{o}_2] = 0
$$
holds as operators $D_\Gamma \to D_\Gamma$. Then $\mathbf{o}_1$ and $\mathbf{o}_2$ are essentially self-adjoint. If we denote by $O_1$ and $O_2$ their unique self-adjoint extensions, then for any unitary functions $f_1,f_2: \mathbb{R} \to \mathbb{C}$ (in the sense of Lem.\ref{lem:unitary_operators_from_functional_calculus}), the unitary operators $f_1(O_1)$ and $f_2(O_2)$ commute:
$$
f_1(O_1) f_2(O_2) = f_2(O_2) f_1(O_1).
$$

\end{lemma}

\begin{remark}
This lemma can be written without using a $\mathcal{D}$-seed; I did so for a notational convenience.
\end{remark}

{\it Proof.} We shall use the special affine shift operators $\mathbf{S}_{(\mathbf{c},\mathbf{0})}$, the Fourier transforms $\mcal{F}_i$ \eqref{eq:mcal_F_i}, $i=1,\ldots,n$, the operators given by multiplication by $e^{\sqrt{-1} \gamma \, a_i^2}$ for $i=1,\ldots,n$ and $\gamma \in \mathbb{R}$, and permutation operators: for each permutation $\sigma$ of $\{1,\ldots,n\}$ define the permutation operator $\wh{\mathbf{P}}_\sigma : \mathscr{H}_\Gamma \to \mathscr{H}_\Gamma$ as $(\wh{\mathbf{P}}_\sigma f)(a_1,\ldots,a_n) = f(a_{\sigma(1)}, \ldots, a_{\sigma(n)})$; let us keep this notation $\wh{\mathbf{P}}_\sigma$ to this proof only, for it may create a confusion with the map \eqref{eq:bf_K_P_sigma}. One can then check that $\wh{\mathbf{P}}_\sigma$ is unitary, $\wh{\mathbf{P}}_\sigma(D_\Gamma) = D_\Gamma$, and the equalities
$$
\wh{\mathbf{P}}_\sigma \, \wh{\mathbf{p}}^\hbar_i \, \wh{\mathbf{P}}_\sigma^{-1} = \wh{\mathbf{p}}^\hbar_{\sigma(i)} \qquad\mbox{and}\qquad
\wh{\mathbf{P}}_\sigma \, \wh{\mathbf{q}}^\hbar_i \, \wh{\mathbf{P}}_\sigma^{-1} = \wh{\mathbf{q}}^\hbar_{\sigma(i)}
$$
hold as operators $D_\Gamma \to D_\Gamma$, for each $i=1,\ldots,n$ and each permutation $\sigma$. So, all the operators mentioned above preserve $D_\Gamma$. Recall e.g. from \eqref{eq:mcal_F_conjugation_on_p_and_q} that we computed the conjugation action of $\mcal{F}_i$:
\begin{align*}
  \mcal{F}_i \left( \frac{1}{(2\pi)^2 \hbar} \, \wh{\mathbf{p}}^\hbar_j \right) \, \mcal{F}_i^{-1} = \, \left\{ \begin{array}{ll} 
- \wh{\mathbf{q}}^\hbar_i & \mbox{if $j=i$}, \\
\frac{1}{(2\pi)^2 \hbar} \, \wh{\mathbf{p}}^\hbar_j & \mbox{if $j\neq i$}, 
\end{array} \right.
\qquad
\mcal{F}_i \, \wh{\mathbf{q}}^\hbar_i \, \mcal{F}_i^{-1} = \left\{ \begin{array}{ll}
\frac{1}{(2\pi)^2 \hbar} \, \wh{\mathbf{p}}^\hbar_i & \mbox{if $j=i$}, \\
\wh{\mathbf{q}}^\hbar_j, & \mbox{if $j\neq i$},
\end{array} \right.
\end{align*}
each of which is an equality of operators $D_\Gamma \to D_\Gamma$. The conjugation action of $\mathbf{S}_{(\mathbf{c},\mathbf{0})}$ on $\wh{\mathbf{p}}^\hbar_i$, $\wh{\mathbf{q}}^\hbar_i$, $i=1,\ldots,n$, is computed in \eqref{eq:conjugation_of_bf_S_on_bf_p} and \eqref{eq:conjugation_of_bf_S_on_bf_q} of Lem.\ref{lem:bf_S_conjugation_on_bf_p_and_bf_q}.

\vs

Let $\mathbf{o}_1 = \sum_{i=1}^n (\alpha_i \wh{\mathbf{p}}^\hbar_i + \beta_i \wh{\mathbf{q}}^\hbar_i)$, for some $\alpha_i,\beta_i \in \mathbb{R}$. We shall first find a unitary $U_1 : D_\Gamma \to D_\Gamma$ such that
\begin{align}
\label{eq:metaplectic_to_accomplish1}
U_1 \mathbf{o}_1 U_1^{-1} = \beta_1' \, \wh{\mathbf{q}}^\hbar_1 \quad\mbox{for some nonzero real $\beta_1'$.}
\end{align}
Not all $\alpha_i$'s and $\beta_i$'s are zero. One can assume that one of the $\alpha_i$'s is not zero; if not, find a nonzero $\beta_i$, so that conjugation by $\mcal{F}_i$ on $\mathbf{o}_1$ results in a nonzero coefficient for $\wh{\mathbf{p}}^\hbar_i$. Then, by applying the conjugation by some permutation operator $\wh{\mathbf{P}}_\sigma$ if necessary, we can make the coefficient for $\wh{\mathbf{p}}^\hbar_1$ to be nonzero. So, assume $\alpha_1\neq 0$ from now on.

\vs

Now define a matrix $\mathbf{c}^{(0)} = (c^{(0)}_{ij})_{i,j\in\{1,\ldots,n\}} \in \mathrm{SL}_\pm(n,\mathbb{R})$ as
$$
c^{(0)}_{ii}=1, \quad \forall i=1,\ldots,n, \qquad c^{(0)}_{1j} = \alpha_j/\alpha_1 \, \quad \forall j \neq 1, \qquad c^{(0)}_{ij} = 0 \quad\mbox{otherwise},
$$
so that its inverse $(\mathbf{c}^{(0)})^{-1} = ( (c^{(0)})^{ij} )_{i,j\in\{1,\ldots,n\}}$ is given by $(c^{(0)})^{ii}=1$, $\forall i=1,\ldots,n$, $(c^{(0)})^{1j} = - \alpha_j/\alpha_1$, $\forall i \neq 2$, and $(c^{(0)})^{ij}=0$ otherwise. From \eqref{eq:conjugation_of_bf_S_on_bf_p} and \eqref{eq:conjugation_of_bf_S_on_bf_q} one can verify that $\mathbf{o}_1' := \mathbf{S}_{(\mathbf{c}^{(0)},\mathbf{0})} \, \mathbf{o}_1 \, \mathbf{S}_{(\mathbf{c}^{(0)},\mathbf{0})}^{-1} = \alpha_1 \wh{\mathbf{p}}^\hbar_1 + (\sum_{i=1}^n \beta_i'' \wh{\mathbf{q}}^\hbar_i)$, for some $\beta_i'' \in \mathbb{R}$. First, suppose that there is no $\ell \neq 1$ with $\beta_\ell''\neq 0$. So
\begin{align}
\label{eq:metaplectic_situation1}
\mathbf{o}_1' = \alpha_1 \wh{\mathbf{p}}^\hbar_1 + \beta_1'' \wh{\mathbf{q}}^\hbar_1 \qquad (\mbox{where $\alpha_1, \beta_1''$ are real numbers, with $\alpha_1\neq 0$}).
\end{align}
Following \eqref{eq:U_gamma}, let the unitary operator $U_2$ on $\mathscr{H}_\Gamma$ be the multiplication by the unitary function $e^{\sqrt{-1} \gamma_1 a_1^2}$ with $\gamma_1 = - \frac{1}{4 \pi\hbar} \cdot \frac{\beta_1''}{\alpha_1}$. As seen in \eqref{eq:U_gamma_conjugation}, one gets $U_2 \mathbf{o}_1' U_2^{-1} = \alpha_1 \, \wh{\mathbf{p}}^\hbar_1$; conjugation by $\mcal{F}_1$ then puts us to the desired situation \eqref{eq:metaplectic_to_accomplish1}.

\vs

Suppose that there exists $\ell \neq 1$ with $\beta_\ell''\neq 0$; after conjugating by a permutation operator, one may assume $\ell=2$. Define a matrix $\mathbf{c}^{(1)} = (c^{(1)}_{ij})_{i,j\in \{1,\ldots,n\}} \in \mathrm{SL}_\pm(n,\mathbb{R})$ as
$$
c^{(1)}_{ii} = 1, \quad \forall i =1,\ldots,n, \qquad
c^{(1)}_{i 2} = - \beta_i''/\beta_2'', \quad \forall i\neq 2, \qquad
c^{(1)}_{ij}=0 \quad\mbox{otherwise}.
$$
Then its inverse $(\mathbf{c}^{(1)})^{-1} = ( (c^{(1)})^{ij} )_{i,j\in \{1,\ldots,n\}}$ is given by $(c^{(1)})^{ii}=1$, $\forall i=1,\ldots,n$, $(c^{(1)})^{i 2} = \beta_i''/\beta_2''$, $\forall i\neq 2$, and $(c^{(1)})^{ij} = 0$ otherwise. From \eqref{eq:conjugation_of_bf_S_on_bf_p} and \eqref{eq:conjugation_of_bf_S_on_bf_q} one can verify that $\mathbf{o}_1'' := \mathbf{S}_{(\mathbf{c}^{(1)},\mathbf{0})} \, \mathbf{o}_1' \, \mathbf{S}_{(\mathbf{c}^{(1)},\mathbf{0})}^{-1} = \alpha_1 \wh{\mathbf{p}}^\hbar_1 + \frac{\beta_1''}{\beta_2''} \alpha_1 \wh{\mathbf{p}}^\hbar_2 + \beta_2'' \wh{\mathbf{q}}^\hbar_2$. First, assume $\beta_1''=0$, so that $\mathbf{o}_1'' = \alpha_1 \wh{\mathbf{p}}^\hbar_1 + \beta_2'' \wh{\mathbf{q}}^\hbar_2$, with nonzero real $\alpha_1,\beta_2''$. Then $\mcal{F}_1 \, \mathbf{o}_1'' \mcal{F}_1^{-1} = \alpha_1' \wh{\mathbf{q}}^\hbar_1 + \beta_2'' \wh{\mathbf{q}}^\hbar_2$, with nonzero real $\alpha_1',\beta_2''$. Using a similar trick we used by $\mathbf{S}_{(\mathbf{c}^{(0)},\mathbf{0})}$, one can find some $\mathbf{c}^{(2)} \in \mathrm{SL}_\pm(n,\mathbb{R})$ such that $\mathbf{S}_{(\mathbf{c}^{(2)},\mathbf{0})} \, (\mcal{F}_1 \, \mathbf{o}_1'' \, \mcal{F}_1^{-1}) \mathbf{S}_{(\mathbf{c}^{(2)},\mathbf{0})}^{-1} = \alpha_1' \wh{\mathbf{q}}^\hbar_1$, as desired in \eqref{eq:metaplectic_to_accomplish1}. Now assume $\beta_1''\neq 0$; again we use similar trick as for $\mathbf{S}_{(\mathbf{c}^{(0)},\mathbf{0})}$. Define a matrix $\mathbf{c}^{(3)} = (c^{(3)}_{ij}) \in \mathrm{SL}_\pm(n,\mathbb{R})$ so that it differs by the identity matrix possibly only at the entries for $i,j\in \{1,2\}$, where it is given by
$$
c^{(3)}_{11} = 1, \qquad c^{(3)}_{12} = 0, \qquad c^{(3)}_{21} = \beta_2''/\beta_1'', \qquad c^{(3)}_{22}=1.
$$
Then the inverse matrix $(\mathbf{c}^{(3)})^{-1} = ((c^{(3)})^{ij})$ differs from the identity matrix possibily only at the entries for $i,j\in \{1,2\}$, where $(c^{(3)})^{11} = (c^{(3)})^{22}=1$, $(c^{(3)})^{12}=0$, and $(c^{(3)})^{21} = -\beta_2''/\beta_1''$. From \eqref{eq:conjugation_of_bf_S_on_bf_p} and \eqref{eq:conjugation_of_bf_S_on_bf_q} one can check that $\mathbf{o}_1''' := \mathbf{S}_{(\mathbf{c}^{(3)},\mathbf{0})} \, \mathbf{o}_1''  \, \mathbf{S}_{(\mathbf{c}^{(3)},\mathbf{0})}^{-1} = \frac{\beta_1''}{\beta_2''}\alpha_1 \wh{\mathbf{p}}^\hbar_2 + \beta_2'' \wh{\mathbf{q}}^\hbar_2$, where both $\frac{\beta_1''}{\beta_2''}\alpha_1$ and $\beta_2''$ are nonzero real. Conjugation by the permutation operator $\wh{\mathbf{P}}_{(1\,2)}$ puts us into the situation \eqref{eq:metaplectic_situation1}, which is already dealt with.

\vs

Let us refresh the notations. If $\mathbf{o}_1$ and $\mathbf{o}_2$ are any non-zero $\mathbb{R}$-linear combinations of $\wh{\mathbf{p}}^\hbar_i$, $\wh{\mathbf{q}}^\hbar_i$, $i=1,\ldots,n$, one can find a unitary operator $U_1$ preserving $D_\Gamma$ such that
$$
\mathbf{o}_1' := U_1 \mathbf{o}_1 U_1^{-1} = \beta \wh{\mathbf{q}}^\hbar_1, \qquad
\mathbf{o}_2' := U_1 \mathbf{o}_2 U_1^{-1} = \sum_{i=1}^n (\alpha_i \wh{\mathbf{p}}^\hbar_i + \beta_i \wh{\mathbf{q}}^\hbar_i)
$$
for some real $\beta,\alpha_i,\beta_i$ with $\beta\neq 0$. Then from \eqref{eq:Heisenberg_relations_of_p_and_q} we get $[\mathbf{o}_1', \mathbf{o}_2'] = - \beta \, \alpha_1 \cdot \mathrm{id}$ on $D_\Gamma$. As we assumed $[\mathbf{o}_1, \mathbf{o}_2]=0$ on $D_\Gamma$, we must have $[\mathbf{o}_1', \mathbf{o}_2'] = U_1 [ \mathbf{o}_1, \mathbf{o}_2] \, U_1^{-1} = 0$ on $D_\Gamma$, hence $\alpha_1=0$.

\vs

 Suppose $\sum_{i=2}^n \alpha_i \wh{\mathbf{p}}^\hbar_i + \beta_i \wh{\mathbf{q}}^\hbar_i$ is zero. Then $\mathbf{o}_2' = \beta_1 \wh{\mathbf{q}}^\hbar_1$. As $\wh{\mathbf{q}}^\hbar_1$ is essentially self-adjoint, so are $\mathbf{o}_1'$ and $\mathbf{o}_2'$. Let $Q_1$, $O_1'$, and $O_2'$ be the unique self-adjoint extensions of $\wh{\mathbf{q}}^\hbar_1$, $\mathbf{o}_1'$, and $\mathbf{o}_2'$. Then we have $O_1' = \beta Q_1$ and $O_2' = \beta_1 Q_1$, as genuine equality of operators.

\vs

Write the unique spectral resolution of $Q_1$ as
\begin{align}
\label{eq:Q1_spectral_resolution}
Q_1 = \int_\mathbb{R} \lambda \, dE(\lambda),
\end{align}
on its domain
$$
D(Q_1) = \{ v\in \mathscr{H}_\Gamma \, | \, \int_\mathbb{R} |\lambda|^2 \, d||E(\lambda) v||_{\mathscr{H}_\Gamma}^2 < \infty\},
$$
where $E(\lambda)$ is some resolution of the identity. We claim that
\begin{align}
\label{eq:our_E_lambda}
E(\lambda) v = \chi_{a_1\in (-\infty,\lambda]} \cdot v, \qquad \forall v\in \mathscr{H}_\Gamma,
\end{align}
where $\chi_{a_1 \in (-\infty, \lambda]}$ is a function on $\mathbb{R}^n$ defined as
$$
\chi_{a_1 \in (-\infty, \lambda]} (a_1,a_2,\ldots,a_n) = \left\{ 
\begin{array}{ll}
1 & \mbox{if $a_1 \le \lambda$}, \\
0 & \mbox{if $a_1 > \lambda$}.
\end{array}
\right.
$$
This is an analog of the $1$-dimensional case: multiplication by $x$ on $L^2(\mathbb{R},dx)$; see e.g. the Example at the end of \cite[Chap.XI.5]{Y}. By an explicit computation one can indeed check that $\int_\mathbb{R} \lambda \, dE(\lambda)$ with this $E(\lambda)$ is multiplication by $a_1$; so the uniqueness of spectral resolution indeed tells us $Q_1$ is written as in \eqref{eq:Q1_spectral_resolution} with $E(\lambda)$ as in \eqref{eq:our_E_lambda}. Then, for any unitary function $f : \mathbb{R} \to \mathbb{C}$, the functional calculus says that the densely defined operator $f(Q_1)$ is defined by the formula
\begin{align}
\label{eq:f_of_Q1}
f(Q_1) = \int_\mathbb{R} f_1(\lambda) \, dE(\lambda)
\end{align}
on its appropriate well-determined domain
$$
D(f(Q_1)) = \{ v\in \mathscr{H}_\Gamma \, | \, \int_\mathbb{R} |f(\lambda)|^2 \, d||E(\lambda) v||_{\mathscr{H}_\Gamma}^2 < \infty\},
$$
on which we can check explicitly by \eqref{eq:f_of_Q1} with \eqref{eq:our_E_lambda} that $f(Q_1)$ is given by multiplication by the unitary function $(a_1,\ldots,a_n) \mapsto f(a_k)$. So we get such a unitary (multiplication) operator; in particular, $D(f(Q_1)) = \mathscr{H}_\Gamma$.

\vs

Similary, one can figure out that
$$
O_1' = \int_\mathbb{R} \lambda \, dE_1(\lambda), \qquad
O_2' = \int_\mathbb{R} \lambda \, dE_2(\lambda),
$$
where
$$
E_1(\lambda) = \chi_{a_1 \in (-\infty, \lambda/\beta]} \cdot v, \qquad E_2(\lambda) = \chi_{a_1 \in (-\infty, \lambda/\beta_1]} \cdot v,
$$
are the unique spectral resolutions of $O_1'$ and $O_2'$ respectively. Observe that $E_1(\lambda) = E(\lambda/\beta)$ and $E_2(\lambda) = E(\lambda/\beta_1)$. So, for any unitary functions $f_1$ and $f_2$, from functional calculus we have
$$
f_1(O_1') = \int_\mathbb{R} f_1(\lambda) \, dE_1(\lambda) = \int_\mathbb{R} f_1(\lambda) \, dE(\lambda/\beta) = \int_\mathbb{R} f_1(\lambda_1 \beta) \, dE(\lambda_1) = \wh{f}_1(Q_1),
$$
where $\wh{f}_1(x) := f_1(x\beta)$. Similarly, $f_2(O_2') = \wh{f}_2(Q_1)$ with $\wh{f}_2(x) := f_2(x\beta_1)$. Note $\wh{f}_1$ and $\wh{f}_2$ are unitary functions, hence the unitary operators $\wh{f}_1(Q_1)$ and $\wh{f}_2(Q_1)$ are multiplication by the unitary functions $\wh{f}_1(a_1)$ and $\wh{f}_2(a_1)$, respectively, hence they commute. So $f_1(O_1')$ and $f_2(O_2')$ commute. As in the statement, recall that $O_1$ and $O_2$ are self-adjoint extensions of $\mathbf{o}_1$ and $\mathbf{o}_2$. So $O_1' = U_1 O_1 U_1^{-1}$ and $O_2' = U_1 O_2 U_1^{-1}$, with $f_1(O_1) = U_1^{-1} f(O_1') U_1$ and $f_2(O_2) = U_1^{-1} f_2(O_2') U_1$. Hence the unitary operators $f_1(O_1)$ and $f_2(O_2)$ commute, as desired.

\vs

Now assume $\sum_{i=2}^n \alpha_i \wh{\mathbf{p}}^\hbar_i + \beta_i \wh{\mathbf{q}}^\hbar_i$ is not zero. We apply the same argument as before to the operator $U_1 \mathbf{o}_2 \, U_1^{-1}$ involving only the indices $i=2,3,\ldots,n$; namely, we can find some unitary operator $U_3 : D_\Gamma \to D_\Gamma$ such that $U_3 ( \sum_{i=2}^n \alpha_i \wh{\mathbf{p}}^\hbar_i + \beta_i \wh{\mathbf{q}}^\hbar_i) U_3^{-1} = \beta' \wh{\mathbf{q}}^\hbar_2$ for some nonzero real $\beta'$, while $U_3 \wh{\mathbf{p}}^\hbar_1 U_3^{-1} = \wh{\mathbf{p}}^\hbar_1$ and $U_3 \wh{\mathbf{q}}^\hbar_1 U_3^{-1} = \wh{\mathbf{q}}^\hbar_1$, all as equalities of operators $D_\Gamma \to D_\Gamma$. Thus $\mathbf{o}_1'' := U_3 \mathbf{o}_1' U_3^{-1} = \beta \wh{\mathbf{q}}^\hbar_1$ and $\mathbf{o}_2'' := U_3 \mathbf{o}_2' U_3^{-1} = \beta_1 \wh{\mathbf{q}}^\hbar_1 + \beta' \wh{\mathbf{q}}^\hbar_2$. Define a matrix $\mathbf{c}^{(4)} = (c^{(4)}_{ij}) \in \mathrm{SL}_\pm(n,\mathbb{R})$ so that it differs by the identity matrix possibly only at the entires for $i,j \in \{1,2\}$, where it is given by
$$
c^{(4)}_{11}=1, \qquad c^{(4)}_{12} = - \beta_1/ \beta', \qquad c^{(4)}_{21} = 0, \qquad c^{(4)}_{22}=1.
$$
Then by \eqref{eq:conjugation_of_bf_S_on_bf_q} we get $\mathbf{o}_1''' := \mathbf{S}_{(\mathbf{c}^{(4)},\mathbf{0})} \, \mathbf{o}_1'' \, \mathbf{S}_{(\mathbf{c}^{(4)},\mathbf{0})}^{-1} = \beta \wh{\mathbf{q}}^\hbar_1$ and $\mathbf{o}_2''' := \mathbf{S}_{(\mathbf{c}^{(4)},\mathbf{0})} \, \mathbf{o}_2'' \, \mathbf{S}_{(\mathbf{c}^{(4)},\mathbf{0})}^{-1} = \beta' \wh{\mathbf{q}}^\hbar_2$. So, for each $j=1,2$, the operator $\mathbf{o}_j'''$, and hence also its self-adjoint extension $O_j'''$, involves only the variable $a_j$; they are acting on different tensor factor in $\bigotimes_{i=1}^n L^2(\mathbb{R}, da_i) \cong L^2(\mathbb{R}^n, da_1 \, \cdots da_n)$. Hence the unitary operators $f_1(O_1''')$ and $f_2(O_2''')$ commute, and therefore so do $f_1(O_1) = U\, f_1(O_1''') U^{-1}$ and $f_2(O_2) = U\, f_2(O_2''') \, U^{-1}$ (where $U$ is some unitary operator we constructed, so that $O_1''' = U^{-1} O_1 U$ and $O_2''' = U^{-1} O_2 U$). \qed

\subsection{Rank $1$ identity: `twice-flip' is trivial}

\begin{proposition}
\label{prop:c_A1_is_1}
$c_{A_1}=1$.
\end{proposition}
{\it Proof.} The LHS of \eqref{eq:twice_flip_to_prove} is
\begin{align*}
& \mathbf{K}^\hbar_{\Gamma\mut{k}\Gamma'} \, \mathbf{K}^\hbar_{\Gamma'\mut{k}\Gamma} = ({\bf K}^{\sharp \hbar}_{\Gamma\mut{k}\Gamma'} \, {\bf K}'_{\Gamma\mut{k}\Gamma'}) \, ({\bf K}^{\sharp \hbar}_{\Gamma'\mut{k}\Gamma} \, {\bf K}'_{\Gamma'\mut{k}\Gamma}) \\
& = \Phi^{\hbar_k}(\wh{\mathbf{x}}^\hbar_{\Gamma;k})  \, \Phi^{\hbar_k}(\wh{\til{\mathbf{x}}}^\hbar_{\Gamma;k})^{-1} \, \ul{ {\bf K}'_{\Gamma\mut{k}\Gamma'} } \, \, \ul{ \Phi^{\hbar'_k}(\wh{\mathbf{x}}^\hbar_{\Gamma';k}) \, \Phi^{\hbar'_k}(\wh{\til{\mathbf{x}}}^\hbar_{\Gamma';k})^{-1}  } \, {\bf K}'_{\Gamma'\mut{k}\Gamma} \\
& =
\Phi^{\hbar_k}(\wh{\mathbf{x}}^\hbar_{\Gamma;k})  \, \ul{ \Phi^{\hbar_k}(\wh{\til{\mathbf{x}}}^\hbar_{\Gamma;k})^{-1} }
\left( \ul{ \Phi^{\hbar_k}(-\wh{\mathbf{x}}^\hbar_{\Gamma;k}) }
\,  \Phi^{\hbar_k}(-\wh{\til{\mathbf{x}}}^\hbar_{\Gamma;k})^{-1} \, {\bf K}'_{\Gamma\mut{k}\Gamma'}  \right) \, {\bf K}'_{\Gamma'\mut{k}\Gamma} \quad (\because {\rm Lem}.\ref{lem:unitary_conjugation_commutes_with_functional_calculus}, \ref{lem:bf_K_prime_conjugation_on_bf_x}, \, \hbar'_k=\hbar_k) \\
& \stackrel{\vee}{=}
\ul{ \Phi^{\hbar_k}(\wh{\mathbf{x}}^\hbar_{\Gamma;k})  \,
 \Phi^{\hbar_k}(-\wh{\mathbf{x}}^\hbar_{\Gamma;k})  }
\,\, \ul{ \Phi^{\hbar_k}(\wh{\til{\mathbf{x}}}^\hbar_{\Gamma;k})^{-1} 
\, \Phi^{\hbar_k}(-\wh{\til{\mathbf{x}}}^\hbar_{\Gamma;k})^{-1} }
\, {\bf K}'_{\Gamma\mut{k}\Gamma'} \, {\bf K}'_{\Gamma'\mut{k}\Gamma} \\
& \hspace{-1mm} \stackrel{\eqref{eq:QD_identity_quadratic}}{=}
\cancel{ c_{\hbar_k} } \, \cancel{ c_{\hbar_k}^{-1} } \, 
\underbrace{ \exp\left( \frac{(\wh{\mathbf{x}}^\hbar_{\Gamma;k})^2}{4\pi\sqrt{-1}\hbar_k} \right) \, 
\exp\left( - \frac{(\wh{\til{\mathbf{x}}}^\hbar_{\Gamma;k})^2}{4\pi\sqrt{-1}\hbar_k} \right) } \, {\bf K}'_{\Gamma\mut{k}\Gamma'} \, {\bf K}'_{\Gamma'\mut{k}\Gamma},
\end{align*}
where for the checked equality we used the fact that the two unitary operators $\Phi^{\hbar_k}(\wh{\til{\mathbf{x}}}^\hbar_{\Gamma;k})^{-1}$ and $\Phi^{\hbar_k}(-\wh{\mathbf{x}}^\hbar_{\Gamma;k})$, where $\wh{\til{\mathbf{x}}}^\hbar_{\Gamma;k}$ and $\wh{\mathbf{x}}^\hbar_{\Gamma;k}$ here stand for the respective unique self-adjoint extensions, commute with each other, as follows from Lem.\ref{lem:commuting} since $[\wh{\til{\mathbf{x}}}^\hbar_{\Gamma;k}, \wh{\mathbf{x}}^\hbar_{\Gamma;k}]=0$.

\vs

We shall prove that the underbraced part is a special affine shift operator. If we identify $\mathscr{H}_\Gamma$ and $\mathscr{H}_{\Gamma'}$ by the natural relabling isomorphisms $\mathbf{I}_{\Gamma\mut{k}\Gamma'}$ and $\mathbf{I}_{\Gamma'\mut{k}\Gamma}$ \eqref{eq:bf_I}, we see from \eqref{eq:bf_K_prime} that ${\bf K}'_{\Gamma\mut{k}\Gamma'}$ and ${\bf K}'_{\Gamma'\mut{k}\Gamma}$ are special affine shift operators. Then one deduces from Cor.\ref{cor:special_affine_shift_operators_form_a_group} that $\mathbf{K}^\hbar_{\Gamma\mut{k}\Gamma'} \, \mathbf{K}^\hbar_{\Gamma'\mut{k}\Gamma}$ is a special affine shift operator as well. As we already know from Prop.\ref{prop:rank_1_constant} that it is a scalar operator, from Lem.\ref{lem:scalar_special_affine_shift_operator_is_identity} we conclude that it is the identity operator. If one would like, one can prove without relying on Prop.\ref{prop:rank_1_constant}. We shall now compute what precise special affine shift operator that $\exp\left( \frac{(\wh{\mathbf{x}}^\hbar_{\Gamma;k})^2}{4\pi\sqrt{-1}\hbar_k} \right) \exp\left( - \frac{(\wh{\til{\mathbf{x}}}^\hbar_{\Gamma;k})^2}{4\pi\sqrt{-1}\hbar_k} \right)$ equals to; this, together with the expressions for ${\bf K}'_{\Gamma\mut{k}\Gamma'}$ and ${\bf K}'_{\Gamma'\mut{k}\Gamma}$ in terms of explicit special affine shift operators, allows one to compute $\mathbf{K}^\hbar_{\Gamma\mut{k}\Gamma'} \, \mathbf{K}^\hbar_{\Gamma'\mut{k}\Gamma}$ as a single special affine shift operator, with the help of \eqref{eq:bf_S_is_group_homomorphism}. One can directly prove by computation that this equals $\mathbf{S}_{(\mathbf{id},\mathbf{0})}$.

\vs

So it remains to prove:
\begin{lemma}
\label{lem:product_of_two_quadratic_exponential_operators}
For any $\mathcal{D}$-seed $\Gamma$, $k\in \{1,\ldots,n\}$, and $\hbar \in \mathbb{R}_{>0}\setminus\mathbb{Q}$, the unitary operator
\begin{align}
\label{eq:product_of_two_quadratic_exponential_operators}
\exp\left( \frac{(\wh{\mathbf{x}}^\hbar_{\Gamma;k})^2}{4\pi\sqrt{-1}\hbar_k} \right) \, 
\exp\left( - \frac{(\wh{\til{\mathbf{x}}}^\hbar_{\Gamma;k})^2}{4\pi\sqrt{-1}\hbar_k} \right)
\end{align}
on $\mathscr{H}_\Gamma$, whose two factors are obtained by applying the functional calculus of \S\ref{subsec:spectral_theorem} to the unique self-adjoint extensions of the operators $\wh{\mathbf{x}}^\hbar_{\Gamma;k}$ and $\wh{\til{\mathbf{x}}}^\hbar_{\Gamma;k}$ on $D_\Gamma$ for the unitary functions $z\mapsto \exp(\frac{z^2}{4\pi\sqrt{-1}\hbar_k})$ and $z\mapsto \exp(-\frac{z^2}{4\pi\sqrt{-1}\hbar_k})$ respectively, coincides with the special affine shift operator $\mathbf{S}_{(\mathbf{c},\mathbf{0})}$, where $\mathbf{c}=(c_{ij})_{i,j\in\{1,\ldots,n\}}$ is given by
$$
c_{ii}=1, \quad \forall i=1,\ldots,n, \qquad
c_{ik} = -\varepsilon_{ki}, \quad \forall i\neq k, \qquad c_{ij}=0 \quad \mbox{otherwise}.
$$
\end{lemma}

{\it Proof.} (beware that there will be lots of notations in this proof). \quad Let $A$ and $\til{A}$ stand respectively for the unique self-adjoint extensions of the essentially self-adjoint operators $\wh{\mathbf{x}}^\hbar_{\Gamma;k}$ and $\wh{\til{\mathbf{x}}}^\hbar_{\Gamma;k}$ on $D_\Gamma$. Let $f^\pm(\lambda) := \exp(\pm \frac{\lambda^2}{4\pi\sqrt{-1} \hbar_k})$ be two functions on $\mathbb{R}$, which are manifestly unitary in the sense of Lem.\ref{lem:unitary_operators_from_functional_calculus}. So the functional calculus in \S\ref{subsec:spectral_theorem} yield two unitary operators $f^+(A)$ and $f^-(\til{A})$; what we mean by the sloppy expression \eqref{eq:product_of_two_quadratic_exponential_operators} is $f^+(A)\, f^-(\til{A})$. Recall the Fourier transform $\mcal{F}_k : D_\Gamma \to D_\Gamma$ \eqref{eq:mcal_F_i}, and its conjugation action written as the equality $\mcal{F}_k ( \frac{1}{(2\pi)^2\hbar} \, \wh{\mathbf{p}}^\hbar_k) \mcal{F}_k^{-1} = - \wh{\mathbf{q}}^\hbar_k$ of operators $D_\Gamma\to D_\Gamma$ (see \eqref{eq:mcal_F_i_conjugation}). We also notice $\mcal{F}_k \, \wh{\mathbf{q}}^\hbar_j \, \mcal{F}_k^{-1} = \wh{\mathbf{q}}^\hbar_j$ for each $j\neq k$. Thus we observe
\begin{align}
\label{eq:mcal_F_k_conjugation_on_wh_bf_x_k}
& \mcal{F}_k \, \wh{\mathbf{x}}^\hbar_k \, \mcal{F}_k^{-1} \,\, \stackrel{\eqref{eq:old_representation}, ~ \varepsilon_{kk}=0}{=} \,\, \mcal{F}_k\left( \frac{d_k^{-1}}{2}\wh{\mathbf{p}}^\hbar_k - \sum_{j\neq k} \varepsilon_{kj} \wh{\mathbf{q}}^\hbar_j \right) \, \mcal{F}_k^{-1}
= - 2\pi^2 \hbar_k \, \wh{\mathbf{q}}^\hbar_k - \sum_{j\neq k} \varepsilon_{kj} \, \wh{\mathbf{q}}^\hbar_j =: Q_k, \\
\label{eq:mcal_F_k_conjugation_on_wh_bf_til_x_k}
& \mcal{F}_k \, \wh{\til{\mathbf{x}}}^\hbar_k \, \mcal{F}_k^{-1} \,\, \stackrel{\eqref{eq:old_representation_tilde}, ~ \varepsilon_{kk}=0}{=} \,\, \mcal{F}_k\left( \frac{d_k^{-1}}{2}\wh{\mathbf{p}}^\hbar_k + \sum_{j\neq k} \varepsilon_{kj} \wh{\mathbf{q}}^\hbar_j \right) \, \mcal{F}_k^{-1}
= - 2\pi^2 \hbar_k \, \wh{\mathbf{q}}^\hbar_k + \sum_{j\neq k} \varepsilon_{kj} \, \wh{\mathbf{q}}^\hbar_j =: \til{Q}_k,
\end{align}
both being equaltities of operators $D_\Gamma \to D_\Gamma$. Since $\wh{\mathbf{x}}^\hbar_k$ and $\wh{\til{\mathbf{x}}}^\hbar_k$ on $D_\Gamma$ are essentially self-adjoint and $\mcal{F}_k$ is a unitary map inducing a bijection $D_\Gamma \to D_\Gamma$, we deduce using Lem.\ref{lem:unitarily_equivalent_essentially_self-adjoint_operators} that the operators $Q_k$ and $\til{Q}_k$ on $D_\Gamma$ are also essentially self-adjoint. If we denote their unique self-adjoint extensions by $B$ and $\til{B}$ respectively, then we have the equalities
\begin{align}
\label{eq:mcal_F_k_conjugation_on_A_and_til_A}
\mcal{F}_k \, A \, \mcal{F}_k^{-1} = B \quad\mbox{and}\quad
\mcal{F}_k \, \til{A} \, \mcal{F}_k^{-1} = \til{B},
\end{align}
each of which holds when applied to elements in the domain of the RHS. Now, Lem.\ref{lem:unitary_conjugation_commutes_with_functional_calculus} tells us that $\mcal{F}_k \, f^+(A) \, f^-(\til{A}) \, \mcal{F}_k^{-1} = \mcal{F}_k \, f^+(A) \, \mcal{F}_k^{-1} \, \mcal{F}_k \, f^-(\til{A}) \, \mcal{F}_k^{-1} = f^+(\mcal{F}_k \, A \, \mcal{F}_k^{-1}) \, f^-(\mcal{F}_k \, \til{A} \, \mcal{F}_k^{-1}) = f^+(B) \, f^-(\til{B})$. Note that $Q_k$ and $\til{Q}_k$ are multiplication by the functions 
$$
f_k := - 2\pi^2 \hbar_k \, a_k - \sum_{j\neq k} \varepsilon_{kj} \, a_j \quad\mbox{and}\quad
\til{f}_k := - 2\pi^2 \hbar_k \, a_k + \sum_{j\neq k} \varepsilon_{kj} \, a_j,
$$
respectively. I claim that, on the whole $\mathscr{H}_\Gamma$, the operator $f^+(B)$ is multiplication by the function $f^+ \circ f_k = f^+(f_k)$, while $f^-(\til{B})$ is multiplication by $f^-\circ \til{f}_k$. A similar such statement is found in \cite[{\S}VIII]{RS70} under the name of `functional calculus', but in my opinion, the well-definedness is not rigorously established by the treatment there; this is why I followed \cite{Y} for the functional calculus. Moreover, the statement in \cite{RS70} does not exactly apply to this particular situation.

\vs

So let us be more careful; here I deal with $f^+(B)$ only. Define a matrix $\til{\mathbf{c}} = (\til{c}_{ij})_{i,j\in \{1,\ldots,n\}}$ by
\begin{align}
\label{eq:our_til_bf_c}
\til{c}_{ii}=1,\quad \forall i=1,\ldots,n, \qquad
\til{c}_{jk} = - \frac{\varepsilon_{kj}}{2\pi^2 \hbar_k}, \quad \forall j\neq k, \qquad \til{c}_{ij}=0 \quad \mbox{otherwise}.
\end{align}
Then it is easy to see $\til{\mathbf{c}} \in \mathrm{SL}_\pm(n,\mathbb{R})$, and therefore $\mathbf{S}_{(\mathbf{c},\mathbf{0})}$ restricted to $D_\Gamma$ is a bijection $D_\Gamma\to D_\Gamma$ (Lem.\ref{lem:bf_S_preserves_D}). One observes
\begin{align}
\label{eq:Q_k_and_wh_Q_k}
\mathbf{S}_{(\til{\mathbf{c}},\mathbf{0})} \, Q_k \, \mathbf{S}_{(\til{\mathbf{c}},\mathbf{0})}^{-1} = \mathbf{S}_{(\til{\mathbf{c}},\mathbf{0})} \left( - 2\pi^2 \hbar_k \, \wh{\mathbf{q}}^\hbar_k - \sum_{j\neq k} \varepsilon_{kj} \, \wh{\mathbf{q}}^\hbar_j \right) \mathbf{S}_{(\til{\mathbf{c}},\mathbf{0})}^{-1}  \,\, \stackrel{\eqref{eq:conjugation_of_bf_S_on_bf_q},\eqref{eq:our_til_bf_c}}{=} \, - 2\pi^2 \hbar_k \, \wh{\mathbf{q}}^\hbar_k =: \wh{Q}_k
\end{align}
as operators $D_\Gamma\to D_\Gamma$. So the operator $\wh{Q}_k$ on $D_\Gamma$ is given by multiplication by $-2\pi^2\hbar_k \, a_k$, which involves only one variable $a_k$. We saw in \S\ref{subsec:positive_representations} when we first defined $D_\Gamma$ that $D_\Gamma$ is canonically isomorphic to the algebraic tensor product of $D_i \subset L^2(\mathbb{R}, da_i)$ \eqref{eq:D_i}, $i=1,\ldots,n$.

\vs

Let $\wh{B}$ be the unique self-adjoint extenion of $\wh{Q}_k$. As in the proof of Lem.\ref{lem:commuting} one can see that the unitary operator $f^+(\wh{B})$ is multiplication by the unitary function $(a_1,\ldots,a_n) \mapsto f^+(-2\pi^2\hbar_k a_k)$ on $\mathbb{R}^n$. The equality \eqref{eq:Q_k_and_wh_Q_k} can be extended to the self-adjoint version $\mathbf{S}_{(\til{\mathbf{c}},\mathbf{0})} \, B \, \mathbf{S}_{(\til{\mathbf{c}},\mathbf{0})}^{-1} = \wh{B}$, which is an equality when applied to elements of the domain of the self-adjoint operator $\wh{B}$; we can write it as $\mathbf{S}_{(\til{\mathbf{c}},\mathbf{0})}^{-1} \, \wh{B} \, \mathbf{S}_{(\til{\mathbf{c}},\mathbf{0})} = B$ which is an equality when applied to elements of the domain of the of the self-adjoint operator $B$. Now Lem.\ref{lem:unitary_conjugation_commutes_with_functional_calculus} tells us 
$$
f^+(B) = f^+(\mathbf{S}_{(\til{\mathbf{c}},\mathbf{0})}^{-1} \, \wh{B} \, \mathbf{S}_{(\til{\mathbf{c}},\mathbf{0})}) = \mathbf{S}_{(\til{\mathbf{c}},\mathbf{0})}^{-1} \, f^+(\wh{B}) \, \mathbf{S}_{(\til{\mathbf{c}},\mathbf{0})},
$$
which is an equality of unitary operators. So, finally, what is $f^+(B) v$ for $v\in \mathscr{H}_\Gamma$? Let us write $v$ as the function $v(a_1,\ldots,a_n)$. Note $\mathbf{S}^{-1}_{(\til{\mathbf{c}},\mathbf{0})} = \mathbf{S}_{(\til{\mathbf{c}}^{-1},\mathbf{0})}$ where $\til{\mathbf{c}}^{-1} = (\til{c}^{\, ij})_{i,j\in \{1,\ldots,n\}}$ with $\til{c}^{\, ii}=1$, $\forall i=1,\ldots,n$, $\til{c}^{\, jk} = \frac{ \varepsilon_{kj} }{2\pi^2 \hbar_k}$, $\forall j\neq k$, $\til{c}^{\, ij}=0$ otherwise. Thus
\begin{align*}
& (f^+(B) v)(a_1,\ldots,a_n) \\
& = 
(\mathbf{S}_{(\til{\mathbf{c}},\mathbf{0})}^{-1} (f^+(\wh{B}) (\mathbf{S}_{(\til{\mathbf{c}},\mathbf{0})} v)))(a_1,\ldots,a_n) \\
& \stackrel{\eqref{eq:def_of_S_c_t}}{=}
(f^+(\wh{B}) (\mathbf{S}_{(\til{\mathbf{c}},\mathbf{0})} v))(a_1,\ldots,\underbrace{a_k+ \sum_{j\neq k} \frac{\varepsilon_{kj}}{2\pi^2\hbar_k} a_j}_{\mbox{\small $k$-th position}}, \ldots, a_n) \\
& =
f^+(-2\pi^2\hbar_k (a_k+ \sum_{j\neq k} \frac{\varepsilon_{kj}}{2\pi^2\hbar_k } a_j)) \cdot (\mathbf{S}_{(\til{\mathbf{c}},\mathbf{0})} v)(a_1,\ldots,a_k+ \sum_{j\neq k} \frac{ \varepsilon_{kj} }{2\pi^2\hbar_k} a_j, \ldots, a_n) \\
& \stackrel{\eqref{eq:def_of_S_c_t}}{=}
f^+( -2\pi^2\hbar_k a_k - \sum_{j\neq k} \varepsilon_{kj} a_j) \cdot v(a_1,\ldots,a_k,\ldots,a_n),
\end{align*}
for each $v\in \mathscr{H}_\Gamma$. Hence indeed $f^+(B)$ is multiplication by the function $f^+ \circ f_k = f^+(f_k)$, as desired. Similar proof shows that $f^-(\til{B})$ is multiplication by $f^-\circ \til{f}_k = f^-(\til{f}_k)$.

\vs

Now that both the unitary operators $f^+(B)$ and $f^-(\til{B})$ are multiplication operators, we easily see that they commute, and their composition is multiplication by the product of $f^+(f_k)$ and $f^-(\til{f}_k)$. This product function is
\begin{align*}
& f^+(f_k) \cdot f^-(\til{f}_k) = \exp( \frac{f_k^2}{4\pi\sqrt{-1} \hbar_k} ) \, \exp( - \frac{\til{f}_k^2}{4\pi \sqrt{-1} \hbar_k})
= \exp ( \frac{ f_k^2 - \til{f}_k^2 }{4\pi\sqrt{-1}\hbar_k} ) \\
& = \exp( \frac{ (f_k+\til{f}_k)(f_k - \til{f}_k)}{4\pi\sqrt{-1}\hbar_k})
= \exp( \frac{ - 4\pi^2 \hbar_k a_k (-2\sum_{j\neq k} \varepsilon_{kj} a_j)}{4\pi \sqrt{-1} \hbar_k} )
= \exp( -2 \pi \sqrt{-1} a_k \sum_{j\neq k} \varepsilon_{kj} a_j ).
\end{align*}
For convenience, define a real-valued function $g_k$ on $\mathbb{R}^n$ as
$$
g_k := - \sum_{j\neq k} \varepsilon_{kj} a_j
\qquad\mbox{so that}\qquad
f^+(f_k) \cdot f^-(\til{f}_k) = \exp(2\pi \sqrt{-1} a_k g_k)
$$
Observe now that for any $v\in \mathscr{H}_\Gamma$,
\begin{align*}
f^+(A) \, (f^-(\til{A}) v))
& = \mcal{F}_k^{-1} \, ( ( f^+(B) \, f^-(B)) \, (\mcal{F}_k v) ).
\end{align*}
In order to use the integral formula for the Fourier transform and its inverse, we restrict to $v\in D_\Gamma$; then $\mcal{F}_K v \in D_\Gamma$, so $( f^+(B) \, f^-(B)) \, (\mcal{F}_k v) \in L^1(\mathbb{R}^n) \cap \mathscr{H}_\Gamma$, hence 
\begin{align*}
(f^+(A) \, (f^-(\til{A}) v)))(a_1,\ldots,a_n)
& = \int_\mathbb{R} e^{2\pi\sqrt{-1} a_k b_k} ( ( f^+(B) \, f^-(B)) \, (\mcal{F}_k v) )(a_1,\ldots,b_k,\ldots,a_n) \, db_k \\
& = \int_\mathbb{R} e^{2\pi\sqrt{-1}a_kb_k} \, e^{2\pi\sqrt{-1}b_kg_k} \, (\mcal{F}_kv)(a_1,\ldots,b_k,\ldots,a_n) \, db_k \\
& = \int_{\mathbb{R}} e^{2\pi\sqrt{-1}(a_k+g_k)b_k} \, (\mcal{F}_kv)(a_1,\ldots,b_k,\ldots,a_n) \, db_k \\
& = ( \mcal{F}_k^{-1} (\mcal{F}_k v) )(a_1,\ldots,a_k+g_k,\ldots,a_n) \\
& = v(a_1,\ldots,\, \underbrace{a_k -\sum_{j\neq k} \varepsilon_{kj} a_j}_{\mbox{\small $k$-th position}}, \, \ldots a_n).
\end{align*}
The final computational result can then be extended to any $v\in \mathscr{H}_\Gamma$, for it is manifestly an expression of a special affine shift operator, which is unitary. \qed

\begin{remark}
A formal heuristic proof can be much shorter, which can be found in \cite{K16}. But here I wanted to have a rigorous proof that is not hand-waiving, as it will be used crucially.
\end{remark}

\subsection{$A_1\times A_1$ identity: the `commuting' relation}

\begin{proposition}
\label{prop:c_A1_times_A1_is_1}
$c_{A_1\times A_1}=1$.
\end{proposition}

{\it Proof.} The LHS of \eqref{eq:A1_times_A1_to_prove} is
\begin{align*}
& \mathbf{K}^\hbar_{\Gamma^{(0)}\mut{i}\Gamma^{(1)}} \, \mathbf{K}^\hbar_{\Gamma^{(1)}\mut{j}\Gamma^{(2)}}
= (\mathbf{K}^{\sharp \hbar}_{\Gamma^{(0)}\mut{i}\Gamma^{(1)}} \, \mathbf{K}'_{\Gamma^{(0)}\mut{i}\Gamma^{(1)}}) (\mathbf{K}^{\sharp \hbar}_{\Gamma^{(1)}\mut{j}\Gamma^{(2)}} \, \mathbf{K}'_{\Gamma^{(1)}\mut{j}\Gamma^{(2)}}) \\
& = \Phi^{\hbar_i}(\wh{\mathbf{x}}^\hbar_{\Gamma^{(0)};i}) \, \Phi^{\hbar_i}(\wh{\til{\mathbf{x}}}^\hbar_{\Gamma^{(0)};i})^{-1} \, \ul{ \mathbf{K}'_{\Gamma^{(0)}\mut{i}\Gamma^{(1)}} }) \, ( \ul{ \Phi^{\hbar_j}(\wh{\mathbf{x}}^\hbar_{\Gamma^{(1)};j}) \, \Phi^{\hbar_j}(\wh{\til{\mathbf{x}}}^\hbar_{\Gamma^{(1)};j})^{-1} } \, \mathbf{K}'_{\Gamma^{(1)}\mut{j}\Gamma^{(2)}}) \\
& = \Phi^{\hbar_i}(\wh{\mathbf{x}}^\hbar_{\Gamma^{(0)};i}) \, \Phi^{\hbar_i}(\wh{\til{\mathbf{x}}}^\hbar_{\Gamma^{(0)};i})^{-1}  \,( \Phi^{\hbar_j}(\wh{\mathbf{x}}^\hbar_{\Gamma^{(0)};j}) \, \Phi^{\hbar_j}(\wh{\til{\mathbf{x}}}^\hbar_{\Gamma^{(0)};j})^{-1} \,  \mathbf{K}'_{\Gamma^{(0)}\mut{i}\Gamma^{(1)}} ) \, \mathbf{K}'_{\Gamma^{(1)}\mut{j}\Gamma^{(2)}}, \quad (\because {\rm Lem}.\ref{lem:unitary_conjugation_commutes_with_functional_calculus}, \ref{lem:bf_K_prime_conjugation_on_bf_x}, ~ \eqref{eq:A1_times_A1_condition})
\end{align*}
where we are implicitly using the fact that each of $d_i$ and $d_j$ is the same for all $\Gamma^{(\ell)}$. Meanwhile, the RHS of \eqref{eq:A1_times_A1_to_prove} without the constant is
\begin{align*}
& \mathbf{K}^\hbar_{\Gamma^{(0)}\mut{j}\Gamma^{(3)}} \, \mathbf{K}^\hbar_{\Gamma^{(3)}\mut{i}\Gamma^{(2)}}
= (\mathbf{K}^{\sharp \hbar}_{\Gamma^{(0)}\mut{j}\Gamma^{(3)}} \, \mathbf{K}'_{\Gamma^{(0)}\mut{j}\Gamma^{(3)}}) (\mathbf{K}^{\sharp \hbar}_{\Gamma^{(3)}\mut{i}\Gamma^{(2)}} \, \mathbf{K}'_{\Gamma^{(3)}\mut{i}\Gamma^{(2)}}) \\
& = \Phi^{\hbar_j}(\wh{\mathbf{x}}^\hbar_{\Gamma^{(0)};j}) \, \Phi^{\hbar_j}(\wh{\til{\mathbf{x}}}^\hbar_{\Gamma^{(0)};j})^{-1} \, \ul{ \mathbf{K}'_{\Gamma^{(0)}\mut{j}\Gamma^{(3)}} }) \, ( \ul{ \Phi^{\hbar_i}(\wh{\mathbf{x}}^\hbar_{\Gamma^{(3)};i}) \, \Phi^{\hbar_i}(\wh{\til{\mathbf{x}}}^\hbar_{\Gamma^{(3)};i})^{-1} } \, \mathbf{K}'_{\Gamma^{(3)}\mut{i}\Gamma^{(2)}}) \\
& = \Phi^{\hbar_j}(\wh{\mathbf{x}}^\hbar_{\Gamma^{(0)};j}) \, \Phi^{\hbar_j}(\wh{\til{\mathbf{x}}}^\hbar_{\Gamma^{(0)};j})^{-1} \,( \Phi^{\hbar_i}(\wh{\mathbf{x}}^\hbar_{\Gamma^{(0)};i}) \, \Phi^{\hbar_i}(\wh{\til{\mathbf{x}}}^\hbar_{\Gamma^{(0)};i})^{-1} \,  \mathbf{K}'_{\Gamma^{(0)}\mut{j}\Gamma^{(3)}} ) \, \mathbf{K}'_{\Gamma^{(3)}\mut{i}\Gamma^{(2)}}. \quad (\because {\rm Lem}.\ref{lem:unitary_conjugation_commutes_with_functional_calculus}, \ref{lem:bf_K_prime_conjugation_on_bf_x}, ~ \eqref{eq:A1_times_A1_condition})
\end{align*}
One can see that each of $\Phi^{\hbar_j}(\wh{\mathbf{x}}^\hbar_{\Gamma^{(0)};j})$ and $\Phi^{\hbar_j}(\wh{\til{\mathbf{x}}}^\hbar_{\Gamma^{(0)};j})^{-1}$ commutes with each of $\Phi^{\hbar_i}(\wh{\mathbf{x}}^\hbar_{\Gamma^{(0)};i})$ and $\Phi^{\hbar_i}(\wh{\til{\mathbf{x}}}^\hbar_{\Gamma^{(0)};i})^{-1}$, from Lem.\ref{lem:commuting} and the Heisenberg commutation relations of the arguments $\wh{\mathbf{x}}^\hbar_{\Gamma^{(0)};j}$, $\wh{\til{\mathbf{x}}}^\hbar_{\Gamma^{(0)};j}$, $\wh{\mathbf{x}}^\hbar_{\Gamma^{(0)};i}$, and $\wh{\til{\mathbf{x}}}^\hbar_{\Gamma^{(0)};i}$. So one has
\begin{align*}
\mathbf{K}^\hbar_{\Gamma^{(0)}\mut{j}\Gamma^{(3)}} \, \mathbf{K}^\hbar_{\Gamma^{(3)}\mut{i}\Gamma^{(2)}} =  \Phi^{\hbar_i}(\wh{\mathbf{x}}^\hbar_{\Gamma^{(0)};i}) \, \Phi^{\hbar_i}(\wh{\til{\mathbf{x}}}^\hbar_{\Gamma^{(0)};i})^{-1} \, \Phi^{\hbar_j}(\wh{\mathbf{x}}^\hbar_{\Gamma^{(0)};j}) \, \Phi^{\hbar_j}(\wh{\til{\mathbf{x}}}^\hbar_{\Gamma^{(0)};j})^{-1}\,  \mathbf{K}'_{\Gamma^{(0)}\mut{j}\Gamma^{(3)}}  \, \mathbf{K}'_{\Gamma^{(3)}\mut{i}\Gamma^{(2)}},
\end{align*}
and hence
\begin{align*}
& (\mathbf{K}^\hbar_{\Gamma^{(0)}\mut{i}\Gamma^{(1)}} \, \mathbf{K}^\hbar_{\Gamma^{(1)}\mut{j}\Gamma^{(2)}})^{-1} \, (\mathbf{K}^\hbar_{\Gamma^{(0)}\mut{j}\Gamma^{(3)}} \, \mathbf{K}^\hbar_{\Gamma^{(3)}\mut{i}\Gamma^{(2)}}) \\
& = (\mathbf{K}'_{\Gamma^{(1)}\mut{j}\Gamma^{(2)}})^{-1} \, ( \mathbf{K}'_{\Gamma^{(0)}\mut{i}\Gamma^{(1)}} )^{-1} \, \Phi^{\hbar_j}(\wh{\til{\mathbf{x}}}^\hbar_{\Gamma^{(0)};j}) \, \Phi^{\hbar_j}(\wh{\mathbf{x}}^\hbar_{\Gamma^{(0)};j})^{-1} \, \Phi^{\hbar_i}(\wh{\til{\mathbf{x}}}^\hbar_{\Gamma^{(0)};i}) \, \Phi^{\hbar_i}(\wh{\mathbf{x}}^\hbar_{\Gamma^{(0)};i})^{-1} \\
& \quad \cdot \Phi^{\hbar_i}(\wh{\mathbf{x}}^\hbar_{\Gamma^{(0)};i}) \, \Phi^{\hbar_i}(\wh{\til{\mathbf{x}}}^\hbar_{\Gamma^{(0)};i})^{-1} \, \Phi^{\hbar_j}(\wh{\mathbf{x}}^\hbar_{\Gamma^{(0)};j}) \, \Phi^{\hbar_j}(\wh{\til{\mathbf{x}}}^\hbar_{\Gamma^{(0)};j})^{-1}\,  \mathbf{K}'_{\Gamma^{(0)}\mut{j}\Gamma^{(3)}}  \, \mathbf{K}'_{\Gamma^{(3)}\mut{i}\Gamma^{(2)}} \\
& = (\mathbf{K}'_{\Gamma^{(1)}\mut{j}\Gamma^{(2)}})^{-1} \, ( \mathbf{K}'_{\Gamma^{(0)}\mut{i}\Gamma^{(1)}} )^{-1} \, \mathbf{K}'_{\Gamma^{(0)}\mut{j}\Gamma^{(3)}}  \, \mathbf{K}'_{\Gamma^{(3)}\mut{i}\Gamma^{(2)}},
\end{align*}
which we know to equal $c_{A_1\times A_1}^{-1} \cdot \mathrm{id}_{\mathscr{H}_{\Gamma^{(2)}}}$ from \eqref{eq:A1_times_A1_to_prove} of Prop.\ref{prop:FG_A1_times_A1_identity}. As each of the four operators $(\mathbf{K}'_{\Gamma^{(1)}\mut{j}\Gamma^{(2)}})^{-1}$, $( \mathbf{K}'_{\Gamma^{(0)}\mut{i}\Gamma^{(1)}} )^{-1}$, $\mathbf{K}'_{\Gamma^{(0)}\mut{j}\Gamma^{(3)}}$, and $\mathbf{K}'_{\Gamma^{(3)}\mut{i}\Gamma^{(2)}}$ are special affine shift operators if we identify the Hilbert spaces by the relabling isomorphisms \eqref{eq:bf_I}, so is their composition, by Cor.\ref{cor:special_affine_shift_operators_form_a_group}. From Lem.\ref{lem:scalar_special_affine_shift_operator_is_identity} we conclude that this is the identity operator. Again, if one would like, one may directly compute $(\mathbf{K}'_{\Gamma^{(1)}\mut{j}\Gamma^{(2)}})^{-1} \, ( \mathbf{K}'_{\Gamma^{(0)}\mut{i}\Gamma^{(1)}} )^{-1} \, \mathbf{K}'_{\Gamma^{(0)}\mut{j}\Gamma^{(3)}}  \, \mathbf{K}'_{\Gamma^{(3)}\mut{i}\Gamma^{(2)}}$ as a special affine shift operator and prove that it is the identity. \qed

\subsection{$A_2$ identity: the `pentagon' relation}

\begin{proposition}
\label{prop:c_A2_is_1}
$c_{A_2}=1$.
\end{proposition}

{\it Proof.} Let us work with the case  $\varepsilon^{(0)}_{ij}=1=-\varepsilon^{(0)}_{ji}$; I claim that a proof for the other case is essentially the same. From \eqref{eq:varepsilon_prime_formula} we get
\begin{align*}
\varepsilon^{(1)}_{ij} = -1 = - \varepsilon^{(1)}_{ji}, \quad
\varepsilon^{(2)}_{ij} = 1 = - \varepsilon^{(2)}_{ji}, \quad
\varepsilon^{(3)}_{ij} = -1 = - \varepsilon^{(3)}_{ji}, \quad
\varepsilon^{(4)}_{ij} = -1 = - \varepsilon^{(4)}_{ji}, \quad
\varepsilon^{(5)}_{ij} = 1 = - \varepsilon^{(5)}_{ji}.
\end{align*}
The LHS of \eqref{eq:A2_to_prove} is
\begin{align*}
& \mathbf{K}^\hbar_{\Gamma^{(0)}\mut{i}\Gamma^{(1)}} \, \mathbf{K}^\hbar_{\Gamma^{(1)}\mut{j}\Gamma^{(2)}} \, \mathbf{K}^\hbar_{\Gamma^{(2)}\mut{i}\Gamma^{(3)}}  \\
& = (\mathbf{K}^{\sharp \hbar}_{\Gamma^{(0)}\mut{i}\Gamma^{(1)}} \, \mathbf{K}'_{\Gamma^{(0)}\mut{i}\Gamma^{(1)}}) (\mathbf{K}^{\sharp \hbar}_{\Gamma^{(1)}\mut{j}\Gamma^{(2)}} \, \mathbf{K}'_{\Gamma^{(1)}\mut{j}\Gamma^{(2)}}) (\mathbf{K}^{\sharp \hbar}_{\Gamma^{(2)}\mut{i}\Gamma^{(3)}} \, \mathbf{K}'_{\Gamma^{(2)}\mut{i}\Gamma^{(3)}}) \\
& = \Phi^{\hbar_i}(\wh{\mathbf{x}}^\hbar_{\Gamma^{(0)};i}) \, \Phi^{\hbar_i}(\wh{\til{\mathbf{x}}}^\hbar_{\Gamma^{(0)};i})^{-1} \, \mathbf{K}'_{\Gamma^{(0)}\mut{i}\Gamma^{(1)}} \, \Phi^{\hbar_j}(\wh{\mathbf{x}}^\hbar_{\Gamma^{(1)};j}) \, \Phi^{\hbar_j}(\wh{\til{\mathbf{x}}}^\hbar_{\Gamma^{(1)};j})^{-1} \, \underbrace{ \mathbf{K}'_{\Gamma^{(1)}\mut{j}\Gamma^{(2)}} }_{\mbox{\tiny move to right}} \,  \\
& \quad \cdot \ul{ \Phi^{\hbar_i}(\wh{\mathbf{x}}^\hbar_{\Gamma^{(2)};i}) \, \Phi^{\hbar_i}(\wh{\til{\mathbf{x}}}^\hbar_{\Gamma^{(2)};i})^{-1} } \, \mathbf{K}'_{\Gamma^{(2)}\mut{i}\Gamma^{(3)}} \\
& = \Phi^{\hbar_i}(\wh{\mathbf{x}}^\hbar_{\Gamma^{(0)};i}) \, \Phi^{\hbar_i}(\wh{\til{\mathbf{x}}}^\hbar_{\Gamma^{(0)};i})^{-1} \, \underbrace{ \mathbf{K}'_{\Gamma^{(0)}\mut{i}\Gamma^{(1)}} }_{\mbox{\tiny move to right}}\, \Phi^{\hbar_j}(\wh{\mathbf{x}}^\hbar_{\Gamma^{(1)};j}) \, \Phi^{\hbar_j}(\wh{\til{\mathbf{x}}}^\hbar_{\Gamma^{(1)};j})^{-1} \,  \\
& \quad \cdot ( \Phi^{\hbar_i}(\wh{\mathbf{x}}^\hbar_{\Gamma^{(1)};i} + \underbrace{ [\varepsilon^{(1)}_{ij}]_+ }_{=0} \, \wh{\mathbf{x}}^\hbar_{\Gamma^{(1)};j}) \, \Phi^{\hbar_i}(\wh{\til{\mathbf{x}}}^\hbar_{\Gamma^{(1)};i} + \underbrace{  [\varepsilon^{(1)}_{ij}]_+ }_{=0} \, \wh{\til{\mathbf{x}}}^\hbar_{\Gamma^{(1)};j})^{-1} \, \, \mathbf{K}'_{\Gamma^{(1)}\mut{j}\Gamma^{(2)}} ) \, \mathbf{K}'_{\Gamma^{(2)}\mut{i}\Gamma^{(3)}} \quad (\because {\rm Lem}.\ref{lem:unitary_conjugation_commutes_with_functional_calculus}, \ref{lem:bf_K_prime_conjugation_on_bf_x}) \\
& = \Phi^{\hbar_i}(\wh{\mathbf{x}}^\hbar_{\Gamma^{(0)};i}) \, \Phi^{\hbar_i}(\wh{\til{\mathbf{x}}}^\hbar_{\Gamma^{(0)};i})^{-1} \, \Phi^{\hbar_j}(\wh{\mathbf{x}}^\hbar_{\Gamma^{(0)};j}+\underbrace{ [\varepsilon^{(0)}_{ji}]_+}_{=0} \, \wh{\mathbf{x}}^\hbar_{\Gamma^{(0)};i} ) \, \Phi^{\hbar_j}(\wh{\til{\mathbf{x}}}^\hbar_{\Gamma^{(0)};j} + \underbrace{ [\varepsilon^{(0)}_{ji}]_+}_{=0} \, \wh{\til{\mathbf{x}}}^\hbar_{\Gamma^{(0)};i} )^{-1} \,  \\
& \quad \cdot \Phi^{\hbar_i}(- \wh{\mathbf{x}}^\hbar_{\Gamma^{(0)};i} ) \, \Phi^{\hbar_i}(-\wh{\til{\mathbf{x}}}^\hbar_{\Gamma^{(0)};i} )^{-1} \,  \mathbf{K}'_{\Gamma^{(0)}\mut{i}\Gamma^{(1)}} \, \mathbf{K}'_{\Gamma^{(1)}\mut{j}\Gamma^{(2)}} \, \mathbf{K}'_{\Gamma^{(2)}\mut{i}\Gamma^{(3)}} \quad (\because {\rm Lem}.\ref{lem:unitary_conjugation_commutes_with_functional_calculus}, \ref{lem:bf_K_prime_conjugation_on_bf_x}) \\
& = \Phi^{\hbar_i}(\wh{\mathbf{x}}^\hbar_{\Gamma^{(0)};i}) \, \Phi^{\hbar_i}(\wh{\til{\mathbf{x}}}^\hbar_{\Gamma^{(0)};i})^{-1} \, \Phi^{\hbar_i}(\wh{\mathbf{x}}^\hbar_{\Gamma^{(0)};j}) \, \Phi^{\hbar_i}(\wh{\til{\mathbf{x}}}^\hbar_{\Gamma^{(0)};j} )^{-1} \,   \Phi^{\hbar_i}(- \wh{\mathbf{x}}^\hbar_{\Gamma^{(0)};i} ) \, \Phi^{\hbar_i}(-\wh{\til{\mathbf{x}}}^\hbar_{\Gamma^{(0)};i} )^{-1} \\
& \quad \cdot \mathbf{K}'_{\Gamma^{(0)}\mut{i}\Gamma^{(1)}} \, \mathbf{K}'_{\Gamma^{(1)}\mut{j}\Gamma^{(2)}} \, \mathbf{K}'_{\Gamma^{(2)}\mut{i}\Gamma^{(3)}}. \qquad (\because \hbar_i = \hbar_j, \mbox{ because } d_i=d_j)
\end{align*}
We would now like to move the factors $\Phi^{\hbar_i}(\wh{\mathbf{x}}^\hbar_{\Gamma^{(0)};j})$ and $\Phi^{\hbar_i}(- \wh{\mathbf{x}}^\hbar_{\Gamma^{(0)};i} )$ to the left, so that we end up in an expression of the form $\Phi( \, ) \Phi( \, ) \Phi( \, ) \Phi( \, )^{-1} \Phi(\, )^{-1} \, \Phi(\, )^{-1} \, \mathbf{K}' \mathbf{K}' \mathbf{K}'$. To do this, we must show for example that $\Phi^{\hbar_i}(\wh{\til{\mathbf{x}}}^\hbar_{\Gamma^{(0)};i})^{-1}$ and $\Phi^{\hbar_i}(\wh{\mathbf{x}}^\hbar_{\Gamma^{(0)};j})$ commute; this follows from Lem.\ref{lem:commuting} and the Heisenberg commutation relation of the arguments $\wh{\til{\mathbf{x}}}^\hbar_{\Gamma^{(0)};i}$ and $\wh{\mathbf{x}}^\hbar_{\Gamma^{(0)};j}$. So the LHS of \eqref{eq:A2_to_prove} becomes
\begin{align*}
& \ul{ \Phi^{\hbar_i}(\wh{\mathbf{x}}^\hbar_{\Gamma^{(0)};i}) \, \Phi^{\hbar_i}(\wh{\mathbf{x}}^\hbar_{\Gamma^{(0)};j}) \,   \Phi^{\hbar_i}(- \wh{\mathbf{x}}^\hbar_{\Gamma^{(0)};i} ) }  \, \, \ul{ \Phi^{\hbar_i}(\wh{\til{\mathbf{x}}}^\hbar_{\Gamma^{(0)};i})^{-1} \, \Phi^{\hbar_i}(\wh{\til{\mathbf{x}}}^\hbar_{\Gamma^{(0)};j} )^{-1} \, \Phi^{\hbar_i}(-\wh{\til{\mathbf{x}}}^\hbar_{\Gamma^{(0)};i} )^{-1} } \\
& \quad \cdot \mathbf{K}'_{\Gamma^{(0)}\mut{i}\Gamma^{(1)}} \, \mathbf{K}'_{\Gamma^{(1)}\mut{j}\Gamma^{(2)}} \, \mathbf{K}'_{\Gamma^{(2)}\mut{i}\Gamma^{(3)}}.
\end{align*}
Note from \eqref{eq:Heisenberg_relations_for_D_q_Gamma}, $\varepsilon^{(0)}_{ij}=1$, and $\hbar_j = \hbar_i$ that $[\wh{\mathbf{x}}^\hbar_{\Gamma^{(0)};i}, \wh{\mathbf{x}}^\hbar_{\Gamma^{(0)};j}]= 2\pi\sqrt{-1} \hbar_i \, \cdot \mathrm{id}$ on $D_\Gamma$, and that the unique self-adjoint extensions of $\wh{\mathbf{x}}^\hbar_{\Gamma^{(0)};i}$ and $\wh{\mathbf{x}}^\hbar_{\Gamma^{(0)};j}$ satisfy the corresponding Weyl relations are satisfied (Lem.\ref{lem:old_representation}). Likewise, we have $[-\wh{\til{\mathbf{x}}}^\hbar_{\Gamma^{(0)};i}, \wh{\til{\mathbf{x}}}^\hbar_{\Gamma^{(0)};j}]=2\pi\sqrt{-1}\hbar_i \cdot \mathrm{id}$, together with the corresponding Weyl relations. Hence Cor.\ref{cor:QD_identity_pentagon2} applies to the two underlined parts above, so that the LHS of \eqref{eq:A2_to_prove} becomes
\begin{align*}
& \cancel{ c_{\hbar_i} } \, \cancel{ c_{\hbar_i}^{-1} } \, \Phi^{\hbar_i}(\wh{\mathbf{x}}^\hbar_{\Gamma^{(0)};j}) \, \Phi^{\hbar_i}(\wh{\mathbf{x}}^\hbar_{\Gamma^{(0)};j} + \wh{\mathbf{x}}^\hbar_{\Gamma^{(0)};i}) \, \underbrace{ \exp\left( \frac{(\wh{\mathbf{x}}^\hbar_{\Gamma^{(0)};i})^2}{4\pi \sqrt{-1} \hbar_i} \right) \, \exp\left( - \frac{(-\wh{\til{\mathbf{x}}}^\hbar_{\Gamma^{(0)};i})^2}{4\pi\sqrt{-1}\hbar_i} \right) } \\
& \quad  \cdot \, \Phi^{\hbar_i}(\wh{\til{\mathbf{x}}}^\hbar_{\Gamma^{(0)};j} - \wh{\til{\mathbf{x}}}^\hbar_{\Gamma^{(0)};i})^{-1} \, \Phi^{\hbar_i}(\wh{\til{\mathbf{x}}}^\hbar_{\Gamma^{(0)};j})^{-1}
\, \mathbf{K}'_{\Gamma^{(0)}\mut{i}\Gamma^{(1)}} \, \mathbf{K}'_{\Gamma^{(1)}\mut{j}\Gamma^{(2)}} \, \mathbf{K}'_{\Gamma^{(2)}\mut{i}\Gamma^{(3)}}.
\end{align*}
Notice the cancellation of the two $c_{\hbar_i}$'s. To the underbraced part we apply Lem.\ref{lem:product_of_two_quadratic_exponential_operators}, to replace it by $\mathbf{S}_{(\mathbf{c},\mathbf{0})}$, where $\mathbf{c}=(c_{\ell m})_{\ell, m\in \{1,\ldots,n\}}$ is given by
$$
\mathbf{c} ~ : ~ \quad c_{\ell\ell}=1, \quad \forall \ell=1,\ldots,n, \qquad
c_{\ell i} = -\varepsilon^{(0)}_{i\ell}, \quad \forall \ell \neq i, \qquad c_{\ell m} = 0 \quad \mbox{otherwise}.
$$
So its inverse $\mathbf{c}^{-1} = (c^{\ell m})_{\ell, m\in \{1,\ldots,n\}}$ is given by $c^{\ell\ell}=1$, $\forall \ell=1,\ldots,n$, $c^{\ell i} = \varepsilon^{(0)}_{i\ell}$, $\forall \ell \neq i$, $c^{\ell m}=0$ otherwise. Hence from \eqref{eq:conjugation_of_bf_S_on_bf_p}, \eqref{eq:conjugation_of_bf_S_on_bf_q}, and $\varepsilon^{(0)}_{ij}=1$ we get 
\begin{align*}
& \mathbf{S}_{(\mathbf{c},\mathbf{0})} \, \wh{\mathbf{p}}^\hbar_i \, \mathbf{S}_{(\mathbf{c},\mathbf{0})}^{-1} = \wh{\mathbf{p}}^\hbar_i, \qquad
\mathbf{S}_{(\mathbf{c},\mathbf{0})} \, \wh{\mathbf{p}}^\hbar_j \, \mathbf{S}_{(\mathbf{c},\mathbf{0})}^{-1} = \wh{\mathbf{p}}^\hbar_j+ \wh{\mathbf{p}}^\hbar_i, \\
& \mathbf{S}_{(\mathbf{c},\mathbf{0})} \, \wh{\mathbf{q}}^\hbar_\ell \, \mathbf{S}_{(\mathbf{c},\mathbf{0})}^{-1} = \wh{\mathbf{q}}^\hbar_\ell + \delta_{\ell,i} \, \sum_{m \neq i} (-\varepsilon^{(0)}_{im}) \wh{\mathbf{q}}^\hbar_m.
\end{align*}
Thus from \eqref{eq:old_representation_tilde} we have
\begin{align*}
\mathbf{S}_{(\mathbf{c},\mathbf{0})} \, \wh{\til{\mathbf{x}}}^\hbar_{\Gamma^{(0)};i} \, \mathbf{S}_{(\mathbf{c},\mathbf{0})}^{-1} & = \mathbf{S}_{(\mathbf{c},\mathbf{0})} ( \frac{d_i^{-1}}{2} \wh{\mathbf{p}}^\hbar_i + \sum_{\ell=1}^n \varepsilon^{(0)}_{i\ell} \wh{\mathbf{q}}^\hbar_\ell) \mathbf{S}_{(\mathbf{c},\mathbf{0})}^{-1} = \frac{d_i^{-1} }{2} \wh{\mathbf{p}}^\hbar_i + \sum_{\ell=1}^n \varepsilon^{(0)}_{i\ell} \wh{\mathbf{q}}^\hbar_\ell = \wh{\til{\mathbf{x}}}^\hbar_{\Gamma^{(0)};i}, \\
\mathbf{S}_{(\mathbf{c},\mathbf{0})} \, \wh{\til{\mathbf{x}}}^\hbar_{\Gamma^{(0)};j} \, \mathbf{S}_{(\mathbf{c},\mathbf{0})}^{-1} & = \mathbf{S}_{(\mathbf{c},\mathbf{0})} ( \frac{d_j^{-1} }{2} \wh{\mathbf{p}}^\hbar_j + \sum_{\ell=1}^n \varepsilon^{(0)}_{j\ell} \wh{\mathbf{q}}^\hbar_\ell) \mathbf{S}_{(\mathbf{c},\mathbf{0})}^{-1} \\
& = \frac{d_j^{-1}}{2}  (\wh{\mathbf{p}}^\hbar_j + \wh{\mathbf{p}}^\hbar_i) + \sum_{\ell=1}^n \varepsilon^{(0)}_{j\ell} \wh{\mathbf{q}}^\hbar_\ell - \sum_{m\neq i} (-\varepsilon^{(0)}_{im}) \wh{\mathbf{q}}^\hbar_m \\
& = (\frac{d_j^{-1}}{2}\wh{\mathbf{p}}^\hbar_j + \sum_{\ell=1}^n \varepsilon^{(0)}_{j\ell} \wh{\mathbf{q}}^\hbar_\ell) 
+ (\frac{d_i^{-1}}{2}\wh{\mathbf{p}}^\hbar_i + \sum_{m=1}^n \varepsilon^{(0)}_{im} \wh{\mathbf{q}}^\hbar_m) = \wh{\til{\mathbf{x}}}^\hbar_{\Gamma^{(0)};j} + \wh{\til{\mathbf{x}}}^\hbar_{\Gamma^{(0)};i},
\end{align*}
as operators $D_\Gamma\to D_\Gamma$, where we used $\varepsilon^{(0)}_{ii}=0$, $\varepsilon^{(0)}_{ji}=-1$, and $d_i=d_j$; such conjugation identities extend to the respective unique self-adjoint extensions. So the LHS of \eqref{eq:A2_to_prove} becomes
\begin{align*}
& \hspace{-7mm} \Phi^{\hbar_i}(\wh{\mathbf{x}}^\hbar_{\Gamma^{(0)};j}) \, 
\Phi^{\hbar_i}(\wh{\mathbf{x}}^\hbar_{\Gamma^{(0)};j} + \wh{\mathbf{x}}^\hbar_{\Gamma^{(0)};i}) \, 
\ul{ \mathbf{S}_{(\mathbf{c},\mathbf{0})} } \, \,
\ul{ \Phi^{\hbar_i}(\wh{\til{\mathbf{x}}}^\hbar_{\Gamma^{(0)};j} - \wh{\til{\mathbf{x}}}^\hbar_{\Gamma^{(0)};i})^{-1} \, 
\Phi^{\hbar_i}(\wh{\til{\mathbf{x}}}^\hbar_{\Gamma^{(0)};j})^{-1} } \,
\mathbf{K}'_{\Gamma^{(0)}\mut{i}\Gamma^{(1)}} \, \mathbf{K}'_{\Gamma^{(1)}\mut{j}\Gamma^{(2)}} \, \mathbf{K}'_{\Gamma^{(2)}\mut{i}\Gamma^{(3)}} \\
& = \Phi^{\hbar_i}(\wh{\mathbf{x}}^\hbar_{\Gamma^{(0)};j}) \, 
\Phi^{\hbar_i}(\wh{\mathbf{x}}^\hbar_{\Gamma^{(0)};j} + \wh{\mathbf{x}}^\hbar_{\Gamma^{(0)};i}) \, 
\Phi^{\hbar_i}(\wh{\til{\mathbf{x}}}^\hbar_{\Gamma^{(0)};j} + \cancel{ \wh{\til{\mathbf{x}}}^\hbar_{\Gamma^{(0)};i} } - \cancel{ \wh{\til{\mathbf{x}}}^\hbar_{\Gamma^{(0)};i} })^{-1} \, 
\Phi^{\hbar_i}(\wh{\til{\mathbf{x}}}^\hbar_{\Gamma^{(0)};j} + \wh{\til{\mathbf{x}}}^\hbar_{\Gamma^{(0)};i})^{-1}
\mathbf{S}_{(\mathbf{c},\mathbf{0})} \\
& \quad \cdot \mathbf{K}'_{\Gamma^{(0)}\mut{i}\Gamma^{(1)}} \, \mathbf{K}'_{\Gamma^{(1)}\mut{j}\Gamma^{(2)}} \, \mathbf{K}'_{\Gamma^{(2)}\mut{i}\Gamma^{(3)}},
\end{align*}
thanks to Lem.\ref{lem:unitary_conjugation_commutes_with_functional_calculus}.

\vs

Meanwhile, we investigate the following factor of the RHS of \eqref{eq:A2_to_prove}:
\begin{align*}
& \mathbf{K}^\hbar_{\Gamma^{(0)}\mut{j}\Gamma^{(4)}} \, \mathbf{K}^\hbar_{\Gamma^{(4)}\mut{i}\Gamma^{(5)}}
= (\mathbf{K}^{\sharp \hbar}_{\Gamma^{(0)}\mut{j}\Gamma^{(4)}} \, \mathbf{K}'_{\Gamma^{(0)}\mut{j}\Gamma^{(4)}}) (\mathbf{K}^{\sharp \hbar}_{\Gamma^{(4)}\mut{i}\Gamma^{(5)}} \, \mathbf{K}'_{\Gamma^{(4)}\mut{i}\Gamma^{(5)}})  \\
& = \Phi^{\hbar_j}(\wh{\mathbf{x}}^\hbar_{\Gamma^{(0)};j}) \Phi^{\hbar_j}(\wh{\til{\mathbf{x}}}^\hbar_{\Gamma^{(0)};j})^{-1} \, \ul{ \mathbf{K}'_{\Gamma^{(0)}\mut{j}\Gamma^{(4)}} } \,\, \ul{ \Phi^{\hbar_i}(\wh{\mathbf{x}}^\hbar_{\Gamma^{(4)};i}) \, \Phi^{\hbar_i} (\wh{\til{\mathbf{x}}}^\hbar_{\Gamma^{(4)};i})^{-1} } \, \mathbf{K}'_{\Gamma^{(4)}\mut{i}\Gamma^{(5)}} \\
& = 
\Phi^{\hbar_j}(\wh{\mathbf{x}}^\hbar_{\Gamma^{(0)};j}) \Phi^{\hbar_j}(\wh{\til{\mathbf{x}}}^\hbar_{\Gamma^{(0)};j})^{-1} \Phi^{\hbar_i}(\wh{\mathbf{x}}^\hbar_{\Gamma^{(0)};i}+[\varepsilon^{(0)}_{ij}]_+\, \wh{\mathbf{x}}^\hbar_{\Gamma^{(0)};j}) \, \Phi^{\hbar_i} (\wh{\til{\mathbf{x}}}^\hbar_{\Gamma^{(0)};i} + [\varepsilon^{(0)}_{ij}]_+ \, \wh{\til{\mathbf{x}}}^\hbar_{\Gamma^{(0)};j} )^{-1} \\
& \quad \cdot \mathbf{K}'_{\Gamma^{(0)}\mut{j}\Gamma^{(4)}} \, \mathbf{K}'_{\Gamma^{(4)}\mut{i}\Gamma^{(5)}}
\quad (\because {\rm Lem}.\ref{lem:unitary_conjugation_commutes_with_functional_calculus}, \ref{lem:bf_K_prime_conjugation_on_bf_x}) \\ 
& = \Phi^{\hbar_i}(\wh{\mathbf{x}}^\hbar_{\Gamma^{(0)};j})  \Phi^{\hbar_i}(\wh{\mathbf{x}}^\hbar_{\Gamma^{(0)};i}+ \wh{\mathbf{x}}^\hbar_{\Gamma^{(0)};j}) \Phi^{\hbar_i}(\wh{\til{\mathbf{x}}}^\hbar_{\Gamma^{(0)};j})^{-1} \Phi^{\hbar_i} (\wh{\til{\mathbf{x}}}^\hbar_{\Gamma^{(0)};i} + \wh{\til{\mathbf{x}}}^\hbar_{\Gamma^{(0)};j} )^{-1} \mathbf{K}'_{\Gamma^{(0)}\mut{j}\Gamma^{(4)}} \, \mathbf{K}'_{\Gamma^{(4)}\mut{i}\Gamma^{(5)}},
\end{align*}
where we used $\hbar_i=\hbar_j$, $\varepsilon^{(0)}_{ij}=1$, and the fact that $\Phi^{\hbar_i}(\wh{\til{\mathbf{x}}}^\hbar_{\Gamma^{(0)};j})^{-1}$ and $\Phi^{\hbar_i}(\wh{\mathbf{x}}^\hbar_{\Gamma^{(0)};i} + \wh{\mathbf{x}}^\hbar_{\Gamma^{(0)};j})$ commute with each other, which holds by Lem.\ref{lem:commuting} and the Heisenberg commutation relation of the arguments. So, by inspection one observes that
\begin{align*}
& (\mathbf{K}^\hbar_{\Gamma^{(0)}\mut{j}\Gamma^{(4)}} \, \mathbf{K}^\hbar_{\Gamma^{(4)}\mut{i}\Gamma^{(5)}} \, \mathbf{P}_{(i\, j)})^{-1} \, (\mathbf{K}^\hbar_{\Gamma^{(0)}\mut{i}\Gamma^{(1)}} \, \mathbf{K}^\hbar_{\Gamma^{(1)}\mut{j}\Gamma^{(2)}} \, \mathbf{K}^\hbar_{\Gamma^{(2)}\mut{i}\Gamma^{(3)}} ) \\
& = \mathbf{P}_{(i\, j)} \, (\mathbf{K}'_{\Gamma^{(4)}\mut{i}\Gamma^{(5)}})^{-1} \, (\mathbf{K}'_{\Gamma^{(0)}\mut{j}\Gamma^{(4)}})^{-1}  \,\, \mathbf{S}_{(\mathbf{c},\mathbf{0})} \, \mathbf{K}'_{\Gamma^{(0)}\mut{i}\Gamma^{(1)}} \, \mathbf{K}'_{\Gamma^{(1)}\mut{j}\Gamma^{(2)}} \, \mathbf{K}'_{\Gamma^{(2)}\mut{i}\Gamma^{(3)}}.
\end{align*}
In this last product, each of the seven factors is a special affine shift operator, hence so is the product (Cor.\ref{cor:special_affine_shift_operators_form_a_group}). By \eqref{eq:A2_to_prove} of Prop.\ref{prop:FG_A2_identity} this equals $c_{A_2} \cdot \mathrm{id}_{\mathscr{H}_{\Gamma^{(3)}}}$, hence by Lem.\ref{lem:scalar_special_affine_shift_operator_is_identity} we deduce that $c_{A_2}=1$ as desired. Again, if one would like, one may compute explicitly the composition of the above seven special affine shift operators and prove that it is the identity, without relying on Prop.\ref{prop:FG_A2_identity}. \qed

\subsection{Heuristic proofs of the hexagon and octagon identities of the non-compact quantum dilogarithm}
\label{subsec:heuristic_proofs_of_hexagon_and_octagon_identities}

We shall find out soon that Fock-Goncharov's $B_2$-type operator identity \eqref{eq:B2_to_prove} of Prop.\ref{prop:FG_B2_identity} is equivalent to the following, up to a constant and modulo some easier operator identities dealt with in \S\ref{subsec:known_operator_identities}.
\begin{conjecture}[the hexagon equation of the non-compact quantum dilogarithm]
\label{conj:hexagon}
Let $P$, $Q$, and $V$ be as in Lem.\ref{lem:linear_combination_of_Heisenberg_operators}. Then
$$
\Phi^{2\hbar}(2P) \Phi^\hbar(Q) = \Phi^\hbar(Q) \Phi^{2\hbar}(2P+2Q) \Phi^\hbar(2P+Q) \Phi^{2\hbar}(2P) .
$$
\end{conjecture}
As in the pentagon case, we will need slight a variation of this.
\begin{corollary}
\label{cor:hexagon}
Let $P$, $Q$, and $V$ be as in Lem.\ref{lem:linear_combination_of_Heisenberg_operators}. Then
$$
\Phi^{2\hbar}(2P) \Phi^\hbar(Q) \Phi^{2\hbar}(-2P) = c_{2\hbar} \, \Phi^\hbar(Q) \Phi^{2\hbar}(2P+2Q) \Phi^\hbar(2P+Q) \exp\left( \frac{(2P)^2}{4\pi\sqrt{-1}(2\hbar)} \right)  .
$$
\end{corollary}

Instead of attempting to give a rigorous proof of Conjecture \ref{conj:hexagon}, here I only give a heuristic argument to `show' this identity by establishing the corresponding identity for the compact quantum dilogarithm. First, let $\mathbf{q}=e^{\pi\sqrt{-1}h}$ be a complex number of modulus strictly less than $1$; then so is $\mathbf{q}^2$ and $1/\mathbf{q}^\vee=e^{-\pi\sqrt{-1}/h}$. We choose a square root of $\mathbf{q}^\vee$ so that $1/\sqrt{\mathbf{q}^\vee} = e^{-\pi\sqrt{-1}/(2h)}$, which is also of modulus less than $1$. From the infinite product expression \eqref{eq:Psi_q} of the compact quantum dilogarithms one immediately observes
\begin{align}
\label{eq:product_of_two_Psi_q_squared}
\Psi^{\mathbf{q}^2}(\mathbf{q}z)\, \Psi^{\mathbf{q}^2}(\mathbf{q}^{-1}z) = \Psi^\mathbf{q}(z),
\end{align}
and likewise $\Psi^{1/\mathbf{q}^\vee}( (\sqrt{\mathbf{q}^\vee})^{-1} z)\, \Psi^{1/\mathbf{q}^\vee}(\sqrt{\mathbf{q}^\vee} z) = \Psi^{1/\sqrt{\mathbf{q}^\vee}}(z)$. We let $\mathbf{P}$ and $\mathbf{Q}$ be formal variables satisfying $[\mathbf{P},\mathbf{Q}]=2\pi\sqrt{-1}h$, so that by e.g. the BCH formula one has $e^\mathbf{P} e^\mathbf{Q} = \mathbf{q}^2 \, e^\mathbf{Q} e^\mathbf{P}$ and $e^{\mathbf{P}/h} e^{\mathbf{Q}/h} = (\mathbf{q}^\vee)^2 \, e^{\mathbf{Q}/h} e^{\mathbf{P}/h}$; at the moment, one may take this only heuristically.

\vs

Observe
\begin{align*}
& \Psi^{\mathbf{q}^2}(e^{2{\bf P}}) \Psi^\mathbf{q}(e^{\bf Q}) \\
& \stackrel{\eqref{eq:product_of_two_Psi_q_squared}}{=} \ul{ \Psi^{\mathbf{q}^2}(e^{2{\bf P}}) \Psi^{\mathbf{q}^2}(\mathbf{q} e^{\bf Q}) } \Psi^{\mathbf{q}^2}(\mathbf{q}^{-1} e^{\bf Q}) \\
& \stackrel{\eqref{eq:compact_QD_pentagon}}{=} \Psi^{\mathbf{q}^2}(\mathbf{q} e^{\bf Q}) \Psi^{\mathbf{q}^2}(\mathbf{q}^3 e^{\bf Q} e^{2{\bf P}}) \ul{ \Psi^{\mathbf{q}^2}(e^{2{\bf P}}) \Psi^{\mathbf{q}^2}(\mathbf{q}^{-1} e^{\bf Q}) } \\
& \stackrel{\eqref{eq:compact_QD_pentagon}}{=} \Psi^{\mathbf{q}^2}(\mathbf{q} e^{\bf Q}) \ul{ \Psi^{\mathbf{q}^2}(\mathbf{q}^3 e^{\bf Q} e^{2{\bf P}}) \Psi^{\mathbf{q}^2}(\mathbf{q}^{-1} e^{\bf Q}) } \Psi^{\mathbf{q}^2}( \mathbf{q} e^{\bf Q} e^{2{\bf P}}) \Psi^{\mathbf{q}^2}(e^{2{\bf P}})  \\
& \stackrel{\eqref{eq:compact_QD_pentagon}}{=} \ul{ \Psi^{\mathbf{q}^2}(\mathbf{q} e^{\bf Q}) \Psi^{\mathbf{q}^2}(\mathbf{q}^{-1} e^{\bf Q}) }\,  \Psi^{\mathbf{q}^2}(\mathbf{q}^4 e^{2{\bf Q}} e^{2{\bf P}}) \, \ul{ \Psi^{\mathbf{q}^2}(\mathbf{q}^3 e^{\bf Q} e^{2{\bf P}}) \Psi^{\mathbf{q}^2}( \mathbf{q} e^{\bf Q} e^{2{\bf P}}) } \, \Psi^{\mathbf{q}^2}(e^{2{\bf P}})  \\
& \stackrel{\eqref{eq:product_of_two_Psi_q_squared}}{=} \Psi^\mathbf{q}(e^{\bf Q}) \Psi^{\mathbf{q}^2}(e^{2{\bf Q} + 2{\bf P}}) \Psi^\mathbf{q}(e^{ {\bf Q} + 2{\bf P} }) \Psi^{\mathbf{q}^2}(e^{2 {\bf P}}),
\end{align*}
where the pentagon equations \eqref{eq:compact_QD_pentagon} for $\Psi^{\mathbf{q}^2}$ are applied to the two arguments $e^{2\mathbf{P}}$ and $\mathbf{q}^{\pm 1}e^{\mathbf{Q}}$ because $(e^{2\mathbf{P}})(\mathbf{q}^{\pm 1}e^\mathbf{Q}) = (\mathbf{q}^2)^2 (\mathbf{q}^{\pm 1}e^\mathbf{Q})(e^{2\mathbf{P}})$, and to $\mathbf{q}^3 e^\mathbf{Q} e^{2\mathbf{P}}$ and $\mathbf{q}^{-1} e^\mathbf{Q}$ because $(\mathbf{q}^3 e^\mathbf{Q} e^{2\mathbf{P}})(\mathbf{q}^{-1} e^\mathbf{Q}) = (\mathbf{q}^2)^2 (\mathbf{q}^{-1} e^\mathbf{Q}) (\mathbf{q}^3 e^\mathbf{Q} e^{2\mathbf{P}})$, while for the last equality we used $e^{2\mathbf{Q}} e^{2\mathbf{P}} = \mathbf{q}^{-4} e^{2\mathbf{Q}+2\mathbf{P}}$ and $e^\mathbf{Q} e^{2\mathbf{P}} = \mathbf{q}^{-2} e^{\mathbf{Q}+2\mathbf{P}}$ which are deduced from the BCH formula. So we just obtained
\begin{align}
\label{eq:we_proved}
\Psi^{\mathbf{q}^2}(e^{2{\bf P}}) \, \Psi^\mathbf{q}(e^{\bf Q})
= \Psi^\mathbf{q}(e^{\bf Q}) \, \Psi^{\mathbf{q}^2}(e^{2{\bf Q} + 2{\bf P}}) \, \Psi^\mathbf{q}(e^{ {\bf Q} + 2{\bf P} }) \, \Psi^{\mathbf{q}^2}(e^{2 {\bf P}}),
\end{align}
where $\mathbf{q}=e^{\pi\sqrt{-1}h}$ and $[2\mathbf{P},\mathbf{Q}] = 4\pi\sqrt{-1}h$.

\vs

Consider the identity \eqref{eq:we_proved} when each $\mathbf{q}$ is replaced by $1/\sqrt{\mathbf{q}^\vee} = e^{\pi\sqrt{-1}/(-2h)}$, while $2\mathbf{P}$ and $\mathbf{Q}$ are replaced respectively by $\mathbf{Q}/h$, and $\mathbf{P}/h$. This is legitimate because $[\mathbf{Q}/h,\mathbf{P}/h] = - 2\pi \sqrt{-1} h / h^2 = 4\pi\sqrt{-1}(-1/(2h))$.
\begin{align}
\label{eq:we_proved2}
\Psi^{1/\mathbf{q}^\vee}(e^{ {\bf Q}/h }) \Psi^{1/\sqrt{\mathbf{q}^\vee}}(e^{{\bf P}/h}) = \Psi^{1/\sqrt{\mathbf{q}^\vee}}(e^{{\bf P}/h}) \Psi^{1/\mathbf{q}^\vee}(e^{(2{\bf P}+{\bf Q})/h}) \Psi^{1/\sqrt{\mathbf{q}^\vee}}(e^{ ({\bf P}+ {\bf Q})/h }) \Psi^{1/\mathbf{q}^\vee}(e^{ {\bf Q}/h }).
\end{align}
Meanwhile, it is natural to regard that each of $e^{2\mathbf{P}}$ and $e^{\mathbf{Q}}$ commutes with each of $e^{\mathbf{Q}/h}$ and $e^{\mathbf{P}/h}$, e.g. from the point of view of the BCH formula, for
$$
[2\mathbf{P}, \mathbf{Q}/h] = 4\pi\sqrt{-1}, \qquad [2\mathbf{P},\mathbf{P}/h]=0, \qquad [\mathbf{Q},\mathbf{Q}/h] = 0, \qquad [\mathbf{Q},\mathbf{P}/h]=-2\pi\sqrt{-1};
$$
further, we assume that any continuous function in each of $e^{2\mathbf{P}}$ and $e^{\mathbf{Q}}$ commutes with any continuous function in each of $e^{\mathbf{Q}/h}$ and $e^{\mathbf{P}/h}$. This is only heuristic, and perhaps the functional calculus for explicit operators is a way to make such an argument rigorous; however, let us not bother doing this here. Also, since $\mathbf{q}^2 = e^{\pi\sqrt{-1}(2h)}$, we have $(\mathbf{q}^2)^\vee = e^{\pi\sqrt{-1}/(2h)} = \sqrt{\mathbf{q}^\vee}$. So
\begin{align*}
\Phi^{2h}(2{\bf P}) \Phi^h({\bf Q}) & \stackrel{\eqref{eq:Phi_h_as_ratio}}{=} \Psi^{\mathbf{q}^2}(e^{2{\bf P}})\, ( \Psi^{1/\sqrt{\mathbf{q}^\vee}}(e^{\cancel{2}{\bf P}/(\cancel{2}h)}) )^{-1}\, \Psi^\mathbf{q}(e^{\bf Q}) \, (\Psi^{1/\mathbf{q}^\vee}(e^{ {\bf Q}/h }) )^{-1} \\
& \stackrel{\vee}{=} \ul{ \Psi^{\mathbf{q}^2}(e^{2{\bf P}}) \, \Psi^\mathbf{q}(e^{\bf Q})  } \,\, \ul{ ( \Psi^{1/\sqrt{\mathbf{q}^\vee}}(e^{{\bf P}/h}) )^{-1} \, (\Psi^{1/\mathbf{q}^\vee}(e^{ {\bf Q}/h }) )^{-1} },
\end{align*}
where for the checked equality we used the commutation relations. Meanwhile,
\begin{align*}
& \Phi^h({\bf Q}) \Phi^{2h}(2{\bf P} + 2{\bf Q}) \Phi^h(2{\bf P} + {\bf Q}) \Phi^{2h}(2{\bf P}) \\
& \stackrel{\eqref{eq:Phi_h_as_ratio}}{=} \Psi^\mathbf{q}(e^{\bf Q}) \, (\Psi^{1/\mathbf{q}^\vee}(e^{ {\bf Q}/h }))^{-1} \, \Psi^{\mathbf{q}^2}(e^{2{\bf P} + 2{\bf Q}}) \, ( \Psi^{1/\sqrt{\mathbf{q}^\vee}}(e^{ (\cancel{2}{\bf P} + \cancel{2}{\bf Q})/(\cancel{2}h)}))^{-1} \\
& \quad \cdot \Psi^\mathbf{q}(e^{2{\bf P}+{\bf Q}}) \, ( \Psi^{1/\mathbf{q}^\vee}(e^{(2{\bf P}+{\bf Q})/h}))^{-1} \, \Psi^{\mathbf{q}^2}(e^{2{\bf P}}) \, (\Psi^{1/\sqrt{\mathbf{q}^\vee}}(e^{\cancel{2}{\bf P}/(\cancel{2}h)}))^{-1} \\
& \stackrel{\vee}{=} \ul{ \Psi^\mathbf{q}(e^{\bf Q}) \, \Psi^{\mathbf{q}^2}(e^{2{\bf P} + 2{\bf Q}}) \, \Psi^\mathbf{q}(e^{2{\bf P}+{\bf Q}}) \, \Psi^{\mathbf{q}^2}(e^{2{\bf P}}) } \\
& \quad \cdot 
\ul{ (\Psi^{1/\mathbf{q}^\vee}(e^{ {\bf Q}/h }))^{-1}
( \Psi^{1/\sqrt{\mathbf{q}^\vee}}(e^{ ({\bf P} + {\bf Q})/h}))^{-1}
( \Psi^{1/\mathbf{q}^\vee}(e^{(2{\bf P}+{\bf Q})/h}))^{-1}
(\Psi^{1/\sqrt{\mathbf{q}^\vee}}(e^{{\bf P}/h}))^{-1} },
\end{align*}
where for the checked equality we used the commutation relations. A reader can immediately verify by inspection that the application of \eqref{eq:we_proved} and \eqref{eq:we_proved2} on the underlined parts yields $\Phi^{2h}(2{\bf P}) \Phi^h({\bf Q})  = \Phi^h({\bf Q}) \Phi^{2h}(2{\bf P} + 2{\bf Q}) \Phi^h(2{\bf P} + {\bf Q}) \Phi^{2h}(2{\bf P})$. We then (heuristically) take a `limit' as $\mathrm{Im}(h) \to 0$ to obtained the desired result.

\vs

Likewise, one may also obtain a heuristic proof of the identity corresponding to Fock-Goncharov's $G_2$-type identity \eqref{eq:G2_to_prove} of Prop.\ref{prop:FG_G2_identity}. This heuristic proof of the $G_2$ identity can be found in \cite{KY}; here let me give a sketch. One first observes $\Psi^{\mathbf{q}^3}(\mathbf{q}^{-2}z) \Psi^{\mathbf{q}^3}(z) \Psi^{\mathbf{q}^3}(\mathbf{q}^2z) = \Psi^\mathbf{q}(z)$; using this and the pentagon equation for $\Psi^{\mathbf{q}^3}$ one can show
$$
\Psi^{\mathbf{q}^3}(e^{3\mathbf{P}})\Psi^\mathbf{q}(e^\mathbf{Q}) = \Psi^\mathbf{q}(e^\mathbf{Q}) \Psi^{\mathbf{q}^3}(e^{3\mathbf{Q}+3\mathbf{P}}) \Psi^\mathbf{q}(e^{2\mathbf{Q}+3\mathbf{P}}) \Psi^{\mathbf{q}^3}(e^{3\mathbf{Q}+6\mathbf{P}}) \Psi^\mathbf{q}(e^{\mathbf{Q}+3\mathbf{P}}) \Psi^{\mathbf{q}^3}(e^{3\mathbf{P}})
$$
where $[\mathbf{P},\mathbf{Q}]=2\pi\sqrt{-1}h$, which suggests
\begin{conjecture}[the octagon equation of the non-compact quantum dilogarithm]
\label{conj:octagon}
Let $P$, $Q$, and $V$ be as in Lem.\ref{lem:linear_combination_of_Heisenberg_operators}. Then
$$
\Phi^{3\hbar}(3P) \Phi^\hbar(Q) = \Phi^\hbar(Q) \Phi^{3\hbar}(3P+3Q) \Phi^\hbar(3P+2Q) \Phi^{3\hbar}(6P+3Q) \Phi^\hbar(3P+Q) \Phi^{3\hbar}(3P).
$$
\end{conjecture}
Ivan Ip told me that the two identites of unitary operators in the above Conjectures \ref{conj:hexagon} and \ref{conj:octagon} would follow from certain modification of the operator identities that he established in \cite[\S3]{Ip}; I myself was not yet able to obtain a rigorous proof this way.

\subsection{$B_2$ identity: the `hexagon' relation}

\begin{proposition}
\label{prop:c_B2_is_1}
$c_{B_2}=1$.
\end{proposition}

{\it Proof.} Let us work with the case  $\varepsilon^{(0)}_{ij}=2$, $\varepsilon^{(0)}_{ji}=-1$; I claim that a proof for the other case is essentially the same. From \eqref{eq:varepsilon_prime_formula} we get
\begin{align}
\label{eq:hexagon_epsilons2}
\varepsilon^{(1)}_{ij}= - \varepsilon^{(0)}_{ij} = - 2,  \qquad \varepsilon^{(4)}_{ji}= - \varepsilon^{(0)}_{ji} = 1, 
\end{align}

The LHS of \eqref{eq:B2_to_prove} is
\begin{align*}
& \mathbf{K}^\hbar_{\Gamma^{(0)}\mut{i}\Gamma^{(1)}} \, \mathbf{K}^\hbar_{\Gamma^{(1)}\mut{j}\Gamma^{(2)}} \, \mathbf{K}^\hbar_{\Gamma^{(2)}\mut{i}\Gamma^{(3)}} \\
& = 
(\mathbf{K}^{\sharp \hbar}_{\Gamma^{(0)}\mut{i}\Gamma^{(1)}} \, \mathbf{K}'_{\Gamma^{(0)}\mut{i}\Gamma^{(1)}})
(\mathbf{K}^{\sharp \hbar}_{\Gamma^{(1)}\mut{j}\Gamma^{(2)}} \, \mathbf{K}'_{\Gamma^{(1)}\mut{j}\Gamma^{(2)}})
(\mathbf{K}^{\sharp \hbar}_{\Gamma^{(2)}\mut{i}\Gamma^{(3)}} \, \mathbf{K}'_{\Gamma^{(2)}\mut{i}\Gamma^{(3)}}) \\
& = \Phi^{\hbar_i}(\wh{\mathbf{x}}^\hbar_{\Gamma^{(0)};i})
\Phi^{\hbar_i}(\wh{\til{\mathbf{x}}}^\hbar_{\Gamma^{(0)};i})^{-1}
\mathbf{K}'_{\Gamma^{(0)}\mut{i}\Gamma^{(1)}}
\Phi^{\hbar_j}(\wh{\mathbf{x}}^\hbar_{\Gamma^{(1)};j})
\Phi^{\hbar_j}(\wh{\til{\mathbf{x}}}^\hbar_{\Gamma^{(1)};j})^{-1}
\underbrace{ \mathbf{K}'_{\Gamma^{(1)}\mut{j}\Gamma^{(2)}} }_{\mbox{\tiny move to right}} \\
& \quad \cdot 
\ul{ \Phi^{\hbar_i}(\wh{\mathbf{x}}^\hbar_{\Gamma^{(2)};i})
\Phi^{\hbar_i}(\wh{\til{\mathbf{x}}}^\hbar_{\Gamma^{(2)};i})^{-1} } \,
\mathbf{K}'_{\Gamma^{(2)}\mut{i}\Gamma^{(3)}} \\
& = \Phi^{\hbar_i}(\wh{\mathbf{x}}^\hbar_{\Gamma^{(0)};i})
\Phi^{\hbar_i}(\wh{\til{\mathbf{x}}}^\hbar_{\Gamma^{(0)};i})^{-1}
\underbrace{ \mathbf{K}'_{\Gamma^{(0)}\mut{i}\Gamma^{(1)}} }_{\mbox{\tiny move to right}}
\Phi^{\hbar_j}(\wh{\mathbf{x}}^\hbar_{\Gamma^{(1)};j})
\Phi^{\hbar_j}(\wh{\til{\mathbf{x}}}^\hbar_{\Gamma^{(1)};j})^{-1}
 \\
& \quad \cdot 
\Phi^{\hbar_i}(\wh{\mathbf{x}}^\hbar_{\Gamma^{(1)};i} + \underbrace{ [\varepsilon^{(1)}_{ij}]_+}_{=0} \, \wh{\mathbf{x}}^\hbar_{\Gamma^{(1)};j})
\Phi^{\hbar_i}(\wh{\til{\mathbf{x}}}^\hbar_{\Gamma^{(1)};i} + \underbrace{ [\varepsilon^{(1)}_{ij}]_+}_{=0} \, \wh{\til{\mathbf{x}}}^\hbar_{\Gamma^{(1)};j} )^{-1} \,
\mathbf{K}'_{\Gamma^{(1)}\mut{j}\Gamma^{(2)}} \,
\mathbf{K}'_{\Gamma^{(2)}\mut{i}\Gamma^{(3)}} \quad (\because {\rm Lem}.\ref{lem:unitary_conjugation_commutes_with_functional_calculus}, \ref{lem:bf_K_prime_conjugation_on_bf_x}) \\
& = \Phi^{\hbar_i}(\wh{\mathbf{x}}^\hbar_{\Gamma^{(0)};i})
\Phi^{\hbar_i}(\wh{\til{\mathbf{x}}}^\hbar_{\Gamma^{(0)};i})^{-1}
\,
\Phi^{\hbar_j}(\wh{\mathbf{x}}^\hbar_{\Gamma^{(0)};j} + \underbrace{ [\varepsilon^{(0)}_{ji}]_+}_{=0} \, \wh{\mathbf{x}}^\hbar_{\Gamma^{(0)};i})
\Phi^{\hbar_j}(\wh{\til{\mathbf{x}}}^\hbar_{\Gamma^{(0)};j} + \underbrace{ [\varepsilon^{(0)}_{ji}]_+}_{=0} \, \wh{\til{\mathbf{x}}}^\hbar_{\Gamma^{(0)};i} )^{-1}
 \\
& \quad \cdot 
\Phi^{\hbar_i}(-\wh{\mathbf{x}}^\hbar_{\Gamma^{(0)};i} )
\Phi^{\hbar_i}(-\wh{\til{\mathbf{x}}}^\hbar_{\Gamma^{(0)};i}  )^{-1} \,
\mathbf{K}'_{\Gamma^{(0)}\mut{i}\Gamma^{(1)}} \,
\mathbf{K}'_{\Gamma^{(1)}\mut{j}\Gamma^{(2)}} \,
\mathbf{K}'_{\Gamma^{(2)}\mut{i}\Gamma^{(3)}} \quad (\because {\rm Lem}.\ref{lem:unitary_conjugation_commutes_with_functional_calculus}, \ref{lem:bf_K_prime_conjugation_on_bf_x} ).
\end{align*}
Like in the proof of Prop.\ref{prop:c_A2_is_1}, we use Lem.\ref{lem:commuting} (and Heisenberg commutation relations) to move some factors around, so that the LHS of \eqref{eq:B2_to_prove} becomes
\begin{align*}
& \ul{ \Phi^{\hbar_i}(\wh{\mathbf{x}}^\hbar_{\Gamma^{(0)};i})
\Phi^{\hbar_j}(\wh{\mathbf{x}}^\hbar_{\Gamma^{(0)};j} )
\Phi^{\hbar_i}(-\wh{\mathbf{x}}^\hbar_{\Gamma^{(0)};i} ) } \,
\,
\ul{ \Phi^{\hbar_i}(\wh{\til{\mathbf{x}}}^\hbar_{\Gamma^{(0)};i})^{-1}
\Phi^{\hbar_j}(\wh{\til{\mathbf{x}}}^\hbar_{\Gamma^{(0)};j} )^{-1}
\Phi^{\hbar_i}(-\wh{\til{\mathbf{x}}}^\hbar_{\Gamma^{(0)};i}  )^{-1} } \\
& \quad \cdot 
 \,
\mathbf{K}'_{\Gamma^{(0)}\mut{i}\Gamma^{(1)}} \,
\mathbf{K}'_{\Gamma^{(1)}\mut{j}\Gamma^{(2)}} \,
\mathbf{K}'_{\Gamma^{(2)}\mut{i}\Gamma^{(3)}}.
\end{align*}
From $\varepsilon^{(0)}_{ij}=2$ and $\varepsilon^{(0)}_{ji}=-1$ note that $d_j = 2 d_i$, so $\hbar_i = 2\hbar_j$. As $[\wh{\mathbf{x}}^\hbar_{\Gamma^{(0)};i}, \wh{\mathbf{x}}^\hbar_{\Gamma^{(0)};j}] \stackrel{\eqref{eq:Heisenberg_relations_for_D_q_Gamma}}{=} 2\pi\sqrt{-1} \, \hbar_j \, \varepsilon^{(0)}_{ij} \cdot \mathrm{id} = 4\pi\sqrt{-1} \hbar_j \cdot \mathrm{id} = [-\wh{\til{\mathbf{x}}}^\hbar_{\Gamma^{(0)};i}, \wh{\til{\mathbf{x}}}^\hbar_{\Gamma^{(0)};j}]$, one may apply Cor.\ref{cor:hexagon} to the two underlined parts above, so that the LHS of \eqref{eq:B2_to_prove} becomes
\begin{align*}
& \cancel{ c_{2\hbar_j} } \, \cancel{ c_{2\hbar_j}^{-1} }
\Phi^{\hbar_j}(\wh{\mathbf{x}}^\hbar_{\Gamma^{(0)};j})
\Phi^{\hbar_i}(\wh{\mathbf{x}}^\hbar_{\Gamma^{(0)};i} + 2\, \wh{\mathbf{x}}^\hbar_{\Gamma^{(0)};j})
\Phi^{\hbar_j}(\wh{\mathbf{x}}^\hbar_{\Gamma^{(0)};i} + \wh{\mathbf{x}}^\hbar_{\Gamma^{(0)};j})
\underbrace{ \exp\left( \frac{(\wh{\mathbf{x}}^\hbar_{\Gamma^{(0)};i})^2}{4\pi\sqrt{-1}\hbar_i} \right) \exp\left( - \frac{(-\wh{\til{\mathbf{x}}}^\hbar_{\Gamma^{(0)};i})^2}{4\pi\sqrt{-1} \hbar_i} \right) } \\
& \quad\cdot
\Phi^{\hbar_j}(-\wh{\til{\mathbf{x}}}^\hbar_{\Gamma^{(0)};i} + \wh{\til{\mathbf{x}}}^\hbar_{\Gamma^{(0)};j})^{-1}
\Phi^{\hbar_i}(-\wh{\til{\mathbf{x}}}^\hbar_{\Gamma^{(0)};i} + 2\, \wh{\til{\mathbf{x}}}^\hbar_{\Gamma^{(0)};j})^{-1}
\Phi^{\hbar_j}(\wh{\til{\mathbf{x}}}^\hbar_{\Gamma^{(0)};j})^{-1} \,
\mathbf{K}'_{\Gamma^{(0)}\mut{i}\Gamma^{(1)}} \,
\mathbf{K}'_{\Gamma^{(1)}\mut{j}\Gamma^{(2)}} \,
\mathbf{K}'_{\Gamma^{(2)}\mut{i}\Gamma^{(3)}},
\end{align*}
where now Lem.\ref{lem:product_of_two_quadratic_exponential_operators} lets us to replace the underbraced part by $\mathbf{S}_{(\mathbf{c},\mathbf{0})}$, where $\mathbf{c}=(c_{\ell m})_{\ell, m\in \{1,\ldots,n\}}$ is given by
$$
\mathbf{c} ~ : ~ \quad c_{\ell\ell}=1, \quad \forall \ell=1,\ldots,n, \qquad
c_{\ell i} = -\varepsilon^{(0)}_{i\ell}, \quad \forall \ell \neq i, \qquad c_{\ell m} = 0 \quad \mbox{otherwise}.
$$
In particular, note the cancellation of the two constants $c_{2\hbar_j}$ and $c_{2\hbar_j}^{-1}$. Note that this $\mathbf{c}$ is the same as the one appearing in the proof of Prop.\ref{prop:c_A2_is_1}, except that we now have $d_j = 2d_i$ and $\varepsilon^{(0)}_{ij}=2$. One can carefully verify that we still have
$$
\mathbf{S}_{(\mathbf{c},\mathbf{0})} \, \wh{\til{\mathbf{x}}}^\hbar_{\Gamma^{(0)};i} \, \mathbf{S}_{(\mathbf{c},\mathbf{0})}^{-1}  = \wh{\til{\mathbf{x}}}^\hbar_{\Gamma^{(0)};i}, \qquad
\mathbf{S}_{(\mathbf{c},\mathbf{0})} \, \wh{\til{\mathbf{x}}}^\hbar_{\Gamma^{(0)};j} \, \mathbf{S}_{(\mathbf{c},\mathbf{0})}^{-1} = \wh{\til{\mathbf{x}}}^\hbar_{\Gamma^{(0)};j} + \wh{\til{\mathbf{x}}}^\hbar_{\Gamma^{(0)};i}.
$$
So, by moving $\mathbf{S}_{(\mathbf{c},\mathbf{0})}$ to the right with the help of Lem.\ref{lem:unitary_conjugation_commutes_with_functional_calculus}, the LHS of \eqref{eq:B2_to_prove} becomes
\begin{align*}
& \Phi^{\hbar_j}(\wh{\mathbf{x}}^\hbar_{\Gamma^{(0)};j})
\Phi^{\hbar_i}(\wh{\mathbf{x}}^\hbar_{\Gamma^{(0)};i} + 2\, \wh{\mathbf{x}}^\hbar_{\Gamma^{(0)};j})
\Phi^{\hbar_j}(\wh{\mathbf{x}}^\hbar_{\Gamma^{(0)};i} + \wh{\mathbf{x}}^\hbar_{\Gamma^{(0)};j}) \,
 \\
& \cdot
\Phi^{\hbar_j}( \wh{\til{\mathbf{x}}}^\hbar_{\Gamma^{(0)};j})^{-1}
\Phi^{\hbar_i}(\wh{\til{\mathbf{x}}}^\hbar_{\Gamma^{(0)};i} + 2\, \wh{\til{\mathbf{x}}}^\hbar_{\Gamma^{(0)};j})^{-1}
\Phi^{\hbar_j}(\wh{\til{\mathbf{x}}}^\hbar_{\Gamma^{(0)};i} + \wh{\til{\mathbf{x}}}^\hbar_{\Gamma^{(0)};j})^{-1} \,
\mathbf{S}_{(\mathbf{c},\mathbf{0})} \,
\mathbf{K}'_{\Gamma^{(0)}\mut{i}\Gamma^{(1)}} \,
\mathbf{K}'_{\Gamma^{(1)}\mut{j}\Gamma^{(2)}} \,
\mathbf{K}'_{\Gamma^{(2)}\mut{i}\Gamma^{(3)}}.
\end{align*}

\vs

Meanwhile, the RHS of \eqref{eq:B2_to_prove} without the constant is
\begin{align*}
& \mathbf{K}^\hbar_{\Gamma^{(0)}\mut{j}\Gamma^{(4)}} \, \mathbf{K}^\hbar_{\Gamma^{(4)}\mut{i}\Gamma^{(5)}} \, \mathbf{K}^\hbar_{\Gamma^{(5)}\mut{j}\Gamma^{(3)}} \\
& = 
(\mathbf{K}^{\sharp \hbar}_{\Gamma^{(0)}\mut{j}\Gamma^{(4)}} \, \mathbf{K}'_{\Gamma^{(0)}\mut{j}\Gamma^{(4)}})
(\mathbf{K}^{\sharp \hbar}_{\Gamma^{(4)}\mut{i}\Gamma^{(5)}} \, \mathbf{K}'_{\Gamma^{(4)}\mut{i}\Gamma^{(5)}})
(\mathbf{K}^{\sharp \hbar}_{\Gamma^{(5)}\mut{j}\Gamma^{(3)}} \, \mathbf{K}'_{\Gamma^{(5)}\mut{j}\Gamma^{(3)}}) \\
& = \Phi^{\hbar_j}(\wh{\mathbf{x}}^\hbar_{\Gamma^{(0)};j})
\Phi^{\hbar_j}(\wh{\til{\mathbf{x}}}^\hbar_{\Gamma^{(0)};j})^{-1}
\mathbf{K}'_{\Gamma^{(0)}\mut{j}\Gamma^{(4)}}
\Phi^{\hbar_i}(\wh{\mathbf{x}}^\hbar_{\Gamma^{(4)};i})
\Phi^{\hbar_i}(\wh{\til{\mathbf{x}}}^\hbar_{\Gamma^{(4)};i})^{-1}
\underbrace{ \mathbf{K}'_{\Gamma^{(4)}\mut{i}\Gamma^{(5)}} }_{\mbox{\tiny move to right}} \\
& \quad \cdot
\ul{ \Phi^{\hbar_j}(\wh{\mathbf{x}}^\hbar_{\Gamma^{(5)};j})
\Phi^{\hbar_j}(\wh{\til{\mathbf{x}}}^\hbar_{\Gamma^{(5)};j})^{-1} } \,
\mathbf{K}'_{\Gamma^{(5)}\mut{j}\Gamma^{(3)}} \\
& = \Phi^{\hbar_j}(\wh{\mathbf{x}}^\hbar_{\Gamma^{(0)};j})
\Phi^{\hbar_j}(\wh{\til{\mathbf{x}}}^\hbar_{\Gamma^{(0)};j})^{-1}
\underbrace{ \mathbf{K}'_{\Gamma^{(0)}\mut{j}\Gamma^{(4)}}  }_{\mbox{\tiny move to right}}
\Phi^{\hbar_i}(\wh{\mathbf{x}}^\hbar_{\Gamma^{(4)};i})
\Phi^{\hbar_i}(\wh{\til{\mathbf{x}}}^\hbar_{\Gamma^{(4)};i})^{-1} \\
& \quad \cdot
\Phi^{\hbar_j}(\wh{\mathbf{x}}^\hbar_{\Gamma^{(4)};j} + \underbrace{ [\varepsilon^{(4)}_{ji}]_+}_{=1} \, \wh{\mathbf{x}}^\hbar_{\Gamma^{(4)};i} ) \, 
\Phi^{\hbar_j}(\wh{\til{\mathbf{x}}}^\hbar_{\Gamma^{(4)};j} + \underbrace{ [\varepsilon^{(4)}_{ji}]_+ }_{=1} \, \wh{\til{\mathbf{x}}}^\hbar_{\Gamma^{(4)};i} )^{-1} \,
\mathbf{K}'_{\Gamma^{(4)}\mut{i}\Gamma^{(5)}}  \,
\mathbf{K}'_{\Gamma^{(5)}\mut{j}\Gamma^{(3)}}  \quad (\because {\rm Lem}.\ref{lem:unitary_conjugation_commutes_with_functional_calculus}, \ref{lem:bf_K_prime_conjugation_on_bf_x} ) \\
& = \Phi^{\hbar_j}(\wh{\mathbf{x}}^\hbar_{\Gamma^{(0)};j})
\Phi^{\hbar_j}(\wh{\til{\mathbf{x}}}^\hbar_{\Gamma^{(0)};j})^{-1}
\,
\Phi^{\hbar_i}(\wh{\mathbf{x}}^\hbar_{\Gamma^{(0)};i} + [\varepsilon^{(0)}_{ij}]_+ \, \wh{\mathbf{x}}^\hbar_{\Gamma^(0);j} ) \,
\Phi^{\hbar_i}(\wh{\til{\mathbf{x}}}^\hbar_{\Gamma^{(0)};i} + [\varepsilon^{(0)}_{ij}]_+ \, \wh{\til{\mathbf{x}}}^\hbar_{\Gamma^{(0)};j} )^{-1} \\
& \quad \cdot
\Phi^{\hbar_j}(-\wh{\mathbf{x}}^\hbar_{\Gamma^{(0)};j} +  \wh{\mathbf{x}}^\hbar_{\Gamma^{(0)};i} + [\varepsilon^{(0)}_{ij}]_+ \, \wh{\mathbf{x}}^\hbar_{\Gamma^{(0)};j} ) \, 
\Phi^{\hbar_j}(-\wh{\til{\mathbf{x}}}^\hbar_{\Gamma^{(0)};j} + \wh{\til{\mathbf{x}}}^\hbar_{\Gamma^{(0)};i} + [\varepsilon^{(0)}_{ij}]_+ \, \wh{\til{\mathbf{x}}}^\hbar_{\Gamma^{(0)};j} )^{-1} \\
& \quad \cdot
\mathbf{K}'_{\Gamma^{(0)}\mut{j}\Gamma^{(4)}} \,
\mathbf{K}'_{\Gamma^{(4)}\mut{i}\Gamma^{(5)}}  \,
\mathbf{K}'_{\Gamma^{(5)}\mut{j}\Gamma^{(3)}}  \quad (\because {\rm Lem}.\ref{lem:unitary_conjugation_commutes_with_functional_calculus}, \ref{lem:bf_K_prime_conjugation_on_bf_x} ) \\
& = \Phi^{\hbar_j}(\wh{\mathbf{x}}^\hbar_{\Gamma^{(0)};j}) \,
\Phi^{\hbar_i}(\wh{\mathbf{x}}^\hbar_{\Gamma^{(0)};i} + 2 \, \wh{\mathbf{x}}^\hbar_{\Gamma^(0);j} ) \,
\Phi^{\hbar_j}(\wh{\mathbf{x}}^\hbar_{\Gamma^{(0)};i} + \wh{\mathbf{x}}^\hbar_{\Gamma^{(0)};j} ) \\
& \quad \cdot
\Phi^{\hbar_j}(\wh{\til{\mathbf{x}}}^\hbar_{\Gamma^{(0)};j})^{-1}
\,
\Phi^{\hbar_i}(\wh{\til{\mathbf{x}}}^\hbar_{\Gamma^{(0)};i} + 2\, \wh{\til{\mathbf{x}}}^\hbar_{\Gamma^{(0)};j} )^{-1} \, 
\Phi^{\hbar_j}( \wh{\til{\mathbf{x}}}^\hbar_{\Gamma^{(0)};i} +  \wh{\til{\mathbf{x}}}^\hbar_{\Gamma^{(0)};j} )^{-1} \,
\mathbf{K}'_{\Gamma^{(0)}\mut{j}\Gamma^{(4)}} \,
\mathbf{K}'_{\Gamma^{(4)}\mut{i}\Gamma^{(5)}}  \,
\mathbf{K}'_{\Gamma^{(5)}\mut{j}\Gamma^{(3)}},
\end{align*}
where we used $\varepsilon^{(0)}_{ij}=2$, and the commutation relations coming from Lem.\ref{lem:commuting} (and Heisenberg commutation relations) to move some factors around. Thus, by inspection one has
\begin{align*}
& (\mathbf{K}^\hbar_{\Gamma^{(0)}\mut{j}\Gamma^{(4)}} \, \mathbf{K}^\hbar_{\Gamma^{(4)}\mut{i}\Gamma^{(5)}} \, \mathbf{K}^\hbar_{\Gamma^{(5)}\mut{j}\Gamma^{(3)}})^{-1} \,
\mathbf{K}^\hbar_{\Gamma^{(0)}\mut{i}\Gamma^{(1)}} \, \mathbf{K}^\hbar_{\Gamma^{(1)}\mut{j}\Gamma^{(2)}} \, \mathbf{K}^\hbar_{\Gamma^{(2)}\mut{i}\Gamma^{(3)}} \\
& = (\mathbf{K}'_{\Gamma^{(5)}\mut{j}\Gamma^{(3)}})^{-1}
(\mathbf{K}'_{\Gamma^{(4)}\mut{i}\Gamma^{(5)}})^{-1}  \,
(\mathbf{K}'_{\Gamma^{(0)}\mut{j}\Gamma^{(4)}})^{-1} \,
\mathbf{S}_{(\mathbf{c},\mathbf{0})} \,
\mathbf{K}'_{\Gamma^{(0)}\mut{i}\Gamma^{(1)}} \,
\mathbf{K}'_{\Gamma^{(1)}\mut{j}\Gamma^{(2)}} \,
\mathbf{K}'_{\Gamma^{(2)}\mut{i}\Gamma^{(3)}}.
\end{align*}
Eq.\eqref{eq:B2_to_prove} of Prop.\ref{prop:FG_B2_identity} tells us that this equals $c_{B_2} \cdot \mathrm{id}_{\mathscr{H}_{\Gamma^{(3)}}}$, while it is a special affine shift operator, because it is a composition of special affine shift operators. Hence Lem.\ref{lem:scalar_special_affine_shift_operator_is_identity} says $c_{B_2}=1$ as desired. Again, this last step can be done directly without resorting to Prop.\ref{prop:FG_B2_identity}. \qed

\subsection{$G_2$ identity: the `octagon' relation}

\begin{proposition}
\label{prop:c_G2_is_1}
$c_{G_2}=1$.
\end{proposition}

I think I presented enough detail of the proof of the other cases, so that a reader can easily construct a proof of this last proposition along a similar line. In particular, it is quite natural to expect that the constants coming from operator identities will cancel each other.

\end{document}